# FUZZY AND NEUTROSOPHIC ANALYSIS OF WOMEN WITH HIV/AIDS

### With Specific Reference To Rural Tamil Nadu in India

## W. B. Vasantha Kandasamy
## Florentin Smarandache

Translation of the Tamil interviews by
**Meena Kandasamy**

**2005**

# FUZZY AND NEUTROSOPHIC ANALYSIS OF WOMEN WITH HIV/AIDS

**With Specific Reference
To Rural Tamil Nadu in India**

**W. B. Vasantha Kandasamy**
Department of Mathematics
Indian Institute of Technology, Madras
Chennai – 600036, India
e-mail: **vasantha@iitm.ac.in**
web: **http://mat.iitm.ac.in/~wbv**

**Florentin Smarandache**
Department of Mathematics
University of New Mexico
Gallup, NM 87301, USA
e-mail: **smarand@gallup.unm.edu**

Translation of the Tamil interviews by
**Meena Kandasamy**

**2005**

**CONTENTS**

**Chapter Four**

## USE OF NEUTROSOPHIC MODELS TO ANALYZE HIV/AIDS WOMEN PATIENTS' PROBLEM

**Chapter Five**

## INTERVIEWS OF 101 WOMEN LIVING WITH HIV/AIDS

**Chapter Six**

## SUGGESTIONS FOR IMPLEMENTATION AND CONCLUSIONS

## APPENDIX

## REFERENCES

## INDEX

## ABOUT THE AUTHORS

# Preface

Fuzzy theory is one of the best tools to analyze data, when the data under study is an unsupervised one, involving uncertainty coupled with imprecision. However, fuzzy theory cannot cater to analyzing the data involved with indeterminacy. The only tool that can involve itself with indeterminacy is the neutrosophic model.

Neutrosophic models are used in the analysis of the socio-economic problems of HIV/AIDS infected women patients living in rural Tamil Nadu. Most of these women are uneducated and live in utter poverty. Till they became seriously ill they worked as daily wagers. When these women got admitted in the hospital they were either terminally ill or the disease has advanced to a chronic stage. We mainly choose this group because they are the most neglected in our male dominated society. Further, after the Aryans entered India, for their comfortable living they introduced and enforced on the natives the *Manusmriti*-based culture that is chiefly responsible for ill treatment of the women in India. The Laws of Manu treat women as a mere object. They can never be independent entities or even policy-makers for themselves. The pitiable impact of Manu on these uneducated rural men has become a part and parcel of them that they treat women very badly. This patriarchy is indirectly one of the major causes of spread of HIV/AIDS among rural women in India.

This book has six chapters. The first chapter gives a brief introduction about HIV/AIDS among rural women in Tamil Nadu. Chapter 2 is an analysis of the situation using fuzzy theory in general and Fuzzy Relational Maps (FRM) in particular. The FRM tool is specially chosen for two reasons: It can give the hidden pattern of the dynamical system and as the socio-economic condition of women infected with HIV/AIDS is interdependent we choose this special fuzzy tool called FRM.

In Chapter 3, we use Fuzzy Associative Memories (FAM) to analyze the problem because FAM is the only fuzzy tool that can give the gradations of each of the nodes/concepts. In Chapter 4 for the first time we use neutrosophic theory in general and neutrosophic relational maps (NRM) in particular to analyze this problem. Interested readers can compare and analyze the two models FRM and NRM. This chapter introduces the notion of Neutrosophic Associative Memories (NAMs) that is an analogous model of Fuzzy Associative Memories. NAMs are applied to this problem and conclusions are based on this analysis.

The sixth chapter gives the translated version of the verbatim interviews of the 101 HIV/AIDS infected, rural, uneducated, and poor women of Tamil Nadu. It is important to record that most of these women talked very openly and that in spite of their sufferings and problems they did not suffer depression. These interviews taken as a part of our research grant, has been published for free distribution by the Tamil Nadu State AIDS Control Society as a book titled, "*Women. Work. Worth. Womb: 101 Real Life Stories of Rural Women with AIDS.*" We have chosen to give these interviews because all the readers may not be from rural Tamil Nadu, so only by knowing the problems faced by these women and their culture it would be easy for the reader to understand why such analysis has been taken. We have given the verbatim opinion of the expert. The final chapter is based not only on our mathematical analysis but also on the interviews, views of experts and observations. This book is not only for mathematicians but also for sociologists, economists, human rights activists, psychologists and to women activists who want to make a difference in the lives of these uneducated, rural, poverty-stricken women. We dedicate this book to women around the world who are bravely battling the AIDS epidemic.

W.B.VASANTHA KANDASAMY
FLORENTIN SMARANDACHE

**Chapter One**

# INTRODUCTION

HIV/AIDS prevention cannot be viewed as a one-time intervention; it must be accepted as a continuous, multigenerational effort that extends well into the lifetimes of our children and their children. The prevention or awareness programmes solely for rural areas – that too among uneducated women – is lacking more in India and especially in villages of Tamil Nadu than in other states. We see that they are the passive victims for no cause of their own. Especially in the case of rural women these prevention / awareness efforts will be successful only if they are coupled with sustained visible efforts to prevent discrimination by their husband's family against those women who are infected and at risk.

At one end of the spectrum, we need to tackle the underlying socio-economic and structural factors that make women vulnerable to HIV/AIDS by giving them few realistic options for self-protection. Empowerment for action is one obvious remedy. The vulnerability could be lessened through higher education, access to credit and equitable rights in the event of divorce or death of the husband or when he is a terminally ill patient. Clearly such action could simultaneously diminish women's vulnerability to other ills and similar socio-economic problems including domestic violence and unwanted pregnancy. To reduce risk options, for example, by raising taxes on risk-related products such as alcohol; policy of obligatory condom use while visiting Commercial Sex Workers (CSWs), which will eliminate the clients risk option of unprotected sex and at the same time empower the CSWs to insist on condom usage.

Another factor specially suffered by the HIV/AIDS affected uneducated women in rural areas is poverty for their day-to-day survival is impossible. Unlike in Africa, we see the reported cases of HIV/AIDS in rural areas in India are higher than in the urban areas.

The reasons we can attribute for this is that several HIV/AIDS cases in urban areas are kept secret. Secondly due to increased awareness, the urban people may be well aware of HIV/AIDS. As in South Africa, Kenya, Tanzania and Uganda, in India also married women are culturally not expected to say 'no' to sex to their husbands. Majority of these rural uneducated woman (this is true of urban educated women also) have little control over either abstinence or condom use at home or their husbands extra marital sexual activities. These cultural inequalities thus increase the vulnerability of women to HIV infection. Here out of 101 HIV/AIDS infected rural women, 91 have been infected by their husbands with whom they had remained monogamous.

Also most of these women being uneducated and from rural areas suffer due to economic inequalities. When a husband dies the relatives are supposed to share or take care of the property of the deceased, women are not supposed to own property, this is not only the case, we see in case of; especially in a specific community if the husband dies of HIV/AIDS and the wife is infected, she is chased out of the home. In many cases she has no one to support her. Being in her early forties she has no parents or her relatives, so many suffer this mental agony. Their (wife's) property, which was acquired, by the husband and his relatives as dowry, are usurped by them and these women land on streets.

Most of these rural women are agricultural labourers and the money earned by them has been used by their husband and his relatives. As in the case of southern Africa, in India also, the high-risk population included long distance truck-drivers, migrant workers, their wives and the CSWs. It is unfortunate to state that even today in rural areas the high rates of untreated sexually transmitted diseases and low rates of use of condom while visiting CSWs (as the men are in the complete drunken state) has facilitated the rapid spread of the infection among married women and their children.

In India and especially in Tamil Nadu sero prevalence rates among sex workers, patients with sexually transmitted disease (STD) clinics and pregnant women were significant and rising steadily. In India although rare among women, premarital and extramarital sexual behaviors are common and tolerated among men even though it is publicly disapproved. Women and men who work as or patronize the commercial sex workers are particularly at a high risk for HIV infection though; CSWs are legally banned in India. Most men have first sexual experience with CSWs. Study of the rate of CSWs affected by HIV is essential.

It is unfortunate that in the case of women, they suffer from varied types of symptoms, which often leads to difficulty in receiving accurate diagnosis. Many suffer for years without any appropriate treatment. For, until 1993, the medical definition of HIV in women did not include the most prevalent indications of HIV in women.

Among women most notably in India and in Tamil Nadu, safe sex negotiation is often embedded in a system of power relations and emotions that can diminish a women's ability to make her partner use condoms for intercourse. Many researchers have studied sexual transmission of HIV in women whose male partners had become infected through CSWs. It is clear that it would be very difficult to implement any effective prevention programs for these women, without putting in place social supports, such as means of independent living and health care.

In many cases in rural Tamil Nadu most men prefer and seek sex from other women and CSWs. So the use of condoms in family needs skilled training of both men and women for mutual negotiation. Unless these support systems are in place, women in relationship with HIV-infected men may often lack the power to translate their knowledge of HIV transmission into actual self-protective behaviors.

Groups such as WORLD (Women Organized to Respond to Life threatening Diseases based in Oakland, California) and Women's AIDS Network of San Francisco and several other advocacy organizations are dedicated specifically to the interest of women and have voiced to protect women from contracting HIV by sexual transmission and the complex issue of maternal transmission of HIV during pregnancy. Still the disease is only a male profile in spite of the fact that women are more vulnerable to HIV/AIDS than men. Further more concern has been expressed about women's roles as 'infectors' than as 'infecteds' ignoring the real risk of infection that many women face. The media has also reinforced imagery of women as infectors, but till date the real problem women face in case of HIV/AIDS has not been analyzed. Concern about women's potential to infect others, also guided discussions of mandatory testing of pregnant women and use of condom while visiting CSWs are on the agenda and in the Indian scenario no effort is even made to test CSWs for HIV.

Although about one third of the estimated 10 millions infected adults worldwide are women, no special attention is given to effective programs to save housewives who are not presently child bearers. Particularly what work has been done often focuses mainly on pregnant women and as a result, little is known about the natural history of HIV/AIDS in women who are not pregnant.

Women with HIV have been almost invisible to the media since the epidemic began. When they have appeared and infection could be traced to partners they still were not portrayed as the innocent victims of HIV/AIDS. In most of the cases they are deserted by their family members. Still we cannot say that there are strong supportive groups for these rural uneducated poor HIV/AIDS infected women victims in Tamil Nadu or in India. Among poor uneducated women these issues are exacerbated by limited access to personal physicians, inadequate health insurance and a general fear and mistrust of health services, thus these women typically enter into treatment only after they are compromised by the disease. Social workers and NGOs must promote the development of resource co-ordination strategies that will improve access and connect these women with health care services, earlier in the process.

In rural Tamil Nadu we see the difference between the awareness levels of AIDS patients by sex has been nil which indicates that there is no difference between men and women here. The only factor is that most women are infected by their partners.

Compared with the male patients, female patients showed more strength of mind than males. Since women in rural areas have the family burden, it is high time strong social groups are formed to empower these women with awareness of HIV/AIDS, a very high level of counseling and monetary help to keep them economically well, so that they get better food and medicine.

It is important to mention that most of these rural women are not women who are drug addicts or injecting drug users etc. They are mostly passive or innocent victims of HIV/AIDS. Hence we in this book have ventured to give the real life history of 101 HIV/AIDS affected women mainly from rural areas with no education and from poor economic strata. 91 of them had been infected by their husbands. In fact 44 of them had tested or come to know that they have been affected by HIV/AIDS only after the death of their husbands. In some cases the disease was chronic and they were seriously very ill and had become invalid and immobile. A few have been deserted by their husband's family and were at the mercy of the hospital. None of them had any help from any of the women's organizations or support groups. The advent of the disease about which they were totally ignorant came as a rude shock to them. Several even at this stage are unaware of the ways of transmission.

Further we were forced to take up this study mainly because of the fact that exclusive research on HIV/AIDS infected women who were past the childbearing age have not been carried out properly. Only the HIV/AIDS study about women have been carried out in case of pregnant women and on CSWs. It is clear that the study of CSWs infected by HIV/AIDS is medically and socially very much different from the HIV/AIDS infected women who are infected by their husbands. Also the study of HIV/AIDS infected pregnant women is different.

Further in case of HIV/AIDS women, age plays a role like pre-menopause or post-menopause and so on and so forth. There is a lack of research analysis and medical insight into this problem; this might be of potential interest to scientists. Besides, we felt a need for analysis of these women's social and psychological status like torture by her husband and his family and are neglected by children. As most of them are dependent on the family their helplessness is magnified and their misery is increased.

Several of them come to know about the disease only after coming to the hospital. But uniformly when compared with men only 3% of them said they were very desperate and felt like committing suicide. Further majority of these women being diseased did not suffer any from depression in spite of being innocent victims. Majority of them were not very much dejected or desperate even though they had no one to provide support to them. Most of the women talked very openly. Some of them who had extra martial affairs, very openly confessed about it and also the cause for seeking such relations.

Thus we give the real life story of these 101 HIV/AIDS infected women interviewed by us using a linguistic questionnaire and most of these interviews are substantiated by additional observations which we have made during the course of the interview.

# FUZZY ANALYSIS OF THE EMOTIONS / FEELINGS OF HIV/AIDS INFECTED WOMEN

Since the study and analysis of the social and psychological problems related with HIV/AIDS affected women are very intricate and can not be attributed to a single cause and varies from women to women we are forced to analyze the problem using fuzzy mathematics. This chapter has two sections. In the first section we introduce the concept of Fuzzy Relational Maps (FRMs), in the second section we analyze the psychological problems of the HIV/AIDS affected women using FRMs. We have particularly chosen fuzzy relational maps because the women's problem of being infected with HIV/AIDS is related, and not an independent one as in the case of men. For majority of the women (91 out of 101) are infected by HIV/AIDS only through their infected husbands so it is pertinent to use only a mathematical model, which is interdependent. So we choose to use Fuzzy Relational Maps (FRMs). Secondly these fuzzy relational maps can be combined in two ways. The first way may give the cumulative or the over-all effect. The second method will interrelate two concepts, which cannot be interrelated directly, but has an indirect relation. Thus, at the outset we are justified in adopting this model for the analysis of capturing the socio-economic problems and emotions of the HIV/AIDS affected rural women patients.

`The second major tool that we use in this study is the Fuzzy Associative Memories (FAMs), because this model alone can give the gradations of each of the attributes or concepts relative to the other concepts. So we have chosen this model, so that by adopting this technique we can say what is the maximum cause for making the women HIV/AIDS patients, is it men, (i.e., their husbands), or is it their lack of awareness or is it their inability to protect themselves from the disease even if they are aware of it and so on. Now we give a brief description of the fuzzy mathematical tools.

## 2.1 Introduction to Fuzzy Relational Maps (FRMs)

We just give a short introduction of learning as encoding. We avoid all types of high and complicated mathematics. As our study is a psychological one, we give more importance to the approach of the mind, of how the human mind correlates, relates and organizes the changes in any system. As our study is about a change in the system that involves the feelings (emotions) of the HIV/AIDS affected women, which is an unsupervised data, we present more of the non-mathematical description of how our model works to make this book readable even by a layman.

D.O. Hebb [41] an expert researcher on behaviour says in his book, "The Organization of Behaviour", that when one cell repeatedly assists in firing another, the axon of the first cell develops synaptic knobs (or enlarges them if they already exist) in contact with the soma of the second cell.

Learning encodes information. A system learns a pattern if the system encodes the pattern in its structure. The system structure changes as the system learns the information. We may use a behavioristic encoding criterion.

Learning involves change. A system learns or adopts or cell organizes when sample data changes system parameters. In neural networks, learning means any change in any synapse. We do not identify learning with change in a neuron, though in some sense a changed neuron has learned its new state.

Synapses change more slowly than neurons change. We learn more slowly than we think. This reflects the conventional interpretation of learning as a semi permanent change. We have learned calculus if our calculus exam taking behaviour has changed from failing to passing, and has stayed that way for some time. In the same sense we can say that the artists canvas learns the pattern of paints the artist smears on its surface, or that the overgrown lawn learns the well-trimmed pattern the lawn mower imparts. In these cases the system learns when pattern simulation changes a memory medium and leaves it changed for some comparatively long stretch of time.

Learning also involves quantization. Usually a system learns only a small proportion of all patterns in the sampled pattern environment. Memory capacity is scarce. An adaptive system must efficiently and continually replace the patterns in its limited memory with patterns that accurately represent the sampled pattern environment. Learning replaces old stored patterns with new patterns. Learning forms internal representation of prototypes of sampled patterns learned prototypes define quantized patterns. Neural network models represent prototype patterns as vectors of real numbers. So we can alternatively view learning as a form of adaptive vector quantization. Thus adaptive vector quantization embeds learning in a dynamic geometry, vector quantization may be optimal according to different criteria. The quantization vectors should estimate the underlying unknown probability distribution of patterns. The distribution of prototype vectors should statistically resemble the unknown distribution of patterns.

Uniform sampling probability, provides an information theoretic criterion for an optimal quantization. We first mention the difference between the supervised and unsupervised learning in neural networks for we often say the data is an unsupervised one. The learning adjectives "supervised" and "unsupervised" stem from pattern recognition theory. The distinction depends on whether the learning algorithm uses pattern-class information. Supervised learning uses pattern class information, whereas unsupervised learning does not.

Unsupervised learning algorithms learn rapidly often on a single pass of noisy data. This makes unsupervised learning practical in many high-speed real time environment where we may not have enough time information or computational precision to use supervised technique. Thus we say the data, which we analyze, are unsupervised for they measure the emotions of the HIV/AIDS patients. Neuroscientist Dr. Hebb [1949] first developed the mechanistic theory of conjunctive (correlation) synaptic learning in his book, "*The Organization of Behaviour*". [41]

We shall now examine four unsupervised learning laws; signal Hebbian, competitive, differential Hebbian and differential competitive.

I. The deterministic signal Hebbian learning law correlates local neuronal signals.

$$\overset{\circ}{m}_{ij} = -m_{ij} + S_i^X(x_i) S_j^Y(y_j)$$

or more simply

$$\overset{\circ}{m}_{ij} = -m_{ij} + S_i(x_i) S_j(y_j)$$

The overdot denotes time differentiation. $m_{ij}$ denotes the synaptic efficacy of the synapse along the directed axonal edge from the $i^{th}$ neuron in the input neuron field $F_X$ to the $j^{th}$ neuron in the field $F_Y$. The $ij^{th}$ synaptic junction is excitatory, if $m_{ij} > 0$, inhibitory if $m_{ij} < 0$.

II. The deterministic competitive learning law modulates the signal-synaptic differences $S_i(x_i) - m_{ij}$ with the zero-one competitive signal $S_j(y_j)$

$$\overset{\circ}{m}_{ij} = S_j(y_j) \, [S_i(x_i) - m_{ij}]$$

$$S_j(y_j) = \frac{1}{1 + e^{-cy_j}}, \, c > 0,$$

a steep logistic signal function $S_j(y_j)$. However the competitive learning law can be reduced to linear competitive learning law

$$\overset{\circ}{m}_{ij} = S_j(y_j) \, (x - m_j),$$

the input pattern vector x represents the output of the $F_X$ neurons; $m_j$ is the synaptic vector with the least distance.

III. The deterministic differential Hebbian learning law correlates signal velocities as well as neuronal signals.

$$\overset{\circ}{m}_{ij} = - m_{ij} + S_i(x_i) + S_j(y_j) + \overset{\circ}{S}_i(x_i) + \overset{\circ}{S}_j(y_j).$$

The signal velocity $S_i$ decomposes with the chain rule as

$$\frac{d \, S_i(x_i)}{dt} = \frac{d \, S_i(x_i)}{d \, x_i} \frac{dx_i}{dt}.$$

IV. The deterministic competitive learning law combines competitive and differential Hebbian learning

$$\overset{\circ}{m}_{ij} = \overset{\circ}{S}_j(y_j)\left[S_i(x_i) - m_{ij}\right],$$

differential competition means learn only if change. Now the deterministic differential Hebbian learning law

$$\overset{\circ}{m}_{ij} = -m_{ij} + S_i S_j + \overset{\circ}{S}_i \overset{\circ}{S}_j.$$

and its simpler version being

$$\overset{\circ}{m}_{ij} = -m_{ij} + \overset{\circ}{S}_i \overset{\circ}{S}_j$$

arose from attempts to dynamically estimate the causal structure of fuzzy cognitive maps from sample data [Kosko, 57-68]. This led to neural interpretations of these signal-velocity learning laws as we shall see. Initiatively Hebbian correlations promote spurious causal associations among concurrently active units. Differential correlations estimate the concurrent and presumably causal variation among active units.

The deterministic competitive learning viz. FRM, which is a special case of the FCM is described now.

In Fuzzy Cognitive Maps (FCMs) causal associations are made among concurrently active units and in all the cases all the causal concepts are taken together and their association marked by edges, which always result in a square matrix called the connection matrix. But in FRMs we don't put all the causal concepts together but we divide them into two disjoint classes. Just for the sake of distinction between the two classes we name them as domain space and range space though in reality there is no difference between them. For the causal relations can be from range space to domain space or from domain space to range space. In some of the cases the concepts in the domain space (or range) space need not be mutually exclusive.

We build relation map using causal relations existing between the domain and range spaces. The FRMs can take edge values from the set $\{-1, 0, 1\}$. If the number of concepts in the domain space is n and that of the range space is m, we get an n × m causal relational matrix, which we call, as the relational matrix. This matrix will have its entries from the set $\{-1, 0, 1\}$. Using this relational matrix we say whether the dynamical system is stable or unstable. If the dynamical system is a stable one we get the hidden pattern of the system, which may a fixed point, or a limit cycle of the FRM.

In FCMs we promote the correlations between causal associations among concurrently active units. But in FRMs we divide the very causal associations into two units, for example, the relation between a teacher and a student or relation between an employee and the employer or a relation between doctor and patient and so on. Thus for us to define a FRM we need a domain space and a range space. We further assume no intermediate relation exists within the domain elements or node and the range space elements. The number of elements in the range space need not in general be equal to the number of elements in the domain space.

Thus throughout this section we assume the elements of the domain space are taken from the real vector space of dimension n and that of the range space are real vectors from the vector space of dimension m (m in general need not be equal to n). We denote by R the set of nodes $R_1,\ldots, R_m$ of the range space, where $R = \{(x_1,\ldots, x_m) \mid x_j = 0 \text{ or } 1 \}$ for $j = 1, 2,\ldots, m$. If $x_i = 1$ it means that the node $R_i$ is in the on state and if $x_i = 0$ it means that the node $R_i$ is in the off state. Similarly D denotes the nodes $D_1, D_2,\ldots, D_n$ of the domain space where $D = \{(x_1,\ldots, x_n) \mid x_j = 0 \text{ or } 1\}$ for $i = 1, 2,\ldots, n$. If $x_i = 1$ it means that the node $D_i$ is in the on state and if $x_i = 0$ it means that the node $D_i$ is in the off state.

Now we proceed on to define a FRM.

**DEFINITION 2.1.1:** *A FRM is a directed graph or a map from D to R with concepts like policies or events etc, as nodes and causalities as edges. It represents causal relations between spaces D and R .*

Let $D_i$ and $R_j$ denote the two nodes of an FRM. The directed edge from $D_i$ to $R_j$ denotes the causality of $D_i$ on $R_j$ called relations. Every edge in the FRM is weighted with a number in the set $\{0, \pm1\}$. Let $e_{ij}$ be the weight of the edge $D_iR_j$, $e_{ij} \in \{0, \pm1\}$. The weight of the edge $D_iR_j$ is positive if increase in $D_i$ implies increase in $R_j$ or decrease in $D_i$ implies decrease in $R_j$ i.e. causality of $D_i$ on $R_j$ is 1. If $e_{ij} = 0$, then $D_i$ does not have any effect on $R_j$. We do not discuss the cases when increase in $D_i$ implies decrease in $R_j$ or decrease in $D_i$ implies increase in $R_j$.

**DEFINITION 2.1.2:** *When the nodes of the FRM are fuzzy sets then they are called fuzzy nodes. FRMs with edge weights $\{0, \pm1\}$ are called simple FRMs.*

**DEFINITION 2.1.3:** *Let $D_1, \ldots, D_n$ be the nodes of the domain space D of an FRM and $R_1, \ldots, R_m$ be the nodes of the range space R of an FRM. Let the matrix E be defined as $E = (e_{ij})$ where $e_{ij}$ is the weight of the directed edge $D_iR_j$ (or $R_jD_i$), E is called the relational matrix of the FRM.*

*Note:* It is pertinent to mention here that unlike the FCMs the FRMs can be a rectangular matrix with rows corresponding to the domain space and columns corresponding to the range space. This is one of the marked differences between FRMs and FCMs.

**DEFINITION 2.1.4:** *Let $D_1, \ldots, D_n$ and $R_1, \ldots, R_m$ denote the nodes of the FRM. Let $A = (a_1, \ldots, a_n)$, $a_i \in \{0, 1\}$. A is called the instantaneous state vector of the domain space and it denotes the on-off position of the nodes at any instant. Similarly let $B = (b_1, \ldots, b_m)$ $b_i \in \{0, 1\}$. B is called the instantaneous state vector of the range space and it denotes the on-off position of the nodes at any instant $a_i = 0$ if $a_i$ is off and $a_i = 1$ if $a_i$ is on for $i = 1, 2, \ldots, n$. Similarly, $b_i = 0$ if $b_i$ is off and $b_i = 1$ if $b_i$ is on, for $i = 1, 2, \ldots, m$.*

**DEFINITION 2.1.5:** *Let $D_1, \ldots, D_n$ and $R_1, \ldots, R_m$ be the nodes of an FRM. Let $D_iR_j$ (or $R_j D_i$) be the edges of an*

*FRM, j = 1, 2,..., m and i= 1, 2,..., n. Let the edges form a directed cycle. An FRM is said to be a cycle if it posses a directed cycle. An FRM is said to be acyclic if it does not posses any directed cycle.*

**DEFINITION 2.1.6:** *An FRM with cycles is said to be an FRM with feedback.*

**DEFINITION 2.1.7:** *When there is a feedback in the FRM, i.e. when the causal relations flow through a cycle in a revolutionary manner, the FRM is called a dynamical system.*

**DEFINITION 2.1.8:** *Let $D_i R_j$ (or $R_j D_i$), $1 \leq j \leq m$, $1 \leq i \leq n$. When $R_i$ (or $D_j$) is switched on and if causality flows through edges of the cycle and if it again causes $R_i$ (or $D_j$), we say that the dynamical system goes round and round. This is true for any node $R_j$ (or $D_i$) for $1 \leq i \leq n$, (or $1 \leq j \leq m$). The equilibrium state of this dynamical system is called the hidden pattern.*

**DEFINITION 2.1.9:** *If the equilibrium state of a dynamical system is a unique state vector, then it is called a fixed point.*

Consider an FRM with $R_1$, $R_2$,…, $R_m$ and $D_1$, $D_2$,…, $D_n$ as nodes. For example, let us start the dynamical system by switching on $R_1$ (or $D_1$). Let us assume that the FRM settles down with $R_1$ and $R_m$ (or $D_1$ and $D_n$) on, i.e. the state vector remains as (1, 0, …, 0, 1) in R (or 1, 0, 0, … , 0, 1) in D), This state vector is called the fixed point.

**DEFINITION 2.1.10:** *If the FRM settles down with a state vector repeating in the form*

*$A_1 \rightarrow A_2 \rightarrow A_3 \rightarrow ... \rightarrow A_i \rightarrow A_1$ (or $B_1 \rightarrow B_2 \rightarrow ... \rightarrow B_i \rightarrow B_1$)*

*then this equilibrium is called a limit cycle.*

Now we give the methods of determining the hidden pattern.

Let $R_1$, $R_2$,…, $R_m$ and $D_1$, $D_2$,…, $D_n$ be the nodes of a FRM with feedback. Let E be the relational matrix. Let us find a hidden pattern when $D_1$ is switched on i.e. when an input is given as vector $A_1 = (1, 0, …, 0)$ in $D_1$, the data should pass through the relational matrix E. This is done by multiplying $A_1$ with the relational matrix E. Let $A_1E = (r_1, r_2,…, r_m)$, after thresholding and updating the resultant vector we get $A_1 E \in R$. Now let $B = A_1E$ we pass on B into $E^T$ and obtain $BE^T$. We update and threshold the vector $BE^T$ so that $BE^T \in D$. This procedure is repeated till we get a limit cycle or a fixed point.

**DEFINITION 2.1.11:** *Finite number of FRMs can be combined together to produce the joint effect of all the FRMs. Let $E_1$,…, $E_p$ be the relational matrices of the FRMs with nodes $R_1$, $R_2$,…, $R_m$ and $D_1$, $D_2$,…, $D_n$, then the combined FRM is represented by the relational matrix $E = E_1 + … + E_p$.*

Now we give a simple illustration of a FRM, for more about FRMs please refer [129, 140, 141].

We just recall yet another new technique using FRMs which is impossible in case of FCMs. This method is more adaptable in case of data when we are not in a position to inter-relate two systems but we know they are inter-related indirectly. Such study is possible when we adopt linked fuzzy relational maps. First we give the definition of pair wise linked FRMs.

**DEFINITION 2.1.12:** *Let us assume that we are analyzing some concepts, which are divided into 3 disjoint units. Suppose we have 3 spaces say P, Q and R we say some m set of nodes in the space P, some n set of nodes in the space Q and r set of nodes in the space R. We can directly find FRMs or directed graphs relating P and Q, FRMs or directed graphs relating Q and R.*

*Suppose we are not in a position to link or get a relation between P and R directly but in fact there exists a hidden link between them, which cannot be easily weighted; in such cases we use linked FRMs. Thus pairwise linked FRMs are those FRMs connecting three distinct spaces P, Q and R in such a way that using the pair of FRMs we obtain a FRM relating P and R.*

If $E_1$ is the connection matrix relating P and Q then $E_1$ is a $m \times n$ matrix and $E_2$ is the connection matrix relating Q and R, which is a $n \times r$ matrix. Now $E_1E_2$ is a $m \times r$ matrix which is the connection matrix relating P and R and $E_2^T E_1^T$ is the matrix relating R and P, when we have such a situation we call it the pair wise linked FRMs.

Thus by this method even if we are not in a position to get directed graph using the expert we can indirectly obtain the FRMs relating them. Now using these three FRMs and their related matrices we derive the conclusion by studying the effect of each state vector.

Now we can define 3-linked FRMs, 4-linked FRMs and so on in general, an n-linked FRM on similar lines and obtain the hidden relation, which are not explicitly calculable.

Now we proceed on to define the hidden fuzzy relational map and the hidden connection matrix in case of pairwise linked FRMs.

**DEFINITION 2.1.13:** *Let P, Q and R be three spaces with nodes related to some problem. Suppose P and Q are related by a FRM and Q and R are related by an FRM then the indirectly calculated FRM between P and R (or R and P) got as the product of the connection matrices related with the FRMs. P and Q, and Q and R is called as the Hidden connection matrix and the directed graph drawn using the Hidden connection matrix is called the Hidden Directed graph of the pairwise linked FRMs.*

## 2.2 Use of FRMs to analyze the socio economic problems of HIV/AIDS infected women patients

In this section we adopt the FRM model to analyze the socio-economic problems of HIV/AIDS affected women patients; which we have obtained from the linguistic questionnaire. This model is specially chosen for the HIV/AIDS infection in women as in majority of the cases women's problems are always inter-related with their husbands and family; unlike men who get infected on their own. We illustrate the model.

The following are taken as the attributes (concepts) related with women, which is taken as the domain space of the FRM.

$W_1$ - Child marriage / widower marriage / child married to men twice or thrice their age.

$W_2$ - Causes of women being infected with HIV/AIDS

$W_3$ - Disease untreated till it is chronic or they are in last stages of life due to full-blown AIDS.

$W_4$ - Women are not traditionally considered as breadwinners for the work they do and the money they earn to protect the family is never given any recognition.

$W_5$ - Free of depression and despair in spite of being deserted by family when they have HIV/AIDS.

$W_6$ - Faith in god / power of god will cure.

$W_7$ - Married women have acquired the disease due to their husbands

The concepts associated with society, men / husband is taken as the range space R of the FRM.

$R_1$ - Female child a burden so the sooner they get her married off, the better relief economically.

$R_2$ - Poverty / not owners of property

$R_3$ - Bad habits of the men / husbands

$R_4$ - Infected women are left uncared by relatives, even by husbands.

$R_5$ - No Guilt or fear of life

$R_6$ - Not changed religion and developed faith in god after the disease.

$R_7$ - No moral responsibility on the part of husbands and they infect their wives willfully

$R_8$ - Frequent, natural abortion / death of born infants

$R_9$ - STD / VD infected husbands

$R_{10}$ - Husbands hide their disease from their family so the wife become HIV/AIDS infected

The related directed graph given by an expert is given in Figure 1.

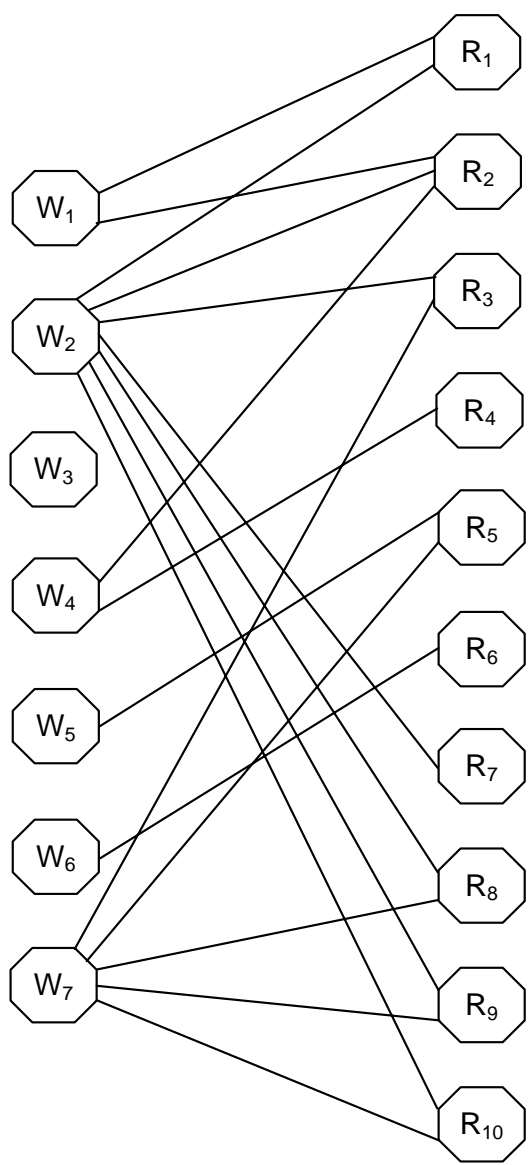

**FIGURE 1**

The associated connection matrix M

$$
\begin{array}{c}
\phantom{M} \\
M \quad = 
\end{array}
\begin{array}{c}
\phantom{W_1}
\end{array}
\begin{array}{c}
\quad R_1\ R_2\ R_3\ R_4\ R_5\ R_6\ R_7\ R_8\ R_9\ R_{10} \\[2pt]
\begin{array}{c}
W_1 \\ W_2 \\ W_3 \\ W_4 \\ W_5 \\ W_6 \\ W_7
\end{array}
\left[
\begin{array}{cccccccccc}
1 & 1 & 0 & 0 & 0 & 0 & 0 & 0 & 0 & 0 \\
1 & 1 & 1 & 0 & 0 & 0 & 1 & 1 & 1 & 1 \\
0 & 0 & 0 & 0 & 0 & 0 & 0 & 0 & 0 & 0 \\
0 & 1 & 0 & 1 & 0 & 0 & 0 & 0 & 0 & 0 \\
0 & 0 & 0 & 0 & 1 & 0 & 0 & 0 & 0 & 0 \\
0 & 0 & 0 & 0 & 0 & 1 & 0 & 0 & 0 & 0 \\
0 & 0 & 1 & 0 & 1 & 0 & 0 & 1 & 1 & 1
\end{array}
\right]
\end{array}
$$

Using the attributes associated with women as the rows of the matrix and the other factors which had made them HIV/AIDS patients as the range space; we using the directed graph we have obtained the connected relational 7 × 10 matrix M.

Now we study the effect of each state vector on the dynamical system M. Let $X = (0\ 1\ 0\ 0\ 0\ 0\ 0)$, a given state vector in which i.e., all states except $W_2$ i.e., causes of women becoming HIV/AIDS patient alone is in the on state and all other attributes are in the off state.

Consider the effect of X on the dynamical system M,,

$XM \quad \hookrightarrow \quad (1\ 1\ 1\ 0\ 0\ 1\ 1\ 1\ 1) \quad = \quad Y\ \text{(say)}$

('$\hookrightarrow$' denotes that the resultant vector has been thresholded and updated).

$YM^T \quad \hookrightarrow \quad (1\ 1\ 0\ 1\ 0\ 0\ 1) \quad = \quad X_1\ \text{(Say)}$

$X_1\ M \quad \hookrightarrow \quad (1\ 1\ 1\ 1\ 1\ 0\ 1\ 1\ 1\ 1) = \quad Y_1\ \text{(say)}$

$Y_1\ M^T \quad \hookrightarrow \quad (1\ 1\ 0\ 1\ 1\ 0\ 1) \quad = \quad X_2\ (= X_1)$

$X_2\ M \quad \hookrightarrow \quad (1\ 1\ 1\ 1\ 1\ 0\ 1\ 1\ 1\ 1) = \quad Y_2\ \text{(say)} = Y_1$

Thus the resultant binary pair is a fixed point given by {(1 1 0 1 1 0 1), (1 1 1 1 1 0 1 1 1 1)} which shows the reason for women becoming HIV/AIDS victims is due to child marriage, (male) widower marriage, child married to an old man, disease untreated till it is chronic or they are in the last stages is in the off state, women not being bread winners in spite of their earning are neglected is on, all other states remain unaffected.

The effect of women becoming HIV/AIDS affected has the following influence on the range space; the causes of women being becoming HIV/AIDS affected patients is due to: female child is burden so sooner they get married off the better the economical relief, poverty, bad habits of men.

Frequent natural abortions, STD/VD of their husbands and most importantly men (husbands) hide their disease from their family become on; $W_3$, $W_6$ and $R_6$ alone remain in the off state which shows religion or its relation to god have no significance to the cause of women becoming a victim of HIV/AIDS. Husbands hide their disease from the family so wife becomes easy victims of HIV/AIDS. We have worked with several vectors. All conclusions and analysis based on our study is given in the last chapter of this book.

Suppose we take the on state of the vector $R_2$ from the range space i.e., the effect of poverty or not owners of property is in the on state and all other vectors are in the off state the effect of Y on the dynamical system M, where Y = (0 1 0 0 0 0 0 0 0 0).

The effect of Y on the dynamical system $M^T$ is given by

$YM^T$     $\hookrightarrow$    (1 1 0 1 0 0 0)       =   X (say)

XM       $\hookrightarrow$    (1 1 1 1 0 0 1 1 1 1)    =   $Y_1$ (say)

$Y_1 M^T$    $\hookrightarrow$    (1 1 0 1 0 0 1)       =   $X_1$ (say)

$X_1 M$     $\hookrightarrow$    (1 1 1 1 1 0 1 1 1 1)    =   $Y_2$ say

$Y_2 M^T$    $\hookrightarrow$    (1 1 0 1 1 0 1)       =   $X_2 (= X_1)$

Thus the resultant vector is a fixed binary pair given by {(1 1 1 1 1 0 1 1 1 1), (1 1 0 1 1 0 1)} which is the same as the one given when the cause of women become HIV/AIDS victims alone is in the on state in the domain space. So from our study poverty and not owning of property is one of the root causes of women becoming HIV/AIDS victims.

We now analyze the second expert's opinion who is a HIV/AIDS patients over 6 years taking treatment in the hospital. The attributes related to women given by him are

$W_1$ - The cause for women becoming HIV/AIDS is marriage of young girls to drunkards, widowers

$W_2$ - Poverty among rural women

$W_3$ - Women have no power to protect themselves from HIV/AIDS infected husbands

$W_4$ - Uneducated rural women in villages do not have any voice or respect or any form of policymaking, regarding family or children.

$W_5$ - Being uneducated and giving birth to many children at a very young age even if they know the value of education or realize the importance of education, they are helpless (i.e. women's work in villages are worthy womb)

$W_6$ - Rural women neglect their health, due to supreme sacrifice for family

$W_7$ - No health centers or hospitals where they can go for medical aid or test; and counseling rural women are entirely dependent on their careless husbands.

The other attributes, which make women, HIV/AIDS infected is taken as the range space.

$R_1$ - HIV/AIDS infected husbands who don't disclose the disease to the family

$R_2$ - Male Ego (machismo / exaggerated masculinity) and bad habits / unconcern about wife

$R_3$ - Most men in the villages don't give any account for the money they earn or they spend

$R_4$ - No higher aim in life; being uneducated they are satisfied by their day-to-day living

$R_5$ - No proper counseling even after they are HIV/AIDS affected. Even if counseled they do not follow instructions/advice as in most cases they hide the news from their family. They do not even disclose to their wife about the disease

$R_6$ - No reformation even after they are infected by HIV/AIDS they continue the same style of living

$R_7$ - They are not bothered about their kith and kin or even their wife

$R_8$ - They (men) suffer from guilt after infection so they cannot be normal at home

$R_9$ - They (men) don't live with courage after being infected most of them (men) are depressed and attempt suicide due to poverty, fear and desperation and continue with all bad habits including visiting to CSWs

$R_{10}$ - They take pride in these bad habits for that is the way they show themselves superior to women and spend money solely on them

Using these attributes $W_1$, …, $W_7$ as the nodes of the domain space and that of attributes $R_1$, …, $R_{10}$ as the nodes of the range space we give the directed graph given by experts is in figure 2.

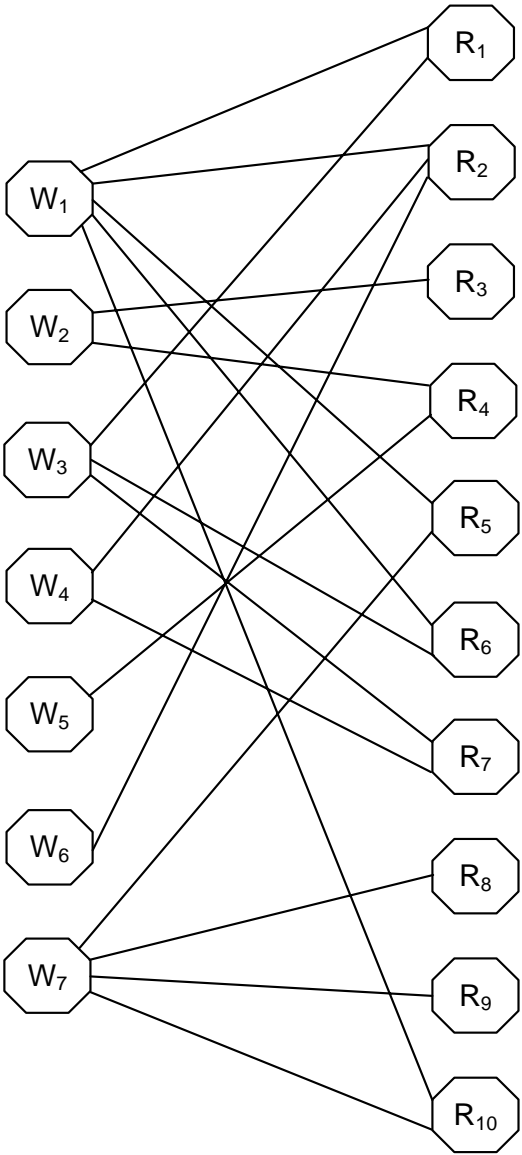

**FIGURE 2**

The connection matrix related with the directed graph is given by $M_1$.

$$M_1 = \begin{array}{c} \\ W_1 \\ W_2 \\ W_3 \\ W_4 \\ W_5 \\ W_6 \\ W_7 \end{array} \begin{array}{c} \begin{array}{cccccccccc} R_1 & R_2 & R_3 & R_4 & R_5 & R_6 & R_7 & R_8 & R_9 & R_{10} \end{array} \\ \left[ \begin{array}{cccccccccc} 1 & 1 & 0 & 0 & 1 & 1 & 0 & 0 & 0 & 1 \\ 0 & 0 & 1 & 1 & 0 & 0 & 0 & 0 & 0 & 0 \\ 1 & 0 & 0 & 0 & 0 & 1 & 1 & 0 & 0 & 0 \\ 0 & 1 & 0 & 0 & 0 & 0 & 1 & 0 & 0 & 0 \\ 0 & 0 & 0 & 1 & 0 & 0 & 0 & 0 & 0 & 0 \\ 0 & 1 & 0 & 0 & 0 & 0 & 0 & 0 & 0 & 0 \\ 0 & 0 & 0 & 0 & 1 & 0 & 0 & 1 & 1 & 1 \end{array} \right] \end{array}$$

$M_1$ is the dynamical system relating the vulnerability of women becoming HIV/AIDS patients and the role-played by men (their husbands) in infecting their wives. Let us consider the state vector $X = (0\ 0\ 1\ 0\ 0\ 0\ 0)$ i.e., rural women have no power to protect themselves from HIV/AIDS is alone in the on state and all other states are in the off state.

The effect of the state vector $X$ on the dynamical system $M_1$ is given by;

$$
\begin{array}{llll}
XM_1 & \hookrightarrow & (1\ 0\ 0\ 0\ 0\ 1\ 1\ 0\ 0\ 0) & = & Y \\
YM_1^T & \hookrightarrow & (1\ 0\ 1\ 1\ 0\ 0\ 0) & = & X_1 \\
X_1\,M_1 & \hookrightarrow & (1\ 1\ 0\ 0\ 1\ 1\ 1\ 0\ 0\ 1) & = & Y_1 \\
Y_1\,M_1^T & \hookrightarrow & (1\ 0\ 1\ 1\ 0\ 1\ 1) & = & X_2 \\
X_2\,M_1 & \hookrightarrow & (1\ 1\ 0\ 0\ 1\ 1\ 1\ 1\ 1\ 1) & = Y_2(=Y_1) \\
Y_2\,M_1^T & \hookrightarrow & (1\ 0\ 1\ 1\ 0\ 1\ 1) & = X_3 = X_2 \\
\end{array}
$$

Thus the resultant vector is a fixed point given by the binary pair $\{(1\ 0\ 1\ 1\ 0\ 1\ 1), (1\ 1\ 0\ 0\ 1\ 1\ 1\ 1\ 1\ 1)\}$.

The node/concept which is in the on state viz. rural women have no power to protect themselves from HIV/AIDS makes all the nodes both in the domain and range space to be in the on state. Thus when these women have no power to protect themselves from HIV/AIDS is in the on state we see that the node; men have all bad habits is in the on state. They lead a careless life even after they are infected, they live in guilt and do not share their problems with the family, in the same manner poverty also plays a vital role and rural women as inferior objects adds fuel to fire by making them easy victims of the disease.

Let us consider the state vector $Y = (0\ 0\ 0\ 0\ 0\ 0\ 1\ 0\ 0\ 0)$ i.e., only the on state of the node least bothered about their kith and kin and all other nodes are in the off state. The effect of Y on the dynamical system $M_1$ is given as follows:

$$
\begin{aligned}
YM^T_1 &\hookrightarrow (0\ 0\ 1\ 1\ 0\ 0\ 0) &&= X \\
XM_1 &\hookrightarrow (1\ 1\ 0\ 0\ 0\ 1\ 1\ 0\ 0\ 0) &&= Y_1 \\
Y_1 M^T_1 &\hookrightarrow (1\ 0\ 1\ 1\ 0\ 1\ 0) &&= X_1 \\
X_1 M_1 &\hookrightarrow (1\ 1\ 0\ 0\ 1\ 1\ 1\ 0\ 0\ 1) &&= Y_2 \\
Y_2 M^T_1 &\hookrightarrow (1\ 0\ 1\ 1\ 0\ 1\ 1) &&= X_2 \\
X_2 M_1 &\hookrightarrow (1\ 1\ 0\ 0\ 1\ 1\ 1\ 1\ 1\ 1) &&= Y_3 \\
Y_3 M^T_1 &\hookrightarrow (1\ 0\ 1\ 1\ 0\ 1\ 1) &&= X_3\ (= X_2)
\end{aligned}
$$

The resultant vector is fixed point given by the binary pair $\{(1\ 0\ 1\ 1\ 0\ 1\ 1), (1\ 1\ 0\ 0\ 1\ 1\ 1\ 1\ 1\ 1)\}$.

Thus this node $R_7$ is also a vital one for it turns on all other nodes in both the state vectors. We see poverty ($W_2$) has no role to play, being uneducated and giving birth to children ($W_5$) etc. has no impact on $W_3$. Likewise $R_3$ and $R_4$ remain in the off state.

We have chosen arbitrarily 10 HIV/AIDS women patients whose expert's opinion was sought on this issue and was converted into a FRM model. Using these 10 FRM models we combined them to get the combined FRM (CFRM) using which the opinion was sought and conclusions were derived. The conclusions are based on the combined opinion only.

Apart from these we had taken opinions from women social workers and general NGOs.

It was strongly felt that women (especially rural uneducated women) have no means to protect themselves from their HIV/AIDS affected husbands for reasons such as they are unaware of the fact their husbands are HIV/AIDS affected, secondly the social set up prevents them protesting against sex even if they suspect their husbands, thirdly women are not so much empowered either socially or economically to protect themselves because they fear that their husband would seek some other women for sex there by they would still become economically exploited.

Also linked FRMs are used to interrelated the nodes which are unrelated.

As the main motivation of the book is to reach non-mathematicians also, we have tried our level best to avoid cumbersome calculations and workings of the problem. In this book we have mainly put forth the conclusions based on our overall study.

Now we work with yet another new set of nine attributes given by this expert. The attributes given by him related to rural uneducated women are

$W_1$ - No basic education
$W_2$ - No wealth / property
$W_3$ - No empowerment
$W_4$ - No recognition of their labour
$W_5$ - No protection from their spouse
$W_6$ - No proper health care center to get proper guidance
$W_7$ - No nutritious food
$W_8$ - Women are not decision makers for women
$W_9$ - No voice basically in the traditional set up

When we arrive at these attributes or work for them as the major reasons we do not at the first instant take this to be true in case of all women. Here our study is only concerned with rural poor uneducated women, so the attributes selected and given by the expert are true for majority of the rural economically poor women.

Now we proceed on to give a brief description of each and every attribute for the main reason for doing so is what the expert precisely means by these attributes.

## W₁ - No basic education

Most of these women have not even completed their primary education i.e. even up to $5^{th}$ standard. Some have never entered school. The rural setup in India is so harsh on the female child because they are either killed at birth as female infanticide or female feticide or by chance they grow up they have to shoulder the burden of the family right from the age of three or so. Girls run errands or if they have a younger brother or sister the child is left under the care of this three year old, in certain cases this child accompanies her mother for work. The main reason which can be attributed for this is that female children in India is looked upon only as a burden and all the more in rural India, where poverty is the only way of life. So they feel female child is a curse upon them so they try to get rid of her. They feel that the food they give her, the money they spend on her is an utter waste so the female children are under-nourished.

The parents of these children do not even think of the psychological feelings of these children they are bred as cows and goats or even worse than that. The unfortunate thing is that they are married of at the age of 10, 12, 13 or so to men double or triple their age and who are invariably drunkards. The marriage, in majority of the cases, breaks down in a year or two. As usual, the husband either gets legally married to some other child or goes behind other women. Thus the torture, problems of stress and strain and above all family burden starts for these females at an early age of ten or twelve i.e., in their very teenage. One of the causes, which can be attributed for this, is that they are not educated. Such things unfortunately cannot take place as everything takes place at a very young age when they are just small children. In several cases their own relatives also sexually abuse them.

## W$_2$ – No wealth / property

In the Indian system of law, that too, after the advent of the Manusmriti (The Laws of Manu), the status of women has become very poor. This system treats women only like an object. She is just an object of the father when she is very young, then the working robot of her husband (includes child-bearing, taking care of family, her in-laws and so on), in the later stages she is ruled by her son. So in their lifetime most rural women do not really know what is freedom. So the Hindu law is so powerful that it is impossible for women to inherit any property that too in rural India. So when she is poor, has no wealth or property by which she can wield power, she is at the mercy of men. Thus women have become powerless because wealth or money alone is the only power uneducated women can possess.

## W$_3$ – No empowerment

The true empowerment of women in India, whether they are rich or poor or uneducated or educated is still a question mark. One should not be hurt or raise their eyebrows if we extend this and say that to some extent this is true all over the world. Women are always considered to be less than men in all walks of life. We do not want to deviate from this issue. We wish to talk about empowerment of women that too about rural, uneducated women in India. Though there are several women organizations yet they do not reach the grassroots level: rural women.

Rural women are naturally empowered with hard labour, natural intelligence, tolerance and above all a strong conditioned feeling that men are superior to them. With these qualities women are more exploited than being given their due. So it has become an easy advantage for men to overpower them, so how can one ever think of easily achieving rural women's empowerment? These women have been slaves to husbands, sons and sons-in-law so they have never known their right or freedom. They work all through their life not caring for themselves. Self-sacrifice is at its supreme in these cases. So they become easy victims of all types of harassment. The main and the most saddening thing is they have become HIV/AIDS patients only because of their husband. Our data shows that their husbands have infected these innocent women willfully, because though they know that they suffer from AIDS they have not confided about it to their wives.

These women start to take treatment only when the disease is chronic, so they are mercilessly thrown out of their homes at this stage with no one around to take care of them. Unless rural women are empowered and at least informed about the risk of this cruel disease it would not be possible to save them from HIV/AIDS. So, the government should not only take appropriate steps to save the uninfected but also steps must be taken to give special attention to these innocent HIV-women and make provision for them to enjoy peaceful last days taking into account that they had been toiling from the childhood till the day they became incapable.

### $W_4$ – No recognition of their labour

When we say "no recognition of their labour" we not only imply the rural women of India but all women; even educated and employed or even very well placed because women are not saved from the family burden even if they are educated and employed. They have to do two jobs at one time: maintaining the family and the other, their office work.

Now the worst thing is that in poor, rural India most of the women are employed as daily labourers (coolies). They cook for the family and then go to work where they labour for over eight hours just for a paltry pay. Their routine is very rigorous, hours of sleep very less yet the male members after a day's labour consume alcohol and invariably in many families they treat the women badly. Thus neither their labour nor the money they earn is ever recognized.

## $W_5$ – No protection from their spouse

This is an important issue in the study of HIV/ AIDS affected rural women and their socio-economic problems. It is the modesty of the married Indian women that they cannot say 'no' to sex with their husbands. This has a very great significance in our study that is why only their husbands have infected 91 of the 101 women interviewed. In most of the cases these husbands fully knew that they were suffering from HIV / AIDS. What can be the explanation given to such act? To be more precise these men statistically have not only infected their wives but also were regularly taking medical treatment for the disease and enjoying all sympathies from the family by saying they were suffering from ulcer or TB or so on. So, on the basis of sympathy they could get the best of food, treatment and attention from the family members. Further majority of these men had all bad habits like visiting CSWs, smoking, illicit relations with other women, alcohol and so on. In fact, they never spent the money they earned on their family but only on themselves. As these women were also employed as daily-wagers they took care of the family. Thus it was impossible to protect themselves from their infected husbands and it lead to their infection.

## $W_6$ -No proper health case center to get proper guidance

We have methodically described the way these rural women lead their life. They have no time either to think of themselves or to take care of themselves. Round the clock they are busy with work. The food they consume is the leftovers; that is the supreme sacrifice a rural woman makes. So she is susceptible to all forms of infection due poor nutrition and poor resistance. Also when she is sick, she does not inform anyone but on the contrary, she continues to carry on her regular routine. Only when she is totally affected and unable to move or do any job she seeks medical aid. By this time the disease becomes invariably chronic. The health-centers in rural India, in majority of the cases have no trained good doctors. They are mainly run by assistants who give some medicine to these women with no proper diagnoses of the real disease. For majority of these women complained that whatever treatment they got from the health center did not cure them or made them feel better. So these women suffer silently. Further by this time in most of the cases they become terminally ill and immobile. These women are then taken to big hospitals (by their mothers or sons, only in some cases by their husbands) for treatment where they are found to suffer from HIV/ AIDS. When they come to know that they suffer from HIV/AIDS, majority of them do not become upset or sad or depressed as men do.

The reasons attributed for this are

    (i)      They are not guilty as men
    (ii)     They are ignorant about the very disease. Their awareness about the disease is nil.
    (iii)    Finally they leave things to destiny.

Thus lack of proper guidance given from the health centers has made them either terminally ill, that is, the disease is chronic.

## $W_7$ - No nutritious food

When we say rural women do not get nutritious food the reader may feel uncomfortable. This title has so much to do with the tradition, culture and the male dominated society in relevance to women. For from the time a female child is born till she dies, she is in all angles taught only about the supreme sacrifice in life. All the India epic role models are Sita, Savithri, Nalayini, Shakuntala, Andal or Radha; so that the Indian women is never allowed to think beyond the concept of self sacrifice. So the idea of sacrifice is in her blood and heart, since her reasoning is never allowed to function.

A special study about women in India with relation with culture, tradition and epic should be made to understand the real setup. These women, even after knowing fully well that only their husbands have infected them do not speak ill or show their anger on them. This will now make the reader know on how they will ever feed on nutritious food.

Just for the sake of more clarity we wish to state how dinner or lunch or breakfast in traditional rural poor family is served. In most of these families the head of the family is first served food: her husband and her father-in-law (i.e. her husband's father) are served food. These men, majority of them are so selfish never even think whether their children or wife or mother has any leftover, they eat to their full.

Next the women folk serve their children and the mother-in-law (mother of her husband), sister-in-law (sister of her husband) after they eat if there is some left over she takes. Mostly her food would be the previous day rice put in water with a chilli. Even this she may not get at times. Those times she would chew tobacco to subside her hunger and go for work. This is the pattern of majority of the families in poor rural India.

## $W_8$ - Not decision Makers

Having said so much about our women, we also need to add that they are not the decision makers. They earn, work and in the meantime they are machines to produce children. They invariably give birth to two or four children. If they do not get a male child they continue to bear children till they get a male child or in some cases the man marries another woman. This is very common in rural India. Thus even to opt for a child is not in the hands of a woman. Likewise, the way in which rural female children are married is very disgusting. In majority of the cases just female children at the age of 12, 14 or so is married to men who are in their thirties or early forties, in many cases knowing fully well that the bridegroom is a drunkard, unemployed or an embodiment of evil. In many cases the children are married to these men, as second wives because the first wives are dead or separated for reasons best known to them.

In our opinion some marriages are very cruelly arranged which ultimately lands them in the Tambaram Sanatorium hospital for HIV/AIDS treatment. The child cannot say no for sex for the fear she would be beaten up or killed or thrashed so whether she likes or not, willing or not should get married at the order of their parent and have sex at the order of the husband. How can poor, rural women ever be decision-makers when they even cannot decide for themselves?

## $W_9$ – No voice in the basic traditional (patriarchal) set up

Already we have stated that they are not decision-makers even for themselves. Now we proceed on to show that the major reason for this is that the traditional patriarchal setup in rural India is so bad. It is a very common occurrence in India that men beat their wife, irrespective of his or her educational status. Further there is no punishment for it, for women fearing that the name/fame of the family would be spoilt do not even speak it out.

The cruelest thing is that after the day's work the rural men return home, in majority of the cases only in the drunken state. It is still more surprising to see that even in their drunken state they only beat their wives. They claim that drunkenness is the reason for beating the wife, but it is very well schemed because their drunken state doesn't make them beat their mother or sister or children, only the wife becomes a victim (This topic needs a very special psychological study). So one can imagine what voice she has in her family. A silent breadwinner, who toils for the family, cooks for the family and takes care of the family, is treated in such a fashion. So, in several cases they fear to voice their feelings. In India it is not a rare occurrence that the wives are easily murdered under the pretext of fake fire accidents (Bursting of stoves or gas cylinders). Thus women have no voice; not even a voice to protect herself from the HIV/AIDS infected husband even if she knows fully well that her husband is infected with HIV/AIDS. This is the voice of the voiceless. Now we proceed on to give the attributes related with men, society and economic status.

The attributes given by the expert making the rural women passive victims of HIV/AIDS

$R_1$ - Male dominated society
$R_2$ - Traditional and cultural values of the nation
$R_3$ - Poverty among villagers
$R_4$ - No proper counseling by doctors to HIV/AIDS infected men
$R_5$ - Media, TV cinema (induces male ego)
$R_6$ - Model Women
$R_7$ - No proper spiritual leaders to channelize the youth
$R_8$ - Men have no responsibility or binding with family (careless living)
$R_9$ - Women in Indian society
$R_{10}$ - Role of women association in the context of these HIV / AIDS affected rural poor women

**$R_1$ - Male dominated Society**

Every nation has patriarchy in a power-wielding manner to some degree or other which cannot be disputed by anyone.

But India is a highly male-dominated society mainly after the incoming of the Aryans, for the sacred book of the East says women are just objects, they are and can be only in possession of men. The poor rural men have imbibed this Aryan culture because it gave lots of food to their ego and above all an unquestioned authority over women. So they took this evil alone and left other factors like education etc. Thus men became more and more powerful not as an independent entity, but only relative to their women. Thus the male dominated society has become the cause of spread of HIV/AIDS in the rural poor women.

The other side of the problem is that these men when they are fully aware of the fact that they are infected by HIV/AIDS not only infect their wife by having unprotected sex but with all the CSWs whom they visit for sex they have only unprotected sex there by infecting the CSWs and indirectly those other men who visit the same CSWs. The migrant labourers openly said this that in spite of the request of the CSWs to have protected sex they did not heed to them. Unless this attitude of domination is psychologically corrected there cannot be any solution to this problem. From childhood, the parents, in-laws and so on, cherish "male ego". Thus unless a change takes place and realization comes among rural men, it is impossible to save women.

## $R_2$ – Traditional and cultural values of the nation

When we talk about cultural and traditional value of the nation in the context of the spread of HIV/AIDS among rural women, we see how these values have negatively worked on our women and has ultimately made them as HIV/AIDS infected.

The newborn Aryan culture that has come to India has several wrong cultural and traditional values. The worst among them is acceptance and treatment of women as an inferior being (treated as object), which is sowed as the seed of several evils like, child marriage, dowry system, denial of education to female child, denial of property right to women, right for the women to make decision for her life and so on.

Our culture that is based on enslaving women and treating them as objects has now become the cause of spread of HIV/AIDS among women. Here is an instance that highlights the point we are making:

In *The New Indian Express* dated 26 April 2005, the first page report carries a photograph of a 12-year-old girl sporting a *mangalsutra* (a symbol of marriage, like the western wedding ring, a sacred thread tied around the neck of the woman at marriage). The photo's caption says that it may be an uncommon sight these days, but not at Varamanaguda village in Krishnagiri in Tamil Nadu. Here child marriage is rampant. While the legal age of marriage for a girl is 18, at Varamanagunda it is 12 years or a month after puberty, whichever comes first. The wretched thing about this is that a 30-year-old tying the knot with a 10-year-old girl is termed routine. The villagers say that it is deeply engrossed in their tradition.

The Gill hamlet in the Tamil Nadu-Andhra Pradesh border has a population of 700 plus most of them belonging to the Valmeeki Nadu community, a Telugu speaking backward caste. Besides Varamanagunda the practice is also rampant in nearby hamlets including Ganamoor, Periyar Nagar, Theertham, Sigaralagiri, Bannapalli and Nachikuppam. And the marriages do not last long, with a majority of the couples separating within a couple of years.

Nagamma (15) of Periyar Nagar was married off at the age of 13. Her husband deceived her after a year to marry again. Radhika (13) and Darini (13) of Periyar Nagar village were also married in a similar fashion. Thimmiah of Varamalagunda who gave his 12-year-old daughter Suguna in marriage this year has no qualms over his action. Averring that he was only abiding by tradition, he does not feel the need to abide by the laws framed by a government that has done nothing to community or village. Thus in India tradition and culture stands above law, that too even above the government rules and laws that is why it is not wrong to say the culture and tradition of our nation is the cause of raise in HIV/AIDS patients among rural uneducated women.

Last month, Varamanagunda villagers organized a mass child marriage programme involving as many as 18 girls between 12 to 14. The village administrative office and other officials were aware of the function but did not report the matter to the government.

Murugan, a social worker in Bargur town says admission of under age pregnant women in Bargar Government hospital is on the rise. While the doctors suppress the information, local politicians are not interested in ending the practice of child marriage fearing that any intervention would cut into their votes.

When the health departments conducted a medical camp here recently they found that 50% of the people were suffering from tuberculosis. Child marriage and inadequate medical facilities were cited as the reason.

Thus we see how culture and tradition has a major role to play in these marriages.

### R₃ – Poverty among villagers

One of the major causes of rapid spread of HIV/AIDS is poverty, for with lack of nutritious good food, the self-sacrificing rural women become not only easy victims of HIV/AIDS but also the time span in which they become seriously infected is short in several cases within six months, in some cases a year or two or three and so on. If they are not very poor or poor, at least they can have good nutritious food and be in a position to fight the disease.

Due to lack of nutrition their immune system is very poor or susceptible to all types of opportunistic diseases and these women become easy victims of full-blown HIV/AIDS. As they are poor they cannot afford any medical treatment when the disease just sets on them, only when they become chronic terminally ill patients some medical aid is given to them. Several of them whom we interviewed suffered from TB and skin problems. Further poverty is one of the causes for them to marry their children early that too in many a cases to men who are twice or thrice their age and at times to divorcees or widowers.

Invariably these men are drunkards and have all bad habits including visiting CSWs. So in some cases within a year of marriage these men die. In some cases their wives become HIV/AIDS infected patients within three to four years. As the African countries which is a poor nation became victims of HIV/AIDS, likewise the poor rural people of India have all chances and in majority of the cases to become easy victims of HIV/AIDS. Thus poverty has to do lot with becoming AIDS patients in a very short span of time.

## $R_4$ – No proper counseling by doctors to HIV/AIDS affected men

This node/concept was specially taken into account because majority of the women whom we interviewed were ignorant about the disease HIV/AIDS, as to how it spreads, etc. Secondly, three of the women said that their husbands committed suicide when the doctors informed they had HIV/AIDS. That is why proper counseling is very important for these women patients to save themselves as well as to be very careful with their spouse once they know their spouse is affected/ infected by HIV /AIDS. Now as majority (91 out of 101) of these women were infected only by their husband, if the doctor would have counseled the husband about the seriousness of the disease and also have forced him to bring his wife and tell her about the modes how this disease is spread certainly some of these women could have been saved. Also the doctors should have made it mandatory to check the blood samples of their wife once the husband comes for HIV/AIDS treatment. By doing so the doctors could have prevented over 66 of them who were terminally ill.

When we say terminally ill in the sanatorium we mean that they are feasted on by lice in their uncared hair and flies eat the sores on their body. Thus these women lead a life worse than animals and no one is by their side. Further they are so ill they are immobile. This is what we mean counseling by the doctors who treat them, if they had counseled the members of their family certainly they will not leave the women like an orphan in the hospital when she has toiled for them for so many years!

## R$_5$ – Media T.V. Cinema (induces male ego)

The media is the tele-serials, skits in T.V. and above all the movies have a very great impact. In almost all these we see a woman with some individuality is always put down and made to surrender to the male ego thus invariably respect for women is lost from the childhood itself.

This is very much seen in villages, especially the male child even from a very young age is beating his mother, ill-treating her in spite of her giving much more importance to the male child than to her female child. So what is ingrained in their mind is the concept that "women are lower than men". Further this ideology is well protected and projected in all the movies, for instance always a women who is original, intelligent, beautiful with a lot of wit will be subdued and made into nothing. So unless the media changes its basic notion and ideology, the male domination would only be cherished in the Indian soil resulting in permanent disfigurement.

## R$_6$ – Model Women

The expert feels that the role models for women in India are either the women who are actresses or the epic women like Sita, Savithri, Andal, Nalayini, Radha, Gandhari and so on. These women were absolute slaves to men they supported their husbands injustice and suffered silently, that is why they were called as role models. For instance, when Sita's husband Rama suspected her chastity she showed that she was pure by walking on fire; but Rama, because a washerman talked scantly of him asked her to leave the palace and live in forest when she was pregnant. Since Sita did not question Rama's purity, she was a role model. Further the Aryan culture and tradition says that women should not mind suffering and that only suffering can make them pure. They should try to suffer and endure it! Is there any logic in this? Thus our role models from epics are so meek. Likewise Nalayini carried her sick husband in a basket to prostitute's house. That was this sacrifice, which made her a role model in our epic. Thus women are asked to encourage their husbands to visit CSWs. The majority of actresses who are taken as easy role models are women who lead flimsy lives and several of them commit suicide

## R₇ - No proper spiritual leader to change the youth

It is a disgusting truth to say that majority of our spiritual leaders are men with no character involved in sex scandals, murder charges, amassing wealth and so on. When these men are our spiritual leaders it is very natural to see our youth following them. So college students murder their lovers and do all types of immoral acts to earn money and use it over bad activities. Thus basically one is very much forced to question the spiritual leaders than the youth, for these spiritual leaders only give a negative impact on them. For they easily quote them in their wrong deeds. As the present day youth is more intelligent certainly their logic happens to be correct and the real problem is with us are the failed moralities of our spiritual leaders.

## R₈ – Men have no responsibility or binding with family so they lead a careless life

This has so much weight for if the rural men whose wives were infected by them had been a little careful in the first instance they would not have become HIV/AIDS patients. Carelessness coupled with irresponsibility had made them HIV/AIDS patients. Only their carelessness and irresponsibility had made their wives also as HIV / AIDS patients. If he had the true responsibility towards his wife or his children would he ever infect his wife after fully knowing well that he is suffering with the disease? This shows that majority of them are not only careless and irresponsible, they more become kind of sadists once this disease affects them.

Here it has become pertinent for us to state that the psychological study of these HIV/AIDS infected men should be made to know why perversion coupled with sadism dominates their character—possibly due to the depression.

For when we interviewed these men, majority of them were just migrant labourers, we saw that these men first took pride in visiting several CSWs.

In the second stage of analysis we come to know that even after fully knowing that they are suffering from HIV/AIDS they have the heart to have unprotected sex with poor CSWs. Not only they had unprotected sex with CSWs they have the mind and heart to have unprotected sex with their wives after being fully aware of their infection.

We are forced to study this when we interviewed the women affected with HIV /AIDS in the T.B sanatorium hospital we found out, from the ward boys and nurses and other women who were assisting the sick women and the women patients themselves, that the male patients who stay in the hospital as inpatients leave the hospital premises in the night to seek sex. Thus one can understand why the study of this node is important. For even after becoming inpatients at a AIDS recovery centre, after suffering and seeing other women and men suffer they go for sex and there by infect several others. This is a prime example of carelessness.

### R₉ – Women in Indian Society

This node was found to be important for our study. Now what is the role played by women in the Indian Society? It is not like any other nation because in this study women mean not only their gender problem, it is dependent so much on caste and also on their economic condition.

To be very precise in India as per the Aryans there are five castes: Brahmins, Kshatriyas, Vaishyas, Sudras and Chandalas or untouchables. The Chandalas or the untouchables would be termed in this book as dalits. The pity is in India even if a dalit happens to be the father of the constitution still he is not recognized as an intellectual but more recognized by his caste says this expert. So the women who have some role to play in India are only from the Brahmin castes and the upper castes.

Dalit women are easy victims of all atrocities by high caste men. So what role can they play in this society? More so is the status of Sudra women. So the role of Indian women or representation of women society in India is just the high caste socially and economically very well off women who just form some less than 5% of India's population.

That is why when we wanted to study HIV/AIDS in relation to women we purposefully choose poor uneducated rural women. For they are the worst neglected in India. To put it properly they are more neglected than animals in India. For there are laws to protect animals! Another point to be kept in mind is that even if rich and high caste women are affected with HIV/AIDS, it would not come in open and it would be only kept as a secret and they take all treatment or best of treatments from very good doctors. So it is only the poor women who are completely neglected by the society come for treatment to these government hospitals. To be more precise these women have no voice in society or even any representation in society, they live and die like plants that have no seeds!

**R$_{10}$ – Role of Women Associations in the context of HIV/AIDS affected rural poor woman in India**

As seen from node R$_9$ when we discussed about women in Indian society we made several points clearly. This node is an essential and an appropriate extension of the node R$_9$. This is the opinion of the expert which he wanted to put on record. Rich and high caste women run most of the women organizations/ associations in India as a past time. Majority of the members are from very rich and high social status. So to make them think of the problems of poor rural uneducated women in the first place is not possible. "This is like making blind persons know the form of an elephant!" he says. At the first instance these elite women do not know anything about the problems or the mode of life led by these rural women.

When this itself is not known how will they ever know the plight of HIV / AIDS infected poor, uneducated rural women! Has any women organization ever visited these HIV/AIDS affected rural women at least in the T.B sanatorium hospital? (We had visited them and after looking at the plight of these women we spent several sleepless nights for as a person we have no power to do anything except give them some money but their close male relatives will take away even the money we give). We had given some suggestions to the government. Women organization or associations in India and Tamil Nadu in particular seek for a ban on some films or some obscene posters simply because such activities would guarantee them some cheap popularity. They are not actively involved in fighting against AIDS, or spreading awareness. Instead such a major, women-concerned theme in left to the HIV/AIDS-infected women victims to form associations of their own and lend support and service.

Now using the experts opinion we give the FRM directed graph in the following.

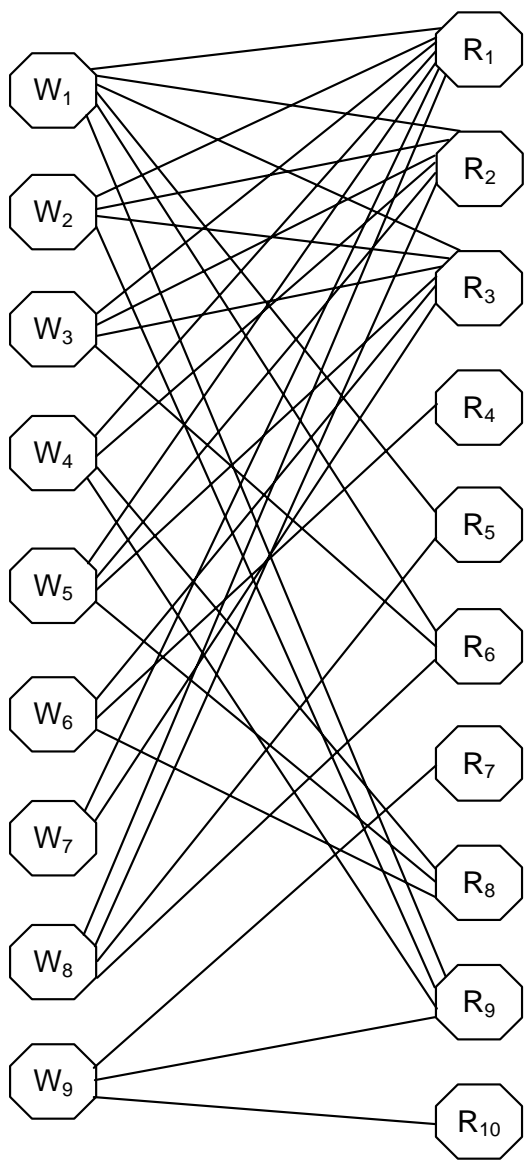

**FIGURE 3**

Using this directed graph we obtain the related connection matrix T which is a $9 \times 10$ matrix.

$$
\begin{array}{c}
\quad\quad R_1\ R_2\ R_3\ R_4\ R_5\ R_6\ R_7\ R_8\ R_9\ R_{10} \\
\begin{array}{c}
W_1 \\ W_2 \\ W_3 \\ W_4 \\ W_5 \\ W_6 \\ W_7 \\ W_8 \\ W_9
\end{array}
\left[
\begin{array}{cccccccccc}
1 & 1 & 1 & 0 & 1 & 1 & 0 & 0 & 1 & 0 \\
1 & 1 & 1 & 0 & 0 & 0 & 0 & 0 & 1 & 0 \\
1 & 1 & 1 & 0 & 0 & 1 & 0 & 0 & 0 & 0 \\
1 & 1 & 0 & 0 & 0 & 0 & 0 & 1 & 1 & 0 \\
1 & 1 & 1 & 0 & 0 & 0 & 0 & 1 & 0 & 0 \\
0 & 0 & 1 & 1 & 0 & 0 & 0 & 1 & 0 & 0 \\
1 & 0 & 1 & 0 & 0 & 0 & 0 & 0 & 0 & 0 \\
1 & 1 & 0 & 0 & 1 & 1 & 0 & 0 & 0 & 0 \\
0 & 0 & 0 & 0 & 0 & 0 & 1 & 0 & 1 & 1
\end{array}
\right]
\end{array}
$$

Now T is the dynamical system of our model. We can study / analyze the effect of any state vector either from the domain space D or from the range space R.

Suppose one is interested to study the effect of $W_9$ on the dynamical system i.e. the rural poor uneducated women have no voice in the traditional set up alone is in the on state and all other nodes are in the off state to study the effect of this state vector.

Let $X_1 = (0\ 0\ 0\ 0\ 0\ 0\ 0\ 0\ 1)$, i.e., $W_9$ alone is in the on state and all other nodes are in the off state. Its impact on the dynamical system T.

$$
\begin{aligned}
X_1T \quad &= \quad (0\ 0\ 0\ 0\ 0\ 0\ 1\ 0\ 1\ 1) \\[6pt]
Y_1T^{\,t} \quad &= \quad Y_1 \\
&\rightarrow \quad (1\ 1\ 0\ 1\ 0\ 0\ 0\ 0\ 1) \\
&= \quad X_2 \text{ (say)} \\[6pt]
X_2T \quad &= \quad (3\ 3\ 2\ 0\ 1\ 1\ 1\ 1\ 4\ 1)
\end{aligned}
$$

after thresholding we get the resultant vector

$$Y_2 = (\ 1\ 1\ 1\ 0\ 1\ 1\ 1\ 1\ 1\ 1\ ).$$

The effect of $Y_2$ on the dynamical system T is given by

$$Y_2 T^t \qquad = \qquad (6\ 4\ 4\ 4\ 4\ 2\ 2\ 4\ 3)$$

after thresholding and updating we get

$$X_3 = (1\ 1\ 1\ 1\ 1\ 1\ 1\ 1\ 1).$$

Now we study the effect of $X_3$ on the dynamical system T.

$$X_3 T = (7\ 6\ 6\ 1\ 2\ 3\ 1\ 3\ 4\ 1).$$

Now thresholding and updating the resultant vector

$$Y_3 = (1\ 1\ 1\ 1\ 1\ 1\ 1\ 1\ 1\ 1).$$

Thus the fixed point of the dynamical system is a pair {(1 1 1 1 1 1 1 1), (1 1 1 1 1 1 1 1 1 1)}.

Now we see when a single node "No voice for women basically in the traditional set up" alone is in the on state all the other nodes both in the domain and range space came to an state which shows that for the spread of HIV/AIDS among rural uneducated women in India is basically due to women voicelessness in the society.

They as a part of nations tradition cannot say no to sex to their HIV/AIDS infected husbands, for male denominated society demands it. So unless the tradition which is basically ruling the Indian society is broken there can be no solution to this problem. Poor rural uneducated women would continue to be the worst victims of HIV/AIDS. Thus the spread of HIV/AIDS in the nation above all in the rural India cannot be controlled.

Now we proceed on to study the effect of a state vector from the range space of the dynamical system T.

Let us consider the on state of the single state vector $R_5$ alone on the dynamical system i.e., "the role of media in inducing male ego" and all other nodes are assumed to be in the off state.

Let $Y = (0\ 0\ 0\ 0\ 1\ 0\ 0\ 0\ 0\ 0)$. Now we study the effect of Y on the dynamical system $T^t$.

$$YT^t \qquad = \qquad (1\ 0\ 0\ 0\ 0\ 0\ 0\ 1\ 0) = X.$$

The effect of X on the dynamical system T is given by

$$XT \qquad = \qquad (2\ 2\ 1\ 0\ 2\ 2\ 0\ 0\ 1\ 0\ ).$$

After thresholding and updating we get the resultant as

$$Y_1 = (\ 1\ 1\ 1\ 0\ 1\ 1\ 0\ 0\ 1\ 0)$$

Now we analyze the effect of $Y_1$ on $T^t$

$$Y_1\ T^t = (\ 6\ 4\ 4\ 3\ 3\ 1\ 2\ 3\ 1).$$

After thresholding and updating we get the resultant as

$$X_1 = (1\ 1\ 1\ 1\ 1\ 1\ 1\ 1\ 1).$$

Now the effect of $X_1$ on T is given by

$$X_1T = (7\ 6\ 7\ 1\ 2\ 3\ 1\ 3\ 4\ 1).$$

After updating and thresholding the resultant vector we get

$$Y_2 = (1\ 1\ 1\ 1\ 1\ 1\ 1\ 1\ 1\ 1\ ).$$

Thus we see the hidden pattern of the dynamical system is a fixed point given by the binary pair.

$$\{(1\ 1\ 1\ 1\ 1\ 1\ 1\ 1\ 1), (1\ 1\ 1\ 1\ 1\ 1\ 1\ 1\ 1\ 1)\}.$$

Thus when only the single fact that the role of media in inbreeding the male ego alone is in the on state, all the nodes both in the domain space and in the range space become on where by it is clearly seen that India is emotionally and culturally ruled by the cinema world. That is why the present youth, have mainly have their role models as actor and actresses, fully for getting themselves and the fact that it is only on the screen and not in real life. Unless politicians take up the issue for the women cause and see to it that movies which put down potent and individual women with a mission should be banned nothing concrete can be achieved. People should come forth to take up movies which give true importance to the role of women and their power. Unless such revolutionary steps are taken it may not be possible to put an end to spread of HIV/AIDS in India that too among poor rural in India.

Not only today's movie inbuilt in males their ego and superiority over women but also invariable in most of all the movies heroes are portrayed with cigar in their hand / or in a pose of drinking alcohol so these bad habits in their real life. Thus the cinema seems to be one of the major cause which inbreeds male ego.

**Chapter Three**

# USE OF FAM TO ANALYZE THE FEELINGS OF HIV/AIDS WOMEN PATIENTS

This chapter has two sections. In section one we introduce FAM and in section two FAM is used in this model.

## 3.1 Introduction to Fuzzy Associative Memories

A fuzzy set is a map $\mu : X \rightarrow [0, 1]$ where X is any set called the domain and [0, 1] the range i.e., $\mu$ is thought of as a membership function i.e., to every element $x \in X$ $\mu$ assigns membership value in the interval [0, 1]. But very few try to visualize the geometry of fuzzy sets. It is not only of interest but is meaningful to see the geometry of fuzzy sets when we discuss fuzziness. Till date researchers over looked such visualization [Kosko, 59-61], instead they have interpreted fuzzy sets as generalized indicator or membership functions mappings $\mu$ from domain X to range [0, 1]. But functions are hard to visualize. Fuzzy theorist often picture membership functions as two-dimensional graphs with the domain X represented as a one-dimensional axis.

The geometry of fuzzy sets involves both domain $X = (x_1,\ldots, x_n)$ and the range [0, 1] of mappings $\mu : X \rightarrow [0, 1]$. The geometry of fuzzy sets aids us when we describe fuzziness, define fuzzy concepts and prove fuzzy theorems. Visualizing this geometry may by itself provide the most powerful argument for fuzziness.

An odd question reveals the geometry of fuzzy sets. What does the fuzzy power set $F(2^X)$, the set of all fuzzy subsets of X, look like? It looks like a cube, What does a fuzzy set look like? A fuzzy subsets equals the unit hyper cube $I^n = [0, 1]^n$. The fuzzy set is a point in the cube $I^n$. Vertices of the cube $I^n$ define a non-fuzzy set. Now with in the unit hyper cube $I^n = [0, 1]^n$ we are interested in a distance between points, which led to measures of size and fuzziness of a fuzzy set and more fundamentally to a measure. Thus within cube theory directly extends to the continuous case when the space X is a subset of $R^n$. The next step is to consider mappings between fuzzy cubes. This level of abstraction provides a surprising and fruitful alternative to the prepositional and predicate calculus reasoning techniques used in artificial intelligence (AI) expert systems. It allows us to reason with sets instead of propositions. The fuzzy set framework is numerical and multidimensional. The AI framework is symbolic and is one dimensional with usually only bivalent expert rules or propositions allowed. Both frameworks can encode structured knowledge in linguistic form. But the fuzzy approach translates the structured knowledge into a flexible numerical framework and processes it in a manner that resembles neural network processing. The numerical framework also allows us to adaptively infer and modify fuzzy systems perhaps with neural or statistical techniques directly from problem domain sample data.

Between cube theory is fuzzy-systems theory. A fuzzy set defines a point in a cube. A fuzzy system defines a mapping between cubes. A fuzzy system S maps fuzzy sets to fuzzy sets. Thus a fuzzy system S is a transformation S: $I^n \to I^p$. The n-dimensional unit hyper cube $I^n$ houses all the fuzzy subsets of the domain space or input universe of discourse $X = \{x_1, \ldots, x_n\}$. $I^p$ houses all the fuzzy subsets of the range space or output universe of discourse, $Y = \{y_1, \ldots, y_p\}$. X and Y can also denote subsets of $R^n$ and $R^p$. Then the fuzzy power sets $F(2^X)$ and $F(2^Y)$ replace $I^n$ and $I^p$.

In general a fuzzy system S maps families of fuzzy sets to families of fuzzy sets thus S:

$I^{n_1} \times \ldots \times I^{n_r} \to I^{p_1} \times \ldots \times I^{p_s}$ Here too we can extend the definition of a fuzzy system to allow arbitrary products or arbitrary mathematical spaces to serve as the domain or range spaces of the fuzzy sets. We shall focus on fuzzy systems S: $I^n \to I^P$ that map balls of fuzzy sets in $I^n$ to balls of fuzzy set in $I^p$. These continuous fuzzy systems behave as associative memories. The map close inputs to close outputs. We shall refer to them as Fuzzy Associative Maps or FAMs.

The simplest FAM encodes the FAM rule or association $(A_i, B_i)$, which associates the p-dimensional fuzzy set $B_i$ with the n-dimensional fuzzy set $A_i$. These minimal FAMs essentially map one ball in $I^n$ to one ball in $I^p$. They are comparable to simple neural networks. But we need not adaptively train the minimal FAMs. As discussed below, we can directly encode structured knowledge of the form, "If traffic is heavy in this direction then keep the stop light green longer" is a Hebbian-style FAM correlation matrix. In practice we sidestep this large numerical matrix with a virtual representation scheme. In the place of the matrix the user encodes the fuzzy set association (Heavy, longer) as a single linguistic entry in a FAM bank linguistic matrix. In general a FAM system F: $I^n \to I^b$ encodes the processes in parallel a FAM bank of m FAM rules $(A_1, B_1)$, …, $(A_m, B_m)$. Each input A to the FAM system activates each stored FAM rule to different degree. The minimal FAM that stores $(A_i, B_i)$ maps input A to $B_i'$ a partly activated version of $B_i$. The more A resembles $A_i$, the more $B_i'$ resembles $B_i$. The corresponding output fuzzy set B combines these partially activated fuzzy sets $B_1^1, B_2^1, \ldots, B_m^1$. B equals a weighted average of the partially activated sets B = $w_1 B_1^1 + \ldots + w_n B_m^{1'}$ where $w_i$ reflects the credibility frequency or strength of fuzzy association $(A_i, B_i)$. In practice we usually defuzzify the output waveform B to a single numerical value $y_j$ in Y by computing the fuzzy centroid of B with respect to the output universe of discourse Y.

More generally a FAM system encodes a bank of compound FAM rules that associate multiple output or consequent fuzzy sets $B_i$, ..., $B_i^s$ with multiple input or antecedent fuzzy sets $A_i^1$, ..., $A_i^r$. We can treat compound FAM rules as compound linguistic conditionals. This allows us to naturally and in many cases easily to obtain structural knowledge. We combine antecedent and consequent sets with logical conjunction, disjunction or negation. For instance, we could interpret the compound association $(A^1, A^2, B)$, linguistically as the compound conditional "IF $X^1$ is $A^1$ AND $X^2$ is $A^2$, THEN Y is B" if the comma is the fuzzy association $(A^1, A^2, B)$ denotes conjunction instead of say disjunction.

We specify in advance the numerical universe of discourse for fuzzy variables $X^1$, $X^2$ and Y. For each universe of discourse or fuzzy variable X, we specify an appropriate library of fuzzy set values $A_1^r$, ..., $A_k^2$ Contiguous fuzzy sets in a library overlap. In principle a neural network can estimate these libraries of fuzzy sets. In practice this is usually unnecessary. The library sets represent a weighted though overlapping quantization of the input space X. They represent the fuzzy set values assumed by a fuzzy variable. A different library of fuzzy sets similarly quantizes the output space Y. Once we define the library of fuzzy sets we construct the FAM by choosing appropriate combinations of input and output fuzzy sets Adaptive techniques can make, assist or modify these choices.

An Adaptive FAM (AFAM) is a time varying FAM system. System parameters gradually change as the FAM system samples and processes data. Here we discuss how natural network algorithms can adaptively infer FAM rules from training data. In principle, learning can modify other FAM system components, such as the libraries of fuzzy sets or the FAM-rule weights $w_i$.

In the following section we propose and illustrate an unsupervised adaptive clustering scheme based on competitive learning to blindly generate and refine the bank of FAM rules. In some cases we can use supervised learning techniques if we have additional information to accurately generate error estimates. Thus Fuzzy Associative Memories (FAMs) are transformation. FAMs map fuzzy sets to fuzzy sets. They map unit cubes to unit cubes. In simplest case the FAM system consists of a single association. In general the FAM system consists of a bank of different FAM association. Each association corresponds to a different numerical FAM matrix or a different entry in a linguistic FAM-bank matrix. We do not combine these matrices as we combine or superimpose neural-network associative memory matrices. We store the matrices and access them in parallel. We begin with single association FAMs. We proceed on to adopt this model to the problem.

## 3.2 Use of FAM model to study the socio economic problem of women with HIV/AIDS

At the outset we first wish to state that these women are mainly from the rural areas, they are economically poor and uneducated and majority of them are infected by their husbands. Second we use FAM because only this will indicate the gradations of the causes, which is the major cause for women being affected by HIV/AIDS followed in order of gradation the causes. Further among the fuzzy tools FAM model alone can give such gradations so we use them in this analysis. Another reason for using FAM is they can be used with the same FRM that is using the attributes of the FRM, FAMs can also be formulated. Already FAMs are very briefly described in section one of this chapter. We now give the sketch of the analysis of this problem with a view that any reader with a high school education will be in a position to follow it. But the conclusions given in chapter four are based on the complete analysis of the problem.

We just illustrate two experts opinion though we have used several experts' opinion for this analysis. Using the problems of women affected with HIV/AIDS along the rows and the causes of it along the column we obtain the related fuzzy vector matrix M. The gradations are given in the form of the fuzzy vector matrix M which is as follows.

$$
\begin{array}{c}
 & R_1 \ \ R_2 \ \ R_3 \ \ \ \ R_4 \ \ \ R_5 \ \ \ \ \ R_6 \ \ R_7 \ \ R_8 \ R_9 \ R_{10} \\
\begin{array}{c} W_1 \\ W_2 \\ W_3 \\ W_4 \\ W_5 \\ W_6 \\ W_7 \end{array}
\begin{bmatrix}
0.9 & 0.8 & 0.7 & 0 & 0 & 0 & 0 & 0 & 0 & 0.7 \\
0.5 & 0.8 & 0.6 & 0 & 0 & 0 & 0 & 0 & 0 & 0 \\
0 & 0.3 & 0.6 & 0 & 0 & 0 & 0 & 0 & 0 & 0 \\
0 & 0 & 0 & 0.6 & 0 & 0 & 0 & 0 & 0 & 0 \\
0 & 0 & 0 & 0 & 0.9 & .6 & .7 & 0 & 0 & 0 \\
0 & 0 & 0 & 0 & 0 & 0.7 & 0.5 & 0 & 0 & 0 \\
0 & 0 & 0 & 0 & 0.6 & 0 & 0 & 0 & 0 & 0
\end{bmatrix}
\end{array}
$$

Now using the expert's opinion we get the fit vectors. Suppose the fit vector B is given as B = (0 1 1 0 0 0 0 0 1 0). Using max-min in backward directional we get FAM described as

$$A = M \circ B$$
that is $\quad a_i = \max \min (m_{ij}, b_j) \quad 1 \le j \le 10.$

Thus A = (0.8, 0.8, 0.6 0 0 0 0) since 0.8 is the largest value in the fit vector in A and it is associated with the two nodes vide $W_1$ and $W_2$ child marriage / widower child marriage etc finds it first place also the vulnerability of rural uneducated women find the same states as that of $W_1$. Further the second place is given to $W_3$, disease untreated till it is chronic or they are in last stages, all other states are in the off state.

Suppose we consider the resultant vector

$$A = (0.8, 0.8, 0.6, 0 0 0 0)$$
$$A \circ M = B. \text{ where}$$
$$B_j = \max (a_i, m_{ij}) \qquad 1 \le j \le 10.$$

A o M = (0.9, 0.8, 0.8, 0.6, 0.9, 0.7, 0.7, 0.0, 0.7) = $B_1$

Thus we see the major cause being infected by HIV/AIDS is $R_1$ and $R_5$ having the maximum value. viz. female child a burden so the sooner they get married off the better relief economically; and women suffer no guilt and fear for life. The next value being given by $R_2$ and $R_3$, poverty and bad habits of men are the causes of women becoming HIV/AIDS victims.

The next value being $R_6$, $R_7$ and $R_{10}$ taking the value 0.7, women have not changed religion and developed faith in god after the disease, lost faith in god after the disease. Husbands hide their disease from family so wife becomes HIV/AIDS infected. Several conclusions are derived from the analysis of similar form. We give the conclusions in the final chapter using several experts' opinion.

Now we proceed on to describe one more model using FAM. Using the attributes related with women affected with HIV/AIDS along the rows and the causes of it along the column we obtain the fuzzy rector matrix $M_1$. The attributes given by the expert are $W_1,\ldots, W_9$ and $R_1,\ldots, R_{10}$ (given in Chapter two).

The gradations are given in the form of the fuzzy vector matrix which is as follows.

|       | $R_1$ | $R_2$ | $R_3$ | $R_4$ | $R_5$ | $R_6$ | $R_7$ | $R_8$ | $R_9$ | $R_{10}$ |
|-------|-----|-----|-----|-----|-----|-----|-----|-----|-----|------|
| $W_1$ | .9  | .8  | .6  | 0   | .5  | .7  | 0   | 0   | .4  | 0    |
| $W_2$ | .8  | .7  | .8  | 0   | 0   | 0   | 0   | 0   | .7  | 0    |
| $W_3$ | .7  | .6  | .5  | 0   | 0   | .8  | 0   | 0   | 0   | 0    |
| $W_4$ | .9  | .7  | 0   | 0   | 0   | 0   | 0   | .8  | .7  | 0    |
| $W_5$ | .8  | .6  | .7  | 0   | 0   | 0   | 0   | .7  | 0   | 0    |
| $W_6$ | 0   | 0   | .4  | 0   | 0   | 0   | 0   | .6  | 0   | 0    |
| $W_7$ | .4  | 0   | .6  | 0   | 0   | 0   | 0   | 0   | 0   | 0    |
| $W_8$ | .6  | .5  | 0   | 0   | .4  | .3  | 0   | 0   | 0   | 0    |
| $W_9$ | 0   | 0   | 0   | 0   | 0   | 0   | .8  | 0   | .4  | .6   |

Now using the experts opinion we get the fit vectors B, given as B = ( 1 1 1 0 0 0 0 0 1), using max-min backward direction we get FAM described as

$$A = M_1 \text{ o } B.$$

That is

$$a_j = \max \min (m_{ij} \, b_j), \quad 1 \le j \le 10$$
$$= (.9, .8, .6, .9, .8, .4, .6, .6, .6).$$

Since .9 is the largest value in the fit vector A and it is associated with $W_1$ and $W_4$, we see the root cause or the major reason of women becoming HIV/AIDS infected is due to no basic education, no recognition of their labour. One may be very much surprised to see why the co ordinate "no recognition of labour" find so much place the main reason for it if they (men i.e., the husband) had any recognition of their labour certainly he would not have the heart to infect his wife when he is fully aware of the fact that he is suffering from HIV/AIDS and by having unprotected sex with his wife certainly she too would soon be an HIV/AIDS patient. Thus $W_1$ and $W_4$ gets the highest value.

The next value is 0.8, which is got by the nodes $W_2$ and $W_5$. $W_2$ – No wealth and property and $W_5$ – No protection from their spouse. If wealth and property was in their hands (i.e., in the name of women) certainly they would fear her untimely death for the property may be taken by her relatives. Further these women cannot protect themselves from their husbands even if they are fully aware of the fact that their husbands are infected by HIV/AIDS.

Next higher value is .6 taken by $W_3$, $W_7$, $W_8$ and $W_9$ that is women are not empowered, they do not get nutritious food and they are not the decision makers or they do not have any voice in this traditional set up. Thus we have worked with several different nodes to arrive at the conclusions.

# USE OF NEUTROSOPHIC MODELS TO ANALYZE HIV/AIDS WOMEN PATIENTS' PROBLEMS

Neutrosophic study / analysis of the problems faced by HIV/AIDS affected rural, poor uneducated women in India, is carried out for the first time in this book.

In this chapter we use two types of neutrosophic tools. Neutrosophic Relational Maps (NRMs) defined in [129] and we to study and analyze this problem have introduced a new neutrosophic model called the Neutrosophic. Associative Memories (NAMs). This chapter has two sections. In the first section we just recall the definition and properties of Neutrosophic relational maps (NRMs) and apply them to study the social problems faced by rural poor uneducated women in India. In section two we for the first time introduced a new neutrosophic model called the Neutrosophic Associative Memories and is applied to this problem. For in this model alone one can compare one attribute / concept with the other in the resultant vector.

## 4.1 Description of a NRM:

Neutrosophic Cognitive Maps (NCMs) promote the causal relationships between concurrently active units or decides the absence of any relation between two units or the indeterminacy of any relation between any two units. But in Neutrosophic Relational Maps (NRMs) we divide the very causal nodes into two disjoint units. Thus for the modeling of a NRM we need a domain space and a range space which are disjoint in the sense of concepts. We further assume no intermediate relations exist within the domain and the range spaces. The number of elements or nodes in the range space need not be equal to the number of elements or nodes in the domain space.

Throughout this section we assume the elements of a domain space are taken from the neutrosophic vector space of dimension n and that of the range space are neutrosophic vector space of dimension m. (m in general need not be equal to n). We denote by R the set of nodes $R_1, \ldots, R_m$ of the range space, where $R = \{(x_1, \ldots, x_m) \mid x_j = 0 \text{ or } 1 \text{ for } j = 1, 2, \ldots, m\}$.

If $x_i = 1$ it means that node $R_i$ is in the on state and if $x_i = 0$ it means that the node $R_i$ is in the off state and if $x_i = I$ in the resultant vector it means the effect of the node $x_i$ is indeterminate or whether it will be off or on cannot be predicted by the neutrosophic dynamical system.

It is very important to note that when we send the state vectors they are always taken as the real state vectors for we know the node or the concept is in the on state or in the off state but when the state vector passes through the Neutrosophic dynamical system some other node may become indeterminate i.e. due to the presence of a node we may not be able to predict the presence or the absence of the other node i.e., it is indeterminate, denoted by the symbol I, thus the resultant vector can be a neutrosophic vector.

**DEFINITION 4.1.1:** *A Neutrosophic Relational Map (NRM) is a Neutrosophic directed graph or a map from D to R with concepts like policies or events etc. as nodes and causalities as edges. (Here by causalities we mean or include the indeterminate causalities also). It represents Neutrosophic Relations and Causal Relations between spaces D and R .*

*Let $D_i$ and $R_j$ denote the nodes of an NRM. The directed edge from $D_i$ to $R_j$ denotes the causality of $D_i$ on $R_j$ called relations. Every edge in the NRM is weighted with a number in the set {0, +1, –1, I}. Let $e_{ij}$ be the weight of the edge $D_i$ $R_j$, $e_{ij} \in$ {0, 1, –1, I}. The weight of the edge $D_i R_j$ is positive if increase in $D_i$ implies increase in $R_j$ or decrease in $D_i$ implies decrease in $R_j$ i.e. causality of $D_i$ on $R_j$ is 1. If $e_{ij}$ = – 1 then increase (or decrease) in $D_i$ implies decrease (or increase) in $R_j$. If $e_{ij}$ = 0 then $D_i$ does not have any effect on $R_j$. If $e_{ij}$ = I it implies we are not in a position to determine the effect of $D_i$ on $R_j$ i.e. the effect of $D_i$ on $R_j$ is an indeterminate so we denote it by I.*

**DEFINITION 4.1.2:** *When the nodes of the NRM take edge values from {0, 1, –1, I} we say the NRMs are simple NRMs.*

**DEFINITION 4.1.3:** *Let $D_1$, …, $D_n$ be the nodes of the domain space D of an NRM and let $R_1$, $R_2$,…, $R_m$ be the nodes of the range space R of the same NRM. Let the matrix $N(E)$ be defined as $N(E) = (e_{ij})$ where $e_{ij}$ is the weight of the directed edge $D_i R_j$ (or $R_j D_i$ ) and $e_{ij} \in$ {0, 1, –1, I}. $N(E)$ is called the Neutrosophic Relational Matrix of the NRM.*

The following remark is important and interesting to find its mention in this book.

**Remark**: Unlike NCMs, NRMs can also be rectangular matrices with rows corresponding to the domain space and columns corresponding to the range space. This is one of the marked difference between NRMs and NCMs. Further the number of entries for a particular model which can be treated as disjoint sets when dealt as a NRM has very much less entries than when the same model is treated as a NCM.

Thus in many cases when the unsupervised data under study or consideration can be spilt as disjoint sets of nodes or concepts; certainly NRMs are a better tool than the NCMs.

**DEFINITION 4.1.4:** *Let $D_1, \ldots, D_n$ and $R_1, \ldots, R_m$ denote the nodes of a NRM. Let $A = (a_1, \ldots, a_n)$, $a_i \in \{0, 1, -1\}$ is called the Neutrosophic instantaneous state vector of the domain space and it denotes the on-off position of the nodes at any instant. Similarly let $B = (b_1, \ldots, b_n)$ $b_i \in \{0, 1, -1\}$, B is called instantaneous state vector of the range space and it denotes the on-off position of the nodes at any instant, $a_i = 0$ if $a_i$ is off and $a_i = 1$ if $a_i$ is on for $i = 1, 2, \ldots, n$. Similarly, $b_i = 0$ if $b_i$ is off and $b_i = 1$ if $b_i$ is on for $i = 1, 2, \ldots, m$.*

**DEFINITION 4.1.5:** *Let $D_1, \ldots, D_n$ and $R_1, R_2, \ldots, R_m$ be the nodes of a NRM. Let $D_i R_j$ (or $R_j D_i$) be the edges of an NRM, $j = 1, 2, \ldots, m$ and $i = 1, 2, \ldots, n$. The edges form a directed cycle. An NRM is said to be a cycle if it possess a directed cycle. An NRM is said to be acyclic if it does not possess any directed cycle.*

**DEFINITION 4.1.6:** *A NRM with cycles is said to be a NRM with feedback.*

**DEFINITION 4.1.7:** *When there is a feedback in the NRM i.e. when the causal relations flow through a cycle in a revolutionary manner the NRM is called a Neutrosophic dynamical system.*

**DEFINITION 4.1.8:** *Let $D_i R_j$ (or $R_j D_i$) $1 \leq j \leq m$, $1 \leq i \leq n$, when $R_j$ (or $D_i$) is switched on and if causality flows through edges of a cycle and if it again causes $R_j$ (or $D_i$) we say that the Neutrosophical dynamical system goes round and round. This is true for any node $R_j$ (or $D_i$) for $1 \leq j \leq m$ (or $1 \leq i \leq n$). The equilibrium state of this Neutrosophical dynamical system is called the Neutrosophic hidden pattern.*

**DEFINITION 4.1.9:** *If the equilibrium state of a Neutrosophical dynamical system is a unique Neutrosophic state vector, then it is called the fixed point. Consider an NRM with $R_1, R_2, \ldots, R_m$ and $D_1, D_2, \ldots, D_n$ as nodes. For example let us start the dynamical system by switching on $R_1$ (or $D_1$). Let us assume that the NRM settles down with $R_1$* *and $R_m$ (or $D_1$ and $D_n$ ) on, or indeterminate on, i.e. the Neutrosophic state vector remains as (1, 0, 0,…, 1) or (1, 0, 0,…I) (or (1, 0, 0,…1) or (1, 0, 0,…I) in D), this state vector is called the fixed point.*

**DEFINITION 4.1.10:** *If the NRM settles down with a state vector repeating in the form $A_1 \rightarrow A_2 \rightarrow A_3 \rightarrow ... \rightarrow A_i \rightarrow A_1$ (or $B_1 \rightarrow B_2 \rightarrow ... \rightarrow B_i \rightarrow B_1$) then this equilibrium is called a limit cycle.*

### METHODS OF DETERMINING THE HIDDEN PATTERN IN A NRM

Let $R_1$, $R_2$,…, $R_m$ and $D_1$, $D_2$,…, $D_n$ be the nodes of a NRM with feedback. Let N(E) be the Neutrosophic Relational Matrix. Let us find the hidden pattern when $D_1$ is switched on i.e. when an input is given as a vector; $A_1$ = (1, 0, …, 0) in D; the data should pass through the relational matrix N(E). This is done by multiplying $A_1$ with the Neutrosophic relational matrix N(E). Let $A_1N(E)$ = ($r_1$, $r_2$,…, $r_m$) after thresholding and updating the resultant vector we get $A_1E \in$ R, Now let B = $A_1E$ we pass on B into the system $(N(E))^T$ and obtain $B(N(E))^T$. We update and threshold the vector $B(N(E))^T$ so that $B(N(E))^T \in$ D.

This procedure is repeated till we get a limit cycle or a fixed point.

**DEFINITION 4.1.11:** *Finite number of NRMs can be combined together to produce the joint effect of all NRMs. Let $N(E_1)$, $N(E_2)$,…, $N(E_r)$ be the Neutrosophic relational matrices of the NRMs with nodes $R_1$,…, $R_m$ and $D_1$,…,$D_n$, then the combined NRM is represented by the neutrosophic relational matrix $N(E) = N(E_1) + N(E_2) +...+ N(E_r)$.*

Now we proceed on to work for the same model for then only comparison would be easy for the reader.

The related directed neutrosophic graph given by an expert, earlier in chapter two.

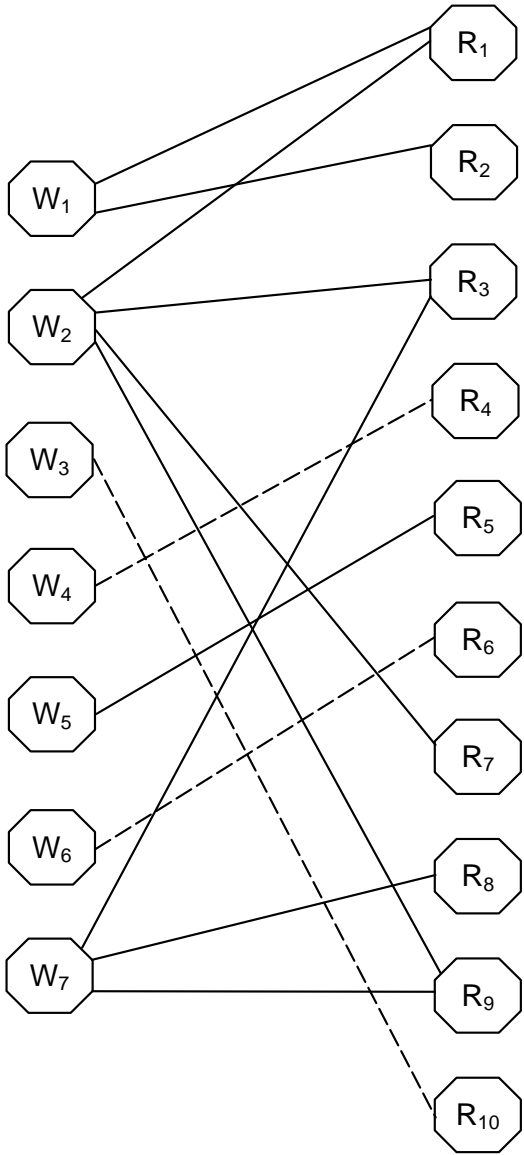

**FIGURE 4**

The domain nodes are taken as $W_1$, $W_2$,…, $W_7$ and that of the range space nodes are $R_1$, $R_2$, …, $R_{10}$.

The related connection neutrosophic matrix N(E) of this neutrosophic graph is as follows.

$$\begin{array}{c c c c c c c c c c c}
& R_1 & R_2 & R_3 & R_4 & R_5 & R_6 & R_7 & R_8 & R_9 & R_{10} \\
W_1 & 1 & 1 & 0 & 0 & 0 & 0 & 0 & 0 & 0 & 0 \\
W_2 & 1 & 0 & 1 & 0 & 0 & 0 & 1 & 0 & 1 & 0 \\
W_3 & 0 & 0 & 0 & 0 & 0 & 0 & 0 & 0 & 0 & I \\
W_4 & 0 & 0 & 0 & I & 0 & 0 & 0 & 0 & 0 & 0 \\
W_5 & 0 & 0 & 0 & 0 & 1 & 0 & 0 & 0 & 0 & 0 \\
W_6 & 0 & 0 & 0 & 0 & 0 & I & 0 & 0 & 0 & 0 \\
W_7 & 0 & 0 & 1 & 0 & 0 & 0 & 0 & 1 & 1 & 0
\end{array}$$

Suppose we consider the on state of the vector $W_1$ i.e.

$$X = (1\ 0\ 0\ 0\ 0\ 0\ 0)$$

that, is the cause for women becoming H*I*V/A*I*DS patient and all other nodes are in the off state. The effect of X on N(E) is given by,

$$\begin{aligned}
XN(E) &= (1\ 1\ 0\ 0\ 0\ 0\ 0\ 0\ 0\ 0\ ) \\
&= Y \text{ (say).}
\end{aligned}$$

Now we study the effect of Y on $(N(E))^T$

$$Y(N(E))^T = (2\ 1\ 0\ 0\ 0\ 0\ 0).$$

After thresholding and updating we get the resultant as $X_1 = (1\ 1\ 0\ 0\ 0\ 0\ 0)$. The effect of $X_1$ on N (E) is given by

$$X_1\ N(E) = (3\ 2\ 1\ 0\ 0\ 0\ 1\ 0\ 1\ 0)$$

Thresholding this vector we get the resultant as

$$\begin{aligned}
Y_1 &= (1\ 1\ 1\ 0\ 0\ 0\ 1\ 0\ 1\ 0) \\
Y_1\ (N(E))^T &= (2\ 5\ 0\ 0\ 0\ 0\ 1).
\end{aligned}$$

After thresholding and updating the resultant vector $Y_1N(E)^T$ we get

$$X_2 = (1\ 1\ 0\ 0\ 0\ 0\ 1).$$

The effect of $X_2$ on the dynamical system is given by

$$X_2\ N(E) \qquad = \qquad (2\ 1\ 2\ 0\ 0\ 0\ 1\ 1\ 2\ 0\ ).$$

After thresholding we get the resultant vector as

$$(1\ 1\ 1\ 0\ 0\ 0\ 1\ 1\ 1\ 0) \qquad = \qquad Y_2\ (\text{say}).$$

The effect of $Y_2$ on the dynamical system $(N(E))^T$ is given by

$$Y_2\ (N(E))^T \qquad = \qquad (2\ 4\ 0\ 0\ 0\ 0\ 3).$$

After thresholding and updating we get the resultant as

$$(1\ 1\ 0\ 0\ 0\ 0\ 1) = X_3 = X_2.$$

Thus the hidden pattern of the dynamical system is a fixed point given by the binary pair.

$$\{(1\ 1\ 0\ 0\ 0\ 0\ 1),\ (1\ 1\ 1\ 0\ 0\ 0\ 1\ 1\ 1\ 0\ )\}.$$

Thus the causes for women becoming HIV/AIDS infected is poverty among rural women and no proper health center are available and they have to depend entirely on their husbands and all other nodes $W_3$, $W_4$ $W_5$ and $W_6$ remain unaffected. In the range space we see the nodes $R_1$, $R_2$ $R_3$ $R_7$, $R_8$, $R_9$ become on i.e., they are affected by the on state on the node $W_1$ however $R_4$ $R_5$ $R_6$ and $R_{10}$ remain only in the off state. The case for women becoming HIV/AIDS infected is due to $R_1$ – HIV/AIDS infected husbands who do not disclose the disease to the family, $R_2$ – Male ego unconcern about wife. $R_3$ – Most men do not give to the family or the account for the amount of money they earn or spend, $R_7$ – They do not bother about the kith or kin or even their wife. $R_8$ – They suffer from guilt so they cannot lead a normal life at home. $R_9$ – They don't live with courage after being infected.

Now we study the effect of the on state the node $R_3$ from the range space and all other nodes remain in the off state. Let $Y = (0\ 0\ 1\ 0\ 0\ 0\ 0\ 0\ 0)$. The effect of Y on the neutrosophical system $(N(E))^T$.

$$Y(N(E))^T \qquad = \qquad (0\ 1\ 0\ 0\ 0\ 0\ 1) = X \text{ say}$$

The effect X on the system N(E) is given by
$XN(E) = (1\ 0\ 2\ 0\ 0\ 0\ 2\ 1\ 2\ 0)$.

After updating and thresholding we get

$$
\begin{aligned}
XN(E) &= & (1\ 0\ 1\ 0\ 0\ 0\ 1\ 1\ 1\ 0) &= Y_1 \\
Y_1\ N(E) &= & (1\ 4\ 0\ 0\ 0\ 0\ 3).
\end{aligned}
$$

After thresholding we get the resultant $X_1 = (1\ 1\ 0\ 0\ 0\ 0\ 1)$. As before we get the same binary pair $\{(1\ 1\ 0\ 0\ 0\ 0\ 1), (1\ 1\ 1\ 0\ 0\ 0\ 1\ 1\ 1\ 0)\}$ which clearly shows, when men do not give account of what they earn or spend they are prone to lose binding with family and infect the innocent dependent women.

Now we make the on state of $W_3$ i.e. women have no power to protect themselves from HIV/AIDS and all other nodes remain in the off state. The effect of $X = (0\ 0\ 1\ 0\ 0\ 0\ 0)$ on the dynamical system N(E) given by

$$
\begin{aligned}
XN(E) &= & (0\ 0\ 0\ 0\ 0\ 0\ 0\ 0\ 0\ I) \\
&= & Y_1 \text{ (say)} \\
Y_1\ (N(E))^T &= & (0\ 0\ I\ 0\ 0\ 0\ 0\ 0\ 0\ 0).
\end{aligned}
$$

After updating we get the resultant as $(0\ 0\ 1\ 0\ 0\ 0\ 0) = X$. Thus the hidden pattern is a fixed point and according to this expert all nodes in the domain space remain off and only $R_{10}$ node becomes an indeterminate.

Now we proceed on to study the FRM model with nodes $R_1, R_2,\ldots, R_{10}$ and $W_1, W_2,\ldots, W_9$ given in chapter two of this book using the NRM.

The directed neutrosophic graph given by expert is

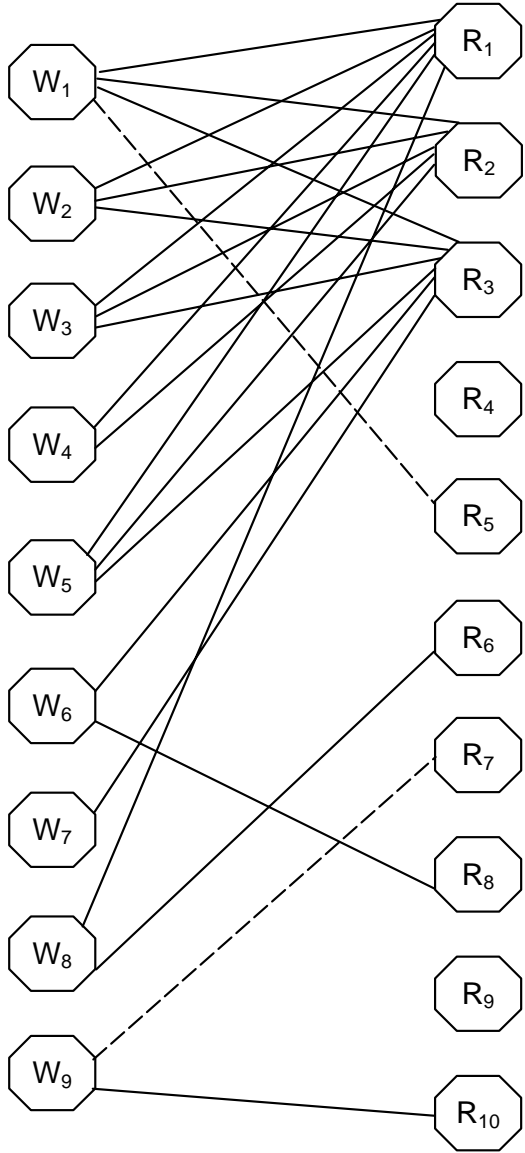

**FIGURE 5**

The related connection neutrosophic matrix N(E) of the NRM is

$$\begin{array}{c c c c c c c c c c c}
 & R_1 & R_2 & R_3 & R_4 & R_5 & R_6 & R_7 & R_8 & R_9 & R_{10} \\
W_1 & 1 & 1 & 1 & 0 & I & 0 & 0 & 0 & 0 & 0 \\
W_2 & 1 & 1 & 1 & 0 & 0 & 0 & 0 & 0 & 0 & 0 \\
W_3 & 1 & 1 & 1 & 0 & 0 & 0 & 0 & 0 & 0 & 0 \\
W_4 & 1 & 1 & 0 & 0 & 0 & 0 & 0 & 0 & 0 & 0 \\
W_5 & 1 & 1 & 1 & 0 & 0 & 0 & 0 & 0 & 0 & 0 \\
W_6 & 0 & 0 & 1 & 0 & 0 & 0 & 0 & 1 & 0 & 0 \\
W_7 & 0 & 1 & 0 & 0 & 0 & 0 & 0 & 0 & 0 & 0 \\
W_8 & 1 & 0 & 0 & 0 & 0 & 1 & 0 & 0 & 0 & 0 \\
W_9 & 0 & 0 & 0 & 0 & 0 & 0 & I & 0 & 0 & 1 \\
\end{array}$$

Now we study the effect of the state vector $X = (1\,0\,0\,0\,0\,0\,0\,0\,0)$ i.e. only the node $W_1$ i.e. No basic education alone is in the on state and all other nodes are in the off state. The effect of X on the dynamical system N(E) is given by

$$XN(E) \quad = \quad (1\,1\,1\,0\,I\,0\,0\,0\,0\,0) \quad = \quad Y.$$

Now the effect of Y on $N(E)^T$

$$Y(N(E))^T \quad = \quad (\,3+I,\,3\,3\,2\,3\,1\,1\,1\,0).$$

After thresholding and updating we get the resultant vector as $X_1 = (1\,1\,1\,1\,1\,1\,1\,1\,0)$.
Now the effect of $X_1$ on the dynamical system N(E) is given by

$$X_1\,N(E) \quad = \quad (6\,6\,5\,0\,I\,1\,0\,1\,0\,1).$$

After thresholding we get $Y_1 = (1\,1\,1\,0\,I\,1\,0\,1\,0\,1)$.

$$Y_1\,N(E)^T \quad = \quad (\,3+I,\,3\,3\,2\,3\,2\,1\,1\,1\,).$$

After thresholding and updating it we get the resultant

$$X_2 \quad = \quad (\,1\,1\,1\,1\,1\,1\,1\,1\,1\,1).$$

Now

$$X_2 N(E) \quad = \quad (6\ 6\ 5\ 0\ I\ 1\ I\ 1\ 0\ 1).$$

After thresholding we get

$$Y_2 \quad = \quad (1\ 1\ 1\ 0\ I\ 1\ I\ 1\ 0\ 1)$$

Now the effect $Y_2$ on $N(E)^T$ is

$$Y_2 (N(E))^t \quad = \quad (3+I,\ 3\ 3\ 2\ 3\ 2\ 1\ 2\ I),$$

after thresholding and updating we get the resultant as

$$X_3 \quad = \quad (1\ 1\ 1\ 1\ 1\ 1\ 1\ 1\ I).$$

The effect of $X_3$ on $N(E)$;

$$X_3 N(E) \quad = \quad (6\ 6\ 5\ 0\ I\ 1\ I\ 1\ 0\ I)$$

after thresholding we get the resultant $Y_3$ as

$$Y_3 \quad = \quad (1\ 1\ 1\ 0\ I\ 1\ I\ 1\ 0\ I).$$

The effect of $Y_3$ on $(N(E))^t$

$$Y_3 (N(E))^t \quad = \quad (3 + I\ 3\ 3\ 2\ 3\ 2\ 1\ 2\ I);$$

after thresholding and updating we get the resultant as

$$X_4 \quad = \quad (1\ 1\ 1\ 1\ 1\ 1\ 1\ 1\ I) \quad = \quad X_3.$$

Thus the hidden pattern is a fixed point given by the binary pair $\{(1\ 1\ 1\ 1\ 1\ 1\ 1\ 1\ I),\ (1\ 1\ 1\ 0\ I\ 1\ I\ 1\ 0\ I)\}$.

Now when the only node no basic education alone in the on state we get all nodes in the domain space on except $W_9$ which is an indeterminate i.e. according to this expert the relation between the concept $W_1$ – no basic education makes the node $W_9$ – No voice for women basically in the traditional set up is an indeterminate.

Further in the range space all nodes $R_1$ $R_2$ $R_3$ $R_6$ and $R_8$ become on. However $R_4$ and $R_9$ are in the off state and $R_5$ $R_7$ and $R_{10}$ become indeterminate. One can work with on state of any other node and analyze the system.

## 4.2 Neutrosophic Associative Memories

In this section we introduce the notion of neutrosophic Associative memories. Now we proceed on to define the new neutrosophic model viz. Neutrosophic Associative memories (NAM) for the first time, mainly to apply to this problem.

Now we proceed on to define essential tools in the introduction of Neutrosophic associative memories.

**DEFINITION 4.2.1:** *Let $<Z \cup \{I\}>$ or $<Q \cup \{I\}>$ or $<R \cup I>$ denote the ring or field generated by Z and I or Q and I or R and I where we use the property I.I = I but I + I = 2I so sum of I taken n times gives n I.*

Now we define min or max function as follows $\min\{x, I\} = x$ where x is real but $\min\{mI, nI\} = mI$ if $m <$ n. The notion of max function is defined by $\max\{x, I\} = I$ where x is real and $\max\{mI, nI\} = nI$ if $m < n$. With this definition of max and min in mind we can define the notion of Neutrosophic Associative Memories (NAM).

The simplest NAM encodes the NAM rule or association which associates a p dimensional neutrosophic set $B_i$ with the n – dimensional neutrosophic set $A_j$.

(Here by a n – dimensional neutrosophic set we mean the unit neutrosophic hyper cube $[0\ 1]^n \cup [0\ I]^n$). These minimal NAMs essentially map one ball in $[0\ 1]^n \cup [0\ I]^n$ to a ball in $[0\ 1]^p \cup [0\ I]^p$. They are comparable to simple neural networks. As FAMs the NAMs need not adaptively be trained. Apart from this the NAM model function on the same rules as the FAM models. Here the fit vector as in case of FAMs are $B_N = (a_1, \ldots, a_n)$ where

$$a_i \quad = \quad 0 \text{ if the node is in the off state}$$
$$= \quad 1 \text{ if the node is in the on state.}$$
$$= \quad I \text{ if it is an indeterminate.}$$

Here every the resultant fit vector $A_N = (a_1, \ldots, a_n)$ can take values in $[0\ 1] \cup I$ as in case of FAMs.

Now we obtain an opinion from the expert which includes the indeterminate values also, so the related matrix of a NAM will be called as a Neutrosophic vector matrix and will be denoted by $N(M)$.

For the $W_1, \ldots, W_7$ and $R_1, \ldots, R_{10}$ given in chapter three of this book the neutrosophic vector matrix $N(M)$ as given by the expert is as follows.

$$
\begin{array}{c}
\quad\ \ R_1\ \ R_2\ \ R_3\ \ R_4\ \ R_5\ \ R_6\ \ R_7\ \ R_8\ \ R_9\ \ R_{10} \\
\begin{array}{c} W_1 \\ W_2 \\ W_3 \\ W_4 \\ W_5 \\ W_6 \\ W_7 \end{array}
\left[
\begin{array}{cccccccccc}
.9 & .8 & .7 & I & 0 & 0 & 0 & 0 & 0 & .7 \\
.5 & .8 & .6 & 0 & 0 & 0 & 0 & 0 & 0 & 0 \\
0 & .3 & .6 & 0 & 0 & 0 & 0 & 0 & 0 & 0 \\
0 & 0 & 0 & .6 & 0 & 0 & 0 & 0 & 0 & 0 \\
0 & 0 & 0 & 0 & .9 & .6 & .7 & 0 & 0 & 0 \\
0 & 0 & 0 & 0 & 0 & .7 & .5 & 0 & 0 & 0 \\
0 & 0 & 0 & 0 & .6 & 0 & 0 & 0 & 0 & 0
\end{array}
\right]
\end{array}
$$

Now using the expert's opinion we get the fit vector. Suppose the fit vector $B_N$ is given as $B_N = (0\ 1\ 1\ I\ 0\ 0\ 0\ 0\ 1\ 0)$. Using max min in backward direction we get NAM described as

$$A_N = N(M) \text{ o } B_N$$

that is

$$a_i = \max \min (m_{ij}, b_j),\ 1 \le j \le 0.$$

Thus $A_N = (I\ .8\ .6\ .6\ 0\ 1\ .6)$ cause of women becoming HIV/AIDS infected becomes indeterminate. The highest value 1 is being taken by the node $W_6$ – Rural women neglect their health due to unconcern and supreme sacrifice for the family. The next higher value is taken by $W_2$ – poverty among rural women followed by $W_3$ $W_4$ and $W_9$ with associated value 0.6

$W_5$ is unaffected for it had taken the value 0; Free of depression and despair in spite of being deserted by family when they have H*I*V/A*I*DS; we have worked with several other coordinates.

We illustrate NAM with yet another model. It is pertinent to mention that the conclusions are derived from working with several experts and many combinations of fit vectors.

Using the model in described in chapter 3 we obtain the neutrosophic vector matrix with row vector $W_1$, $W_2$,…, $W_9$ and column entries as $R_1$, $R_2$,…, $R_{10}$. N(M) the neutrosophic vector matrix is a $9 \times 10$ matrix with entries from [ 0 1 ] $\cup$ *I*.

|     | $R_1$ | $R_2$ | $R_3$ | $R_4$ | $R_5$ | $R_6$ | $R_7$ | $R_8$ | $R_9$ | $R_{10}$ |
|-----|----|----|----|----|----|----|----|----|----|-----|
| $W_1$ | .6 | .8 | .5 | 0  | I  | 0  | 0  | 0  | 0  | 0   |
| $W_2$ | .3 | .7 | .8 | 0  | 0  | 0  | 0  | 0  | 0  | 0   |
| $W_3$ | .6 | .5 | 0  | 0  | 0  | 0  | 0  | 0  | 0  | 0   |
| $W_4$ | .7 | .8 | 0  | 0  | 0  | 0  | 0  | 0  | 0  | 0   |
| $W_5$ | .8 | .6 | .9 | 0  | 0  | 0  | 0  | 0  | 0  | 0   |
| $W_6$ | 0  | 0  | .7 | 0  | 0  | 0  | 0  | I  | 0  | 0   |
| $W_7$ | 0  | .6 | 0  | 0  | 0  | 0  | 0  | 0  | 0  | 0   |
| $W_8$ | .7 | 0  | 0  | 0  | 0  | 0  | 0  | 0  | 0  | 0   |
| $W_9$ | 0  | 0  | 0  | 0  | 0  | 0  | 0  | 0  | 0  | .5  |

Now we study the effect of a fit vector $B_N = $ (0 0 1 0 *I* 0 0 0 0) on N(M).

$A_N = $ N(M) o $B_N = $ (*I* .8 .0 .0 .9 .7 0 0 0). The resultant shows $R_3$ – when poverty in the villages is in the on state, $R_5$ – influence of media in an indeterminate state we get the resultant fit vector as $W_1$ – No basic education becomes an indeterminate. $W_5$ – takes the highest value 0.9 i.e., No protection from their spouse, followed by $W_2$.

No wealth or property is in their possession, with $W_6$ – taking value 0.7 i.e. No proper health care center to get proper guidance. However $W_3$, $W_4$ $W_7$ $W_8$ and $W_9$ remain in off state as they take the value zero, in the fit vector, $A_N$.

Now we find the effect of the fit vector $A_N$ on $N(M)$ defined by

$$A_N N(M) = \max \{\min (a_i \, m_{ij})\}$$
$$B'_N = (.8, .8, .9, 0 \, I, 0 \, 0 \, I \, 0 \, 0).$$

$R_5$ and $R_8$ becomes indeterminate i.e. impact of media and men have no responsibility, …, $R_4$, $R_6$, $R_7$ $R_9$ and $R_{10}$ are only zero, for they are not influenced by the fit vector, $B_N = (0 \, 0 \, 1, 0 \, I \, 0 \, 0 \, 0 \, 0 \, 0)$ by its resultant $A_N = (I \, .8, 0 \, 0 \, .9, .7 \, 0 \, 0 \, 0)$. $R_3$ takes the highest value i.e. no women empowerment, followed by $W_1$ and $W_2$ holding the next higher value, no basic education and no wealth / property in women's possession.

Thus one can work with any of the attributes and study the situation. All the conclusions based on our study have been included in the chapter on conclusions. Further we have not given several models and their analysis for our motive is that it should be read by non-mathematicians also.

# INTERVIEWS OF 101 WOMEN LIVING WITH HIV/AIDS

This chapter is entirely devoted to giving the condensed data of 101 rural women living with HIV/AIDS. It is important to note that these women were interviewed at the Government hospital of thoracic medicine at Tambaram, Chennai (herein after referred to as the Tambaram Sanatorium); where they were inpatients or taking treatment as outpatients in this hospital. Another point which is worth mentioning is that these women were mainly from the rural areas and belonged to poor or lower-middle class strata of the society, with no education or were educated up to the primary level. Further these women were housewives or agricultural coolies or worked in beedi factory. Also we wish to bring to the notice of the concerned reader that some of the interviewed patients did not talk with us freely or were not willing to answer our questions. Might be they were sad / shocked / depressed / dejected. So in certain cases we could understand only a section of their feelings. We have totally avoided technical terms with an intention that it can be read by a layman with a high school education.

At some places in the interviews, a reader may find certain contradicting answers or statements, yet we have given it as it is, for, our main motivation is to give the verbatim interview and not mould or manipulate it. So a shrewd reader may understand that the sequence of answers and questions are also not stereotypes.

It was left for the patient to speak out her mind for when they are upset or disturbed or angered by our interference they may not talk to us their real feelings. In some cases they refuse to answer us, some are so well prepared that their answer is "do not know". Some argue with us "what is the use of interviewing", when one is fully aware of the fact that there is no cure for this disease.

The predominant feeling, in 90 of the women 101 patients, is the expression of their inability to protect themselves and take care of their children. In case of older women who had married off their children, there was desperateness and dejection because they are not able to fondle their grandchildren. The overall feeling of these patients is helplessness. Most of them in the hospital accepted that the doctors and nurses in the Tambaram Sanatorium look after them properly. At least 60 of them were unaware of how the disease spreads. Uniformly all of them wanted the government to give them medicine. Majority of them felt that the government must help the HIV/AIDS infected children with free medicine, free food and free education. Those children who were orphaned must be provided with free hostel and good food.

As they were frightened of the social stigma, majority of them did not give the name of their husbands or the exact place of their stay or their native place.

Some of the women were very frank and openly acknowledged their extra-marital sex with other men. 37 of them had lost their husbands due to HIV/AIDS, 54 had their husbands living with HIV/AIDS; 6 had not yet requested their husbands to test for HIV/AIDS. The answer for this may be well known to their wives. A few (3 of them) of the husbands were not affected by HIV/AIDS.

96 of the women have been married when they were in the age group 12 to 19. It is worth mentioning that they may not be able to handle their children with maturity. Further even if they had grown fully or biologically, their mental state would be that of a child. In several cases they are married to men twice their age or widowers who were old men, which certainly is wrong both medically, biologically and psychologically, for there is an explicit absence of equality and blatant differences in opinion, view points.

Unless the Government implements the Sarada Act, the social status and the mental status especially of the rural women is at stake. This is more important in the case of rural woman. Also they were unaware of any other vocation. Only, they work in the fields and some of them in beedi factories so the, lack knowledge of any other occupation or self-support. This is one of the reasons for their early marriage, for the parents feel it as a burden to feed and support them. Some of them were married to old men and so these ladies for better sexual satisfaction approached other men secretly for sex. Most of their husbands had the habit of consuming alcohol. 68 of them were only agriculturist coolies. Now we give the brief summaries of the interviews along with our additional observation wherever possible.

**1**

On 01-08-2002 we met a 45 year old, uneducated woman. She hails from Tiruvannamalai. She is an agricultural labourer. Her husband is deceased. She has two children: a male and a female. Her son is an auto driver who has studied up to 10th standard and taken up this profession as the wages he gets as an agricultural coolie is very little in comparison with his profession as an auto driver.
She had chest pain and cough so they took X-ray of her chest. She says doctors look after her well. She is completely unaware of even the name of the disease HIV/AIDS; she says that she does not know; so she asked us what to think about it. Doctors do not say anything to her they only note down on a paper.

*Additional Observation: She looked very weak and was not able to stand and sit on her own. Nobody was around her.*

*She does not know the cause of her husband's death. She refused to disclose any symptoms with which her husband suffered or why he died. She did not answer about the position of her daughter, from her talks we could conclude that her daughter must be elder to her son. Further her voice was very feeble. She looked like a tuberculosis patient and she was coughing all the time. As she belonged to very poor strata of the society, she could not spend her time in seeing TV. She looked very old beyond her age and tired. She was not depressed because of the disease for she is totally unaware of anything about HIV/AIDS.*

## 2

We interviewed on 01-08-2002 a lady who appeared to be in her late thirties. She refused to answer the questions on educational qualification, number of children, whether married, her native village, parents' name. Later by slip of tongue she said that she has no children and she is not educated. But she very openly acknowledged, when she was in the age group 19 to 25 years she had the habit of having sex with many men. She said nobody has spoilt her. She went for sex on her own. She said she had not taken any other treatment like Siddha or indigenous medicine. She says she never believes in god. She came to know that she had HIV/AIDS in 1994. She feels that doctors take care of her well. She says her husband does not know about her activities and she has so far not disclosed them. She also knows weaving. Then we asked the question have you ever used condoms? To this question she answered in negative. When we asked why people do not speak openly about HIV/AIDS she says that she would boldly say that she has AIDS.

*Additional Observation: In the beginning of the interview she was reluctant to speak out, even did not give her age or name of her village. But as the interview advanced she gave particulars only about her education and her marital status. Till the end she did not disclose the name of her village or*

*her parent's name. She was looking happy and was not frightened or dejected. Might be her disbelief in god has made her feel free of any karma or such fictitious things. She looked contented, did not bother about the fact that she is suffering from HIV/AIDS. This is one of the qualities of her that need to be taken into consideration. She was willing to come out in public and say that she has HIV/AIDS. A very rare patient who was outspoken on certain occasions when it mainly concerned her and not her family. Such women can be used to spread awareness about HIV/AIDS.*

## 3

We interviewed the lady in her late thirties or early forties who did not disclose her age or education or hometown or number of children or her marital status on 10-8-2002. She refused to answer all these questions. She was having cold, could not walk, she lacked energy and could not even eat. She says every one in her village knows about her sickness. She acknowledges that she does not know any work. She refused to give the type of work her father was doing. She said that from the age of 21 she had the habit of having sex with several men. However she says that after her marriage she did not have sex with other men. She has never traveled. She says no one has ever forced her to have sex with other men. She has not taken any indigenous  medicine. She believes in God. She happily acknowledges that doctors take care of her well. She further says that her husband is fully aware of her past activities.

She says she came to know that she had HIV/AIDS in the year 1991. She further said her children have not acquired this disease. She says she does not need any help from anyone. She is uneducated. People in her village are frightened to move with her. However she says her parents did not say anything about her.

***Additional Observation****: Though she really looked sick, was unable to walk, she did not look dejected or speak with*

*dejection. She looked and talked quite contentedly. She did not feel sad because she is HIV/AIDS affected. She did not speak about her children in the first round of interview, she did not answer the question whether she has children. Even now we were not able to make out how many children she has. She was all alone; no body was around to take care of her. She had been in the hospital, as an inpatient for over 10 years. She says she feels better now. She never gave the name of her hometown or about the very existence of her brothers or sisters. Though she acknowledged that she is married she did not give the profession, name or address of her husband. From our conversation we come to understand that he is living. She did not answer the question whether her husband was tested for HIV/AIDS.*

## 4

On 01-03-2002 we met a 46 year old woman staying as an in-patient of the Tambaram Sanatorium hospital. She had studied up to 8th standard. She hails from Malliayakkarai. She is married and both of them work in the forest. She was married when she was 15 years, that is, she has 31 years of married life. Her husband aged 48 has studied up to the 7th standard. They have three children. All her children are married. She said she got this disease only after marriage. She says 8 years ago she had an attack of jaundice and was given 3 bottles of blood. She has cough and she further says doctors look after her well and give her medicine only for cough. She says she would have got jaundice because she ate a type of stone called *Sukkan* stone. When she went to the doctors for treatment they said she had no strength and it was due to weakness and gave her an injection. She said her main symptom was fever, cold, loose motion and lack of appetite. She was asked to take blood test. She says the doctor informed her son that she had HIV/AIDS and he informed his father, that is her husband. She says all of them take care of her well and treat her with affection.

When we told her that there is no cure for HIV/AIDS, she so causally said "what is there? Till the day I will be able to live, I will live." She said the disease spreads by blood transmission and injection, she said to prevent this, one should lead a proper life.

When questioned about her faith in god, she said some say this disease can be cured by faith in god and some say it cannot be cured by God. Whatever be the case she says it is the wish of god. She lives in a nuclear family and not in a joint family. She lives only in a hut as she is very poor. She is unaware about diseases like diabetes, cancer etc. She says they give medicine, so it is okay; if she becomes better and well. She says how can the children affected by HIV/AIDS be ill treated or discriminated. She says that is wrong. All her family members move well and affectionately with her.

*Additional Observation: She is least sad or upset about the fact that she has HIV/AIDS. She further gets all help from her family. From the conversation we understand that her husband is healthy. Being a woman she is not in a position to advice nor does she know that her husband must be tested for HIV/AIDS. Further it is a surprise to us that the doctors have not advised her husband to go for HIV/AIDS test. At times from her actions, we feel that she might have acquired the disease from blood transmission or injection. In this case the possibility she had acquired the disease by blood transmission looks more relevant for, she is able to say that the disease spreads by blood transmission and by injection. She does not suspect her husband. Even in her subconscious she does not think of him to be having extra martial affairs. We are not able to find out from her conversation whether her husband drinks or smokes. She too does not indulge in any other activities; only blood transmission seems to have been the cause of her disease. She is composed and free of any fear or frustration.*

**5**

On 1-3-2002 we interviewed a women aged 35 years who has studied up to the 5th standard. She hails from Salem. She is married for the past ten years and has a son. Her husband is a bore pump worker who usually visits home once or twice a year. She was with her mother and only her mother answered most of the questions. Her mother says she is not aware of how her daughter has acquired this disease. She was suffering from palpitation and cough. They have not tested their son for HIV/AIDS. She says she does not think anything about HIV/AIDS, but she says people in this hospital say that after taking medicine in this hospital they feel better. She thinks that this disease spreads by 'women' and through injection needles. Further she said that the doctors say that HIV/AIDS is not like tuberculosis, which spreads in air and can easily infect anyone, so she argues that children affected with HIV/AIDS can be admitted in any school without any discrimination for it does not spread by air or personal contact like several other diseases like tuberculosis. She also said that people in the hospital have said the previous doctor used to inspect all the wards regularly, the hygiene and the quality of food, but however now the food they get is not up to the standard as in those days. Several patients who are in the hospital for over 5 years have recorded these statements.

She was sick, so she had an X-ray taken and at first the doctor said that her intestines were shrunk. She had high palpitation in spite of taking medicine, so only she was advised to take treatment in this hospital. The doctor in this hospital however said it is not a very dangerous disease. Only people were talking very fiercely about it for which she said she felt very sad. She advises three things for the youth:

1. Should accept transmission of blood only after testing the blood sample.
2. Should not have sex with many women and strangers (CSWs).
3. One should not have more than one wife.

When asked can this disease be cured by faith in god? She said the disease cannot be cured by faith for she says that this particular disease HIV/AIDS there is no god. For she says to cure small pox or chicken pox Muthumariamma / Karumariamma is there. For cholera, Kali Amman cures it, but for HIV/AIDS there is no special god to cure it. She too acknowledged that doctors in this hospital take proper care of her. When we asked her what she expects from the government, she counter-questioned us by asking, will the government do what we expect? Then she said the Government should cure us. She says she has seen about HIV/AIDS advertisements in TV and they say this disease is a cruel one.

*Additional Observation: She was a very frank person. She even gave her photo to us because she does not feel bad about the disease. She was the one who insisted that while transmission of blood we should check the blood sample for HIV/AIDS. She further adds with so much of Hindu philosophy that so far people have not discovered the god, which could cure HIV/AIDS; like gods for small pox and so on. She did not blame her husband though her husband was a bore pipe labourer who used to visit her only; once or twice a year, she has not asked her husband to test for the disease; neither have the doctors advised him to go for blood test. From her conversation we feel that he is not very ill. Her special trait was when questioned she always puts out a counter question to us, which made the interview more interesting. She also shared with us the view of the other patients in several questions.*

**6**

We interviewed on 07-03-2002 an uneducated 21 year HIV/AIDS patient in the Tambaram Sanatorium. Her husband is also uneducated and works as coolie. Her parents in laws are just construction labourers. She has two female children and both of them are not infected by HIV/AIDS.

She is from Kadaya Nalloor near Tenkasi. They are very poor, they live only in a hut. She says she has got this disease from her husband. She says her husband is also an HIV/AIDS patient. Her husband too was with her when we interviewed her. He has had premarital sex twice at the age of 20, however he has not confided about this to her. He says he has not ill treated or tortured his wife. He is in this hospital for the past 2 years. She says she really feels sad for being infected by HIV/AIDS, by her husband. She says he had sex with CSWs in the platform of Kutralam, a tourist spot where a lot of CSWs are affordable on platforms, bus-stands etc. He spent Rs.100/-. After this he had scabies and boils in private parts for which he had taken treatment from a doctor in Tirunelveli. He was married at 28. He says he had told this to his wife only in the hospital and not before.

We also learnt from her that he used to beat her in anger when she used to ask him in which lady's house he stayed during night? When we asked, "you have infected your wife only by staying / having sex with CSWs?", "What is wrong with it?" he asked. He was infected by a CSW only, but when she asked "her house" he became angry and beat her. They recommend free HIV/AIDS test for all because by this scheme the poor like them will be benefited. They say they have faith in Siddha medicine but felt that Siddha cannot cure HIV/AIDS. He says he came to know about condom only 2 years ago. His wife was weeping for they have no one who is going to take care of them. She says she knows about HIV/AIDS.

She says he had no education that is why he is ignorant about HIV / AIDS. She said her parents-in-law blame her saying that their son is a good person and he would never do any wrong, only she infected him. They scold her of being an immoral and a characterless women. She wept more when she said about this. We asked, "Is your husband interested in second marriage?". She said her husband is uninterested, whereas her parents-in-law are interested in it. Since she has faith in him she is living with him, holding the double stigma of being a HIV/AIDS patient and being abused by her husband's family as a characterless women. They said they have taken indigenous medicine and spent Rs.8000 for it in 50 days, but it did not cure them. They only became more unhealthy. Her husband says that they took this treatment by the force of his parents. No one has so far come to see them. She says her parents know about her disease and they wanted to give her special medical treatment, she did not agree for it so they are little upset with it.

She says her husband has not wasted her property. She had some problems and had taken medicine for it also. She says often her husband was sick, and now she also became sick so they came to this hospital and have confirmed that they have HIV/AIDS. Her parent-in-laws, said by consulting a   fortuneteller, we can appease gods and one would be cured. She wanted to commit suicide. But she thought of her two daughters and so only she is living. They said people with social service mind must adopt HIV/AIDS children. Regarding providing job to HIV/AIDS patient by government they said that they are not even able to sit, how can ever they be able to do 8 hours of work. They felt their sorrow was lessened by our interview.

*Additional Observation: She was in tears several times during the course of our interview. They were discussing with us for over 30 minutes. They were happy to be with us. Both were very weak and were not even able to sit or stand easily. They are from very poor social strata. She says in the presence of her husband that he used to beat her. Only he has infected her, but her husband's family has abused her and called her an immoral lady. Further they have a selfish wish to marry him again. They do not think he would be infecting one more woman and that lady's life too will be ruined.*

*The actions of parents-in law are very cruel and treacherous. But her husband happens to be kind and considerate with her. She was very emotional and wept several times. She was worried about her two daughters. She feels sad that no one visits her. She recommends free HIV/AIDS test by government for she feels the poor people will certainly be benefited. Both were not able to sit, that was their health condition at the age of 21 and 28. It was saddening for us to witness their health conditions in comparison with their age.*

On 7-3-02 we met a married muslim woman aged 27 years with no educational qualifications presently in the Tambaram Sanatorium as an in-patient. Her husband aged 30 years is alive but does not go for any work. He is educated up to 12th standard; her father-in-law is no more. She hails from Ariyamangalam in Tiruchi. She has 3 children; 2 boys and a girl. They are very small children just studying 6th, 5th and 2nd standard respectively. Both were working as daily labourers in the beedi factory. She says she was living happily and it was a happy family. Only now after she contracted the disease she is away from her husband. She claims that she has no hereditary property. They earned Rs. 45 each. She says when they lived as a joint family they used to have only petty quarrels and it used to get over in due course of time. She however acknowledges that she never suspected her husband. She has seen advertisements about HIV/AIDS in T.V. She said she had dots all over her face and body. She took treatment for the same from private doctors and has spent lots of money over it.

After blood test they found that she had HIV/AIDS, so she came to this hospital and got admitted. She denies any relationship with any man before marriage. She says she does not know about her husband's relation with other women and he has so far not confided to her of any such relationship. She says that AIDS spreads through mosquito bite, and sex relationship with strangers; and only in these six months she has been infected.

When we asked whether they have tested her husband she says he always said he had ulcer. Now, she has realized that might be he is hiding the disease from her for the fear that the family name would be ruined. However she added that when he comes to see her she would ask and verify from him. She says her husband does not drink but has stayed away from home on several nights in the pretext that he was seeing movies or playing cards. It was not known to her that he would have had other habits like extra marital sex. No one in the family knows about it and only her mother-in-law admitted her in this hospital. Her mother also knows about this but she is very kind and affectionate with her. Her mother feels very upset that because of her son-in-law, her daughter is suffering like this. She confides that none of her friends know about this.

She says she is not aware of condoms but she has used contraceptive methods to stop pregnancy. Only now she has come to know about HIV/AIDS. She has seen about condoms advertisement in TV. She says since it was a relationship between herself and her husband they have never used it. She said she does not suspect her husband but on certain days he never used to come home in the nights. When asked why he had not come he would say they did not give him salary. But she says that he took medicine regularly; now she feels that he should have taken medicine only for HIV/AIDS, of which she was not aware of. She is of the opinion that more advertisement about HIV/AIDS must be given. School children should be taught about HIV/AIDS in 10th to 12th standard; only then, the youth can be saved.

She is of the opinion that government jobs can be given to HIV/AIDS patients provided they are healthy and can carry out their duties without any difficulty. The government must protect children affected with HIV/AIDS. She says government should give good medicine so that they can be cured. If they ask for rice or provisions from the government and give details of their real pay, money they earn the hospital workers and other patients look down upon them. She believes in god and she likes Jesus Christ very much and has lots of faith in him. She also says that every day before she goes to sleep, she prays to Jesus. For youth, she advises that they should never make the mistake of having sex with strangers. Her only symptom was dots on her body; at first doctors mistook it for mosquito bite, for when she scratched it, it became black. Apart from this she had no other symptoms. She is hopefully waiting to go home to see her children after she is cured.

*Additional Observation: This is an instance were the husband has not disclosed about his disease to his wife but has been taking treatment secretly and without any conscience, he had infected her even after he knew he was suffering from HIV/AIDS. Further he has no strength even to go for job. His stayed away from the family in the nights alone has made his wife feel very uncomfortable and sad. He had said he was taking medicine for ulcer when he was really taking medicine for HIV/AIDS. They have 3 children: all of them below the age of twelve. They have not been tested for HIV/AIDS.*

*She is cheerful, open and happy, she has complete hopes that she will be cured soon. Though she is a muslim she has enormous faith in Jesus Christ. Only her mother-in-law has admitted her to the hospital. This woman's mother was very sad to see her daughter infected by HIV/AIDS by her son-in-law.*

## 8

On 12-03-02 we met a HIV/AIDS patient from the Tambaram Sanatorium who did not want to give any personal information about her. She declined to answer the questions like her age, martial status, educational qualification, native place, whether she has children and her parents' name. Her first symptoms were lack of energy, cold, inability to walk and eat. All friends and known persons ran away with fear from her. She knows to repair cycles; apart from this, she does not know any other work. She started having sex at the age of 20. She says even before her marriage she had sex with several men. However she has not gone out of her native place. She says no body has ever forced her and she did all these acts of her own choice. She feels better after taking medicine. She says she has not taken indigenous medicine.

She has no belief in god. She had come to know that she was affected by HIV/AIDS from 1996. She acknowledges that the food given in this hospital previously was of good quality but now the food is not up to the standard. She further says that the rice is so much spoilt and stinks of bugs. She further adds that it would be of very great help if they provide good food, at least good rice in the hospital. She knows about condoms but has never used it. She acknowledges that doctors take care of her, further she confirms her husband knows all about her activities. She says she is not frightened about HIV/AIDS she would come openly and say in public that she has HIV/AIDS. She refused to give any opinion about CSWs.

*Additional Observation: She looked contented and happy and did not show anxiety or sadness because she has HIV/AIDS. She did not blame anybody. Doctors should take steps to check the blood samples of her husband for HIV/AIDS. She did not expect anything from the government but only good rice and good medicine. She must be in her late thirties (we could only guess) she was very thin and was not able even to sit or walk.*

**9**

On 04-04-2002 we interviewed an in-patient from the Tambaram Sanatorium. She must in the age group of 27-29 years. She is from Chennai. She is married and has two children. She has so far not tested her children for HIV/AIDS. Her husband is deceased. She does not know of what he died. When we asked her, "Did he die of HIV/AIDS?" she did not answer. However she said that for the past 5 or 6 months she had been in this hospital taking treatment. She says she is affected by headache and vomiting. Now she wants to be discharged from this hospital to go home. She says she was married when she was very young, around 13 years and she never liked this married life. She says all the more she did not like her husband. Without taking her full permission her parents had married her to him and now he is dead. Her relatives particularly her brother, is very kind and affectionate with her. She is very worried and does not know how she is going to support her children. She is not able to comprehend the inability of the scientists all over the world who are not in a position to find any cure for this disease.

*Additional Observation: This is a tragedy of forced marriage. When she was hardly 13 years old she was forced into marriage by her parents and relatives. Now her husband is dead and she has to carry both the social stigma of the disease and widowhood for no fault of hers. She has the added burden of supporting the two children who are in their early teens. She knows both Tamil and English. She is more worried about the future of her children than her disease.*

## 10

We met an uneducated 40-year-old woman on 08-03-2002 at the Tambaram Sanatorium. She is from Salem. She knows both Tamil and Telugu. Her husband is also uneducated. Their main means of livelihood is breaking stones, and working as mason. She has three children. She was married when she was 18 or 19 years old. She says her marriage was an arranged one; her husband was quite old and she was his second wife. She did not say anything about his first wife. We could not even know whether his first wife was living or dead. They earn about Rs.2000 to Rs.3000 a year. They have a own house at Salem. She has no other property. She lives in a joint family. Her marriage was forced on her without her consent. She did not like her husband. She says her husband had sex with other women also. She is affected with HIV/AIDS for the past four years. She has tuberculosis, asthma and heart disease. Only after check-up she came to know that she has HIV/AIDS. She says her doctor in K.K.Nagar said that she would have got the disease through injections. Only he advised her to admit herself in the Tambaram hospital.

She says her husband had only asthma because he is very old they (doctors) say he cannot get HIV/AIDS. He is normal but they have not tested his blood for HIV/AIDS. She says they have not tested their children for HIV/AIDS.

She says only doctors advised that it was not essential for them to undergo this test. She is unaware of the presence of this disease among her friends . She says happily that they go for movies and return home. Her husband does not drink everyday, he drinks occasionally and has not wasted any property of hers by the drinking habit. She says none of her relatives know that she is suffering from this disease. Only her brother and mother know that she is infected with HIV/AIDS and they only have admitted her in this hospital. She does not say this  to others.

She says she is not sad or depressed, she is courageous and says she does not bother about the disease. She says that she does not know anything about HIV/AIDS. Only after the blood test, the doctor said that she had HIV/AIDS. She says she has seen advertisements but forgets about them immediately. She says only when doctors or nurses say or discuss about HIV/AIDS she understands about it deeply, thinks about it and worries how she has got it. By thinking and researching of this, she says she is losing weight. She first had faith in the deity Periyasamy but for the past 3 years she has become a faithful follower of Jesus with a fond hope that he would cure her. She says soon she will be baptized.

She does not know about providing HIV/AIDS awareness education to school children. When asked how to take care of the HIV/AIDS affected children she says they can be kept in the hostel. She said if someone informs her of some work she would be very happy to do it. The doctors take good care of her. Since her mother was very old she did not give any dowry at the time of her marriage. All her things are intact at home and her husband has not wasted anything. She feels advertisements about HIV/AIDS can be given in villages. She also welcomes free HIV/AIDS test for one and all; she feels it is good. She says the youth should not make any mistakes. She was very happy to be interviewed.

*Additional Observation: This is an instance of forced marriage. But it is not known how the wife is affected but*

*the husband is free from HIV/AIDS. She is happy,*
*philosophical and does not feel openly sad about HIV/AIDS,*
*but from her statement that everyday she loses weight,*
*worrying about how she got HIV/AIDS made us think twice*
*and confirm that she is a little upset over the disease but*
*does not show it out openly.*

## 11

On 07-03-2002 we met a 40-year-old woman from the
Aalancheri village in Uthiramerur town, Kancheepuram
district. She is an agricultural labourer and an in-patient of
Tambaram Sanatorium. Her husband is 60 years he is also
an uneducated agricultural labourer. She has a small land,
which will yield 6 to 7 thousand rupees a year if the yield is
good. She was married at 13 years and has 3 children, her
eldest son is 26 years old now working in a private
company. Her second child is a daughter and is given in
marriage at the age of 14. She had incurable cold and she
went to a specialist and took treatment, he referred her to
the Tambaram Sanatorium. That time they said she had only
tuberculosis but later when she began to get fever and cold
they suspected her of having HIV/AIDS so they carried out
blood tests for the same. She is from a nuclear family.

None of her relatives know about it, she wants to check
her children and husband for HIV/AIDS. She has not seen
any advertisement about condoms. She is least bothered
about what should be done to small children who are
affected with HIV/AIDS, she says let the government do
whatever it likes. She is totally supportive of the idea of the
government giving job, employment to the HIV/AIDS
patients who are healthy and can do work.

She feels doctors take good care of her. She says sex
education to children can be given from the age of six. She
does not know about HIV/AIDS and does not know how
one gets HIV/AIDS. She prays to Muneeswaran once in a
year. She advocates that a man must have only one wife.
She asks; "will a husband ever say to his wife about his
extramarital affair?". She says she suspects her husband.

She says government should find means to cure the patients from this disease. However she acknowledges that the medicine they give saves her from headache and fever.

*Additional Observation: She is a victim of child marriage. That too she was married to a man aged 33 years when she was barely 12 or 13 years old. She has followed this convention by marrying off her daughter at the age of 13 or 14 years. She has become a mother at the age of 14. That is why she says she has a son aged 26 years when she is hardly 40 years. This is one of the major drawbacks of our society. She is cheerful and does not look dejected or worried. She does not know about the fierceness of this disease or the social stigma it carries.*

*However she feels it is a 'bad disease' so she has not disclosed this to her relatives. She also advocates that extra-marital sex is bad for youth and men. In her innermost mind she suspects her husband. It is unfortunate that many men take treatment for HIV/AIDS without even informing of it to their wives. That is why by the time their wives come to know of the disease it is in the chronic stage.*

## 12

A Gounder women aged 30 years, hailing from Cheyyar was interviewed on 03-09-2002. She and her husband (who is around 40 years) are both uneducated. They live in a nuclear family. She is poor and she says her husband drinks occasionally. She says she does not know how she has got the disease. She was married at the age of 18. It was an arranged marriage. She checked about this disease in Cheyyar. She went to the doctor for chronic cold and she was losing body weight; so only she was admitted in this hospital. Now she says she feels little better. She says even her state of mind is little better only in this hospital. She says she had no connection with any other man before marriage. When asked whether her husband has extra marital affairs, she says "how can I know that?" She does not know whether the blood sample of her husband was tested for HIV/AIDS. She does not know anything about HIV/AIDS and she has not seen any advertisement on it.

She is an illiterate so she says she could not read the posters behind the buses or autos. She is for more awareness programs in villages. Doctors take proper care of her. She says the children affected with HIV/AIDS must be protected. She does not feel that HIV/AIDS is a stigmatized disease. She expects good medicine and injection from the government. She is not aware of condoms. Her major symptom was cough; cold and at times she coughed blood.

*Additional Observation: She did not answer about children. She says her husband has so far not been tested for HIV/AIDS. She is totally unaware of how the disease spreads. She and her husband have never entered school. They are fully ignorant about the disease. She does not even suspect her husband. We were not in a position to educate her about how the disease HIV/AIDS spreads. She is happy at the hospital for she feels better.*

## 13

On 04-03-2002 we met the 5th wife of a man belonging to the Chettiar caste. He was in the age group 40 to 50. She was just 22 years old. She had been married at the age of 12 years. This is a picture perfect example of child marriage. She has two children both male. She proudly acknowledges that she has studied up to the 2nd standard and her husband has no education. He is a merchant and she is a construction labourer. She had loose motion, so she took treatment for a month as an outpatient. Now her health became worse so she has admitted as an inpatient in this hospital. HIV/AIDS was confirmed only after the blood test. She does not know about her disease fully for she says "the doctors have written in the chit that she has AIDS". She says that AIDS will make one lose weight ultimately causing death.

She says there is no one who has discovered any medicinal cure for HIV/AIDS. She says she does not know the ways in which HIV/AIDS spreads. She proudly says that her husband has five wives and she is the fifth one. He has all bad habits including chewing and smoking tobacco.

When asked "does your husband visit CSWs", she says that once a while he has said that he has visited CSWs but she did not believe it. She hails from Tiruvannamalai. To prevent this disease she says one should only ask the doctors. When asked "what does our society feel about HIV/AIDS, she says they think the patients will die soon. She says she does not need money from government if government cures her it is enough. When asked, "Can god cure HIV/AIDS?" she says "He gives more and more diseases, that too one after another." She did not give any direct answer!

*Additional Observation: At the first place this is yet another instance of child marriage. That too, to a man twice her age, who already has four wives, and she takes the place of the 5th wife. Secondly she has no voice to ask her husband to go for HIV/AIDS test, neither has she any knowledge or idea to test her two sons for HIV/AIDS. It is a pity that the status of women especially in rural areas is very poor. Further she is really proud to state that she is educated up to second standard whereas her husband is uneducated. She is a child marriage victim who does not feel the sufferings or the social stigma of this disease. She is calm and contented because of her ignorance. She is innocent to know the consequences of a possible HIV/AIDS infection of her sons.*

## 14

On 11-03-2002 we met a woman who says she does not know her age. She belongs to the Mudaliar caste and her husband is dead. She proudly acknowledges that her husband has studied up to 1st standard but she does not know his age. Both her father-in-law and mother-in-law are dead. She hails from the Authur in Salem district. She had a son who died at the age of 8. She has a daughter who has never been to school, and is now married. She says she prays to Mariamma and Murugan. She is an agricultural labourer and has no hereditary property. She belongs to nuclear family. She knows Tamil and little Telugu.

When she was 13 years old she got married to her uncle that too just 3 months after attaining puberty. She lives in her father-in-laws house. She does not suspect her husband. She says two years have passed since her husband's death. She further says only after his death this disease started manifesting in her and in a month she became very sick. He was also very sick and died due to the over swelling of the prostate glands. She says she sold most of her property to give treatment to her ailing husband. She says, before marriage and after marriage she had no relation with any other man.

When her husband had this disease she had suspicion. Further he used to drink everyday. He also had ulcers in the private parts. It seems that once she has asked him, whether he has visited CSWs, he said he has not. She does not know the year when her husband got infected. When we asked her how her son died she says that because her husband scolded him he hanged himself in the house. She says she came to know about HIV/AIDS during the past one year. In a year she was very ill four times. All relatives know that she has HIV/AIDS but they are kind with her. Even doctors are kind with her. Everyone in the hospital is friendly to her. She is uneducated, has not read the advertisement about HIV/AIDS but she has seen it on TV.

When asked why people do not talk about this disease openly, she says only God knows. She adds that it spreads from husband to wife and wife to husband and is a 'bad' disease. As he was ill for a very long time, most of the property she had was sold for giving him medical treatment. Now she is all alone and there is no one to even give her a cup of tea. She has changed from her native religion to Christianity and prays only to Jesus Christ. When asked about what should be done to the HIV/AIDS affected children she says only God knows the answer. She feels those affected with HIV/AIDS can be given government job. She says if she is cured it is enough for her. She also said that when she first heard that she was affected with HIV/AIDS she was heartbroken.

*Additional Observation: This is another instance of child marriage. This woman proudly acknowledges that her husband has studied up to 1st standard. From this we indirectly learn that women, that too from the remote villages, value education in their inner minds and they want to get educated and aspire for educated husbands. From our interview we are not in a position to know whether her husband died of HIV/AIDS, one thing we can say is that he had STD/VD and she says that it must have been the cause of his death.*

*She has spent a lot of money on giving treatment to her husband because of which now she is so poor that she cannot even afford a cup of tea. This makes us feel sympathetic towards these women who are from villages and live like orphans in the hospital. They must be given some counselling and material and physical support. If there are support groups it may help them. She is little upset and does not wish to inform her son-in-law for the fear that he may torture her daughter. She looked weak and her clothes showed her poverty. The Government should help such patients with free sarees and proper bedding. Like giving of free dhotis and sarees for the poor they can distribute clothes to HIV/AIDS patients especially the in-patients in the hospital.*

## 15

On 08-03-2002 we met a woman aged 31 years hailing from Vasudevanpattu village in Thiruvannamalai district. She has studied up to the 8th standard. Her husband aged 45 years has studied up to 6th standard. They have three children. Her first son is studying in the 9th standard. One daughter is studying in the 6th standard and the last daughter is in 2nd standard. She was married when she was just thirteen. She says her marriage was a love marriage. Her parents-in-law are dead. She lives in her own house. She has one-cent of land inherited from her parents. Her job is brick-cutting. Her husband also does the same work. They earn around Rs.15,000 per year. They had lived happily. She says she does not know how she got this disease. Two years ago her husband was very ill. She took treatment from every doctor in Thiruvannamalai but it was in vain. No one said anything about his cure. He had blisters all over the mouth. For two years he was in a hospital in Thiruvannamalai and as he was not cured she took him home. She did not think she would get this disease, also she says she does not know how he got this disease. He had fever and was very ill so they took blood test in the Tambaram Sanatorium where it was diagnosed that he has HIV/AIDS. He is taking treatment for the past three years. She says her husband did not say to her that he was suffering from HIV/AIDS.

Recently, she had fever, headache and pain in hands and legs so she got admitted in the hospital and they say she has HIV/AIDS. She further adds that for the past three years she had no sexual relation with her husband. She is surprised then how she has got the infection. Her parents are not aware of the disease from which she suffers. She further says she was aware of HIV/AIDS, as she has seen advertisements in TV, posters behind buses, etc. She knows it is a dangerous disease, that is why she did not have sex with her husband. Therefore she is surprised. When she first heard she was infected with HIV/AIDS she was sad. Sometimes she felt that she should not live. When asked what the government should do for small children who are affected with HIV/AIDS she said "they must be given free medicine and the government should protect them." She expects free medicine from the government. She says doctors take proper care of her and give her good medicine.

She says she was aware of HIV/AIDS for the past five years. She said that at that time her husband was hale and healthy. She says they are the first ones in their village to get HIV/AIDS. If government gives employment quota to HIV/AIDS patients she says it would be of immense help to them. Since this disease is transmitted sexually, she says people think very badly of this disease. She says that she has not had any extramarital affair with anyone, but her husband has extramarital affairs and he has confided that to her. She says he drinks everyday. He has squandered a lot of money on drinking. She has earned a lot of money on her own and purchased several things. When we asked her whether everyone can be given free HIV/AIDS check-up by the government she said it would be nice for everyone will definitely take up free check-up.

She says sex education can be given from the age of 10 i.e., from 5th standard. She advises youth not to have any kind of sexual relationship with strangers. She says "if you take injection buy new injection needles." She is happy to be interviewed. She feels no one should get this disease. She has not mentioned about the disease to anyone. She says she will carry out blood test for her children for HIV/AIDS only in the month of April when they are having the summer vacation.

*Additional Observation: She is being treated well by the doctors. In case of her husband, only the doctors made the disease chronic because they did not diagnose it for a very long time. She has become a victim of HIV/AIDS only because of her husband. She is not aware of the fact that the incubation period for HIV/AIDS takes over a period of several years depending on the immunity of the body and the food habits. If they take good food and are strong they can postpone the disease. This is another case of child marriage, she was married at 13 years and now she has become a HIV/AIDS patient at 31 years, after giving birth to 3 children. She is the only one who advises others to buy new needles for injection and not to use the used needles, which the doctors have used on others.*

## 16

We interviewed (on 09-04-2002) a construction labourer aged 35 years from Ennore who is presently in the Tambaram Sanatorium. She is uneducated and has never even entered school. She was married when she was just 15 years old and has four children. She says 20 years have passed since he ran away from her.

This information is contradictory; perhaps the age, which she has given is wrong. She says he is a drunkard who drank everyday. She says she had taken treatment in the Chennai General Hospital and Stanley Hospital. She answered in the negative to questions like knowledge about HIV/AIDS, about the disease suffered by the patients in the neighbouring bed and the results of the blood test. She says doctors look after her properly. She says her husband had all bad habits including visiting CSWs.

*Additional Observation: Her age at marriage and present age, the number of children and the fact that her husband had left her 20 years ago do not have any proper correlation. But it was not possible to re-question her. She does not say unlike the other women, the age of her four children and what they are doing. As interviews with HIV/AIDS patients is a delicate matter we could not cross-question her. She was not looking depressed. She is composed and neither agitated nor dejected by HIV/AIDS. At no time has she given any words of concern about her children. She did not answer questions related to her parents and her husband's family. She has not answered the question "Who is taking care of her and her children." We could not ask her about her husband, as he was not living with her now. She is very reserved with the neighbouring patients also.*

## 17

On 14-04-2002 we interviewed a women aged 25 years educated up to 3rd standard, a native of Otteri, Chennai. She is married and is a housewife. She was married when she was just 19 years old. Her husband is aged about 30 years and used to run a canteen. She is in the hospital for the past 20 days, she came to this hospital for treatment for tuberculosis but now she is suffering from HIV/AIDS also. She had fever so she took treatment in a Government

Hospital for tuberculosis and she became very ill so they referred her to this hospital. She lost weight, has no hunger, could not eat; vomiting and fever were the symptoms she was suffering. She says doctors look after her well. Now she feels a little better.

She says that HIV/AIDS is contracted because "women go to men", she does not know any other means by which this disease spreads. She says she did not suspect her husband after marriage and she did not know about his practices before marriage. She says she lives in a rented house paying a rent of Rs. 400/- a month. She has no children; the fetus got naturally aborted four times. Doctors said she has a weak womb. Her main symptom was ulcers / thrush in the mouth and itching all over the body. She further adds that people think of this disease as a very dangerous disease and people affected with this disease will not live long. She feels that by faith in Jesus Christ the disease can be cured.

*Additional Observation: By talking with her and observing her activities we feel that she would not have been be infected by her husband for he is hale and healthy. Also when asked "did he go for test?" she did not answer the question.*

*She does not suspect him; we are forced to think whether she was infected in the hospital while she got admitted four times due to natural abortion and blood transmission. This can be the only cause for her being infected by HIV/AIDS.*

## 18

A woman from Kosavapalayam in Pondicherry was interviewed on 08-04-2002. She was 40 years old and has studied up to 6th standard. She had been employed in an export industry. She is married for the past 22 years. She was married when she was just 18 years. She has 3 children, 2 male and one female.

Her husband aged 45 years, used to work in a private firm in Guindy. He had studied up to SSLC. For the past 2 years she is taking treatment in this hospital as an in-patient. She says she has contracted the disease from her husband. Her husband died of AIDS just two days back. He had tuberculosis, elephantiasis and Jaundice apart from HIV/AIDS. He used to smoke and drink cheap liquor. She and her husband used to quarrel often. She says there is no method to stop the spread of HIV/AIDS. She says she is not able to eat properly after the disease. She says God has not saved her. None of her friends or neighbours know about her being infected with HIV/AIDS. She knew about HIV/AIDS only 2 years ago; that too only after getting this disease.

*Additional Observation: Her children do not give their salary to her, after they came to know that she has HIV/AIDS. She said that her husband, after 22 years of family life has given her 3 children and the stigmatized disease HIV/AIDS. She too used to be employed. Her husband had been very ill for the past 3 years and has died of AIDS. She is sure that she has got the disease only from her husband. She is not very dejected or bothered about her children. She is composed and does not feel worried of the disease for she feels one day everyone is going to die. All her children are grownups. She did not answer the question as to whether she has given her daughter in marriage.*

## 19

An educated woman aged 30 years who had finished schooling and an in-patient of the Tambaram Sanatorium was interviewed on 13-03-2002. She is from Thuraiyur, is married and has two children: one boy and one girl. Her husband aged 35 is educated up to the 10th standard. She lives in a joint family. She is married to her relative. She has both HIV and piles. She said she came to know about HIV/AIDS in Coimbatore. She says HIV/AIDS is a life killing disease. She says she should have got the disease only by injection needles for she has never had done any wrong. Her husband has the bad habit of drinking and he is also a drug addict. She is not aware of the disease like STD or VD.

She says they had taken blood test for her husband. She lives in a rented house paying a rent of Rs.1000/- per month. She says she has not informed about the disease to the neighbours due to fear of social stigma. She says this disease is a self-respect issue so it cannot be taken lightly. She says that to be free from the disease one should use only new needles for injection, while taking blood transmission one should test the blood and then only accept it.

Finally she advises that one should lead a self-controlled and moral life. She had learnt about condoms at the age of 20 from the radio and TV advertisements. She says all the CSWs must be killed. She is of the opinion that the HIV/AIDS affected children must be educated and given admission in all schools. She says the government should first give proper treatment to HIV/AIDS affected patients in hospitals. Free medicine must be given to them without delay.

She suffered from the symptoms of headache and pain in limbs. When she came to know she had HIV/AIDS she was feeling very bad, thinking of who will take care of her children. She says when a HIV/AIDS patient sees the HIV/AIDS advertisement; he would become more and more depressed and reach the dying state. When they are with other HIV/AIDS patients they feel better. Only by thinking of the children she becomes depressed. She has not used indigenous medicine for treating the disease.

*Addition Observation*: Two major observations from this interview were:

- *She says when a HIV/AIDS patient sees the advertisement about HIV/AIDS he becomes more and more depressed reaching the level of wishing for death.*

- *When these patients are with other HIV/AIDS patients their tension is lessened and they feel emotionally better.*

*She is very much worried about her children; for she does not know who would take care of them.*

## 20

On 13-03-2002 we interviewed a 23 year-old-woman who has studied up to the 6th standard and now is an in-patient of the Tambaram Sanatorium. Her husband, aged 26 years, has studied up to the 8th standard. She is a housewife and married just one and half a years ago. She has a 6 month old child. He is a carpenter and they live in Madras in a rented house paying a rent of Rs.800/-p.m. They did not give the name of their native village. She was suffering from cold and cough. The blood test confirmed that she has HIV/AIDS. She says she does not know about HIV/AIDS. She feels better after getting treatment from this hospital.

*Additional Observation: It is unfortunate that she continues to breast-feed the child and has not yet tested the child for HIV/AIDS. She looks lean. Further she has no voice to ask her husband to go for a blood test for HIV/AIDS. Further within in a year of marriage, she has become a HIV/AIDS patient. No one really knows the cause.*

*A major question is: Can one be infected with HIV/AIDS with in one and half a year of marriage? What is the incubation period of HIV/AIDS? Or do women easily get affected even within a very short period? This woman did not look depressed nor does she know about the seriousness of the disease. She was weak and pale, did not know anything about HIV/AIDS. We did not feel like putting more questions to her for it might depress or make her sad when she was already suffering. She fears even to mentally suspect her husband.*

## 21

On 14-03-2002 we met a women in the Tambaram Sanatorium who was 35 years old and was married at 27 years. She has two girl children. Her husband was a M.Com degree holder. She did not give his job description, but she said he died two years after marriage i.e., in the year 1996. However she did not answer the question of how he died or of what disease he died. She hails from Tirunelvenli. The exact place is Vannarapettai in Tirunelveli. Her father aged 70 and educated up to 9th standard is still living. Her mother aged 60 and educated up to 10th standard is also living. Her husband had no habits like smokes or alcohol but he used to be very free with women. She suspects his character, thinks that he had sexual relationship with women before and after marriage. She further says that she should have got the disease from him when she lived with him as his wife for two years. She was first admitted to an hospital in Tirunelveli for treatment. She had chest pain, tiredness, fever, breathlessness and swelling of legs as major symptoms. Only for the past 50 days she in taking treatment in the Tambaram Sanatorium and she feels a little better after treatment.

*Additional Observation: She does not say about the children's health and is not even willing to carry out blood test for them. Further she does not say to us whether the children are living with her parents or with her husband's family. But she gave us her photo and further she gave the real names of her parents. She says to prevent the spread of this disease only doctors should help. She further adds she does not know anything about this disease. She says "AIDS is the disease got by going with women" She lived only in a joint family. Her only wish is that she should be cured.*

### 22

On 14-03-2002 we met a woman from Musiri in Tiruchi. She is an in-patient of Tambaram Sanatorium. She is 32 years old and has studied up to the 5th standard. She worked as a daily wage labourer in arranging stones. She is married and has two male children aged 16 and 12 years. Her husband has studied up to the 10th standard, is aged 42 years and has worked as a temporary staff in a post office.

They belong to the Gounder community. She says she was married at the age of 15. She was often ill and those doctors gave her injections and she felt better. After a stage the doctors advised her to see doctors in Keeraiyur. Only in Keeraiyur the doctors said that her blood is infected. So she was advised to get her blood tested for HIV/AIDS.

She says she does not know about HIV/AIDS. She says she should have got the disease because she often took injections given by the doctors. However she said her husband has given his blood for HIV/AIDS test. When asked "do you think your children will have the disease?" she said "no chance". Her major symptom was fever. When asked whether her husband has any extra marital relationship she says "no". She says to keep oneself away from the disease one should be very careful. Even if one has extra marital sex or sex with CSWs one should use condoms. She says when we go for treatment we should buy new needles and should not accept the needles used by the doctors. She lives in a joint family.

About her faith in God, she says she cannot say precisely. Asked to advise youngsters she said she cannot advice, for youth will only think very cheap of her because she is affected with HIV/AIDS. She says even if one gets HIV/AIDS through injections or blood transmission, people in general think that they should have got the disease only through sex. She says she must be infected only after marriage; that too recently. She has no idea of how to spread the awareness of HIV/AIDS among people. She had come to know about condoms only after she came to know about HIV/AIDS. She welcomes free HIV/AIDS test by the government.

*Additional Observation: She is of the view that even if one gets HIV/AIDS by using infected needles or by blood transmission, the public views them only as persons who are immoral or have acquired it by sex with CSWs or through extramarital affairs. Further she says that if a HIV/AIDS affected person advices the youth to be careful*

*the youth will only think cheaply of her. She is for the government giving free test for HIV/AIDS.*

## 23

We met a woman who had acquired HIV/AIDS after 28 years of married life. She is 45 years old and was married at the young age of 17. She has studied up to the first standard. She is a housewife. She has given birth to four children, 3 male children died just after birth. Her husband has studied up to the 6th standard and is 50 years old. She says that she suffers from cough, lack of appetite, vomiting and several other problems which she did not wish to state. She is affected both by tuberculosis and HIV/AIDS. She does not know how she has got this disease. Her husband has all bad habits; he was addicted to tobacco, alcohol and used to visit CSWs. She says her husband had both premarital and extramarital sex. She feels that her children are not affected by this disease. She says that when she came to know she had HIV/AIDS, she started researching on how she has acquired this disease. She is living in a joint family, in her own house and does not know anything about HIV/AIDS. She says that she feels a little better after taking treatment here.

She says the disease should be cured. More than that she feels it is better to find medicines to prevent it. She is unaware of how the disease spreads. When asked whether one will get HIV/AIDS by touching the HIV/AIDS patients she says, "only you will know". She was unaware of condoms; only after coming to this hospital she has come to know about it. She says, if one speaks about this disease, they think ill of that person. She believes in God and thinks she would be cured by God.

*Additional Observation: She was the first one among those interviewed to talk of preventive injection or medicine for HIV/AIDS. Just like polio drops for polio and vaccination for small pox one should find some preventive medicine for HIV/AIDS. She also proudly says she has studied up to 1st*

*standard. She is the one who says she is doing research on how she has acquired this disease. Though she is uneducated, a villager and a housewife she seems to be very original and creative in her answers. She is not dejected or sad.*

## 24

On 14 -03-2002 we met a woman aged about 40 years from Neidhavoil in Meenjur. She has studied up to 6th standard. Her husband is a 50 years old agricultural coolie. They have no children. She suffered from cold, breathlessness, cough and so she went for blood test on the advice of doctors. For the past one month she is staying in this hospital. She belongs to the Udayar community. When questioned how she got this disease, she asks, "How can I know about how I got this disease?" We could not ask her more on HIV/AIDS disease. When we asked her, "What kind of habits does her husband have?" she says he is alcoholic and he also smokes and visits CSWs. She says "Suppose I questioned him about having sex with other women he used to beat me." She owns a house. She says he had affairs with many women before and after marriage.

*Additional Observation: She is from a lower middle class or poor economic strata. She looked tired, sick and sad. This is an instance of a woman being affected because of her husband's weak character. All the more he was violent with her whenever she spoke of his bad habits. He was brutal. This sort of action cannot be attributed to any cause except his male ego and arrogance. She looked feeble and worn out. She is fortunate that she has no children so she is free from the worries of supporting her children. She is frightened to ask her husband to go for blood test for HIV/AIDS. She is fully unaware of the social stigma this disease holds. She is sick and wants to be cured. It is the duty of the doctors to advice her husband to go for HIV/AIDS test for there is an absolute danger that he would be infecting the CSWs and other women.*

On 13-3-2002 we interviewed a woman aged 40 years who was an in-patient of Tambaram Sanatorium. She belongs to the Mudaliar community. She has studied up to 8th standard. She was married to a coolie aged 45 years who has studied up to SSLC. They are native of Maduravoil. Her mother-in-law was taking care of her. She was very silent before her mother-in-law. Her mother-in-law answered all questions directed towards her.

She (mother-in-law) said 'it has no children', where 'it' refers to her daughter-in-law. The patient was admitted in this hospital for cold, asthma and cough. The patient's mother in law has visited the hospital 3 times. She says she is unaware of HIV/AIDS.

*Additional Observation: This shows that women enslave women. The mother-in-law never allowed the daughter-in-law to speak. Even at the age of 40 the daughter-in-law was literally frightened to talk. So we were only able to question the mother-in-law. In fact we were frightened to put forward questions like "How she got HIV/AIDS? , The bad habits of her (patient's) husband or has he taken up the HIV/AIDS test?" because the mother-in-law of the patient was unapproachable. In fact from the mother-in-law's talk we were able to conclude that her son might have been admitted in the same hospital for she says she has come to the hospital only three times and one time may be to admit her son. The patient said she has stayed only for four days in the hospital and this is the first time her blood was tested.*

**26**

On 14-03-2002 we interviewed a 38 year-old woman from Madhuranthakam. She is uneducated, with two sons. Her first husband died ten years ago. She used to work in a rice mill. Since her husband died when she was just 28 years she says she illegally 'kept' her neighbour and she says that this is the cause for her disease. She gave his name too. Now she is not with him for he is married. She is in this hospital taking treatment for the past 3 months. Her neighbour is 3 years younger to her and has studied up to the 8th standard and is an agricultural labourer. He is Gounder by caste. She is poor, and lives in a hut. She says if a woman has sex with any man she will get this disease. She says people don't say anything about the disease because she is uneducated. She asks, "How will I know how one gets this disease?" She says, "Even the plate and tumbler in which we eat and drink tea are thought to be so unclean and infectious so they don't touch it."

*Additional Observation: She is open and has accepted that she got the disease due to her sexual relationship with her neighbour. She however did not disclose how her husband died ten years ago. Now what worries us more is the possible plight of her neighbour with whom she was involved all these years and above all what is the status of the woman who has just married him. As she is illiterate she does not understand any of the awareness posters on the autos or buses. She did not even say the age of her sons. Further she did not answer the question, "Who is taking care of her?" We are not able to find out when she got married. She is not dejected but is contended and happy. She feels a little sad only when people ill treat her. She seems least bothered about her sons for she never spoke about them.*

### 27

A labourer from an oil-mill hailing from Padappai was interviewed on 14-03-2002. She is 31 years old, has studied up to 6th standard and is now an in-patient of the Tambaram Sanatorium. She is married, but her husband got married again to another woman. She has 3 children: one male and two female. Her marriage was a love marriage.

Her husband is uneducated and is a driver. She says that when he went on trips he used to have sex with CSWs in Andhra, Kerala etc., he used to return home only after 20 days or so. She says she suspected him. She has visited the hospital for the third time. She says before she took treatment here she had taken treatment in the Ramachandra hospital also. She first went for blood test in the Ramachandra hospital and they only said that she has HIV/AIDS. Her father takes care of her. He admitted her to the Ramachandra hospital. She says her mouth is full of blisters; she is not able to eat anything without vomiting. She says she does not know how the disease spreads or how to control this disease. She has come to know about condoms recently by seeing T.V. She says every day hundreds of people come for treatment and some also get admitted. Her main symptom was fever, for three months, she had continuous fever. She said her father had informed the neighbour that she is suffering with tuberculosis only for otherwise they would talk ill of her.

*Additional Observation: She is ill and is very ignorant about this disease. Her father who is 60 years old takes care of her with kindness. We are not able to conclude about the position of the three children; she does not even know that she should test her children for HIV/AIDS. The bigger cruelty is that her husband a driver by profession after infecting her is now married to another woman. Soon she too will became an HIV/AIDS patient. As she is fully ignorant about the disease no further questions regarding awareness program or spread of HIV/AIDS could be asked. It is all the more pitiful that even her father is as ignorant as her regarding HIV/AIDS. They only know this disease is a social stigma.*

## 28

On 12-03-2002 we met a woman from Gudiatham who had just got admitted as in-patient of the Tambaram Sanatorium. She is just 25 years. She, her husband and their one child are infected with HIV/AIDS. She is uneducated and does not know any occupation while her husband has studied up to plus two. He is 40 years old. He is a painter, works in a mill and also knows the job of laying tiles. Her mother is uneducated and her father is no more. Her father-in-law is a retired teacher and mother in-law is a flower vendor. Their native village is Kallapadi. She says that she wept when she first came to know that she was HIV/AIDS infected. Now she says that she has picked up courage. Her village people do not know about this, for if they know, they will chase them from the village. Her family is not affectionate with her and they often speak offensively about her. Her husband's family does not know about the disease. She has two children of which only one child is infected. She has one male and one female child.

She is of the opinion that from the fifth standard onwards sex education should be given to school children, which in her opinion would create awareness and prevent the spread of HIV/AIDS. She feels highly depressed and sad when she sees the HIV/AIDS patients who in her opinion are nearing the grave faster than others. She acknowledges that she has contracted the disease from her husband. Only 12 days ago she came to know that her husband is also HIV/AIDS infected. She says her husband was in love with a very rich lady, who said she wanted to marry him, but since he was already married she was not able to marry him. She says, "might be, my husband would have got the disease from this rich lady." According to this patient the rich woman was supposed to have a 'very bad character'.

She says she had blisters, fever, shivering and she coughed blood. She says her husband rarely earned and wasted all her property. She says he used to sell them whenever she went to her mother's place. Very often he would not leave for work. She further adds that both her parents-in-law have tortured her. Her parents-in-law do not know that their son is affected by HIV/AIDS. She welcomes free HIV/AIDS test by the government. She is of the opinion that people do not speak about this disease because scientists have not yet discovered the medicine. For HIV/AIDS affected children she expects free education, medicine and food.

*Additional Observation: She is depressed and sad. She did not talk much about her children for she did not even answer the question, "Who takes care of your children." She is of the opinion that the disease is caught by her husband from a very rich lady. She talked very vulgarly in anger about that woman for infecting her husband and thereby infecting her and her child. She welcomes free HIV/AIDS test for all by the government. However she do not wish to give any advice to youth. She feels that if she was educated she would have been very careful and would have saved herself from this dreaded disease.*

## 29

We met a lady who was an inpatient of the Tambaram Sanatorium on 08-03-02. She did not say her age, job or native place or her educational qualifications. She said she is married for the past 4 to 5 years and has two children one male and one female, but she did not disclose their age. Asked whether her husband is affected with HIV/AIDS, she says she does not know since he looks well. She tested her blood only in the hospital in her village. She was upset when she came to know that she was affected with HIV/AIDS. She does not know how the disease spreads. She is not worried of her children, does not know the name of her parents-in-law and refused to say her husband's name. She on the whole did not want to talk anything personal.

*Additional Observation: Nobody is taking care of her. She is affected by HIV/AIDS externally. She had blisters in the lips. She had scabies in her legs and hands. She did not divulge her age, however she looked like in her late thirties or early forties and was from the poor strata. Practically speaking, she did not like us asking questions. Somehow or the other she wished to keep her sickness with herself. She*

*did not even say the symptoms she is suffering from that had lead her to do the blood test.*

### 30

We met a woman aged 37, on 08-03-2002, who was uneducated and had never stepped into school. She is an in-patient of the Tambaram Sanatorium. Her husband who is deceased was also uneducated. He died at the age of 40. They belong to the Gounder community. Their native place is Chengal in Karur. She has two sons. The eldest son aged 18 years is working in a shop. The second son aged 15 years is studying in the 8th standard. She used to work as an agricultural coolie. Their basic occupation is agriculture. She has only one younger brother who also is an agricultural coolie. She was married at the age of 15; the marriage was an arranged marriage. She believes in Perumal.

She was affected with continuous fever. She went to Madurai for treatment, people in that hospital said that they would not treat her, but advised her to go to the Tambaram Sanatorium. When asked whether her husband had this disease she says no. Only she is affected with it, her husband used to drink cheap liquor everyday and was very sick when he died. What he earned would only be spent on drinks and he never used to give any money to the family. She had educated the children only with her earnings. She has one acre of land, which is her sole property. Till now she does not know how she got this disease.

*Additional Observation: The cause of her husband's death is not known but from the symptoms, which she elaborately described, we can conclude that he should have died of HIV/AIDS. Her husband must have infected her for she says that she never had any extra marital affair. She however did not give the name of her parents. She spoke about her only brother. She is unaware of how she has contracted the disease and she is not in a position to suspect her husband's extramarital affairs. She had been an earning member of*

*the family and with her earnings she had educated both her sons. Further the money her husband earned as a daily wager (as an agricultural coolie) was spent on alcohol. She is from a upper middle class background and does not suffer in poverty. This shows her commitment towards the family and her husband's irresponsibility. She is totally ignorant about HIV/AIDS.*

## 31

We interviewed a Hindu Udayar woman hailing from Cuddalore. She is 35 years old and is uneducated. Her husband who had studied up to 12th standard is deceased. Both her parents-in-laws are still living. She has one 18-year old daughter who has studied up to the 10th standard. She has not yet married her off. This woman was married when she was just 14 years. Her marriage was a love marriage. He was living in the neighbouring, next street and they used to visit each other and were married with the support of her family. However the husband's family was not willing for the marriage. Hers was a nuclear family. Her parents paid good dowry for this marriage though it was a love marriage.

She lives only in a hut. Her husband was working as a cycle mechanic. They were able to make both ends meet with difficulty. She has some hereditary property of about 1 acre. She says even before marriage they had sex. She says her husband was very sick for a year with fever, headache and inability to move, she did not take care of him. So he died after a year of suffering without any medical attention. Since she was feeling very ill doctors advised her to go for a blood test and they said she is infected with HIV/AIDS.

The doctors said to her "if your husband has the disease then you would have it." When she heard that she had HIV/AIDS she felt sad for she was worried about; who is going to take care of her children. She says the disease spreads by using the same needle and from sexual intercourse. She now says that her husband was tested in Salem and the doctors in Salem said he is infected with HIV/AIDS. She says her husband had sex with several women, whoever came to work he would have relationship with them. She said he used to consume lot of liquor and go with women, if she ever dared to question him he would beat her.

Further he was also suffering from ulcers and scabies in the private parts. She says she too was affected in the same way and in spite of medical treatment it did not heal. She further says she suffers from breathlessness. She says in villages no one knows about condoms and her husband has never used it.

She says she never had any sexual relation with any other man. Here, the doctors are kind and consoling. She says all affected persons with HIV/AIDS, all over Tamil Nadu and nearby states like Andhra Pradesh come here (Tambaram Sanatorium) for treatment. When we asked "Can the government give jobs to persons like you?" She says, "They will not give me any work or job, as I am HIV/AIDS infected." When we asked, "What can be done to small children, who are HIV/AIDS infected?" She says, "As they are very small nothing can be done." She expects nothing from the government except cure. She says she feels better after taking treatment in this hospital. She says this is an 'indecent disease' so people do not talk about it openly. When she sees other HIV/AIDS patients she feels sad because God has created all these diseases.

She is unaware of the difference between HIV and AIDS or the full name of the disease. She says she has seen advertisements about HIV/AIDS in T.V. She feels more advertisements can be given to rural uneducated people. She says her relatives look down upon her, and don't visit her.

*Additional Observation: She is a victim of child marriage, she was married at 14 and widowed at 34. Her husband who used to beat her badly when she questioned him of his bad activities infected her. She says at first that he did not say about his disease and died untreated after suffering for a year, but later she says she came to know her husband had HIV/AIDS when his blood was tested in the Salem*

*hospital. From her talk it is clear that both she and her husband suffered from STD. Her old parents-in-law, aged over 70, are still living. Only bad habits had driven her husband to HIV/AIDS and finally to death. It is very depressing to hear that even after so much of scientific advancement a man suffers with a disease for a year and then dies unattended without any medical aid.*

## 32

A woman from Renigunda, Tirupathi was interviewed on 07-03-02. She is uneducated and did not wish to say her age. Her husband was also taking treatment in the same Tambaram Sanatorium but was in a different ward. She has a son and she says he is not infected by HIV/AIDS. She has not tested him yet she says he is not infected by HIV/AIDS. She does not disclose his age also. Her husband does the job of weaving baskets from bamboo. She says she is affected because of her husband. She does not say how her husband was infected. Asked how she came all the way from Tirupathi to take treatment here she did not answer. She belongs to a nuclear family. Her husband has studied up to the 10th standard. She says that she has not taken any indigenous treatment for this disease.

*Additional Observation: She looked very thin, had scabies on her hands and legs. She did not say the age of her son or his schooling or his work. She only says that she is infected because of her husband. She does not tell how she was first tested for HIV/AIDS. She does not disclose the symptoms she suffers or suffered from. But she says life feels better after getting admitted in this hospital. She seems fully free of any tension or worry. She does not say about his bad habits. Her teeth reveal that she was a continuous tobacco chewer for it looked brown, but her lips were bleeding for which she was applying some ointment. She looked happy and comforted.*

On 07-03-2002 we met a native of Tamil Nadu settled in Andhra Pradesh and now taking treatment in the Tambaram Sanatorium, she is about 40 years old. She is taking treatment for the past two years. She is a housewife and her husband is a coolie. She is a Christian. She was born in Perunthurai and has three children: two daughters and one son. She says her husband infected her. She hails from a joint family and lives in her own house. She says she was infected by HIV/AIDS only after her marriage. Her husband too is living with HIV/AIDS. She prays to Jesus Christ. Her husband has all addictive habits. She suffered from headache, breathlessness, loss of appetite, stomach pain, ulcers, scabies and loose motion. She says she has taken treatment in Chitoor, Chennai, Saligrahmam and Tambaram. Doctors in T.Nagar did not admit her in the hospital. She says now she has taken treatment in Tirupathur, CMC and Apollo. She says she has spent up to Rs.50, 000 for treatment. Each time at least she has to spend Rs.2,000/-. She says social service groups do not help her.

It is her own money that she spends on treatment. She says she owned a flower shop and had earned a lot of money and also has saved some money. She ate well, before she had HIV/AIDS. She says she has lost her weight, from 50 kgs she became 42 kgs in a short span of 3 months. She used to earn from Rs.40 to Rs.100 a day when she was running the flower shop. Four years have past since she had stopped working as a flower seller. She has never heard about HIV/AIDS.

She says she first took treatment form doctor Shah Jahan in the Chitoor bus stand. She says only she had asked the doctor to perform HIV/AIDS test because she suspected that she would have HIV/AIDS. She said that once again to confirm this disease since she felt the doctor was speaking lies she paid Rs.1200/- and took up another blood test in Thirupathi. Doctors felt sad and she was confirmed to have HIV/AIDS. She says doctors here look after her nicely. She says she has never taken any indigenous medicine. She also feels sad that nurses and ayahs do not behave well. She says she has tipped them ranging from Rs.20 to Rs.200.

She first informed about the disease to her daughter; then to her husband. She advises the youth to perform blood test before marriage and not commit the mistake of having sex with strangers or CSWs.

*Additional Observation: She is not very sad in spite of spending over Rs.50,000 for treatment. She did not speak anything about her husband expect that he is the cause for her to get the disease. She is also calm and composed. She did not speak much about her children. She did not answer where her husband took treatment. She has not performed HIV/AIDS test on her children, she does not even think it relevant. She is the first case reported in our interviews, who has taken HIV/AIDS test second time from a different place to confirm that she has HIV/AIDS. Might be her monetary position has made her self sufficient and self-reliant. She does not think HIV/AIDS to be a social stigma. She is at one place a little hurt when she was not admitted in the T-Nagar hospital. Both private and government doctors must be advised to respect their profession and thereby respect the patients whatever be the disease they suffer from; then only the social stigma associated with this disease can be wiped out.*

## 34

On 05-03-2002 we met a plus two educated women. She was 27 years old. She is a tailor by profession. The total interview was answered by her brother-in-law i.e., her husband's younger brother as she was not in a position to talk. She is from Karaikuruchi, Puddur of Nammakkal district. Both her parents are alive. Both her in laws are alive. She looked sick, shocked and her brother-in-law said that she did not talk. Her husband died 6 months back. Since her husband had AIDS and he died of it, she asked her brother-in-law to admit her in the hospital and carry out blood test for HIV/AIDS. He says nothing was wrong with her, she was feeling all right. Her husband died when he was 34.

*Additional Observation: They are from a middle class family. She was married at 25 years and widowed at 26 and half a years, she lived with her husband only for one and half a years. Only a week ago her brother-in-law confided to her that her husband had died of HIV/AIDS so only she was very sad and dumb-struck and wanted to take a test. No one in her village knows about it. They have not even informed their parents. As her brother-in-law answered we could not ask about the habits like smoking and alcoholism, which his brother possessed.*

### 35

We interviewed on 05-03-2002 a young girl, 16 years old, an inpatient of the hospital. As she was in a very serious state only her mother answered all the questions on her behalf. The sixteen year old was married when she was just 11 years to a 32 year old man as a second wife since his first wife died; they say of some strange disease. She was married even before she attained her puberty. She lived with him only for a month and he died. Now this girl is HIV/AIDS infected due to him. He died of swelling of legs and several other complications, which is prevalent among HIV/AIDS patients who do not take treatment or come for treatment even at the last stages. Her mother says she was fully aware of the character and bad habits of her son-in-law. In the Indian society it often happens that the son-in-law is older than the mother-in-law. For the mother-in-law is just 28 years while her son-in-law was 32 years. She says poverty was the reason to perform child marriage, she has five children all of them younger to this girl.

*Additional Observation: Her mother does not show any expression of sadness to see her child suffer. This girl was widowed at 11 years even before attaining puberty and was seen to observe strict widowhood though she has lived with*

*him only for one month. This moved us a lot. But as we are in a public place, a hospital we could not express any of our feelings. Her mother is kind with the child but we do not see any concern on her part. She takes it lightly to see a small girl who was just 11 years old being used sexually by a 32-year-old man. We wonder if she ever thinks of the child's body. When the child has not reached puberty what will be her state of mind when a 32-year-old has sex with her. Even in appearance she is thin, very short and does not even look like sixteen. We were not able to comprehend how her mother permitted this in the first place. In the second place we wonder how the mother is treating such a small child as a widow, when the married life lasted only for a month. Does any one ever think of giving this child counseling or consolation or at least make the child speak out its fear and inhibited feelings? Now she suffers from the disease. I actually wept on seeing her. I feel in the 21st century such barbaric actions must never take place. A girl without reaching her puberty is married and widowed in a month and now after 5 years is suffering with full blown AIDS because of that one month of marriage. Which women organization will take up these issues? They say no child marriages but it takes place in practice and child widowhood is observed, who is going to take care of her? What is her life? What is her future?*

## 36

This is an exceptional case. This woman never talked about her husband or even about his job, his native place or even the answer to the question of whether he is living or dead. She is a Hindu Gounder. She is 25 year old and has studied up to 6th standard. She was married 7 years ago that is at the age of 18. She has two children: one male and one female. She does not know anything about HIV/AIDS. She is from Paramathrveli in Neyveli. She has told every body about the disease. She suffered the symptoms of headache, fever and inability to walk.

*Additional Observation: She does not say anything about her husband. We are not even in a position to know whether he is living or dead. She did not say how she got it. She never answered any question pertaining to him. We felt very uneasy, for whenever any question was put to her regarding her husband her answer was a sigh or a silence of sadness, she never even wishes to recollect memories of him. How she got this disease is not known. Who is taking care of her is not known. Who takes care of her children is not known. She did not look sad or depressed, whereas she was composed. She was not able to walk, she could just sit. She looked very ill.*

## 37

On 05-03-2002 we met a Vanniar lady educated up to 5th standard in the Tambaram Sanatorium. She is 25 years old, hails from a remote village Chibichakkaram of Namakkal district. She was married when she was just 19 years. She gave the address of her parents. She is being looked after by a social service volunteer who is an old lady. Her husband is 30 years old and he is uneducated. Her marriage was an arranged one. She also accepted and liked the arranged marriage. Her parents-in-law are alive. Everyone knows about her disease. She lives in a hut and he is from a poor family.

She said that for a year she did not have any child, and then a child was born in the hospital. Only in the hospital they tested her and she was found to have HIV/AIDS. When she informed this to her husband he ran away to Bangalore and never returned. He also said that he was not affected by the disease only she is infected. Both her husband and her mother-in-law beat her badly; she said this with tears. So she went to her mother's place. Her husband did not want to live with her so he has obtained divorce from her. It is 2 years since he got her divorced. She says two children were born to her and they died soon after birth. When we asked her why she suspected her husband she says, often he was sick, with loose motion, loss of appetite and general weakness. She suffers from headache, inability to walk and loose motion. She says she has only one brother, he ran away from home with other boys. Her people take care of her. Asked if she knew about HIV/AIDS, she says surely she would have been careful and would not have been a victim if she had been made aware about HIV/AIDS earlier. She says sex education can be given at the age of 10. She says even after they knew they had AIDS she had lived with her husband for 4 years. She says she cannot question him for he would not listen to her. He was so characterless that even he used to have sex with his brother's wife.

She says CSWs cannot be stopped from practicing their professions because according to her police have relationship with them. How can they be changed? She says as a woman she has come to a low state because of her marriage. She at times felt like committing suicide. When asked about her state of mind; she said that each one would be affected in their own way. She says some think badly of her. She feels very sad about the disease. She says now they are trying to marry her divorced husband to another lady. She says he and his mother are least bothered about anyone being affected by this disease.

She advocates compulsory HIV/AIDS test before marriage. She says she has so far not taken indigenous medicine for treating her disease. She is not able to eat well. Further she is contended with the advertisement and awareness given by the government in recent days. She mentioned that her husband has given back the entire dowry taken by him.

*Additional Observation: This is the worst case we have observed where she is very much dejected and wanted to commit suicide. Several things are to be observed. It is highly questionable and unethical on the part of her husband and his family to beat her for being infected with HIV/AIDS. Secondly it is unfortunate to note that when men are infected, their wives take good care of them but when the wife is infected the husband obtains divorce. This has caused this woman to be totally wrecked. She is dejected,*

*depressed and very sad. However her only consolation is that her mother and father are very kind and considerate with her.*

## 38

We interviewed a lady from Kurujipadi in Cuddalore district on 10-02-02, now an in-patient of the Tambaram Sanatorium. She is 30 years old and has studied up to 6th standard. Her husband has studied up to 10th standard, she says that he knows all jobs. He is 40 years old. They have three children, one girl and two boys. She was married at sixteen. She is in the hospital for the past one and half months. She says she is affected by HIV/AIDS. She says one gets HIV/AIDS by using same injection needles, transmission of blood and sex with strangers affected by HIV/AIDS.

She says that when she came to know that she had HIV/AIDS she thought she has to die. She further adds this is a very bad disease and no one should be affected by it. Later she accepts that she should have got the disease from her husband. When asked about his habits she said he is a smoker but whether he visits CSWs is not known to her. She says only now she is infected and that no one else knows about this. She also says that if her neighbours come to know about it they cannot live in that locality because of the ill treatment and the torture that would be meted out to them.

She says more than the public, the people with whom she is staying now, think badly about her. She says the fluid oozing from the blisters in her body is also infectious for it will spread HIV/AIDS. To prevent this disease the only way is the discovery of medicine. Nowadays doctors come and visit her weekly once and go away. She says she suffers from fever, headache, itching and pain in her limbs as the main symptom.

*Additional Observation: She did not answer the question whether her husband was tested for HIV/AIDS. Women are*

*frightened even to suggest their husbands to take up the HIV/AIDS test. She has not tested her three children, all of them are below 13 years. She could not ask for it. Her body was full of blisters. These looked like blisters that came by pouring hot water. All these blisters had full water in it. No treatment was given for it. It looked as if they are frightened to treat it. She says the fluid from these blisters is so infectious that if children touch this fluid they would get HIV/AIDS. She is not well attended by the doctors. Whether her husband takes care of her or who visits her or who attends on her is not known, for she did not speak anything about her parents or other relatives. She never even answered whether they are living or dead. Nor did she talk about her husband. She looked very sick and was suffering from pain because of these blisters in several places. It looks pathetic for no medicine is even applied externally for her blisters.*

### 39

On 10-02-2002 we met a woman aged 28 years, an in-patient of the hospital. She was married at the early age of 17 and widowed at 27 years. She has two children, one boy and one girl. She has studied up to 7th standard. She is from Puliathur, Tindivanam. Her father has studied up to 8th standard. She says she had come to this hospital for treatment as people in her place said one could get the best treatment for HIV/AIDS only from this hospital. She also adds that she had first taken treatment in Kodambakkam. In that hospital they asked her to bring her children and also test them. Her husband died some 9 months ago. He died of HIV/AIDS. She says she is not aware of how her husband got this disease, but she says that the doctors told her that her husband had sex with CSWs and other women even before marriage and that is why he was so badly affected by HIV/AIDS and died. She says that she must have got this disease from her husband. She does not know anything about HIV/AIDS. Her marriage was a love marriage. She says she is unable to eat and suffers from fever, cough and cold. She suffers often from shivering, fever and giddiness.

She adds, "To control HIV/AIDS, we should be self-controlled." Her only advice to the youth is "lead a self-controlled life". She says her husband very often used to say that he would live only for two months or so and that he is going to die. Her only sorrow and tension is "how to protect these children who are just below ten years? Who is going to bring them up?" She belongs to Mudaliar caste and is from a middle class family. She says some doctors are kind, while some fear to come near the patients, they just stand at a distance and look at them.

*Additional Observation: This is an instance of tragedy. Her husband who has not lived a controlled moral life has ruined her life. She is just 28 years old. She does not say whether her children are also affected with HIV/AIDS. It looks as though even the parents of this lady are not coming forward to take care of her. She is all alone in the hospital. Her parents-in-law too are not visiting her. Her husband had not disclosed that he suffered from HIV/AIDS; neither did he ask her to take the HIV/AIDS test, but often he said he would die soon in a month or two. She could not guess that he was suffering with HIV/AIDS. The doctors too who treated her husband had not advised her to go for HIV/AIDS test. Now the same story follows. She has no one to advise her to carry out HIV/AIDS test for her two children. The doctors told that if her husband had sex even before marriage there is every chance that these children were infected with HIV/AIDS.*

### 40

We interviewed on 11-04-2002 an uneducated Chettiar woman 35-years-old, now in Tambaram Sanatorium as an in-patient. Her parents and her husband are deceased. He had studied up to 9th standard. She is in this hospital for the past four years. She has taken treatment earlier from a hospital in Dharampuri and has spent over Rs.70, 000. They promised to cure her but did not cure her. She lost weight very rapidly from 85 kgs to 40 kgs, in a short span of three months. She has two children, one is married and the other is studying some computer course in Tiruppur.

She says that after her husband died the whole village knew that their family is infected with HIV/AIDS. Even her close relatives like bothers and sisters have advised her not to visit them. She says they send money to her by money order. They are very much frightened that their children will also get this disease. Her husband was an agriculture labourer in the beginning but as there was acute failure of agriculture and a very poor yield he took up driving as a profession. Fearing the ill-treatment that could be meted out by her relatives she said she has joined her son in the hostel.

She says her husband infected her. She further adds that her husband was a good person who used to take care of the children and family very well and used to bring home from Rs.6, 000 to Rs.7, 000 per trip; but he must have been infected by the CSWs. She says for the past 3 years he has been taking treatment from several doctors and also taking lot of medicines. She says he had hidden the disease from her. One day he was very ill as if in dying state, so she not able to take him to the doctor, instead she, went to the doctor who was giving her husband treatment. Only that doctor told her that her husband was suffering from HIV/AIDS. Since she was very much aware of the disease HIV/AIDS and all the more aware of the social stigma it carried with it; she came home and asked her husband what was his disease. She only got the answer from him that he did not suffer from HIV/AIDS. Suspecting both the doctor and her husband she carried out a blood test for her husband with the help of her brothers and it was confirmed that he had HIV/AIDS. As he was not able to bear the shame he ran away from the house. Soon after a year she came to know that he had died in a nearby village and thus everyone in village came to know that he has died of HIV/AIDS.

She says to control this disease only men i.e., husbands should be careful; if they wish; then alone the disease can be prevented. To have pleasure for some minutes the men forget their family, wife, children and above all the family prestige. "Can men do this?" She advises youth to be honest, truthful and not to spoil other families. She feels that those affected with HIV/AIDS should be killed by an injection for they hold a larger stigma in the society. She says "even after their death the stigma on the family never subsides. We are not in a position to live in the same village, we are not in a position to even buy a cup of tea from public tea shop or provisions. The only solution for such families is to leave the village." She describes the social stigma her family faces even today among relatives and among the villagers with whom they live.

She wants to beat with a slipper the doctors who had cheated her by promising to cure her and taken a lump sum of over Rs.70, 000. She says that after her husband's death she was not even able to buy provisions from any shop in the village. She further says doctors stand a mile apart from her and look at her. They don't come near to check or even hear her health problems. She feels a little better now. She feels if she complains about the ill-treatment meted out to her, would the government take any step to punish the guilty even after a year. She says only if party workers say something they may take some action about it. She further philosophizes that we pray to God not for cure for God cannot cure and he has no capacity to do it, we pray only for our mental satisfaction and peace. She is a Kannada Chettiar. She says HIV/AIDS is a disease in which germs get into the blood and spoil the blood. She says she is an orphan in spite of being born with brothers and sisters because she is affected by this disease.

She says that CSWs must be burnt alive. Some reside near her house; they come in cars and taxies and live in pomp and prosperity. She says because of them only she is now so much affected, suffers the social stigma and disrespected by one and all.

***Additional Observation****: Though she is just 35 years she had lots of information to share. First and foremost she is well aware of the disease and the social stigma it holds. She* *says because of this disease her children have become orphans and she is not in a position to even fondle her grandchildren. She is sick and depressed not because of the disease but because of the ill treatment she faced from her kith and kin and the villagers. She says a prostitute lives with respect, her clients come in posh and big cars but because of them she suffers the social stigma. She observes, "Unless men correct themselves it is impossible to prevent this disease."*

*She is drastic for she says kill the person who gets HIV/AIDS so that the family is saved from the social stigma and harassment. She is of the opinion that the government does not take proper action to punish people who ill-treat HIV/AIDS affected and their close relatives. She is rich and has spent lots of money on a private doctor who promised to cure her. It is unfortunate to see the husband not disclosing his disease to his wife and family and running out of the home when his illness was discovered. The only solution is that, if there was a support group they could have counseled them and their relatives, so that such types of discrimination do not take place. Further this man after a month was brought dead from a neighbouring village; by some persons known to the family. This has a double effect on the affected individual and the stigmatized family.*

### 41

An agricultural labourer from Thiruvannamalai was interviewed on 05-04-2002. She belongs to the Gounder community, is 24 years with no education. Her husband aged 28 years has studied up to 7th standard. He too is an agricultural labourer. She was married 2 years ago and she had a child, which died. She is comfortably well off, she owns 3 acres of land. She says her marriage was a love marriage. Even before her marriage she had sex with him. She said she has never used condoms. She feels she was infected by her husband. For she says her husband had some relationship with another lady even before marriage. She says she worships Ezhumalaiyan. She says she used to pray to god only in the festive season. She says after she got this disease she is not even praying to god. The doctors take good care of her and they give her good food. When we asked her whether she had sex with other people she said she has never had sex with anyone except her husband. She was first checked for HIV/AIDS in Aadaiyur. They said she had HIV/AIDS. She felt sad and she knows how this disease spreads. She says that the  government must give free medical test for one and all. She says sex education can be given to children from 10th standard onwards. She says she has not seen any HIV/AIDS advertisement but advocates that giving such advertisements are good. In villages also it is important to give advertisements about awareness of HIV/AIDS. She expects only medicine from the government. She says people in general think very badly about this disease.

*Additional Observation: She is composed and does not feel depressed or dejected. However she did not speak about her husband being tested for HIV/AIDS. She feels it is impossible to advice her husband to go for HIV/AIDS test. It is unfortunate to note and impossible to comprehend that with married life of two years she has become a victim of HIV/AIDS. She is proud to state that she loved him and got married. It is clear from her statements that her husband is not visiting her. Even her parents do not visit her. We are not able to know what is her stand on this disease. Her husband's whereabouts is not known. He does not support her materially or physically. He is least bothered to visit her. As it is a love marriage she has no voice with her husband's family or with her parents. She says she is composed but does not like to talk about her husband. She gave only the name of her parents-in-law but did not give the name or address of her parents.*

**42**

We met a woman about 45 years ailing with HIV/AIDS on 05-04-2002 in the Tambaram Sanatorium. She is a native of Maranelli village in the Karaikudi district. She is a Hindu. She is uneducated, and used to earn by selling sweet potatoes and other petty eatables. Her husband is also uneducated and he is a coolie. She has three children, of which two are no more, and one is living. Her son is studying 10th standard. She is poor; she says government has given them some house. She was married when she was 25 years; she had 20 years of married life. She says her husband used to drink once in a while and he used to beat her but now he is kind and does not beat her. Doctors in the hospital are kind with her. She says that if persons with HIV/AIDS are able to do work, they can be given employment by the government. She says the HIV/AIDS affected must be taken care of. She is also for the government giving free HIV/AIDS test for one and all. She says she does not know how she got the disease.

When we asked her whether she had relationship with any other man she said yes. She had relationship with other men when she was 20 years. When we asked her how long, she said for a year she was associated with some other man. She further adds that her husband does not know about this. She says her husband is well, however they have given his blood for the HIV/AIDS test. They have not yet got the result, they are awaiting the same. She said she felt sad when she came to know that she was affected by HIV/AIDS. However she says this disease holds with it a strong social stigma. She says all her relatives stand at a distance from her. She is unaware of how the disease spreads.

She says God only questions of the sins everybody have committed. She has now changed her faith, for she says she believes in Jesus Christ. She prays twice a day and goes once a week to the church. She says people cannot speak out about this disease, because it is a social disgrace.

*Additional Observation: She was the first interviewed woman to advice men not to believe all women and get married only after blood test. She is herself a person who had premarital sex. She says openly that from the age 20 to*

*25 she has led a careless life. She never blames her husband as other women did. She is being taken care of by him. This happens to be a very rare case.*

## 43

This 25 year lady was married when she was just 21. She has studied up to 8th standard. She belongs to the Kallar caste. Her husband is uneducated, an agricultural labourer by profession. Her marriage was an arranged one. She does not know how she has acquired the disease. She is from Kallaperambur in Tanjore. She too did not disclose her parents' name. She is at home and does not go for any job. She had paid some dowry to get married. She says her husband does not go for any job, only her father is supporting the family. Her husband is not yet tested. Further she says doctors don't see them at all, only the nurses come and give medicine and injection. She says she does not know how the disease spreads. She did not want to say about the habits of her husband. She says he does not go for work and he stays at home and takes medicine regularly. She does not even answer the question about his health status. She has no children, this is revealed from her answers.

*Additional Observation: She cannot even question her husband when he does not go for work. It is unfortunate that for the past two years her father is supporting the family. This girl says her husband does not go for a job, might be her husband is sick and weak so he is not in a position to work. She did not answer the question whether her husband is healthy. He has not taken the HIV/AIDS test. Her father-in-law is very old and she says her husband had gone to Mumbai on a pleasure trip and he has told her how he enjoyed life with friends. She feels from then he had become inflicted with some sexually transmitted disease. She is least bothered about her disease. She is not bothered about her husband. She looked composed, might be she is*

*internally upset and dejected. One thing we could see is that she is not very happy with her marriage and her husband.*

## 44

We interviewed a women from Kuppasamuthiram in Vellore, now an in-patient of Tambaram Sanatorium on 08-03-2002. She is 35 years old, uneducated and Hindu Gounder by caste. Her husband had studied up to 8th standard, and died at the age of 38. She has two female children aged 14 years and 8 years studying in 9th and 3rd standard respectively. She works in a beedi factory. Her mother-in-law is an agricultural coolie.

She was married at the early age of 16. Her marriage was an arranged one. She says her husband was always trying to run away from her. She says only after she gave birth to the second child she became very ill, she had fever so she was admitted in the CMC hospital. Her brother-in-law (that is her husband's brother) joined her in this hospital, that too after the death of her husband. She says she is taking treatment in the Tambaram Sanatorium for the past 3 years.

Till his death she did not know about the fact that he had HIV/AIDS. Then only her brother-in-law informed her that her husband died of HIV/AIDS so there is all probability that she too is infected by it. Her husband has confided his sex with other women to her. She says she never had any relation with any man other than her husband.

She says her husband had lots of scabies and ulcers in private parts and when she asked him about it, he said it was due to heat and dried fish that he ate. It is 7 years since his death. She was taking treatment for asthma, only now she knew she had HIV/AIDS and so she joined as an in-patient. She says she does not know much about HIV/AIDS but, has seen advertisements on T.V. She further knows that there is no medicine for HIV/AIDS. She has no other property, she lives in her own hut and is from a nuclear family. She has not checked her first daughter for HIV/AIDS as she is free from fever and cold but when she checked her second daughter it is said that she is infected with HIV/AIDS and the second child is taking treatment as an outpatient. She says one cannot disclose about HIV/AIDS to anyone because the stigma it carries. She says this disease has spoiled her whole family.

She says she prays on Saturdays, on Pongal and during other festival days but after she is infected by HIV/AIDS she does not pray at all. The god did not take care of her children and her husband so she says she has no faith and it is of no use praying to him when he does not do anything good to her family. She further says she is not even able to sit or walk. For HIV/AIDS affected children she says government can give them free food, free education, free medical treatment and free medicine. She feels that government can carryout free medical test for one and all. She says the disease spreads because of sex relationships with men and women. She says no one is suffering from this disease in her village.

She argues that because her husband had gone to wrong places he is infected and the whole family is ruined. She says only her sisters are kind and considerate but others are cold and ill-treat her. The doctors are good and take care of her in a nice way. She first thought she was the only one affected with this disease, but after admitting in the hospital she found that many are infected with HIV/AIDS. She is uneducated and does not even know to read Tamil. She says her husband used to sell her jewellery and drink everyday. She has no property. Her house is built in poromboke land. She says she has reduced to one fourth of her weight.

*Additional Observation: She is a victim of child marriage, who lost her husband at 28 years. Her husband had hidden the disease from her. Only after his death she came to know from her husbands brother that her husband was affected by HIV/AIDS. Her second child is also affected by HIV/AIDS. Both her parents-in-law are dead. She feels that after 4 months of married life, he started to have relationships with other women and tried to elope with them and so on. With*

*the birth of a second child she became a HIV/AIDS victim and he died after a year of the birth of the second child. She is from a very poor family. She says that in spite of praying to God he has not protected her family so she says she would not worship him.*

*Further she is much worried because the disease has affected the whole family and the stigma attached to this disease is totally ruining her family.*

**45**

A 25-year-old woman was interviewed on 05-04-2002. She has studied up to 7th standard. She was married at the young age of 15. She is from Tiruchenkodu District in Thillai Nagar. They used to work in the forest for living. Her marriage was an arranged one. Her husband is 27 years old. They have two children: one male and one female. She said she is suffering from HIV/AIDS, but she does not know about HIV/AIDS. She says she does not know how she got this disease or how one gets HIV/AIDS.

She says her husband is alcoholic, addicted to smoking and had sex with other ladies. She says that at present her husband is in Mumbai. She feels that if neighbours come to know of the disease they will think badly about her family. She says that she has not seen or heard HIV/AIDS advertisements. One of her children, is also infected with HIV.

*Additional Observation: She is yet to test her husband for HIV/AIDS. She is unaware of how she got this disease. She says he is currently in Mumbai. She did not answer the question of why her child was tested for HIV/AIDS, before she was tested. There are high probabilities that her husband was infected in Mumbai. She is unaware of how she is infected. She does not even know about the activities of her husband. She does not blame her husband. She does not believe that prayer can cure her. However she is yet to know that this disease has no cure.*

On 05-04-2002 we interviewed a woman aged 22 years, she has studied up to 5th standard. She is married and has a child. Her husband was a driver. He is deceased. She is from Tiruchi. Her marriage was arranged. Her father had two wives. She works as a domestic help to earn a living. HIV/AIDS test was first carried out for her child who was continuously suffering from vomiting, loose motion and so on. She is taking treatment for the past three years. When the doctors found her child to be affected with HIV/AIDS they said that she must undergo the blood test for HIV/AIDS, and it was confirmed that she has it. She says it is her destiny that she has got this disease.

She says she does not know how her husband died. In Tiruchi she is paying a rent of Rs.150/- and living in a hut. She is from very poor strata of the society. She is fully unaware about the disease. She says she does not know if HIV/AIDS can be prevented; but however she says people think very cheaply of the persons who are affected with HIV/AIDS. She says she has so far not taken any indigenous treatment. She belongs to the Naidu community. She thinks if God wishes she could certainly be cured. She says when she saw advertisements about HIV/AIDS she felt sad. She says she does not know about condoms.

*Additional Observation: Her husband who is a truck driver has died of an unknown illness. She says she did not know what exactly caused his death. He was sick and for some time he did not go for trips. This is a unique case she was not aware of the HIV/AIDS and it was her sick child which made the doctors to advice her to go for HIV/AIDS test. When the child was affected by HIV/AIDS, doctors asked this woman to take the HIV/AIDS test and it was found that she too was infected with HIV/AIDS. She is not affected by any illness and she is healthy; only the child is very sick. She is unaware of the fact that her husband must have died of HIV/AIDS or that he has infected her and consequently*

*she infected their child (MTCT – Mother-To-Child-Transmission).*

## 47

We interviewed a 35-year-old Gounder woman from Chendamangalam. She has studied up to 7th standard. Her husband who was a lorry driver, died two months ago. She has only one child who is studying in 11th standard. She lives only in a rented house. She pays a rent of Rs.250/- p.m. She says her husband's blood was tested and it was found that he was HIV positive, so only she came here for her husband's treatment. She was at first a little sick, and then she became very sick, so she did her blood test in Namakkal. She was feeling very sad when she came to know about the disease.

Now she says she feels a little better. First she had scabies in the private parts and she says she spent over Rs.5000/- for the treatment. She says she has also taken indigenous medicine for cure. Some said it is better if she goes for treatment to the Tambaram Sanatorium. She was not cured by any medicine from any other doctor. So only she got admitted in the Tambaram Sanatorium. She says she was his only wife and he was keeping her well. She says he was always kind with her and he died. She says she prayed to all the Gods. She says now she does not pray to any God as Gods have not taken any care of her but have left her in sorrow. She says even if she spends over several lakhs she cannot be cured. She stopped going to temples for the past one and half years. She says she has no T.V. at home because she wanted her son to study well in plus one (11th standard). She says that her son knows that she is affected by HIV/AIDS. He is kind to her. She says most of the friends are kind to her. She says she is unaware of how HIV/AIDS spreads. She says she has tested her son and her son is free from this disease.

***Additional Observation:*** *She is very devoted to her son's study that she denying herself even the TV because she*

*wishes her son must study well. Her husband who was a lorry driver had become infected with HIV/AIDS and he died, after infecting her. This woman after spending lot of money has at last come to the Tambaram Sanatorium for treatment. She is different from many other women because she does not blame or find fault with her husband for infecting her. She is also suffering from STD. She is taken care of and loved by all friends and neighbours even after she was affected by HIV/AIDS. She has renounced the worship of gods for according to her, the Gods have not been kind or considerate. She has spent lot of money for treatment of STD and HIV/AIDS.*

## 48

We interviewed a 26-year-old woman who studied up to the 10th standard. Her husband is a driver by profession and is uneducated. She has only one child. Her parents are living, they are uneducated, they are framers. Her husband's family is also uneducated. She says she got the disease from her husband. She has not tested the child for HIV/AIDS.

*Additional Observation: She did not answer the question whether her husband is affected or not but she emphatically says that she has got the disease from her husband. He is a driver by profession. She did not give her native place or even the place of her living. She did not say the particulars of when she got married or the age of her child. She is unconcerned about the disease, might be this is due to her ignorance. She does not even disclose the health status of her husband. So this interview does not give us much of news or views.*

## 49

We met a woman aged 30 years, who was married at 11 years; even before she reached her puberty. She is uneducated. Her husband is 35 years old and he is by profession a farmer. She has two children: one male and one female. She says she has not yet tested her children. She did not give her native place. She says her blood was first tested in her native village. She lives only in a hut. She did not disclose the health status of her husband.

*ADDITIONAL OBSERVATION: She got married at 11 years. Both her children are not yet tested. She does not disclose anything about her husband's health status. She did not answer whether her husband has been tested for the disease. She is poor. We are not able to make out whether he is living or dead. She does not know how she has got the disease. She does not know how one gets infected by HIV/AIDS. She says she was very ill and after admitting in this hospital she feels little better. We were not able to get any information from her.*

### 50

On 05-02-2003 we interviewed an inpatient in the Tambaram Sanatorium. She is just 30 years old and has studied up to 8th standard. Her native place is Tanjore but she is now living in Tiruchi. Her husband is a science graduate who owns a butcher shop. She says her husband had habits like smoking, visiting CSWs and he was also alcoholic. She says she does not know how she got the disease. Doctors have confirmed after blood test that she is HIV positive. She says she does not know anything about HIV/AIDS but knows it is a very cruel disease. She has two children: one male and one female. She says doctors take proper care of her.

When asked "can the government give free test for HIV/AIDS?" She says "it is difficult, and government cannot directly give free test for HIV/AIDS." But when we asked her why, she says, "it is very difficult practically." She feels that HIV/AIDS affected children can be joined in regular schools and need not have separate schools. She says one gets this disease by taking a wrong path. When we asked what she means by 'wrong path' she says men have sex with strangers / CSWs. However she said she had not taken any treatment from other doctors. She says she found that she has HIV/AIDS in June 2000 and from then on she is taking treatment in this hospital. She says now there are tablets and contraceptive methods so that one can 'be better' by using them. She clarified to us that she meant one can have sex by using these preventive and protective methods, thus they can be benefited.

*Additional Observation: She says that she is not well aware of HIV/AIDS, but she knows it is a cruel disease and people feel shy to speak out about this disease. She did not answer the question whether she has tested her children for this disease. She does not openly blame her husband for being the cause for this disease. Further she did not answer any question concerning her husband; except the fact that he is a science graduate and is a butcher by profession.*

*She did not answer even about his whereabouts, when we asked whether he comes to see her she did not answer that question. She said that no one ever visits her. She did not even answer the question whether her husband underwent blood test for HIV/AIDS. At this stage we feel that it is duty of the doctors to force the spouse to undergo HIV/AIDS test when one of husband or wife is infected by HIV/AIDS. Further we see she is not dejected she looked composed. She did not talk about her children, and when we asked her who takes care of her children she did not answer anything. She said her native place is Tanjore but now she lives in Tiruchi. But she did not give the name of the native place of her husband.*

## 51

An illiterate woman aged 32 years was interviewed on 09-04-2002. She is from Perambur. Her husband is deceased. She has two children, a 17 year old girl who is uneducated (never entered school) and a 15 year old boy who studied up to 8th standard and is now at home. She feels bad, because she has not married off her daughter.

She does not know when she was married or even the age at which she was married. She has some property, which is with her mother. She says her husband does not drink, but is a regular gambler. He used to fight and beat her whenever she questioned him. Her marriage was an arranged one. She says from the day of her marriage she did not lead a happy life. She says her son is in the village without any work. No one visits her, she is all alone in this hospital. When asked how she got the disease she says she must have got it from her husband. She said her husband was ailing for a few years but he did not say the disease he was suffering from. But she says in his last days he was only taking treatment in this hospital.

She says she had headache and fever that is why she tested her blood on the advice of the doctors. They said she has HIV/AIDS. She is ignorant of the ways by which the disease spreads. She is very ignorant of her husband's relationship with others. She says he did not say anything. But she says she got the disease from her husband; that is what the doctors who give her treatment said. She is not willing to give her opinion about awareness programmes or about care of children with HIV/AIDS or about free HIV/AIDS test. She says, "I don't know anything about these things."

*Additional Observation: She must have been married at the age of 13 or 14 years. It is once again an instance of child marriage. She was unhappy and felt it was difficult to get along well with him. Her husband was cruel, used to beat her even for small things as he wanted to lead a unquestioned life. He felt that he can waste money on gambling or on CSWs but she should not question him. Till his death she did not know the disease with which he was suffering. He would never confide to her but she recollects that he has taken treatment in Tambaram Sanatorium. She looks strong and healthy: one needs to thank the doctors who advised her to go for blood test for HIV/AIDS as soon as her husband died. As a result the disease in her is not chronic, but only in an initial stage.*

On 09-04-2002 we met the grandmother and grandfather of an HIV/AIDS patient. The patient (grandchild) has married her uncle, thereby becoming their daughter-in-law also. The grandmother answered on the behalf of the patient, who is 26 years old. This is also an instance of child marriage because this patient was married at the age of 13. She had studied up to the 4th standard and her husband has studied up to 9th standard. The grandmother of this lady says that her son was a scoundrel, had all bad habits and above all he was a ganja addict, a drunkard and visited CSWs. When we asked why she made such an alliance for her grandchild she says only her granddaughter loved him and got married to him. The grandmother says that even before and after marriage they used to fear that their son would infect their granddaughter. She says her fears have come true. All the questions were only answered by the grandmother as the patient (granddaughter) was so ill, she could not sit or talk.

The major symptom the women suffered was from boils and scabies all over her body. They said that her blood was polluted so they went for blood test, which confirmed that she was infected by HIV/AIDS. She is taking treatment in this hospital for the past 2 years. She has been taking treatment in many hospitals for the problems of boils and scabies. The grandmother feels that in their family such a thing must have never happened and sadly accepts that all these were due to her son. She is from Mettur in Salem. Both her grandparents are pretty sure that their son has infected their granddaughter. In fact they openly feel for it. Her grandfather joined the discussions.

He says that he will not even claim him as his son, this was the feeling of the grandmother also. He too complained of all the bad activities of his son. The patient is more bothered about her children because she hopes soon she will be cured and would go and look after her children. When they came to know she was affected with HIV/AIDS they were shocked because they only expected to hear the blood is spoilt and it had to be purified, thus they say they did not expect her to be affected by HIV/AIDS. Her grandfather felt very sad to know his granddaughter is affected with HIV/AIDS when he saw even small children visiting this hospital he says he consoled himself.

He also accepts the free HIV/AIDS test by the government for one and all. Further he says this disease carries a social stigma with it that is why people do not speak openly about it. He is of the opinion that CSWs should not be allowed to practice their profession. He says only in villages more awareness programmes about HIV/AIDS should be given. They are unaware of condoms. She says they pray to god once a year. After they are affected with HIV/AIDS they have started to hate the gods. However they have not changed their religion.

*Additional Observation: This interview was answered on the behalf of the 26 year old patient by her grandparents. Part of the interview was answered by the grandmother and part by grandfather both of whom were taking care of this woman. She married him on her own accord. He had all bad habits. We are not able to understand from our interview whether he is living or dead. His father personified him as a mountain of evil and he has only infected this innocent girl. They have started to hate god and they don't believe in god, for he has given them more problems than any good. They answered around 60 questions. The whole family was interested in discussing with us. They were of the opinion that only the husband has infected the lady. It is unfortunate that even his present place of living is not known to them. They feel more awareness program must be given in the villages. It should be such that even uneducated must be in a position to follow and understand the advertisement and its implied meaning.*

### 53

We met a native of Elankottam from Villupuram on 09-04-2002. She is 35 years old and uneducated. She has only one daughter who is 13 years old. She is studying in the 8th standard. Her husband died nine years ago. She became a widow at 25; she was married at the age of 15. Her husband was uneducated. She did not give the names of her parents or of her parents-in-law. Her marriage was an arranged one. She says her husband was a drunkard. She is unaware of his other activities. She says her husband died of HIV/AIDS. She says she must have got the disease only from her husband. When asked whether she had sex with any other men after her husband's death, she said no.

Her relatives do not know that she has this disease; only her mother knows and her mother is very kind and affectionate with her. She says her mother is the sole person who takes care of her. She says the whole family would have died without her mother. The doctors are very kind, they asked her how she was infected. She says she narrated her story completely. She was first tested for the disease by one doctor named Parasuraman in Kallakuruchi. She had cold, fever and cough. As she had no money she did not go once again to the doctor. She went to work. She says when she sees all the patients she is consoled for she used to think that she was the only one who is ill, but she observes all of them are seriously ill like her. She advocates employment for HIV/AIDS patients who are healthy enough to do work.

She says children must be cured. She says only in this hospital people say HIV/AIDS is curable. She admits that she has not seen HIV/AIDS advertisements; HIV/AIDS advertisement is given only for the benefit of people that too for people to be good and proper. She is also supportive of the idea of free HIV/AIDS test by government. She expects that if she is cured she can go for a job and earn her living; that is her sole expectation.

She does not know for how long she has been ailing with HIV/AIDS. She says she had cold, fever, cough, rapid heartbeat and pain in her limbs, she says now she feels better after taking treatment in this hospital. Asked about people in general she says, if they know someone has HIV/AIDS, they do not go near, they become very aloof. This disease affects the social status of the people who are infected by HIV/AIDS that is why they do not talk about this disease. She says from 7th standard onwards sex education can be given to children. She advocates youth to live with morality and self control, she says only that can save people from HIV/AIDS. She says she used to visit temples twice a year. But now she has no belief in god. She says CSWs are arrogant women and only because of them she has got HIV/AIDS. She says they are women who like to covet the property of others. She says they expect money and have sex with men, which is not proper. She says they should not practice such trade.

*Additional Observation: She says her husband has died 9 years ago, but she does not know whether he died of HIV/AIDS. She feels she has got the disease only from him. She is against CSWs but at the same time she does not answer why men seek them. When we told her that there is demand for them and that is why they practice, she refuted by saying that CSWs are bad since they covet the money of others. She is not aware of how her husband died. Here it is important to mention that more than sex education to children, education about HIV/AIDS should be given to these totally uneducated men and women from the age of 15; for they cannot read and do not care or afford to see TV. They may see TV only from public centers for in their poor areas, they do not possess TV in every home.*

## 54

We interviewed a woman, on 09-04-2002; she is a native of Kullangulam from Perambalur district. They are Chettiars, who do the job of stone cutting. She is 35 years old and her husband is 40 years old. She did not say the names of her parents or parents-in-law. Both she and her husband are uneducated. She wants to see her children. She complains that she is unable to eat well. She says she has 3 children, all of them are females. She was not in a position to talk so her husband answered most of the questions directed towards her. She was married at the age of 13. Her husband says it is since 22 years since they are married. He has married off one daughter; another daughter is 8 years old and is studying in the 3rd standard and another one is just 3 years. He says all the three children have been tested and they are free from HIV. The doctors and nurses take good care of them and do not show any sign of dislike.

Their marriage was an arranged one. Since his wife had no parents he did not get any dowry. He says now he does not have the habit of drinking, earlier he used to drink. He says after a hard day's work he would drink, then eat and sleep. He however clarifies he has never beat his wife. He says he was tested at Salem. The doctors said, "If you want to live for some more time please go to Tambaram hospital and take treatment." He says he did not have sex with other women. When asked whether his wife had sex with any man other than him he said no. He says only god knows how the disease spreads. He says he knows well that he has not made any mistake. He is not aware of condoms. Only after coming to the hospital he says he knows about HIV/AIDS. He says only four or five days has past since he got admitted in this hospital. They have been given medicine.

He says only after taking treatment he is suffering from itching. Earlier he used to get only fever and headache, once in two months. None of their relatives know about his illness. Only both he and his wife know, even their children do not know about it. He says he does not know why people do not speak about this disease. When asked if HIV/AIDS is a cruel disease, he says he does not know about it. He looked strong and healthy. When asked whether government can give employment to people like him, he says are ready to do any type of work.

He says children affected with HIV/AIDS must be educated. He says he has never seen a TV set. He says for his status he cannot afford a TV. He says his village is very remote and there is only one TV in the Panchayat office which he has never seen as he had never been to that office. He says if they give advertisement in their village about HIV/AIDS, it is good. He further adds that only due to ignorance of the villagers they become victims of HIV/AIDS. He further feels that if government announces a free check up for HIV/AIDS it would be good. When he heard that he had HIV/AIDS he felt sad for in his life he has never been admitted to any hospital. He is however little satisfied when the doctors gave him some medicine and discharged him. He feels very sad for all the HIV/AIDS affected patients. He says when they shout of pain he feels very depressed. Since he has not seen the soft porn magazines he does not know about it. He says it is wrong on the part of CSWs to sell sex for money. He says by believing in the government there is no use or help; because the government does not do anything in time or for the poor people like him. He is a stonecutter; when he asked help from the trade union it turned out to be of no use. He earns a pay from Rs.1000/- to Rs.1500/- per month. He says he has no property.

*Additional Observation: This is the first case were both the husband and wife are affected with HIV/AIDS. The wife is very serious. The husband is affected but he is only an outpatient. He is hale and healthy. Till the end we did not get to know how he has been affected by HIV/AIDS. He is of the opinion that the government is too slow to help poor people like him. He is affectionate to her and takes care of her in the hospital. He is very poor. He is willing to take up any job within his capacity if he is offered the same. At the same time he feels the government will not do that. He says awareness programs in very remote villages where over 90% of them are uneducated must be given compulsorily; otherwise he says these poor uneducated people become easy victims of HIV/AIDS due to ignorance. It is unfortunate that he is well aware of all other things but unaware of any information about HIV/AIDS.*

## 55

On 08-04-2002 we met a woman, who is an inpatient of the Tambaram hospital. She had studied up to the 6th standard.

Her husband is deceased. She is from Rishivandiyam in Vizhupuram. She has only one daughter. She was married at the young age of 17. Her daughter is studying 3rd standard; she has not tested her in spite of the doctor requesting her to do so. Her marriage was an arranged one. She says her husband was a drunkard and a smoker. He has also beat her and made scenes under the influence of alcohol.

She has hereditary property of 30 to 40 acres of land given by her mother. She is ill for the past one and half a months. She says doctors give her injection. She has cold, fever and cough. She advocates employment must be given by the government to HIV/AIDS patients. She does not know how to deal with infants who are affected with HIV/AIDS.

She however feels that such children must be with their parents and relatives only. She says, she has not seen advertisements about HIV/AIDS in T.V. or on buses or autos.

She further adds such things are not properly shown in remote villages. She says more awareness programs to reach the uneducated adults must be given in these villages. She wants the government to give free HIV/AIDS test to one and all. She says sex education can be given from the age of eight. She says she was first tested for the disease in Kallakuruchi. They immediately advised her to go to Tambaram Sanatorium.

She says she feels very sad to look at other HIV/AIDS patients. She says she has spent over Rs.10, 000/- to cure the persistent cough but she has not been cured, only in this hospital it has been cured. She acknowledges that her husband had several extra marital affairs but he has not disclosed any of them to her. She says she never had any such relationship with any man. She says she knows about condoms but her husband had not used it. When we asked her "did you check your husband for this disease," she replied that he was no more. (She used singular words to address her husband). She says her husband was killed by some of his enemies. She advises youth not to go to CSWs for it is due to them that people get infected. She expects only free medicine from the government. She says pictures depicting HIV/AIDS are good to be fixed as posters in public places because only then people will come to know that when they have sex with CSWs they will get HIV/AIDS.

She speaks very derogatorily about CSWs. She says their very existence spoils the village. She says only her mother knows and no one in the village or relatives know about her illness.

*Additional Observation: From her talk we understand that her husband was possibly a rowdy and his enemies murdered him four years ago. She says he had many bad habits and some of his bad habits caused, him to develop enmity, which resulted in his murder. She did not say the names of her parents or the names of her parents in law. She says her husband lived like an unquestioned 'daadha', drunk all the time, whenever she questioned him or even asked simple things she was beaten up.*

*She says the disease must have been spread to her only from him since he had no moral standing. The way she addresses him reveals his character and her unhappy married life. She must have been married young for she speaks of 12 years of married life with her husband and she is 30 now and her husband has died four year ago. She did not say about the belief in god or her offering prayers after disease.*

*From our conversation with her we find that she used to pray previously, after she came to know that she was infected by HIV/AIDS, she does not wish to pray. She was very weak and was just able to sit. However she says it is not possible for her to even walk 10 feet for she says it makes her breathless. She is happy for she is cured of tuberculosis and cough. She is not sad or depressed. She is composed and very calm. She does not think this disease to be a big problem. She is in bad spirits only when we talk about her husband. She became very angry and said (Avan seththei puttan – He is dead).*

We met a woman from Manali who is an in-patient in Tambaram Sanatorium on 10-06-02. She never revealed her age, education or native place, or her job. She gave her present address to be in Manali. She just nodded her head in affirmative for the question "Are you married?" She says she has no children. She adds that for the past six years she is married. She does not know how she has got this disease. She says that they have not yet tested her husband for the disease. He may be infected, she adds.

She first took treatment and test in the Stanley hospital in Chennai. When she was confirmed to be HIV +, some of the doctors called her separately and advised her to eat well and asked her not to be depressed over it. She says no one visits her since they came to know that she was HIV/AIDS affected. She says even her husband is not visiting her. She says he has confided to her that several times he had visited CSWs before marriage but had promised to her that he had not visited them after marriage. She suffered from vomiting and lack of appetite.

*Additional Observation: She looked sick and was not able to talk or even sit. She looked uncared for; she had not combed or oiled her matted hair. It looked as if she was not able to do them for herself. Her parents, her parent-in-laws are living but we are not able to understand why they do not come and take care of her. She says as soon as they knew the news that she has HIV/AIDS, they don't visit her. It is unfortunate that even her husband does not come to see her. She has no children, she may be only in her late twenties. Doctors who looked after her were kind with her and advice her not to be worried or sad. She feels that only her husband must have infected her. He has so far not taken the HIV/AIDS test. She feels she can do nothing about it. This is a special case where even the kith and kin including parents, parents-in-laws and husband desert a woman patient with HIV/AIDS.*

A woman now an inpatient of the Tambaram Sanatorium was interviewed on 10-06-02. She did not give her age, place of birth, native place, present address or even her years of married life. She must be in her early thirties or late twenties. She was married to a much older man, as her parents were very poor. She says she is very unhappy over the marriage. As he was a very old man she had no satisfaction in her life.

She says her husband was not a moral man. He used to have affairs with several other women in his younger days. She says his neighbours told her about his character. She says only he had brought her here and had admitted her in this hospital. She says there is no one to buy her medicine or even a cup of tea. She further says her husband always suspected her of some relationship or the other with other men. When asked to give advice to the youth she said, "my health is spoilt for ever". She had nothing to advice the youth. She is worried that her husband is not testing their son for HIV/AIDS. She feels since money will be spent on it, her husband is not carrying out the blood test for their son. She says the other patients are happy and kind with her, some cry, some suffer of pain, for some their husbands are also infected. They have not carried out the blood test for HIV/AIDS for her husband. Further she does not know whether her husband is already undergoing treatment for it.

The son she mentions here is not her child. It is the child of his first wife; she has no child of her own. She does not know how his first wife died. The son of his first wife is older than this woman so her husband is constantly suspecting her of having sexual relationship with him. Based on this issue they used to quarrel. She says her husband never accepts her chastity and always suspects her. Now she says that if she asks him to test his son he will suspect her doubly.

*Additional Observation: This is a very strange and unusual situation. The old man has married a woman who is younger than his own son (from his first wife); and he tortures her by suspecting her fidelity. She is not bothered about the disease as much as she is bothered about her suspecting husband. She does not know how she is infected; her husband has so far not disclosed that he has HIV/AIDS. She is more bothered about testing this son than her husband for HIV/AIDS.*

*One is not able to understand anything, but one thing is clear that the old man really is bothered about her that is why he has taken all steps to admit her in this hospital and give treatment. Whether he is infected with HIV/AIDS is yet to be determined; for she herself acknowledges her husband takes care of her. She is least bothered about the old man or his health.*

## 58

We interviewed a 5th standard studied woman from Tirukoviloor who is now an inpatient of the Tambaram Sanatorium. She did not give her present address. She says she is married for the past 13 years and has 3 children: one female and two males. She says she always suspected her husband of having illegal relationships with other ladies. So they used to quarrel over this matter. She says she used to keep him at a distance for the fear of being infected. But she says she could not always protect herself. She says she was first very dejected. She knows it spreads by sex. She says all her village people go away from her. However her parents in law are kind to her.

*Additional Observation: She could not sit or walk, she looked very weak. Her husband is deceased and she has not tested her children for HIV/AIDS. She however did not disclose how her husband died. She says her husband must have died of HIV/AIDS. He visited CSWs before and after marriage. It is a surprise she did not say whether her parents are living or dead, and none of her blood relatives*

*visit her. As she did not say her age and how many years have passed after the death of her husband we could not estimate her age even approximately.*

## 59

We met a 45-year-old woman with HIV/AIDS in the Tambaram hospital. She is a daily wager. She is from Pondicherry, and has studied up to 7th standard. She says, after seeing other patients in the hospital, she felt little calm, as so many people are affected by HIV/AIDS. She was interviewed on 17-06-2002. She says she is also affected with STD/VD. Whenever she saw advertisements in the T.V about HIV/AIDS, she said she was feeling guilty. She knows that this disease is a disgraceful one and people with this disease hold a social stigma.

*Additional Observation: She looked ill with scabies and blisters all over her body. She was able to just sit with great difficulty. She did not talk about her children or even about the work or living conditions of her husband. What her husband is doing does not come to light from this interview. She never talked about them (her husband and children) or never even showed any concern about them nor was she willing to advice them to go for HIV/AIDS test. There was no one around her.*

*She says no one sees her or visits her. She does not disclose whether she had affair with any other men. She says if people come to know that one has HIV/AIDS they cannot even get into a teashop and buy a cup of tea. In any public place they cannot move if anyone around knows that they are HIV/AIDS patients. However till the end she did not answer the question whether she got this disease from her husband. She never blamed him. She did not answer the question whether her husband has any other bad habits. So we are unable to conclude of how she had been infected by HIV/AIDS.*

We interviewed a woman from Guduvancheri who must be in her late fifties or early sixties on 21-06-02. She is presently residing in Tambaram Sanatorium. She says her husband is dead, but does not remember exactly when he died. She has three children, two sons and a daughter. All of them are married and she also has grandchildren. She says she was working in a shoe company in Guduvancheri. Her husband was a carpenter. They used to pray to all the gods. She was first tested in Chennai as her legs had swelling. She took treatment in a private hospital and spent a lot of money over her disease. Her stomach became hard and she could not be cured. They took several X-rays. Then she was advised by some other patients in the hospital to go for treatment to the Tambaram Sanatorium. She says after she got admitted in Tambaram Sanatorium she is well. She says she does not know how she got HIV/AIDS. She says she did not have any extramarital affair with others. "If one does not have such habits how will one ever get HIV/AIDS?" She says the doctors and the nurses also asked her the same question. She says she does not know, we asked whether any one in the shoe company had HIV/AIDS; she says all the persons who work in the company are only women.

She is very poor, for she said her sister visited her and said the thatched roof of the house was fully torn and the children wanted her monetary support. She says she has forgotten the names of her parents-in-law. She showed us the photo of her elder sister. She says her sister visited her 10 days ago. She says she has seen HIV/AIDS advertisements. She says she was very sad that her eldest son did not visit her yet. She feels her relatives are kind and considerate with her. She happily recollects her daughter and her sister who visited her last week and gave her good food. She says she does not know how she got the disease. She is having a complete faith that she would be cured soon. She says the body has lost strength and she adds that now she feels better. She says her husband died 30 years ago and at that time her children were studying in high school. Only one daughter had attained puberty. Only then she went and took up the job in the shoe company, and she says she got the disease only from the company.

She says her owner does not know that she has HIV/AIDS. She says she had swelling of legs and could not walk and could not eat, for 2 months she was eating only bread and milk. The owner on seeing her health conditions gave her a week's leave and then a month's leave. She has not yet returned to her work place. She has blisters in her mouth and says she could not eat well. But however she feels she is better now. Before she could not walk or even fold her legs but now she is able to do so.

Only for the past two months she says she is ill. Her husband had only one wife. He used to drink a lot. Her in laws have never troubled her because they are natives of Tanjore so, she says they visited her very rarely. She says her marriage was an arranged one.

She has never donated blood or received blood. When we put the question "this disease comes through sex? How did you ever get it," She says the doctors said the disease has started just like and they said she could be cured. Her husband has educated the three children untill his death to the classes 8th, 7th and 5th respectively

*Additional Observation: She looked like she was in her late forties. So when she says 30 years have past after her husband's death it implies she must be in her early sixties, for she says the children at the time of his death were studying in high school and her daughter has attained puberty. By all possibilities her husband should have died some 10 or 12 years ago. Might be he could have had HIV/AIDS and she when was in the delicate time of menopause should have been affected by the virus which was in her. The fact that she is morally good is evident, from the respect she holds from her close relatives and friends who visit her with good food and fruits to comfort her knowing fully well she is HIV/AIDS affected.*

We interviewed a woman of the Reddy caste on 01-08-02 in Tambaram Sanatorium. She is 50 years old and uneducated, from Kuppakkadu, Sholinger. She is married and her husband died recently, must be of HIV/AIDS. He was illiterate and worked as a coolie. Her parents are no more, she lost them when she was very young. Only her grandmother brought her up. She has just been admitted in this hospital. She says that her husband was a drunkard, he was very ill when he died. She says she was troubled with continuous cough. She has four children: 3 males and one female. She says she has not tested them as they are keeping good health; her son is working in a cloth company. She is from middle class background living in her own house. She feels that she will be cured and then go home. She says the government will never do anything good.

*Additional Observation: It is evident from the interview that her husband must have infected her for she says he was very ill at death. They have not given him any medical treatment. He was a drunkard and had many bad habits. She had the burden of looking after four children and they could not have monitored her husband for he was just a coolie and drank all the time. She has brought up her children with lots of difficulty.*

**62**

On 1-08-02 we interviewed a 42-year-old uneducated women who is an inpatient of the Tambaram Sanatorium. She is married, to a man who is unemployed at present. Her husband was a mill worker and currently he has left the job. Her parents from Kothagiri had come to see her. She says her son takes care of her. She feels better after getting treatment from here. At first doctors said she was affected with tuberculosis but later they informed that she is affected with HIV/AIDS. She said she would have got the disease from her husband: people in the hospital had said so to her.

She says her husband has good character and she does not suspect him. But however her husband stayed away from the home and at that time he could have visited CSWs. She says he does not have any disease.

So far her husband has not been tested for HIV/AIDS. She says she has not had extramarital sex with anyone. She says she used to have petty quarrels at home. Nowadays he visits home only twice a month and stays in Coimbatore. She prays to God to cure her soon, forgive her and send her home. Only few relatives know about it and her sister's husband said she would not become cured so she is very depressed. She wants the government to find medicine to cure HIV/AIDS. She does not know how one gets this disease. She supports the proposal of free HIV/AIDS test for one and all by government.

*Additional Observation: She is affected by HIV/AIDS but does not know how HIV/AIDS spreads or how one is infected. She is very ignorant about the disease. Her husband has stopped going for job for the past two years. Her children support her. She was very sick before and now she is able sit and talk. She has blind faith in God, prays to God for complete cure. Her husband visits the family only twice a month makes us think - whether he is married to someone. She is not able to question him on his erratic behaviour. She has not tested her husband. It is not known whether he is taking treatment secretly. We are not in a position to find out how she was infected.*

### 63

We met a woman aged 22 years, a native of Vaniyampadi on 13-08-02. She has studied up to S.S.L.C. She is now an inpatient of the Tambaram Sanatorium. She is married to a 35-year-old man, who is uneducated. He works as a coolie in Mumbai. She has two children: 1 male and 1 female. She was tested positive for HIV/AIDS in the year 1999. She had no doubt of any kind when she got married to him in Mumbai for at that time she did not know about Mumbai.

Her husband had boils and scabies, which never healed, she says even at that time she did not suspect him. She does not know when he got HIV/AIDS but both came to know that they are infected at the same time. Even after marriage her husband did not take her to Mumbai. He would come and stay with her in Vaniyampadi for one month in a year.

She said she had symptoms of STD soon after she had sex with him. She however said she never dreamt of being infected by HIV/AIDS. The doctors here look after them well. She says it is a must that school children are given sex education. She says people shun talking about HIV/AIDS because of shame and the social stigma it holds. Her husband has got the disease only by extramarital affairs. She says her husband has never acknowledged his bad activities. Even after being confirmed as an HIV/AIDS patient he did not disclose his bad habits. She is upset over this disease. She recommends that the government can give jobs for HIV/AIDS affected people. Both her children are unfortunately affected by HIV/AIDS. She first disclosed her disease to her mother.

Only in this hospital she came to know how HIV/AIDS spreads. She says even in Tambaram Sanatorium a nurse used the same needle for several persons. When she said to the nurse that the disease will spread, the nurse scolded her and asked her to mind her own business! She has also taken indigenous medicine for treatment. She also advocates free HIV/AIDS test for one and all by government. She has changed from her native religion to Christianity. She loves Jesus and hopes he will cure her. She belongs to the Chettiar caste. She is poor. Her husband is her close blood relative. She lives in a nuclear family.

Her parents-in-law look after her well. Her mother-in-law does not feel sad about her son. Only her brother and mother know about the disease. She says advertisements are not sufficient in villages to create awareness among the uneducated. She says advertisements should be also given about the disease after one gets the disease. She came to know about HIV/AIDS only after she was infected by it. She knew about HIV/AIDS but did not know it was such a cruel disease. She says, her marriage was an arranged one. She was not for marriage, only her parents forced her into this marriage. She says that for HIV/AIDS children government should help in education and give them free medicine and medical aid. Only sometimes she feels worried about the disease. She says to people not to commit mistakes. She says that boils, blisters, the cruelty of the disease and the affected parts should be photographed and publicly displayed so that by looking at it people will fear to do mistakes. This will not only create awareness but also fear in the minds of people so that people will not make repeated mistakes and become HIV/AIDS victim. This was her line of thinking.

*Additional Observation: Her husband infected her. She is affected by HIV/AIDS when she is merely 20 years old. She was married off at sixteen and at 20 she has become a full-blown HIV/AIDS patient. Her husband has married her in Vaniyampadi, left her there itself and he went to Mumbai. He is seriously affected with HIV/AIDS, STD and VD and visits her only once a month in a year. He did not only make his wife a HIV/AIDS patient but has also ruined his two children who were born HIV+. Thus the whole family is affected with HIV/AIDS. She says that he is unconcerned and does not repent about having ruined her life. She says some awareness programme must be given about the HIV/AIDS patients, their sufferings and the types of symptoms and problems so that when they have these symptoms they can go for medical aid and take treatment before the disease becomes chronic. This will help not only a HIV/AIDS infected patient but also every one to check themselves about symptoms and know their health problems. She sadly expressed that in the Tambaram Sanatorium some nurses use the same injection needle for several persons and when she pointed this out she was ill treated. Further she advocates more awareness programmes in red light areas so that both the CSWs and the customers would be careful. She further feels that the*

*government should take care of the education of HIV/AIDS children and give them free medical treatment.*

## 64

We interviewed a woman on 20-09-02 who did not wish to reveal her age to us. She is now an inpatient of the Tambaram Sanatorium. She is married for the past 8 years. She has two children. She was tested for HIV/AIDS only in the Tambaram Sanatorium. She says her husband has infected her. She further says her husband visits CSWs whenever he goes on trips. He is a lorry driver. She says she was tortured by her husband's family very often for more dowry and more jewellery. They were always dissatisfied with the dowry. She has two children below the age of ten: one male and one female. However she says that her relatives are as usual with her. She says she is fully aware of condoms, for she had worked before marriage in a social service center from where she has learnt about this. She says she does not know whether her husband had extramarital sex.

*Additional Observation: She must be in her late twenties. She is a victim of dowry harassment. She did not give the name of her native place, age, education or her parents' name. She did not even answer the question whether her parents are living or dead. She looked very thin and ill. She talked more about the cruelty meted out to her by her parents-in-laws than about HIV/AIDS. She says she cannot and does not have the courage to ask her husband to go for HIV/AIDS test because it would make her undergo further torture by her husband's family. We come to know that he was already taking medicines. When asked why he was consuming them he said to her that he was suffering from ulcer. Further she has scabies and blisters all over the body. Her husband also had scabies and blisters all over the body, which has made it impossible for him to go regularly for trips regularly. She did not confide to us of who had admitted her in the hospital or who was taking*

*care of her. She was all alone. She was not very frank but she did not feel dejected for having HIV/AIDS.*

## 65

On 14-08-02 we met an inpatient of the Tambaram Sanatorium who was seriously ill with palpitation. She was not even able to sit. She has no one attending on her. She said she was just 40 years old. She is a coolie, a daily wager. She is from Choliang Karai in Arakkonam. She is uneducated. She has four children: 3 males and one female. She proudly says "I live in a tiled house and it is my own house." She had cough and cold, blisters all over the body and could not fold or move her hands or legs. She is not in a position to walk. She says she does not know about HIV/AIDS. She says her husband had died after a prolonged sickness.

*Additional Observation: Unfortunately no one was by her side. She did not answer the question of who was taking care of her. She looked very ill and was panting for breath all the time. She made many answers by just nodding her head. We felt very sad to look at her in that state. We guess that she had got the disease from her husband because she said that her husband also had scabies and blisters all over the body, that he was lying like a log for several days, unable to sit or move. She says that she used to give him food as he did not like to see any doctors or take treatment he died after some months.. She says she is affected just like her husband and there is no one to give even a cup of water to her. We can only assume that he could have died of HIV/AIDS. As she was very ill we could not put forth more questions.*

## 66

On 11-9-02 we interviewed a 26-year-old woman who had studied up to 10th standard. She was married at 17 years to a 28 year old man who was her distant relative. Her husband is dead. She knows tailoring. Since her relatives were pestering, her mother gave this woman in marriage to a TB affected man. He was a goldsmith by profession. They lived in a joint family. Her husband had never disclosed any truth, however she does not suspect her husband. She says her husband was a tuberculosis patient from the day of marriage. He said that he has not made any mistake. She says her parents-in-law have cheated her because they knew fully well that their son was affected with HIV/AIDS. They tested her blood very secretly but did not disclose to her that she had HIV/AIDS. She has a son studying in the 4th standard who is just 9 years old. She says no one in her family are affected with HIV/AIDS. If she had known about the disease she would have been certainly careful she adds. She feels that HIV/AIDS affected children must be given free medical aid and must be educated. She lives in her own house.

She gave her native place and address of residence. She is from the Cheyyar District. Her mother-in-law said that he was good for the past 7 years but later he changed his habits. He went as a tourist and returned after a week or so. He was happy, she says only now she suspects him. The interviewed woman says that her husband had sudden itching with black spots all over the body and even to look at him she felt very uneasy. She was frightened even to look at him at that stage. She feared to sleep with him, because he became black and thin with a shrinking of the skin all over the body.

People think it is a disgrace to talk about this disease openly. She says one should discover the medicine to cure HIV/AIDS. She says only after his death she came to know about HIV/AIDS. She also recommends free HIV/AIDS test for every one, because the poor will benefit by that. She says self control and leading a life of morality alone can prevent the spread of HIV/AIDS.

She says she will not come out to talk about HIV/AIDS because it will spoil the self respect of her and her family. She says she is very frightened to look at certain HIV/AIDS patients because she fears she would become like them. She says that she has no expectations from the government but the government can provide jobs for HIV/AIDS patients who are capable of working. She says she is very sad and she weeps at times and is afraid that her child may get HIV/AIDS. She says she does not know how the disease spreads. She gives statistics that out of her husband's 5 friends only one is dead.

Now the interview is answered fully by the mother of the patient. Her daughter lived with her husband only for 5 years. He was first infected with tuberculosis, later on they found out that he was suffering from HIV/AIDS. Even knowing fully well that he had HIV/AIDS he had sex with his wife. Her mother feels very sad about this act of her son-in-law. Her daughter was well for four years then she became sick. Very often she was ill. Only on the basis of suspicion only her mother tested the daughter for HIV/AIDS as her son-in-law was a HIV/AIDS patient. She first took treatment in Erode. Then she saw an advertisement in TV. about Rajas Siddha Cure on the way from Tiruchi to Karur. There were no patients in the hospital so they thought only the name of the hospital was big but was not popular. She says even to just see the doctor one has to pay Rs.50. When she said she would first see the doctor and then pay the money that said no. They gave a rate like this

| | | |
|---|---|---|
| Registration charge | - | 500.00 |
| Testing charge | - | 7000.00 |
| Medicine (monthly cost) | - | 6000.00 |
| Consultation | - | 1500.00 |
| Total | - | 15000.00 |

Only when someone gave Rs.15,000 fully the doctors would see them. When she requested that she would pay after some time they did not agree. So, she had no other option expect to come back. She consulted her sons and they said on advice from their friends to go to Chennai Tambaram Sanatorium for treatment. She says her daughter is little sick only after consuming medicine. The mother says her daughter's health has become worse after admitting in this hospital.

***Additional Observation:*** *The marriage was an arranged one and she was married at the age of 17. She was widowed at 21. She says their parents-in-law were fully aware of the fact that their son was infected with HIV/AIDS but they cheated her. From the day of his marriage he was ill with tuberculosis, on checking they found him to be affected with HIV/AIDS. He had sex with her even after he knew he had HIV/AIDS.*

*She was taking treatment for the disease from several doctors in several places. The symptoms, which her husband suffered, were fearful. His body shrunk and his skin became black in colour. He died of blisters, unbearable itching and scabies all over the body. She says if she should have been aware of this disease and the way it spreads she would have been very careful.*

## 67

A 50-year-old uneducated woman, from Neyveli was interviewed on 4-9-02. She is married for the past 30 years. She has 3 sons. Her husband has died. She is taking treatment in this hospital for the past two years. She says she has also taken treatment in different hospitals. The doctors said she has all diseases in her blood. She became very weak so she had gone for treatment. She says her husband is a "gem of a man" and so she does not know anything about him. However, she also says her husband had all the bad habits like smoking and drinking alcohol. She says the doctors here come give the injections and go. She says she feels a little better now. She is a vanniyar by caste. She does not know much about HIV/AIDS or the way it spreads. When asked about her belief in god she said she does not believe in god. She could not be questioned more as she was very sick.

*Additional Observation: This woman is taking treatment for the past two years in the Tambaram Sanatorium. She is very sick, yet she talked with us. Her husband had died and she does not know the number of years that have passed since his death. She looked uncared for. Though she has four sons there was no one with her. She was not even able to sit on her own. She said it is one and a half years since she ever walked. Earlier she was able to walk to attend the call of nature, now she is not even in a position to do so.*

*She looked grave. She thinks she is admitted to this hospital for lack of strength and for tuberculosis, the medicines have not made her well but she feels somewhat better. When asked how one gets HIV/AIDS she says that she does not know how it comes or how it goes. She absolutely has no faith in god. She does not show her depression. She lives uncared for, all alone.*

## 68

On 18-9-02 we interviewed a Muslim woman who was just 24 years old and now an inpatient of the hospital. Her husband has died. She says she is educated up to the 5th standard and says her husband is uneducated. She is from Tirupathur in the Vellore district. She has no child. She was married when she was just 19 years old. Her husband died just one year after marriage. Only her sister informed her that her husband had died of HIV/AIDS. He was from Mumbai. Her marriage was an arranged one. Often he used to get fever. He was an auto-driver in Mumbai. She feels bad that he has concealed his disease from her. She was just a housewife. She came to know about HIV/AIDS only an year ago. Her husband was tested positive for HIV sero prevalence in a hospital in Tirupathur within a few months of marriage. She said that the doctors informed her that he was affected with HIV/AIDS. Then he took treatment in Erode. She was sick and she was also tested at Erode and the doctors there said she has HIV/AIDS. She took only indigenous medicine in Erode, but it did not improve her health. She was taking treatment in Erode for nearly two years spending lot of money on the medicine. She lost weight from 51 kg to 40 kg. She did not become well or feel better. She has not seen advertisements about HIV/AIDS, so she had never suspected him. All her friends and relatives are kind with her. She has no property of her own. She was totally unaware of HIV/AIDS. Her husband had at times used condoms. She now feels that because he has the disease the doctors would have requested him to use condoms. Only now she thinks about it and the cause of her husband using condoms. In her family she is the first one to be infected.

She lived only in a joint family and her mother-in-law was working as a servant maid in several houses. She only thinks of when she would go to home. She says it is a 'wrong disease' that is why people do not talk about it. She says the HIV/AIDS affected children must be given good food, good medical aid and education. She feels it is better to stay in hospital. She says doctors here are very capable for in one day they cured her of the blisters in her mouth. She says making HIV/AIDS test free for all is a good move and that everybody would be happy to accept it. She says more advertisements must be given. She advices youth to lead a life of self control coupled with morality.

She prays 5 times a day. Now she lives only in a rented house and her parents are very kind and considerate with her. The awareness programmes must be given to one and all. She suffered from vomiting, loss of appetite, stomach pain and other ailments.

*Additional Observation: She was married at 19 years and her husband died within a year due to HIV/AIDS. The relatives of her husband did not disclose his disease to her family. She feels very sad for being infected. Further she feels bad that just after two years of her marriage she has become a victim and is now suffering with HIV/AIDS. Her desperation about her husband cheating her and infecting her when he knew fully well about the disease remains unforgiven in her heart. This case clearly shows the poor status of women in our society and their untold suffering not*

*only of the body and soul but also of the social stigma with which they are held in the society.*

## 69

We met a 35-year-old woman who had completed 7th standard on 24-9-02; She was an inpatient in the Tambaram Sanatorium. She has only one child, he is studying in the eleventh standard. She lives in a rented house. She was married at 17. Her husband was a lorry driver. She had taken treatment from Kerala for over a year and she says it was of no use. She was first tested in Namakkal for HIV/AIDS. After she contracted this disease, she has lost faith in god. She says when she sees the patients, she feels very courageous and comforted for not only her but also hundreds of people are affected by HIV/AIDS. She does not disclose whether her husband had died of HIV/AIDS.

*Additional Observation: She looked healthy, contented and composed. But she was not in a position to walk or fold her legs. Her husband, a lorry driver had been living loosely. He was very ill and had taken medicine and had never disclosed to her that he was suffering from HIV/AIDS. She did not give the name of her hometown or her present address. It is unfortunate that no body is taking care of her. Some co-patients who are healthy help her. She did not wish to talk about her husband's death. She has however lost completely the faith in God. She was dejected and desperate when she came to know she had HIV/AIDS but after seeing hundreds of HIV/AIDS patients she says her desperation is gone. She is from the poor strata of the society. She does not blame his misdemeanor, but she feels sad for not informing about his disease to her and not explaining why he was taking medicine. She says symptoms she suffers are similar to those her husband had suffered. First he too had scabies all over his body and was suffering from lack of appetite and vomiting. She too suffered with the same symptoms. Many a time she thought her husband might be suffering from Jaundice, because he is a driver traveling*

*long distance but she never dreamt he would be infected by HIV/AIDS.*

## 70

We interviewed a 27-year-old Gounder woman on 25- 9-02. She was an inpatient of the Tambaram Sanatorium. She has no educational qualification and she does not know any way to support herself. She is from Kudumbathaur in Dharmapuri district. Her husband was a lorry driver. He died of HIV/AIDS. She says she has three male children. None of her relatives know about her sickness. Only her mother in law knows about it, she scolds her by saying that her son is a good man, only she is of 'loose character' that is why her son died. Now she has changed her religion to Christianity. She has no one to take care of her or her children. She says she was in the streets. Her husband's parents have taken all the dowry and other properties given by her family and have forced her and her children out of the house saying she is a characterless woman. Only a priest (father) from Dharmapuri has taken pity on her and he is taking care of her and her children. She was weeping bitterly when she said about it and after this she could not talk.

*Additional Observation: It is to be noted that several Gounder families (that we could learn about due to this case study) are very strict in getting the dowry for marriage. They are very concerned about the property. They do not wish to divide the property into many parts. So when the husband dies the wife is chased out of the house empty handed. Now she suffers from multiple sufferings. One her husband is dead and as she is uneducated and does not know any occupation to support her and her 3 small children, and in this case the last is the biggest problem for both her parents are dead. Secondly, the mother-in-law knowing well that her son has died of HIV/AIDS and has also infected her daughter-in-law; scolds her daughter-in-law as an 'immoral' woman. They say the daughter-in-law*

*was of 'loose character' and their son was good. But ultimately, chasing the daughter-in-law out of the house is unjust. Thirdly the woman patient's mental agony of being infected and her physical suffering of this socially stigmatized cruel disease.*

*Thus this case was very pathetic except for the father (a priest), she would have been driven to suicide. She says only this Dharmapuri priest took pity on her, when she was chased away from her home with her three children. She says only this priest visits her in the hospital and is educating and protecting her children, all the three of whom are just less than ten years. No one knows whether these kids are infected by HIV/AIDS.*

## 71

On 05-03-2002 we met a 41 years old woman in the Tambaram hospital. Her husband is 47 years and uneducated. He is from Manaparai in the Tiruchi district. She is also uneducated, she has never stepped into school. She has three children: two female and one male. She says she has married off both her daughters. Her son is just 18 years. Her husband was a lorry-driver. She said her husband was sick for six months and he had been taking medicine but had not informed her that he was affected with HIV/AIDS. She says she went for treatment for tuberculosis.

In Chengleput she went for an operation, there the doctors said that she was infected with HIV/AIDS. She does not know how she was affected but she felt very bitter that her husband had cheated her. She had asked her husband for and explanation to which he replied that he had made the mistake only once but she says that from his young age he has been continuously doing such actions. She says she need not test all her children for they are grown-ups but she wants to check only for her 18-year-old son. She says none of her relatives know about this. But her sons-in-law knows about this, and sometimes they show faces of disregard and discontentment.

She says when she sees her small grandchildren she feels she should live for some more years. Further she says she has to live for some more time and perform her son's marriage. She says she had seen HIV/AIDS advertisements and recommends more advertisement about HIV/AIDS to rural people. She feels the doctors look after her well. She advices the youth not to become slaves to sex. She expects the government to take care of all the HIV/AIDS patients in this hospital. The government should mainly take care of the education and medical care of the HIV/AIDS affected children.

*Additional Observation: This lady is more vexed by her husband's activities than by the disease. She was married at the young age of 15 years. She has married off her two daughters and has grandchildren but the pity is at this age she is suffering a social stigma, that too in front of her sons-in-law. She is inconsolable by her husband's activities. She feels more awareness program must be given specially to villagers for in certain villages over 50% of population would have not even entered school premises. Thus she feels some sort of education about HIV/AIDS must be given to uneducated adults in villages.*

## 72

We met a 31-year woman on 20-9-02. She is from a village named Pallappalayam in the Karoor district. Her husband died when she was 29 years old. She has had no education. Her husband was first working as a driver. Then he changed his profession as an agricultural labourer. She was checked for HIV/AIDS in Karoor, there the doctors advised her to go to Tambaram sanatorium hospital for treatment. Her husband had HIV/AIDS for several years. She says her sister-in-law told her that when he was injured in the head due to a quarrel, the hospital people refused to stitch the wound on his head citing the excuse that he was a HIV/AIDS patient. Two days later, her husband died of a heart attack as he used to drink everyday. He has never told her about his relationship with CSWs.

She has only a 3 year old son. She has left him in her sister-in-law's house. She has performed a HIV/AIDS test for the child and he is free from this disease. She says none of her relatives visit her, only her mother stays with her. She is now in her mother's house and her native place is Mohanoor. She says the doctors here are kind to the patients, give them good medicines but at times they make faces, which pains us. She says she expects only free medicine from the government.

Only a self controlled life can prevent HIV/AIDS**.** The sex should stop with husband and wife and they should not seek sex with others. She is sure that she got the disease only from her husband. She says she has not seen advertisements about HIV/AIDS. She had come to know about HIV/AIDS only recently. She was not upset when she came to know that she was affected with HIV/AIDS. Her early symptoms were fever, lack of appetite, tiredness and blisters in private parts. She was feeling very ill over all.

***Additional Observation****: This is a instance that proves the inhuman behaviour of the doctors. When a man has suffered from head injury and comes for first aid it is the basic duty of the doctor to help the patient survive and it is highly unethical on the part of the doctors to refuse first aid to the person on the basis of his being a HIV/AIDS patient. It is very unfortunate that the man died after a few days. The only thing she feels bad about her husband is his keeping the disease as a secret and at the same time infecting her. She is also affected by STD/VD, her body is full of blisters and scabies, she is not able to sit or walk. The doctors do not properly attend on her and the sores are frequented by flies and the wounds are open. She is calm and composed with a fond hope that she would be cured though she says that AIDS is incurable. One positive point is her mother who is good enough, takes care of her daughter very sincerely. This shows the pure love of her mother.*

A 10th standard educated Naidu woman from Andhra Pradesh was interviewed on 4-10-02. She is just 22 years old. She is married and she has no children. She does not know about the educational status of her husband. She is married for the past one year. She says her marriage was an arranged one. He is a car driver. She says she knows about condoms and has used it. She had no relationship with any other man before marriage. She says she has not asked her husband if he had any sex before marriage. She does not know how one gets this disease.

She says her husband is free from HIV/AIDS. He was tested and doctors said that he is free from HIV/AIDS. She says she does not know how she alone has got this disease. She wept when she came to know she was affected by it. She asks how will the government ever give jobs for those affected with HIV/AIDS? She says she would stay in this hospital for a month and then go home if she becomes better. Children affected with HIV/AIDS must be given free education, medicine, food and shelter. She adds that the doctors here take good care of her. She says she does not know why people don't talk about this disease. It is better if the government gives advertisements to uneducated villagers about HIV/AIDS. She has seen advertisement on TV, which state if you don't use condoms you may get HIV/AIDS. She had an attack of STD with blisters and scabies all over the body. After this she had continuous bleeding for 10 to 20 days.

*Additional Observation: She did not answer several of the questions, and till the end refused to give her parents name or parents-in-law names. Her native place was broadly given as Andhra. However she is sick with STD/VD. Her body is covered with scabies and blisters. Her husband has been tested for HIV/AIDS and he is found to be free from it. It is not known who has admitted her here. She was all alone. From her talk it can be guessed that none of parents or parents-in-law visit her. She was composed. It is*

*courageous that she has come all alone from another state to take treatment in Tamil Nadu. No more information can be obtained. She is the only woman who says that it is none of her concern about the existence or non-existence of CSWs. As she was a Telugu woman, we were unable to get several points across, we suffered from the language barrier. She only understood little Tamil, and replied to us in broken Tamil.*

## 74

We met a woman on 06-10-2002, who has studied up to the 4th standard but she says she does not know how to read and write. She is 37 years old. Her husband had a college education but she does not exactly know the degree he had obtained. She is from Manali, Chennai. Her husband is a mechanic and she is a housewife. They have no children. They live only in a rented house and as a nuclear family. She says her marriage was a love marriage so her parents-in-law do not even entertain their son, so her husband never visits. She is a native of Tirunelveli. She says she has no voice to even ask her share of property from her parents, since hers was a love marriage – even though she has a share of the big ancestral house in Tirunelveli. She was born with 10 children, all sisters are married and her last sister is a nurse. Her brothers work in a hotel, her last brother works as a cleaner. She is the eldest among women and she has two elder brothers and no elder sister.

She has come to Chennai as a maid servant, was asked to stay in a house and work for the residents. She was paid every month and was given food. She says she knows only Tamil but since her husband has visited several places he knows Tamil, Telugu and Hindi. Their only property is the house. She says she was married at 25 years and seven years have past since then. She says she did not inform about her marriage to her family. After some years of marriage she went home and informed them that she was married and her grandfather accepted it. After two years, her husband took her to his home and she was there for ten days to perform and participate in the marriage of her husband's brother.

She fainted one day, some six months ago and became very sick thereafter. She went to see a doctor. Doctors said that she has no blood in her system. All her relatives advised her to get treatment in the Stanley hospital. They tested her and admitted her as an inpatient of Stanley. After 10 days of performing all tests they called her husband and said some thing about her disease and advised him to admit her immediately in the Tambaram Sanatorium. She made arrangements for everything and was admitted in the 30th ward. She was so depressed in the hospital that she ran away from there and came home. After she went home she had fever, blisters in her mouth and was unable to eat. Every day she suffered from fever with shivering. She took tablets from the pharmacy and consumed them but she could not be cured. Her husband told her "if you stay here certainly you will die soon" and he brought her once again to the Tambaram Sanatorium and admitted her. She says she has not had sex from the time she has been infected. She does not know how she got this disease.

But her husband told her that he had sex with three or four women before marriage. Her husband is healthy. She says he has lot of blood so he is well and she is ill because she has undergone two natural abortions. Her husband was first tested and was found to be free from HIV/AIDS. She says she never had any relationship with any man before marriage. Her husband used to drink before marriage but after marriage she says he does not even drink. She says she held a lot of respect in his family because she had made him leave the bad habit of drinking.

At times her mother-in-law was kind, at times she was rude. No one knows that she is suffering with HIV/AIDS, to everyone they have said that she suffers from tuberculosis. She says she has seen HIV/AIDS advertisements. She says the depression she suffers is so great that she wants to be killed by some injection.

Her only expectation from the government is good medicine. She says that before the government finds the cure for this disease she may die. She advises youth not to make the mistake of seeking sex, as they are the ones who are affected to the maximum. In some places, women get spoilt and in many other places, men get spoilt. She is a Yadava by caste. Her husband is a Christian and she has changed to the faith of her husband. Her husband had always taken care of her.

She is unaware of condoms. After she was affected by HIV/AIDS he does not have sex with her. She says the food is not proper in this hospital. She has seen HIV/AIDS advertisements. She says she did not know about HIV/AIDS, only after she got the disease she has come to know about HIV/AIDS.

*Additional Observation: This patient's husband seems to be very good, kind and considerate with her inspite of the fact that she alone is affected with HIV/AIDS whereas he is free from it. She does not know how she has acquired the disease. We doubt that when she was treated for the two natural abortions she might have been infected. This is the only possibility as she says she has not had any sex before or after marriage with any one other than her husband. She requests for better food especially good rice since the rice severed in this hospital stinks of bugs and is very old.*

## 75

A 28-year old woman from Andhra Pradesh was taking treatment in the Tambaram Sanatorium. We interviewed her on 14 -10 -2002. She is from Kothampatta in Andhra Pradesh. She was married at 18 years. She did not disclose the name of her parents but she gave the names of her parents-in-law. She is uneducated; her husband who is also uneducated works as an auto driver. She has two children, one has died and the other is just 9 years old. She says the younger child died of fever, vomiting and loose motion. She belongs to the Reddy community. Her son is studying in the 4th standard. She has no property. Her husband used to earn from Rs.50 to Rs.100 a day. She worked as a tailor. She was affected with HIV/AIDS at the age of 26, which was confirmed by blood test.

She says she has never had premarital sex. She says she is unaware of the ways by which HIV/AIDS spreads. She does not know whether her husband had sex with other woman. Her husband was tested to have HIV/AIDS and he too is in the same hospital but in a different ward. She does not know about condoms. Her husband does not drink liquor everyday. He drinks only occasionally. She and her husband were first tested to have HIV/AIDS in a hospital in Andhra Pradesh. They were in the Jarla hospital for 11 days. Some friends in that hospital advised them to go and get treatment from the Tambaram Sanatorium. For the past two years she is taking treatment here.

She has not seen advertisement about HIV/AIDS. She concludes that CSWs should not be allowed to practice for the disease spreads only because of them.

*Additional Observation: She was calm and composed. She advocates people to accept the disease with laugher and smiles. She says that when one is affected by it, becoming sad or worried about it cannot change the destiny. She very curtly told us that she had not made any mistakes and only her husband had infected her. Now her second child had died of HIV/AIDS. Only after the birth of the second child, she too had several health problems like STD/VD. She advocates free treatment for one and all and free medicines for all HIV/AIDS patients. She says that for people like her family where there is no means for living, government should give some help for their survival.*

*She says her husband has become very weak and for the past four years he was not been able to earn any money because he used to feel so tired and never went for driving. She supported the family by tailoring but this tailoring job was only seasonal. Only in festive seasons she got some work, after people came to know that her husband had HIV/AIDS they never gave her cloth to stitch. She changed her very residence because she felt it was impossible for them to live in that area. She said her husband had first lied*

*to her that he had only tuberculosis; so he used to have sex with her. She says that her husband even after knowing that he was affected with HIV/AIDS had sex with her, despite being aware that he would be infecting her. This act of her husband lies very fresh in her memories and makes her very sad and depressed.*

## 76

We met a woman in a near-by teashop on 9-10-02, who is a bonded labourer from Coimbatore. She is just 32 years old, and is taking treatment in this hospital for the past 2 years. She is taking treatment only as an outpatient. Both she and her husband are coolies and work as the construction labourers, wherever they are assigned to work they would go to that place. Both of them are uneducated. Her parents are alive. She does not know the age of her husband. Her parents-in-law have died. She says they died of HIV/AIDS and tuberculosis. She does not know how they were affected by it. However only her neighbours have told her about all this; her husband has not yet disclosed these things to her. She also knows the work of stone cutting. They live in a rented house. She says she had a own house, when her husband was very ill, he sold it for Rs.20, 000. For the past 5 or 6 years her husband is only at home taking treatment.

Further he had an affair with a co-worker on whom he used to spend a lot of money. She says that this coworker had relationships with several other men and had a loose character. She says often her husband used to run away from home, live with that woman for several days and came back only for the money. Just within a year of her marriage he ran away and returned after many months. She was married at the age of 22. It was an arranged marriage. From the day of marriage he was thin and looked tired. Her only child is also affected by HIV/AIDS, he used to have white patches in the mouth suffer from fever and could not walk or stand or talk. It was just an alive mass of flesh. She then took the child for treatment, in the hospital they said it is HIV+. The doctors then tested her for HIV/AIDS. She was well aware of HIV/AIDS from her young age because in her neighbouring house a whole family had died of HIV/AIDS.

She says she was the second wife of this man for the first wife of her husband had already died of HIV/AIDS. Till date he had not disclosed this to her. Even her parents-in-law have not told about this to her, when she once visited her parents-in-law in search of her husband, their neighbours told her all this background. Even her parents do not know that her husband was married once and his first wife had died of HIV/AIDS. Even today she says her parents aged 75 and 65 years are working as coolies in Coimbatore. They are hale and healthy.

*Additional Observation: Her parents take care of her. They visit her regularly. Her husband has not visited her till date for the past 2 years in spite of knowing that she is in this hospital. He is least bothered to take care of their child who is also a HIV/AIDS patient. Her husband takes treatment as an outpatient in the same hospital. She says she has seen him coming to the hospital to collect medicine. However, even on those occasions he had not visited her or their child.*

*This is an instance of helplessness because she was well aware of HIV/AIDS and could not protect herself from the infection. He has not only infected her but he has also infected his child and his first wife. She was married at the age of 22. She says he had taken up such activities only 10 years ago. We are not able to analyze his flimsy relationship with yet another woman a construction labourer. She says both her parents-in-law have died of HIV/AIDS. The whole family is affected with HIV/AIDS for the second generation also, which is very shocking.*

## 77

We met a woman from Nathangudi in Thiruvannamalai district, she is 25 years old and has studied up to the 5th standard. Her husband has studied up to the 10th standard. She was married at 21 years. They have no children. Her husband is a barber by profession. He was her close relative, her mother's brother. The marriage was an arranged one.

He used to drink. He was first tested for HIV/AIDS in Pondicherry and; when he was told that he was HIV + he committed suicide. One year has passed since his death. Just after two years of marriage, he came to know that he had HIV/AIDS. She said a lot of their property was wasted because of the treatment. Several doctors in the beginning stages had cheated them by saying that they can be completely cured of HIV/AIDS. She was also tested only at Pondicherry.

She had a petty quarrel with her mother-in-law but never with her husband. She says she has some property. He was also an Nathaswaram artist. He had a barbershop of his own. First she lived in a nuclear family, when the money they earned was not sufficient they started to live in a joint family that is with her parents-in-law. She suffered from headache and stomach pain. After the death of her husband the doctors in Pondicherry advised her to take treatment in Tambaram Sanatorium. Here, the doctors take proper care of her. The nurses are also kind to her.

She says that in some families when they come to know somebody is suffering from HIV/AIDS they forcefully send them out of their homes. She feels the government can give them some job. She says she has seen HIV/AIDS advertisement on T.V. She does not know about condoms. She says children affected with HIV/AIDS must be given free medicine. She is supportive of the government giving free HIV/AIDS test to one and all. People are very scared about the disease that is they don't talk about it. She always thought her husband was good, only now she is feeling very bad, since he has not told her anything about his affairs with other woman.

She denied having any extramarital sex. She says she did not have premarital sex with her husband. She says only her husband has infected her. She says, "we have become infected, what to do? we should die with the disease." She does not have scabies or blisters in her body. " She says when the CSWs have a mindset to practice that trade, what can we do?"

*Additional Observation: This is a very tragic case. She lost her husband within 3 years of their marriage. She was married at 21 years, he became very sick within a year of marriage. She spent a lot of money to cure him. When the doctors in Pondicherry disclosed that he had HIV/AIDS he committed suicide. It is unfortunate that such things happen. Doctors must be sensitized to treat the patients not only physically but also prepare them mentally before they disclose the disease. There should be support groups to counsel the patients so that such mishaps can be prevented. She is a victim of not only the disease, loss, but also a symbol of total suffering. She did not give the names of her parents or her parents-in-law that is her grandparents.*

## 78

On 21-10-02, we interviewed a woman hailing from the Perambalur district aged 28 years and educated up to the 7th standard. She was an inpatient of the Tambaram Sanatorium. She was married at 15 years. Her husband is uneducated and does combustion work. She has two children; her daughter is studying in the 8th standard and her son in 5th standard. She has tested both her children for HIV/AIDS and the test revealed that they are free from it. She has no hereditary occupation. Her parents-in-law are agricultural labourers. Her marriage was an arranged one. Her husband has never earned any money; only from the earning of her father-in-law he was eating. She belongs to a nuclear family. All the dowry money has been spent for their living. She says her husband was not a alcoholic. She says she prays weekly once i.e., on Fridays. She denies having sex with any man other then her husband.

She is unaware of her husband's habits. She says, "How can I ever ask my husband whether he had sex with other women before marriage?" Her husband was bedridden for some time, so they took blood test and the found he was to be HIV +. She says he used to get big boils, which took a very long time to be healed. He used to say he got these boils because of heat. Doctors take good care of her. She said the life of HIV/AIDS children must be saved. She had symptom of white discharge and STD/VD. She took medicine and was cured. Only after her husband's death the doctors suspected a HIV/AIDS infection and advised her for a blood test.

She says more and more advertisements can be given in the remote villages about awareness of HIV/AIDS. Every one can be tested for HIV/AIDS freely by government. She says as HIV/AIDS comes to a person due to improper behaviour and because it is a killer disease, people feel bad to talk about it.

She has not seen advertisement about HIV/AIDS in TV or behind buses or autos. She says CSWs are "very bad and their activities spoil the lives of family women." She says that she has no right to say that they cannot be in this country because the world has both good and bad people. Her relatives know about her suffering from the disease and they are kind and considerate with her.

*Additional Observation: Her husband infected her and he recently died of HIV/AIDS. She says the doctors said she was affected with STD/VD. She is very sad that her husband has infected her. She has no respect for him; what he earned he spent on his own accord and he did not support his own family. It is her father-in-law and mother-in-law who fed the family and their children. She feels from that he had been spending money on improper things. From her talk it was seen clearly that she feared to question him. He did not support the family with what he earned and he was the root cause of selling all the dowry articles and her jewellery. She thinks he had spent the money on CSWs or he had some other wife and family unknown to her and her parents-in-law. She philosophizes about the profession of CSWs. The patient is worried about the future of her two young children. She is otherwise calm and composed.*

We interviewed a woman from Walahjahpet in Vellore district on 27-10-2002. She has studied up to the 8th standard. She is 38 years old and is living as an inpatient in the Tambaram Sanatorium. Her husband died four months ago. He was an agricultural labourer but took to lorry driving due to lack of employment opportunities in his village. She did not reveal the names of her parents or parents-in-laws. Her husband is uneducated. They have three children: one male and two female. Their son has studied up to 10th standard and now he is employed. Her second daughter is studying 9th standard and her third daughter is studying 6th standard.

The patient was married when she was 20 and nearly 18 years have passed since her marriage. She says her husband used to go on long trips as a driver even up to Delhi. Her parents arranged her marriage. Her husband was a drunkard. If she ever questioned, he used to beat her very brutally. Her husband had sex with CSWs and other women. She says once she was infected very badly by him and when she fought over it, he said he had been to CSWs only once. Both of them were highly infected with STD and they became very sick. They had to take treatment for which they were given injections and tablets. They became cured soon. After that he was hospitalized with HIV/AIDS**.** He took a year long treatment in this hospital. She took care of him but it was in vain. The doctors who came on rounds forced her to undergo blood test and gave her medicines to be taken daily. The doctors take good care of her but do not give glucose or give injections. She feels giddy very often. She knows only Tamil. Her husband had suffered for over 3 years with fever and loss of appetite. She says he has used condoms with her after he was affected with the disease.

She says she never had either pre-martial or extra-martial sex. She says, "How will one ever give jobs to an HIV/AIDS infected patient?" She says several children with HIV/AIDS are treated in this hospital. She says no private hospital will treat HIV/AIDS; to prove her point she said she had seen several specialist doctors but none of them agreed to treat her. She has never seen awareness advertisements about HIV/AIDS on T.V. She has seen advertisements on buses and autos that AIDS is a killer disease. She is of the opinion that more advertisements must be given that too in remote villages.

First she said that sex education must be started from the 1st standard, when we asked her whether children in the age group of 5 could understand about sex, she said one could start the sex education from 5th standard. When the government gives free HIV/AIDS test she thinks people should come forward to take it. They think HIV/AIDS is a 'dangerous and indecent disease' so people do not talk about it.

If government takes care of them and provide them with good medicine it is sufficient for her. She did not have any symptoms of disease. But now she is unable to walk, has vomiting and no appetite. She says if all of them think of god, once in a year or in two years then they visit temples. She has faith in god. All her relatives know that she suffers from HIV/AIDS for her husband had suffered from the disease. All her relatives are kind to her. She says she lives only in a hut and she owns no property. She advices youth to lead a self controlled life coupled with morality. She says HIV/AIDS spreads due to unprotected sex, and use of same injection needles. She says she does not know much about this disease. She says only because of the easy availability of CSWs, men are spoilt.

*Additional Observation: She was composed and spoke very frankly. Her husband had died in the same hospital with HIV/AIDS some four months ago and she felt little upset over it. She is of the opinion that only because of her husband's improper living she has become a HIV/AIDS victim. Being a lorry driver her husband went on very long trips, he became a victim of HIV/AIDS by CSWs. She suffers from STD/VD. She is not able to sit or walk freely. All her relatives take care of her; does not suffer the social stigma that is usually attached with HIV/AIDS, so she is free of*

*many social tensions. Her husband died in the hospital, the body could not be transported to Vellore as it was in a very decayed state, so the cremation took place in Tambaram. So all the relatives had come to the hospital and that is why everyone knows that he died of HIV/AIDS and she is also suffering from it. The patient takes the disease on an easy stride. She is not much concerned about her children, for they are being taken care of by her relatives.*

## 80

We interviewed a Muslim woman on 16-08-2002, who has studied upto the 5th standard. She was a housewife and her husband, was a driver by profession. He had studied up to 8th standard. Her husband died of HIV/AIDS. She is just 38 years old. She says she is a Muslim Rowther. She has two children: one male and one female. Her son is 23 years old and has studied up to 10th standard. Her daughter is just 9 years old studying in 5th standard. She was married at 13, she said her first husband's child is the 23 years old son and her second husband's child is the 9 year old. She married second time because her first husband divorced her. He too got married to some other lady. Her second husband has only infected her because even now her first husband is hale and healthy. Her second husband had boils, which did not heal for 7 years. She says just two years of her married life she was looking after her ailing husband for 7 years. He was taking treatment in GH in Chennai. Then he died of HIV/AIDS. She says he had several other women apart from the CSWs whom he used to frequent. Her second husband has told her about his sexual relation he had with other women when he was in Mumbai. He was an agent for cine actresses. When he returned from Saudi Arabia he took to business and was a popular man in Mumbai. She says most of the cine actresses are the reason for the spread of HIV/AIDS. Their very living is like CSWs. The CSWs and these actress lead a life of dignity but we housewife after becoming infected by our husbands lead a very miserable life full of social stigma. That was the main reason why he was affected with HIV/AIDS. She was first affected by headache, and then by fever, her doctors in the village said she was suffering from typhoid and had given treatment for typhoid.

A doctor named Kamala from Aranthagi said that this disease seemed difficult to be diagnosed so she referred her to the CRC hospital in the Tiruchi Bus stand. There an American doctor checked her and he said she should go for blood test. She spent Rs.1350/- for a blood test; they confirmed it to be HIV/AIDS. He advised her to go to the General Hospital in Madras. But in the GH they said that they have shifted the HIV/AIDS check up / treatment to the Tambaram Sanatorium. She was in Tambaram Sanatorium for 12 days, and then she went home. Then she decided to go for treatment in a big private hospital so she came and got admitted in the ORGS hospital at Adyar and stayed there for 13 days and spent Rs. 60,000. The doctors in that hospital said that they can treat only tuberculosis so they asked her to take treatment for HIV/AIDS in the Tambaram Sanatorium. She took treatment as an outpatient of Tambaram Sanatorium. They gave her medicines and fully cured her of tuberculosis. They told her that she had HIV/AIDS, which cannot be cured but only controlled from increasing so they gave medicine for the same. She took that medicine and as a result she became erratic and mentally distributed and her mind was affected so badly that she was talking and acting incoherently. On her own accord she stopped it. But as she was still weak she went to a hospital in Ratoor for treatment. She was in that hospital for 10 days and then returned home. Still she did not feel well and could not eat. So she came once again to Tambaram. She feels this is a very cruel disease. She says the doctor had made medical observations only once. The nurses come and give medicines. She says there is no medicine in any other hospital, however much they spend for it. She says small children affected with HIV/AIDS can be given free education and free medical treatment. She said she fainted as soon as she came to know that she was infected with HIV/AIDS. Even her late husband wept for her. The whole family was in tears. She has not yet tested her children for HIV/AIDS. All her relatives are kind and considerate with her, they do not discriminate her. They wash her clothes and stay with her and take care of her very well. Only her mother-in-law used to scold her and thought her son was very good and he has died only of jaundice. She has seen HIV/AIDS advertisement on TV, at first when she saw it she was very much frightened. Now she says she is used to it, some advertisements are good some are degrading. She has a house given by her parents. Her husband educated his brother up to engineering and married off his two sisters, and he did not give anything to her except HIV/AIDS disease. She advocates free HIV/AIDS test for one and all. She says by looking at pornographic books a good person will correct himself but a bad person will become worse. She says CSWs are always close to death. She says it is very wrong to entertain men, using sex for money. She says that HIV/AIDS can spread by touch, by kiss, by saliva, by food etc. So only she has kept herself separate brush, soap, towel etc. But she says her relatives do not believe that. She is not praying to god for the past 1 year.

*Additional Observation: This woman spoke very frankly. She is from Aranthangi in the Pudukottai District. She is married for the second time. Her first husband divorced her. Her second husband was a HIV/AIDS patient. He has contracted the disease due to his careless, high-risk behaviours a result of his pompous living. He has connections with rich cine actresses. He has had sex with several women and CSWs at Mumbai. Just after two years of marriage her husband developed boil that did not heal for 7 years and he had taken treatment for the same for 7 years in GH. Only in GH they confirmed that he had HIV/AIDS. So they also warned her to take test for HIV/AIDS. She found herself to be affected with HIV/AIDS. She is contended. She says she lost control of mind when she first took medicines from the Tambaram Sanatorium. This throws light on the observation that some people when treated for HIV/AIDS become very depressed may be due to*

*the medicine. She had many wrong notions about HIV/AIDS. We felt it was our only duty to interview so we did not intervene or argue with her on any matter. She is against half-clad women acting in cinema.*

## 81

On 25-10-2002 we interviewed a 25 year old woman who has studied up to 12th standard. She has no job. She is married. Her husband is a barber who has studied up to 8th standard. He has a shop of his own in Villivakkam. He was her close relative and their marriage was an arranged one. She was married at 19. She says doctors just come on rounds and the hospital premises are very unclean. They don't get enough water in the toilet and bathrooms. She says soon after her marriage, that is, within 3 months of marriage he was very sick. They have a son who is free of HIV/AIDS. The main symptoms that she suffered was fever and headache. Her husband was a smoker, drunkard. He used to give his pay to the family only once in a way. Her husband when he was affected with HIV/AIDS he said "I don't know somehow I have caught this disease". He has not disclosed to her how he was infected with HIV/AIDS.

*Additional Observation: She suspects him because he never used to give his pay to her or say how his business was. He was unquestionable. He once said that he has got the disease by injection needles, she was not able to know how he was really infected. She is very calm and composed. She did not talk about her parents or her parents-in-law. It looks as if she is ignorant of the ill effects of this disease. She has not seen any HIV/AIDS advertisement on TV or buses. She is very aloof. She has no respect or regards for her husband and at the same time she does not suspect him of visiting CSWs. Her husband was not bothered that his wife would get HIV from him. That is why even after he knew that he had the disease he carelessly had sex with her. That was the reason why she became seriously ill. She feels*

*bad about her husband who despite knowing the consequences had infected her.*

**82**

On 07-11-2002 we interviewed a woman from Mumbai, affected with HIV/AIDS. She is 38 years old and studied up to 10th standard. Her husband is uneducated. For the past 7 to 8 years he is working as a cook in Mumbai. She has been in this hospital for a year. Doctors have confirmed that she has HIV/AIDS. He was first tested for HIV/AIDS in Mumbai and he was confirmed to have HIV/AIDS. She says he was so upset that he even tried to commit suicide.

She says it is impossible to come to know from one how they become infected. Only very few say the truth. She says her husband has all bad habits. He used to spend around Rs.500/- a day very carelessly. His parents know that he is suffering with some disease but they do not know that he is suffering from HIV/AIDS. Her parents know that she and her husband are suffering from HIV/AIDS. They could do nothing, because they are Malaysian citizens. He has taken treatment from Kerala for cure but it had no effect. Further when he was tested in Kerala the result was negative, but when the test was carried out in Tambaram he was declared HIV +. The symptoms of the disease that they suffered from was fever, cold, scabies in legs and itching in genital areas. She says, "if someone says he has HIV/AIDS, he is removed from the job. How will anyone ever provide jobs for HIV/AIDS affected persons?" She says she is from a middle class family. She advices that only that a self controlled life can save people from HIV/AIDS.

*Additional Observation: She did not give answer to the question whether she has children. The fact that her husband spends Rs.500/- a day on filthy things like cigarette, CSWs etc., means that he is very well off. She says he has all bad habits. She is sad about that. Her husband is secretive since he has not disclosed the disease even to his parents. She consoled her husband who was about to*

*commit suicide when he came to know that he was infected by the disease. She is courageous and contended. As her parents are in Malaysia she is not taken care of by any one. We are not able to know whether her husband is an inpatient or an outpatient of this hospital for she was all-alone and there was no one by her side to help her. She is affected by STD/VD. She had lot of scabies over her feet and she found it difficult to walk with it.*

### 83

We met a 38-year woman now an inpatient of the Tambaram Sanatorium on 2-11-02. She has studied up to 8th standard. She is from Vadalur in Cuddalore district. She belongs to the Kavara Naidu community. Her husband is also educated up to 8th standard. Her parents are alive but her parents-in-law are dead. She gave birth to 6 children. The first two were born and died after some time. Two are living. The last two children also died as infants. She has married off her daughter who is just 19 years. This patient was married at 13. She attained her puberty six months only after her marriage. Her husband was a drunkard, a driver by profession. After the birth of the last child he never came home regularly. He used to visit her only once in a while. She does not know where he is now. For the past six years he has not visited her. She says she has her own doubts on whether he is living or dead.

She was unable to keep her daughter so she married her off. She has a grandson. She says for the past 1½ months she is very ill with cold and fever. So she took treatment from one doctor named Saravanan. He thought she must be suffering from typhoid. She was treated for 15 to 20 days for typhoid. When she did not become better he tested her blood and he was very sad to inform her that she was affected by HIV/AIDS.

The doctor said "a disease one should never get has infected you." He advised her not to worry and remain courageous. He wrote a letter to the Tambaram Sanatorium with her blood test report and sent her to this hospital. Her marriage was an arranged one. She never wanted to marry because she was just studying 8th standard. She used to enjoy life by going in a bicycle. But they married her by force. Her mother wanted to appease her brother so she married this girl to her uncle (mother's bother).

She used to sell bangles, beads and other plastic goods and earn a living. When her husband used to visit her used to give Rs.200 to Rs.300 a month. He was a drunkard, he used to drink all the time. He used to carry loads of Bamboo to different places. She has no hereditary property. She did not doubt him at first. Later her neighbours used to say that her husband was having illicit affairs. That is why he used to visit her once a month then once in three months then once a year after wards for the past 6 to 7 years he is never seen or heard. She has never asked him about how many women he had for sex.

She says even if one asks such questions the husband will not reply. She says he always suffered from STD/VD for he had blisters and scabies in the genital areas for which he used to take treatment very often. He has beaten her many a time. He used to get drunk fully and say the food was not proper and throw away the food and fling things in the house.

She says she never had sex with any other man. Now she is living alone with her daughter. She has so far not tested her daughter for HIV/AIDS. She lives in a rented house. She knows both Tamil and Telugu. She used to pray only on *Ammavasai* (new moon day). She says she has seen all the advertisements regarding HIV/AIDS. When she came to know that she had HIV/AIDS she was very much affected. She told the doctor that she had lead such a frugal life, she had even ran the family with Rs.10/-, so she asked the doctor how could she be affected by it. The doctor consoled her by saying that her husband should have lead a characterless life since he was a driver. The doctor said this disease does not show its symptoms immediately it takes time to spread. He just consoled her and sent her to Tambaram Sanatorium.

She has no relatives or friends. She does not know how the disease spreads. However she is happy that the doctors are kind with her. She says everyone can be tested for HIV/AIDS by the government. She says she is very frightened to see the HIV/AIDS affected people. She says she would stay in this hospital for ten or fifteen days and go home soon.

She says if the government cures her that is enough for her. She says she feels relaxed after talking with us. She says only self-control and morality can save the present day youth. She says her son-in-law is a wonderful man. None of her relatives know about her sickness. Her son-in-law visits her but she has not told him that she is affected with HIV/AIDS.

*Additional Observation: Yet another instance of child marriage occurring in the 20th century. Even before she attained her puberty she was married to her uncle without her consent. He was drunk all the time and he subjected her to domestic violence. He did not come home regularly. Now the pity is that she does not know whether he is living or dead, for he has not returned home for the past 6 to 7 years. She adds that only three years ago her daughter's marriage was performed even for the marriage he did not come. Their friends and relatives have not seen him yet. No one knows where he is.*

*Her husband has infected her. Except for the doctors, no one could have consoled her. She is sad and worried and wants to go home soon for she says she becomes very much dejected when she sees the sufferings of the other HIV/AIDS patients. This feeling varies from patient to patient for some feel consoled and some dejected. Her parents too do not know about her disease; she says they are very old and she does not want to trouble them. Only after he ran away from the family she became self-employed by selling bangles, beads and small plastic things to maintain herself and her two daughters. Soon she wants to get her other daughter also married off. She says her only visitor is her son-in-law.*

We met a woman hailing from Apparkulam in Tirunelveli on 2-11-02. She has studied up to 4th standard and is married. She says her husband is illiterate. They work in sand quarries. They belong to the Christian Nadar Community. They have two children.

She did not give the names of her parents or their age but she gave the names of her parents-in-law and their ages. Both of them are alive. They had been to Mumbai as sand quarry workers only. Both her parents are based in Mumbai and she was born there. Her first child is 3 years old and her second child is just 6 months. She does not know the cause of her getting the disease. However she says she suspects that her husband should have infected her. The red light area was a curse for most men get infected from that place.

She had fever with shivering, inability to eat, bloating of stomach and itching all over the body. She is taking treatment in this hospital for the past four months. She tested her blood first in Pondu in Mumbai. The news was broken to her by her mother; she was very unhappy to see her daughter affected with HIV/AIDS. All her relatives know about the disease. Both her parents are very kind with her. Her friends are a little distant they talk very little, they do not move closely with her as before. Her mother is more upset than her. Her father-in-law is dead. Her mother-in-law causally said "go to your parents' house, they alone can take care of you".

She was married at 21 years. Her marriage was a love marriage. He is not her relative, he is her neighbour. He became her friend when she was doing the sand quarrying work. He too used to do the same job during the time of marriage. Her husband suffered only from cold and cough. She says she knew about HIV/AIDS. She knows how it spreads. She says the HIV/AIDS affected children must be given good education so that they can stand on their own feet. The world will not become spoilt by the discovery of HIV/AIDS medicine; for she says those who commit the mistake will continue to do so, but it will not affect the good or bad. She thought of committing suicide as soon as she came to know about the disease. All her neighbours told her to be courageous and never to fear for it. She says government should find a cure for this disease.

She does not know how to spread awareness about HIV/AIDS among villagers. She says she has never been having sex before marriage. She says she had heard about her husband's immorality. She says he promised her that he would correct himself after marriage. She married him in spite of knowing his bad qualities and weaknesses like the fact that he drinks, smokes and goes for cinema with CSWs. Before marriage she had worked in an export company. She said she used to earn Rs.1200 per month. She lives in a tiled house, which is her mothers'. After her marriage they lived in a rented house paying a rent of Rs.600 p.m. If he works for a day he would get Rs.100/- He was a daily wager labourer. She had never taken indigenous medicine. She welcomes free HIV/AIDS test for one and all. She believes that god will cure her.

*Additional Observation: She must be in the age group 25 to 30 years. She is thin and weak. Both her children are affected with HIV/AIDS. Her husband is infected with HIV/AIDS but has he concealed to her about his taking the treatment. Her mother-in-law's causal behaviour has pained this girl. However it is nice to note that her parents are very kind and considerate. She says she has committed the mistake of marrying him in spite of knowing very well about his lack of good character. It is very sad to see that all the four in the family are affected with HIV/AIDS. It could have been better if her husband and mother-in-law gave them at least moral support.*

*This is the first case among the interviewed women who expressed their desire to commit suicide. Her parents are supporting her both morally and economically. She says Tambaram is the best hospital for HIV/AIDS treatment. She concluded the interview by saying that she could not get this type of treatment from any hospital in Mumbai.*

We met a married woman, who is a construction labour on 16-01-2003. She is uneducated. She is from Chengam in Tiruvannamalai. She has only one daughter. Her husband is also uneducated. He was a *maestry*. He died of HIV/AIDS. She says her husband used to smoke, drink and go to CSWs for sex. He was a domineering man. He had the habit of visiting CSWs and other women, both before and after marriage. She hopes she will be cured of HIV/AIDS. She is a Gounder by community, Christian by faith.

*Additional Observation: Her age must be around 32-34 she but looks very much older, weak and withered. She says it is because she works in cement and sand all the time that too in the hot sun so only she has lost her complexion and skin texture. She has scabies in her hands and legs. She looked unattended with matted hair. Her mother-in-law or father-in-law are not interested in supporting her after the death of her husband. She is not taken care of even by her parents. She is left all alone. We are not able to understand from the interview who is taking care of her only daughter. Her husband was a maestry who died at the age of 40-42, leaving behind this woman and their child. From her talks, it is revealed that he had many wives for she claims herself to be one among them. She had married him despite knowing his faulty habits. They (these supervisors) do not marry the co-workers who work under them, but misuse most of them. She too is a victim of this kind.*

**86**

We met a 21-year-old agricultural labourer on 16-1-03. She has studied up to 10th standard. Her husband is also an agricultural labourer. They are from Royapettah, Chennai. She refused to give the names of her parents or her parents-in-law. She was married at 16 years. Her first child died soon after birth. Now the second child is just three years old. The child is free from HIV/AIDS. Her husband was found to be infected with HIV/AIDS. The symptoms that she suffers from are swelling of the knee joints, boils in the legs and swelling in the feet due to which she was unable to walk. Doctors, however suspected her to have HIV/AIDS and tested her blood and confirmed her to be infected with HIV/AIDS. Both her parents are with her. Both of them are very kind and affectionate with her. Her husband is hale and healthy though he is HIV +. He gets medicines from the hospital as an outpatient and takes them regularly.

*Additional Observation: She looked very thin and had a big boil on her swollen feet. She was suffering with pain. She was not even able to walk or fold her leg. There are scabies all over her body. She was not depressed or worried about her child for she does not understand the seriousness of the illness of her child. She has taken it lightly. Her husband was aloof. He does not visit her now. She suffers the pain all the time, as the boils have not healed for a long time.*

### 87

On 15-04-2003 we met a 40-year woman. She is uneducated and sells books for a living. She is married. He is also uneducated and sells books. Her parents are dead. She is a Christian. She is from the Tirupathi district. She has no children, she prays to Jesus and has belief in Christianity alone. She does not have any property. They earn from Rs.100/- to Rs.150/- a day. She was married when she was just fifteen years old. She has been married for 25 years. Since her parents were dead she got married on her own. Her mother died because of the constant quarrels with her husband who was a drunkard (Her father). Her father's brother brought her up, He was a weaver by profession.

Her husband belongs to Padappa. She says doctors are rude with her. She says some job should be provided to HIV/AIDS patients. She says government should provide good medicine to them. She says she had sex with her husband even before marriage. She says she has this disease from very small age. She says she had heart pain, fever and typhoid. From typhoid only, it has become HIV/AIDS. She has pain in the bones. She was also tested in Tirupathi. They did not say properly to her about the disease. She is affected but does not bother about it. Doctors say she should stay for a month in the hospital but her husband asks her to come home with him. She had boils in her private parts that did not heal for a long time.

She says she gave birth to 3 children and all of them died soon after birth.

She says she has not seen HIV/AIDS advertisements in her village but she has seen however HIV/AIDS advertisements in Tirupathi. She says she does not know how HIV/AIDS spreads.

She says people will mock at her if she speaks openly that she has HIV/AIDS. She says 50% of those infected are cured and 50% cannot be cured of HIV/AIDS. She says HIV/AIDS first comes as small blisters, then it becomes scabies and she knows a woman who died of it. She feels sad, since she is not able to sit on her own.

***Additional Observation****: Only she is infected by HIV/AIDS, her husband was tested to be negative. However this is yet another instance where the husband has taken care of her and admitted her in the hospital and continues to take care of her. She said that even when she was very young she had this disease. She is very thin. She suffered with boils in her genitals also. She had a very unhappy childhood because she lost her parents when she was very young and she had to choose her husband on her own. From her talks it is revealed that she has been abused even in her childhood. But she cannot be HIV infected from her childhood for she had lived up to 40 years. She looked lean, weak and frustrated.*

### 88

On 26-04-2003 we met a 36 year woman, an agricultural coolie by profession. She is in the Tambaram Sanatorium as an inpatient. She is uneducated. Her husband is also an illiterate agricultural coolie. She is from the Tiruchi district. She did not give the names of her parents; she however said that both her parents-in-law have died. She further said that she was married at 13 years. She has two sons and one daughter.

One son is married and her daughter has passed 10th standard. Her eldest son studied only up to 8th standard. Now her youngest son is studying in 9th standard. Her husband died of heart attack 10 years ago.

She says her marriage was an arranged one. She has no hereditary property. Her husband was an agricultural labourer. She says all her dowry property was sold away to take medical treatment. She says her husband was a drunkard. She often suffered from loose motion. A major hospital in Tiruchi referred to Tambaram Sanatorium. She says as a widow she had never been out of the house. If at all she went, it was only to work as agricultural coolie. She has not seen TV. She says even if government gives a free test for HIV/AIDS, people will not accept it. She said she was bedridden with chest pain.

*Additional Observation: She is affected by HIV/AIDS but her husband had died of cardiac arrest . This is also an instance of child marriage. Her husband was a drunkard. She has spent a lot of money for his treatment. She is totally ignorant of HIV/AIDS. She said she had no feelings of sorrow when doctors said she had HIV/AIDS. She says "live till we can live and die when we have to die" is her philosophy.*

*She suffers only from loose motion, she says she does not have any other disease or symptom. She refused to advice the youth because she says she has not lead any bad life. On the whole we see that she is highly innocent, ignorant of the disease, of how it spreads and how it has a social stigma associated with it, so she is fully free of all tensions and does not think anything about it. So is contented and happy, though she looked tired and could not sit for a long time.*

We interviewed a 27 year old illiterate woman on 26-04-03, from Uraiyur in Tiruchi District. She is a Mudaliar by caste. She has no job. Her husband has studied up to the 12th standard. He is a painter by profession. She did not give her parent's names or parents-in-law names. She has only one daughter who is four years old. A blister on her tongue did not heal for over 90 days. She says doctors don't give treatment to her. She says she simply eats in this hospital. Doctors come on rounds only once a week. She says she lives in a rented house because she has no hereditary property. She has no father only her mother is taking care of her. She was married at 23 years. Her marriage was an arranged one. She says her husband does not have the habit of smoking or drinking alcohol or visiting CSWs. She says her husband is a very good person. She says she doesn't know whether he did any mistake without her knowledge. She says she has no extramarital affair with anyone.

She says very often she used to go to the GH to get injection for general weakness etc. She suspects if she could have been infected in the hospital. She has checked the child and it is free from HIV/AIDS. Because her husband had HIV/AIDS so only she was tested for it. She was not even able to sit, only now she is a little better. She lives in a rented house by paying Rs.500 per month as rent. Only if he goes out to earn he will get around Rs.50 per day, otherwise no pay. According to her, doctors don't visit them only nurses come and give medicine.

She feels that if the doctors come and test them every day it would be better. She says they come once in a month (she is contradictory for earlier she said they visit once a week.) She says she has seen the doctor only once so far. She says most of the HIV/AIDS patients cannot work continuously, if they work for 10 days they will become bedridden for 10 days. She is of the opinion that even if the government gives free check up for HIV/AIDS she doubts whether the people would accept it. But however she says that by checking one can take treatment before making the disease chronic. She says men don't correct themselves by seeing advertisement; even after getting HIV/AIDS most men do not correct themselves.

She says several men in this hospital who are HIV/AIDS patients go out in the night to visit CSWs. She is not aware of how one gets HIV/AIDS. She suffered from lack of appetite and she lost her long hair. Then her husband went for a free HIV/AIDS test in a Government hospital and doctors in that hospital said he was HIV positive. The next day she paid Rs.350/- to check her blood for HIV/AIDS. She too was affected with HIV/AIDS. The doctors referred them to the Tambaram Sanatorium. She says she is taking medicine for the past two years. She says she felt sad when she came to know that she had HIV/AIDS. She says that after coming to this hospital and after seeing hundreds of HIV/AIDS patients she is least worried about HIV/AIDS. People fear this disease, and they fear even to talk with a HIV/AIDS patient. She says that she has once seen a pornographic book just after marriage. She says she feels if one views the book in a wrong perspective it is wrong but if one views it in the right perspective it is right. What she feels is "we can see and know and learn to avoid bad things."

When asked about the difference between HIV and AIDS she said "In Tamil they call it as AIDS and in English they call it as HIV". She says she does not know about these. She says one should get married only after the HIV/AIDS test. She advises the youth not to get married for money, high dowry and many sovereigns of jewellery. If they seek for these, they have to be admitted in the hospital, what they should seek is HIV/AIDS test before marriage. She says the CSWs work as CSWs when they are young after that they are the most affected ones because their hands and legs shrink and they lose all their beauty and suffer unbearable diseases. She says all CSWs should be called and injected with poison because they serve for money. She has used condoms after getting this disease. She advocates more awareness programmes in remote villages.

She says if they give advertisements once in a while people forget it.

She says on Sundays' a doctor or an expert should talk about HIV/AIDS for half an hour on television. Only when people hear and see for a long duration it will stay in their mind. She says only the patient should be responsible, take the medicine properly and live courageously. She says others cannot do anything for us. Advertisements in the daily papers about HIV/AIDS must be given. She says "bring all the college students and show them the HIV/AIDS patients in this hospital only then they will understand, till then they will not follow anything about HIV/AIDS." She finally adds that she paid lots of dowry and 10 sovereigns of gold but all of it has been spent on medical treatment.

*Additional Observation: We have interviewed over 500 HIV/AIDS patients but this interview was the best in terms of suggestions especially to youth and about the awareness program. She talked with us for over 50 minutes. She gave good ideas about HIV/AIDS infected children. She gave consultation on what a HIV/AIDS patients must do.*

*She is just 27 years old has studied upto 5th standard but her approach to solving, studying and analyzing of the HIV/AIDS problem is very analytical and realistic. About awareness program she feels that the advertisement which is flashed in T.V. for a minute or two about HIV/AIDS cannot certainly stay in their minds for they would forget when they see several other pictures, songs and dances; so she says on Sundays when everyone is free they should talk or display a movie about HIV/AIDS for half an hour or more everyday so that people see it and start to think about it.*

*She advocates blood test for HIV/AIDS for both bride and bridegroom before marriage. For the HIV/AIDS affected patients she says, "courage is important".*

*Secondly they must properly take all the medicines given by the doctors and third, she says, "Several men who are inpatients of this hospital go out in nights to have sex,*

*by doing so they are not only infecting others but they become more and more vulnerable to the disease." Thus this lady's view points moved us.*

## 90

We interviewed a woman from Kathanaddappu near Nagari in Tirunelveli district on 30-05-03. She is just 27 years and she has studied up to 5th standard. She is a Hindu Nadar. Her husband is a car driver and has studied up to 10th standard. She has two female children, aged 9 and 5. She has three brothers. They have not studied much. She was married when she was just 17. Her marriage was an arranged one. Her husband used to consume alcohol everyday. She does not know whether he has extramarital sex with other women.

She feels that she cannot ask such questions and he too has not said anything regarding this. She says God has given her this disease, already her husband is suffering from it. He was first tested for HIV/AIDS in Nagerkoil. She says that some people when they know they are HIV/AIDS infected, attempt suicide and some suffer from fever with shivering. She did not do any such thing. She says she does not know from where she got HIV/AIDS from Mumbai or from Tamil Nadu. We do not know because both of them were in Mumbai sometime. She was only a housewife, only her husband was a cab driver in Mumbai.

When asked whether her husband has ever been with other women she says, "Who can ever control a man?" She says that when she came to her mothers' residence for delivery he should have been leading a lose life. Because when she asked him to go for a HIV/AIDS blood test he shouted no to her but when he got high fever he took a blood test but did not inform her. She says he is in the same hospital taking treatment in ward number 20. He is now feeling a little better.

She said she had fever with shivering, and a strange kind of skin disease in her genital areas and cough. Nothing could be cured so she came to this hospital. She says after visiting this hospital, she is feeling better. She is fully unaware of what had happened to the dowry and gift articles given by her mother. She is from an upper middle class background and is not very rich.

She says her parents-in-law are very bad and demanding, a house was to be given as dowry. Now her parents have pledged that house and are paying for her treatment. She says both the children are affected with HIV/AIDS. Her parents are going to take care of them.

Now the related story is that this lady's parents are also infected with HIV/AIDS. Her mother is taking treatment for the past 6 years and her father is taking for the past 8 years. Her mother who intervened in between our interview was also a HIV/AIDS patient. She says that the advertisement about HIV/AIDS given in villages is very meager and has a little or no impact on viewers.

She says several men go behind any women not knowing about HIV/AIDS. So if proper HIV/AIDS advertisement is given it would be better and it could prevent men from doing this type of mistake. She supports the proposal of free HIV/AIDS check up that will be conducted by the government. She says, "it is not wrong, but will the public accept it?" She says, "it is a cruel disease and is incurable." She says she has not seen any pornographic books. Further, she adds, "it is only because of the CSWs that this disease is created."

She doubts whether the government will give jobs to HIV/AIDS affected persons. She says from the government she expects money, medicine and cure.

*Additional Observation: This is a very unusual case. We visited this 27-year woman in the hospital to know her problems. To our shock not only she and her husband and their two children were HIV/AIDS infected but her mother and father are also HIV/AIDS patients. All the six are in the hospital in different wards. The son in law at first did not accept that he had HIV/AIDS but when he became very ill with fever and stomach pain, he was advised by the doctors to take an HIV/AIDS test and was found to be HIV positive.*

*It is a cruel and a very sad thing to note that their two children are also affected by HIV/AIDS.*

*Three generations are affected by HIV/AIDS in that family; the daughter is very ignorant of HIV/AIDS and her mother who had been suffering with HIV/AIDS had not informed it to her due to the fear that her daughter's husband and his parents would ill treat them. She said, "After spending so much of money, including a house as a dowry, I did not want the marriage to be a broken one. Only on that basis, I did not inform my daughter." She said during the delivery of the second grandchild itself she had suspicion because her daughter suffered a strange type of skin disease, which remained uncured for many years in spite of treatment from many hospitals. Thus she says that when she saw her daughter sick with HIV/AIDS she was totally broken.*

## 91

We met a 23 years old woman in the Tambaram Sanatorium on 30-05-03. She is married. Her mother's sister takes care of her, because she has no parents. Four years ago they checked her in the Gori Meru hospital and found her to be HIV +. So far she has spent 1½ lakhs for cure. She denies taking indigenous medical treatment. The patient has a son. She was married at 16 years and infected at 19 years. Her marriage was a love marriage. Her husband died of HIV/AIDS recently. She came to know that he had HIV/AIDS; only after that she found herself to be infected. No one in her village knows about this. She said they married the patient as a child because of poverty. They belong to the Vanniyar community. They live in a hut. Only if they go for jobs, they can eat. There were no symptoms. She says she has no faith in indigenous medicine. Before coming to the hospital they had spent Rs.60, 000. She has not changed her religion. She says jobs can be given to HIV/AIDS patients both in the private and government concern. According to her they don't treat a HIV/AIDS patient properly because it is a disease associated with immorality so it is considered to be indecent, so they attach a social stigma with it.

*__Additional Observation__: This woman was interviewed on 16-04-02. She is a nurse by profession. They are from Tanjore. Her marriage was a mishap. She has a male child; nothing was said about this child's health status. We could not understand whether they have tested the child for HIV or not? She avoided the answer to this question. Later she said the marriage was a love marriage and her husband had died of HIV/AIDS. She has spent over 1 ½ lakhs for cure. She says that only in Dr Deivanayagam's time they had regular medicine and good food. Now everything is topsy-turvy. She said this with tears and complained that in the present day the doctors are giving them rice that smells of bugs and is old and cannot even eat a mouthful. The milk they offer is also very watery. Earlier they used to get buttermilk up to three tumblers now they provide only a small cup of buttermilk. She concludes that on the whole the hospital is unclean now and food is of very poor quality.*

## 92

On 01-06-2003 we met a 32-year old woman in the Tambaram Sanatorium. She is a flower vendor by profession. Her husband is also a flower vendor. She is married for the past 20 years. She has 3 children. She says after the birth of her last child, she contracted HIV/AIDS. She said that her husband has infected her with HIV/AIDS and the children are free from HIV/AIDS. They live in a rented hut paying a rent of Rs.150/- p.m. They underwent blood test in the local hospital. She is a native of Sindhoor pettai. When she heard that she was infected with HIV/AIDS, she said she began to recollect how she would have got it. She says this disease was first discovered in America. She says she does not know much about the disease. She says that none of her relatives know about her husband suffering from HIV/AIDS.

*Additional Observation: She is only thinking of how she was infected. Her husband too is infected. She does not blame her husband. She is from very poor strata of the society. Both are flower vendors. She was married at an early age of twelve, even before she attained puberty. Her marriage was an arranged one. She is very ill and is unable to sit or walk. Her husband is a little better now than before and he takes medicine regularly from this hospital. He is an outpatient of this hospital. She and her husband are uneducated, we also come to know that her children have not stepped into the school premises. Both her parents and parents-in-law are alive; they are hale and healthy. No one is there to take care of her, whenever her husband visits the Tambaram Sanatorium to collect medicine for himself he just visits her with some food because he still sells flowers to earn a livelihood.*

## 93

We met a native of Tanipadi in Tiruvannamalai on 16-06-03 at the Tambaram Sanatorium. She is a Vanniyar woman who has studied up to 9th standard. Her husband is uneducated, he is a lorry driver by profession. She has three children aged between 6 and 12. She is just 27 years old. She says she was married at the age of fourteen when she was just studying 9th standard. She discontinued her studies soon after marriage. Her marriage was arranged. Her husband is from Natanthal near Vettavalam. She says that before marriage she never had sex since she was studying. Her husband had many illicit affairs and he left her in lurch; 5 years have passed since he last visited her. He abandoned her. The children are now living with her husband. When he left her, she was living with his brother-in-law for one and a half years. Now he too has got married. She says she was very ill due to fever and for the past 6 to 8 months. She says after she has been admitted in this hospital she feels a little better, though she is often laid with fever.

When asked of how she would have contracted this disease she says she had sex with an auto-driver who was very young, she thinks she might have got the disease from him. She says she had sex with this man only after her husband left her. She suspects that he should have infected her. Her husband and children are hale and healthy: none of them are infected by HIV/AIDS. But her friend the auto-driver is very sick. He met with an accident and his legs and hands have shrunken and are immobile. He is suffering from fever and is bed-ridden.

She says that she felt sad when she came to know she was affected with HIV/AIDS. Her husband does not allow her to see their children; but news has reached her that they are doing well. Nobody knows that she is affected by HIV/AIDS. Only her mother knows about it. Her mother takes care of her but she discriminates her, i.e., does not eat or drink in the same vessels but washes them before use. Her mother never drinks the water, which she drinks. She says she is not even able to sit. She has seen HIV/AIDS advertisements on TV and new papers. She says she has seen and heard about

1. Transmission of HIV/AIDS from the mother to child in the embryo stages can be prevented.
2. No medicine has been so far discovered to cure HIV/AIDS completely.

She says once she was hurt in the head and in CMC they gave blood. She wonders whether that could have infected her. She says HIV/AIDS advertisements must be given in the villages because the villagers keep the disease as a secret and thus it is left untreated resulting in death. They fear the social stigma and the consequences if others know they suffer this disease. She says free compulsory HIV/AIDS test may make people to take treatment and thus lessen the social stigma. She says by reading pornographic books the reader becomes a victim of illicit sex and thereby there are chances of his getting HIV/AIDS. She says the CSWs are persons who take up the profession because of the easy money and she says she cannot do anything about them. She says if this information is passed onto others at least they will be careful.

She says that this morning she had taken a larger dose of medicine so she was feeling giddy; only now she just opened her eyes and has started talking with us. She says people don't talk about HIV/AIDS since it is disgraceful to have got the disease from CSWs due to ones own immoral character. She has not made daughter's her husband check for HIV/AIDS.

The patient's mother visits her daughter occasionally. Her mother says first her son-in-law broke her head and it was stitched in the CMC Hospital. Her daughter had a fight with her husband so she drank the paste of some poisonous seed and they spent over Rs.2000/- in Thiruvannamalai and saved her. So her mother says that she has nothing to comment upon. She says that her son-in-law drinks only beer and brandy but she has only heard this and has not seen because she is 20 miles away from her daughter's family. She says that her daughter has not said anything because she visits her daughter once, or maximum, twice a year. She says her son-in-law has taken everything and has left behind only his wife. He is working as a driver. She sometimes sees her grand children.

Now the stranger part of this life history beings when we causally asked the patient's mother of her husband is she began to narrate another interesting story. She says that her husband has left her and he is living with another woman. He has sold all her property and she says that he has taken TV, Mixie, Grinder and all other things to the other woman's house. She is also in Thiruvannamalai. But he pays her (the patient's mother) Rs.2000/- a month. That lady has a daughter who has attained puberty. She has questioned him several times and he has also beaten her after which he ran away from home. She says that her husband's current lover was by profession a CSW, and she suffers from HIV/AIDS. So he is supporting that woman and her child. His real wife says she does not see him. Recently when he was driving he met with an accident so she went to see him; as soon as he saw her he closed his face with the blanket. She opened his face and looked at him and cried for some time. That woman with whom he is currently living also did not say anything so she came off.

*Additional Observation: This is a very strange case. The woman is infected by HIV/AIDS, but husband is free from it. In fact she is the mother of three children. He is supporting the children and educating them. However he does not allow her to visit the children. He has not married any other woman.*

*On the other hand, this woman has alone lived for 2 years with her brother-in-law. Only after he too got married she was living alone. She lived scandalously only for a month or so, then she became a friend of an auto-driver who was much younger to her. Only when she was living with him she became very sick and was infected with HIV/AIDS. She says that now her (auto driver) friend is also seriously ill with immobility of the legs and hands, he has been hospitalized with fever. Thus it is learnt that her husband and children are free from any disease.*

## 94

On 18-6-03 we met a woman who was educated up to 8th standard; she is now an inpatient at the Tambaram Sanatorium. She is 43 years old, taking treatment for HIV/AIDS for the past 14 years. She is an agricultural labourer by profession. She was married at the age of 17 years to a 6th standard educated man.

For the past 9 years she has been staying in Guindy, Chennai. She is a Naidu. She has three children one son and two daughters. The son has studied upto 10th standard. She says her son is 22 years; the first daughter 18 years and the second daughter is just 16 years. None of them are married. They are idle at home. She knows both Tamil and Telugu.

She goes more to the church than to the temple because she feels that she gets more consolation from the Christian god than from the native god. She says her marriage was arranged by her parents. Her husband is her mother's brother.

Even before marriage he used to drink. The family elders and her grandparents said that if he was married he would change and become reformed. It is her destiny that she had to marry him.

She says that she knew her husband had sex with other women even before she married him. She undertook a hysterectomy (family planning operation) at the Vellore Government hospital. They said at that time (14 years ago) that her blood was infected with HIV/AIDS. She has not told about this to anyone including her husband.

She was not able to understand what the doctors meant by the statement "blood is spoilt". She has kept this in her mind.

They (the doctors at CMC hospital in Vellore) tested the blood of all her children, they said that only for the last daughter they were not able to finalize whether she is infected or not. She wants to bring her here and carry out the blood test for her. Her husband was staying only in Vellore only now he has come to Guindy.

The doctors at the Tambaram Sanatorium wanted her husband to be tested for HIV/AIDS. When she first went for family planning operation 15 years ago, the doctors at the government hospital said that her blood was infected. That is why soon after a year she had her blood tested at the CMC hospital in Vellore, where it was confirmed that she has HIV/AIDS. From then onwards she is taking treatment. Her husband had some sickness for which the doctors had given him some injection on his spinal cord; now he was affected by lots of scabies at the private places so the doctors have advised him to take blood test.

She too had the same kind of attack of scabies and was later cured of it. She was first affected by tuberculosis and she took medicines for the same. She says that she does not know what medical treatment they give her but she is regularly taking medicine for HIV/AIDS. Her husband was a drunkard, he started to drink at the age of 18. She says she has been beaten by him many a time. He also used to beat their children so the children did not like him. She says she is poor and has no property of her own. She does not know about government policies.

She says that when she thinks that this disease is incurable she feels that has suffered throughout her life because of marrying such a man. When she sees her children also suffer, she is very dejected because she thinks of their future. She says that she has seen a lot of HIV/AIDS advertisements, such advertisements are good. She feels that there are all possibilities that a person who commits mistake will reform himself by seeing these advertisements. If the government gives free HIV/AIDS test, the poor people will be benefited. She proclaims that this is a fierceful, indecent disease so no one dares to talk about it openly.

She says no one knows about her disease except her younger sister. She has not even informed her children about this. She says her sister is kind and affectionate with her. Her sister knows that their uncle is a very bad man and she is fully aware of the problems the patient has faced with him. The patient further says that her husband is a jealously fellow and would not like her sisters to be well off. She says she feels well only when she takes medicines otherwise she becomes very ill. She does not know much about the commercial sex workers.

She says that only because of CSWs this disease spreads at a very rapid rate. Some of the HIV/AIDS patients are happy, they laugh and have taken up the disease positively whereas some of them are always thinking of the disease and crying, when she sees them she feels very sad.

She says she should question god for having given her this disease. She says the youth should be moral, they should not seek sex before marriage, they should lead a self-controlled life then certainly the disease will not spread. She says that for the sins men commit, women become passive victims. She says that when she thinks of this she becomes very upset and is angered by the attitudes of men. She says her son visited her and wept and also asked her when she was going to be discharged.

*Additional Observation: This is a very rare incident where the patient has been HIV/AIDS infected for the past 14 years and has never informed her husband or children or her parents. She looked hale and healthy. She did not look weak or sick like other HIV/AIDS patients in that hospital. She however did not give her parents' name or parents-in-law names. She is from an agricultural family. She did not name her native place but she is in Chennai (Guindy) for the past 9 years. She is composed and has a very strong will-power to live with this disease for the past 14 years without sharing her sorrow with anyone.*

## 95

We interviewed a 30-year old woman on 20-12-2002 in the Tambaram Sanatorium. She has studied up to 8th standard. She is married and her husband is 40 years. He is studied up to the 6th standard and he works as a lorry driver. They are from Velavalampatti. Her husband is also infected with HIV/AIDS and is taking treatment in the same hospital. She has never scolded her husband though she says only he has infected her. She has two children, aged 9 and 6. Her native place is Chengaipalayam. No one in her family know about this. Even her brothers do not know that she is infected by HIV/AIDS.

She says people should correct themselves. She has many friends and this disease has not affected any of her friends. She says that a person can get any disease but no one should get this disease. Her husband came to see her from the neighbouring ward he also joined our interview.

When we asked him what was his first reaction when he came to know that he had HIV/AIDS, he said that he wanted to commit suicide. Only people around him had consoled him and asked him not to do so. We asked whether her husband has so far asked her pardon. She said he has not asked her pardon. When we asked him why he did not confide about this disease to his wife, he said that he did not want her also to suffer. They said they are very poor. He joined work as a driver in 1982. They are married for the past 10 years. He says he has had premarital sex. He says he has visited several CSWs but he has not talked with them. Those CSWs had offered him condoms, but he had not used it. He spoke all these things before his wife. He proudly says that he has gone to all states and has visited CSWs of all states. He said that only government should take steps to stop CSWs from practicing their trade.

*Additional Observation: She looked very sick. She was married at 20 years. Her husband came when we were interviewing her. All things that he had not confided to her so far he talked freely to us right in front of her, while she was attentively listening. But she looked cheerful and composed. He seemed to be more impulsive because he thought of suicide as soon as he came to know that he was infected with HIV/AIDS. Now he says that no one in this world should be infected by HIV/AIDS. They have not yet tested their children for HIV/AIDS. They do not show any concern in doing it. She is like Sita / Savithri / Nalayani of the epics and does not scold him or find fault with him. She does not blame him; she just says that her husband has infected her. She asks us, "What is the use of quarreling with man who has erred?"*

**96**

On 05-01-2003 we interviewed an outpatient from Kallakuruchi who had just come to take treatment at the Tambaram Sanatorium. She is aged 31 to 32 years, is married and has 3 children aged 16, 13 and 9. Her husband is a taxi driver aged 39 years and is an HIV/AIDS patient. They are from Karunapuram village in Villupuram. She had attacks of cold and cough and lost her hair. Later she became very dark. She says that with these external changes she became more desperate. She was taking some medicine from Kerala and spent a total amount of Rs.6, 200 that too only for 10 days. All these medicines had no effect on her; but she became very sick with prolonged fever and diarrhoea.

Recently, as a follow-up we came to know that she died on 01-05-2004. Then the following data was collected from her home. Her sister gave this information! She says the diseased husband used to visit her only once a year and stay in the house only for 2 months. He worked as a taxi driver in Kerala, Mysore, and other neighboring states. She says that for the past two years the diseased did not go for any job because he too was sick. The family was supported by parents. He is suffering from boils, fever etc but refuses to take any treatment. Further we come to know reliably that there are about 35 HIV/AIDS affected patients in Alathur in Kallakuruchi. They are all taking treatment. She says even last week one man died of HIV/AIDS, he was a driver.

When we asked her reason for the prevalence of so many HIV/AIDS patients in that village she says that some five years ago there was a lodge, where they supplied CSWs. So she is very certain that there are many more in that village and neighboring places who are HIV/AIDS affected. In fact she says that she knew a family of five: husband, wife and three children, all of them are affected by HIV/AIDS. She says only last month that man died and now the wife is taking treatment. Several of these patients have taken treatment from Kerala spending from Rs.5000 to Rs.15000.

*Additional Observation: Of the 101 woman interviewed, we have come to know that as the following information this woman has died of HIV / AIDS passed away. She was very sick and she was sent home by the hospital authorities. She died within a day or two.*

## 97

We interviewed a Chettiar lady from Pachayur in Krishnagiri on 27-06-2003. She is now an inpatient of the Tambaram Sanatorium. She is 45 years old and has studied up to the 3rd standard. She is married and is a housewife. Her husband is uneducated and is a cloth merchant. He too is 45 years old.

They have 3 children, two sons and one daughter. Their children are healthy. Her parents and parents-in-laws have died due to age. First she was very sick and suffered with continuous fever, so the doctors advised her to go for a blood test. She was tested and found to be HIV positive. The doctors in that hospital immediately wanted to test her husband.

They tested him and he was found to be infected with HIV/AIDS. They carried out the blood test first in Tiruppathur hospital. Both were found infected in 2001. Her husband had tuberculosis and sore throat. She says that she feels better now but her husband has become very sick. He is not eating anything and is not even taking medicines. They are from very poor strata. She says that 4 years have passed since her husband left his job. She does not know how he got this disease. She says when he went on business and stayed away from home for weeks he should have been infected by HIV/AIDS. She says that her mother tongue is Tamil but they can also talk Telugu.

No one else knows that they are infected by HIV/AIDS. They have spent over Rs.10,000 for the treatment. She says that certainly there are some people who will become good after seeing advertisements about HIV/AIDS. They live in a rented house paying a rent of Rs.4, 000 per year. Their marriage was an arranged one.

*Additional Observation: She was very sick so she went for treatment where the doctors compelled her to take blood test and found her to be HIV/AIDS infected. Then they forced her husband to take blood test. It is unfortunate to note that her husband is now very serious and she is by his side taking care of him. He is not even in a position to sit or walk. He does not take any medicine properly. He was very sick unable to talk. He looked at us just once or twice and for over 45 minutes kept his eyes closed. She answered over 70 questions. She is composed and does not blame her husband. On the other hand she says it is her destiny that she should suffer.*

We met a Gounder woman aged 22, who has studied up to 6th standard. She is now an inpatient of the Tambaram Sanatorium on 18-6-03. She is married and works as an agricultural coolie like her husband. Her husband is 23 years old and has studied up to 5th standard. They have two male children: one is 5 years old and the other six months old. They have not carried out HIV/AIDS test for the children. Now after taking treatment here, they feel better.

On the advice of her neighbour she came to this hospital for checkup. She says no one knows about her disease, she felt very sad when she came to know that she had HIV/AIDS. She says that when she sees other patients she feels that she is comparatively better. As her husband had a very bad cold, which could not be cured, they admitted him in the hospital for the HIV/AIDS test. She wants her husband to be cured. She says day-by-day there is an increase in the number of HIV/AIDS patients since she sees there are many new patients visiting this hospital daily.

She says the main reason is people have extramarital sex with anyone and that they don't lead a controlled life. She is well aware of how HIV/AIDS spreads. She is from Karoor. She says her husband might have had sex with many other women.

He used to say the truth sometimes. Whenever she scolded him, he used to say that he would die soon. They earn around Rs.1000/- per month. They live in a hut and they are very poor. They get jobs seasonally, not always. Her parents-in-laws are unaware of the fact that this woman and her husband are affected with HIV/AIDS. Her marriage was an arranged one; she married her neighbour. He was a chain-smoker. He has never seen any HIV/AIDS advertisement. She says that her husband had extramarital affairs mainly in Andhra. Further she says that he has another wife who is in the same locality and is a Chettiar by caste. They belong to the Gounder community.

She says that the Chettiar woman was good looking. She says that when her husband had the first attack of HIV/AIDS his private parts became swollen and he could not be cured for over a month. That time itself doctors had suspected that he had HIV/AIDS. Her husband's affair with Andhra women began when one of his friends from Erode asked him to accompany. Her husband has acknowledged that he had visited the Andhra CSWs thrice. They used to leave home, go for cinema and have sex. The Chettiar lady has no children. She says her husband had sex with other ladies both before and after marriage. She says her husband suffers from STD/VD also.

*Additional Observation: She has become a victim of AIDS at the young age of 22. She is a Gounder by caste. She has come to know about her husband's activities. She often suspected him of having HIV/AIDS because of the symptoms he suffered. She says that he has made mistakes several times. She says he is very sick, with scabies all over the body. His private parts are swollen forever. He is not in a position to walk. She was married at 16, she is a AIDS patient at 22. They are poorly dressed. She says there is no one to take care of her. She further says that after her husband got admitted to this hospital, most of the time she is with him and helping him to eat, sit and answer the call of nature. They are not even in a position to test her children.*

### 99

On 18-4-03 we met a 35-year old woman from Thirukoyilur from Villupuram district in the Tambaram Sanatorium. She and her husband are uneducated and they work as agricultural coolies. She did not give her parents or parent-in-law names. She has four children: one daughter and 3 sons.

Her daughter is married off. Her youngest son is studying in the 7th standard, the other two children work as agricultural coolies. Her son and her brother-in-law have admitted her to this hospital a fortnight ago. She is now all alone in the hospital. It is since 10 months her husband ran away from the house. Her husband is a drunkard. Once before in a quarrel the conductor pushed him out of the bus. He was very much injured.

They kept him in the Ooty hospital for three months, when he was able to walk a little once again he ran away from the house. Her husband is 50 years old. She says that her children left him free, she does not know whether he is living or dead but some say he is still alive. She says that she never had any sex before marriage. She says she had sex with her brother-in-law. Now she is living with him. She thinks she has got this disease from her brother-in law. She says that he also drinks, he works in TNEB in Kezhiyur. No one in her family knows of this disease.

They think she suffers only from fever. She has not even informed her children because they would feel bad about it. She does not know much about HIV/AIDS. She has not even seen advertisements about HIV/AIDS.

She says more advertisements can be given about HIV/AIDS in rural areas, it would do the villagers a lot of good. In villages, if someone, says and explains they would follow otherwise they won't follow.

*Additional Observation: This is another rare case. She is 35 years old, she lives with her brother-in-law, her husband has run away from the family. We see that she is frank when she says she lives with her brother-in-law. She is taking treatment for the past two years. Her brother-in-law infected her.*

*She has no proper knowledge of HIV/AIDS and is ignorant of other non-sexual modes of transmission. She has not informed anyone in the family about her disease. She says they would be upset. She did not look very unhealthy.*

### 100

We met an inpatient from the Tiruchi district. She is 31 years old. She is an agricultural labourer during the season and idli vendor when she has no employment as a coolie.

Her husband studied up to the 10th standard and became a lorry driver. Earlier he was working as an agriculture labourer. She has also studied up to 10th standard. She did not give the names of her parents but she gave the names of her parents-in-laws. She has three children. She says weekly once he goes on trips. Every month he used to bring home Rs.3500/-. Even now she lives only in the joint family. She has been confirmed of HIV/AIDS last year. She was tested for this disease at Tiruchi. In Katoor she has spent over Rs.5000/- for treatment. The doctors there wanted more money from her; when she refused they used abusive language against her (like *chi chi, you dog,* etc.,) so she came and admitted herself in this hospital.

She has no hereditary property; she is very poor. She says she knows both English and Tamil. Her husband also knows both English and Tamil. He used to drink alcohol once a while. Her husband is still earning a livelihood by driving lorries. She says she never had sex with any man other than her husband. She says she does not know about her husband's affairs with other women or CSWs. She does not suspect her husband. She says that she might have got this disease due to her monthly periods. She says that she would have got this disease due to heat. She admits that she does not know how one gets HIV/AIDS. She says that she had fever for over ten days, the doctors tested for blood and declared she was HIV +. Her husband said that he had never been to CSWs. She was sad when she came to know that she had HIV/AIDS but when she saw so many HIV/AIDS patients she is not at all sad.

Her major symptoms were prolonged fever with complete lack of appetite. She says that now her husband has been tested and he is affected with HIV/AIDS. Her first two children are free from the disease because they were born before she got infected. Her 3rd daughter who is just 5 years is affected with HIV/AIDS and she is also taking medicine for it. She says she has seen a HIV/AIDS advertisement in which they say one should not have sex with strangers. She says this disease has lots of social disgrace and stigma associated with it, that is why people do not talk about HIV/AIDS publicly or openly. Everyone knows that she has HIV/AIDS. Relations don't talk much with her, they are aloof. Her husband, her children and her parents-in-law don't mingle with the neighbours. She wishes to be cured from this disease. She says that no one can ever know who has this disease so it would be better if one does have sex with strangers. She says that CSWs take up this profession only due to poverty because they know if they get the disease they have to suffer life long. They go mainly to satisfy their stomach so she does not want to blame them. She says only in the villages people should give awareness program, for the uneducated become easy victims of HIV/AIDS. She knows how it spreads. She says this sort of interviews can also spread awareness among people.

*Additional Observation: She was first of the opinion that she had got HIV/AIDS due to heat and due to menstrual problems. She has very soft disposition, first women to say that the CSWs do that work only due to poverty. Nobody on earth would risk their life with HIV/AIDS for the sake of pleasure, it is the acute poverty, which has forced these woman to take up this profession. Certainly her views will make people think twice before they talk ill of CSWs. She says everyone in the family and everyone in her village knows of the disease that she is suffering from. She is reasonable and rational.*

We have kept this as the last and 101st interview with a special significance because this woman is a social worker. After being affected by HIV she built the Positive Women's Network (PWN).
Now we proceed on to give her brief interview with us.

## 101

We met a 28-year old woman on 24-1-2002. She has studied upto the 12th standard. She is from Namakkal. Now she is in Chennai. She has established the Positive Women's Network (PWN) and she is its president. She says she gives training to students in colleges about HIV/AIDS awareness. She has acquired the disease from her husband. Her husband died, at that time when she was affected with HIV/AIDS she was sad and felt guilty. She says she had the first fear that she would die soon. She was infected by this disease just three months after marriage. First, all her relatives were frightened to move with her. She feels that people generally view HIV/AIDS infected as 'bad' persons and have obtained the disease only by their bad deeds. She says the spread of HIV/AIDS can be prevented by advertisements about the training courses in TV, radio and in the college campus. She says that if the uneducated are given proper explanation they understand well. She is of the opinion that only 10% of the doctors take care of patients properly. She says that all foods that are rich in proteins are good for the health of the HIV/AIDS affected patients.

*Additional Observation: We see that she is a living light who have taken up the cause to fight the discrimation against HIV/AIDS that too in societies like ours where widowhood is a large stigma. This brave woman bearing two sorrows has come out boldly to fight for a good cause. In India, that too in south India at that time there was practically no one to take up this cause. In this male dominated society we see that the women are more vulnerable to HIV/AIDS due to the faults of their husbands. It is to be noted that she has given several training programs at the college level to spread awareness.*

Now the emotions of these 101 women patients will be analyzed using fuzzy techniques like Fuzzy Relational Maps (FRMs) and Fuzzy Associative Memories (FAMs). We feel that as the data is an unsupervised one, full of transitory human feelings and (some of which is expressed, and several, suppressed, due to fear or status or sex) we felt it deem fit to adopt fuzzy model for the mathematical analysis.

**Chapter Six**

# SUGGESTIONS FOR IMPLEMENTATION AND CONCLUSIONS

The suggestions and conclusions are based on our study of 101 HIV/AIDS patients interviewed using linguistic questionnaire, they were either inpatients or outpatients of Tambaram Sanatorium whom we had met either in the tea shop or in Sanatorium and from the mathematical analysis of data using fuzzy theory and neutrosophic theory. It is also pertinent to mention here that most of these women i.e., 96 of them are from rural areas with no proper education and from poor strata. Of the 101 patients, 91 of them were infected by their husbands, 36 of these women are widows because their husband has died due HIV/AIDS after being sick for a long time. 55 of the men (infected women's husband) were either taking treatment in the hospital as inpatients or outpatients. Forty-four of the interviewed women had never entered school and 16 of them had studied primary classes from $1^{st}$ standard to $5^{th}$ standard, 25 had studied from the class $6^{th}$ to $8^{th}$ standard, only one has studied $12^{th}$ standard and only one had undergone a nursing course.

It is important to note that among women patients, if the disease is found at the early stage then they live longer than men patients, the only reason is that even after becoming

HIV/AIDS patients men continue the same style of life, which makes them more and more vulnerable to infection. So most of them die within 5 years.

Also we wish to state that several of the interviews we held were freewheeling discussions, allowing the patients to speak their worries and trauma, we shared all their suggestions.

Since most of the women we interviewed were from rural areas with no education and were from economically poor strata, our research is restricted only to these poor uneducated rural women and our data or study do not pertain to rich/urban women or working women with good educational background. So our study, conclusions and suggestions are only for poor rural uneducated women who are passive and innocent victims of HIV/AIDS.

This chapter is divided into three sections as follows:

1. Observations from our study
2. Suggestions given by the HIV/AIDS patients
3. Conclusions from our mathematical analysis

## 6.1 Observations from Our Study

1.   80% of the women patients interviewed complained that they were married very young so they could not find a meaningful life; further they became mothers and their life was more of a misery and the responsibility in their teenage itself was very high and they found life more a dreary routine than a happiness. This is mainly attributed to the fact that the parents of the girl children feel it a big burden to feed them, they think that if they marry their daughter at 14 and not at 24 years, they can gain the 10 years money to be spent on her; apart from this, the bitter truth is that they also deny basic schooling and education to their female children. Several of them also said they do not know any other occupation to sustain them economically.

2. Most of the women patients not only talked more and freely but they also gave the general opinion, a consolidated problems/feelings/sufferings of other patients. This was absolutely absent when we interviewed the male patients to even speak about them was itself not from the heart.

3. The female patients were so submissive to fate or to natural course of due action that most of them did not feel dejected or depressed. The suicidal tendency in them was comparatively less than that of men or completely absent in majority of the cases.

4. Doctors should strongly insist the male partner (husband) to take the HIV/AIDS test when the wife is infected and when the husband is infected they should demand test to be done for the wives.

5. Husbands who know that they are infected with HIV/AIDS hide this information from their wives. Even if the wives are infected and sick they do not disclose to them that they are HIV +. This mostly expresses their irresponsible social behavior and an uncaring attitude they show towards their wives. On the other hand, any woman who is infected by her husband does not blame him badly.

Women have no voice to force their husband to go for a blood test for HIV/AIDS. At this stage it is important that the doctors force these men to take up HIV/AIDS test once their wives have been found infected. The two reasons are one, they can take proper treatment before the disease becomes chronic and advanced. The two of them can lead a proper life.

6. From our study and analysis we see that the women take so much pride of being educated, for, when we ask them about educational qualifications they proudly say they have studied up to $2^{nd}$ standard. Some of them take even pride in saying that their husband has studied up to $1^{st}$ standard. Thus these rural women long for education and they are happy to say that they have studied up to $1^{st}$ or $2^{nd}$ standard, shows their love for education. When we say uneducated it implies they have not stepped into any school premises. These interviewed women even if they have studied $1^{st}$ standard, they proudly say that they have studied up to $1^{st}$ standard.

7.  We are forced to think if there is a possibility of getting infected by HIV/AIDS for those women who undergo natural abortion or from the hospitals in which they are admitted where they take the immediate treatment. In some case we feel or get a feedback from our study that there are possibilities when women are treated for natural abortion become victims of HIV/AIDS. This study needs some attention in the medical field, for we see some women who suffered natured abortion had become infected with HIV/AIDS whereas their husband remain unaffected by HIV/AIDS, which was confirmed in these cases by blood test, or more relevently, "Is their any relation between HIV/AIDS infection and natural abortion?".

This is left for the doctors to carry out further research. Or does frequent natural abortion indicate a symptom that the patient is infected with HIV/AIDS.

8.  We have observed that men patients immediately say they have got HIV/AIDS from women (not CSWs) but we see the vast difference that women do not blame even their husbands, most say they do not know how they got this disease. Some suspect only when we put forth questions like "Did your husband have sex with any other women?" Most attribute it to injection or blood transmission, but actually from the data 91 women out of 101 women were infected by their husbands.

9.  We feel it is the professional ethics of doctors to give treatment to any patient for when they under take the, Hippocrates oath. They should not say that we don't give treatment for such a disease. In all cases it is mandatory for the doctors to give treatment to the patient especially when an alternative is totally absent. We hear from patients that doctors don't give treatment even at the time of emergency. It is grave injustice, which they do to their profession.

10. In several cases the husband and the in laws do not say even to the wives that they suffer from HIV/AIDS.

11. Several people from Andhra Pradesh come to take treatment in the Tambaram Sanatorium. They say they get the best treatment in South India for HIV/AIDS only in Tambaram Sanatorium.

12. Several men knowing fully well that they are infected with HIV/AIDS do not disclose the disease to their wife but on the other hand have unprotected sex with them which makes them infected in a very short span of time say within a year. Within a year these women become seriously ill and suffer with all sorts of opportunistic diseases apart from STD/VD and HIV/AIDS. This attitude of men must be changed.

13. If there can be some women support groups and if the doctors who treat the infected husbands call their infected wife and advice them to avoid sex these women can be saved from HIV/AIDS. Not only do these women become infected, at times the children born to them are also affected and these children die as infants in many cases.

14. Yet another important observation is that most of these women who have got HIV/AIDS not through their legal spouse clearly name their friends from whom they have acquired it. By no means do they hide the fact. Also they give the cause of having such affairs. In some cases it is sheer insecurity because their husband had run away leaving them alone, others did it for money. They have illicit sex to ensure monetary support for their family. Apart from this we do not see any woman committing mistakes to enjoy life etc. Further they do this mistake only after marriage and after they are cheated or left uncared by their husband. It is important to mention here that we give only the cause; we give no justification of their acts.

15. We wish to state that of the 101 women we interviewed only one woman did not say about her husband's character prior to marriage whereas the other 100 of them openly said that their husband had sex even before marriage. This is one of the main reasons for HIV/AIDS among infants and mother to child transmission of HIV/AIDS. Also this habit of the men to have illicit affairs continued even after marriage in 90 of the cases, which is really questionable, and this is one of the major reasons for women being infected with HIV/AIDS and the spread of HIV/AIDS at large.

16. We observed that among the 101 women only 5 had premarital sex, of which 4 of them later married that man with whom they had sex.

17. It is still important to note that of the 101 women we interviewed only 7 women had sex with men other than their husbands and lead a loose life in their twenties and late teens.

18. It is also important to mention here that some of the women whom we interviewed said that there was difference between CSWs and prostitutes; most of the CSWs do their profession due to poverty or by some forced circumstances unlike prostitutes who sell their body due to their sexual urge and also at the same time make maximum money and thereby exploit man who comes across them. A few women distinguish these two terms in Tamil as *Vesi* and *Vibachari*. There is a marked distinction between the two.

19. Further women do not have multi-partner sex. Of the 101 women, the men had for sex was only three by one single woman. Three women had sex with their husband's brother also. Thus we see the sex behaviors of rural women and rural men are incomparable.

20. Most of the women, even if they are aware of HIV/AIDS, are not in a position to save themselves from their HIV/AIDS infected husband for the following reasons.

    i.   Husbands do not disclose their HIV + infection disease even if they are taking treatment for it.

    ii.   They don't restrain themselves from sex with their wife, thereby infecting them. Also the vulnerability of a woman to contract HIV is 24 times greater than a man's vulnerability.

    iii.   Even if they know their husband has HIV/AIDS, most of the women are ignorant of how HIV/AIDS spreads and thus become victims.

21. On the other hand if a husband knows his wife has HIV/AIDS (Of the 101 women only one woman has not told her husband or children that she is HIV/AIDS infected), the husband fears and does not have sex with his infected wife. Thus we see the careful way men try to save themselves, but that carefulness is totally absent when they think of their wife being infected by them.

22. Also when the husband is infected by HIV/AIDS invariably the wife serves him and takes maximum care. On the other hand if a woman becomes infected with HIV/AIDS we see most of the men desert their wife, or they are already dead and in these cases the women are either left in the streets or all alone suffering in the hospital or their parents take care of them and support them.

23. Another sad thing to be noted is that these women are very badly abused by their parents-in-law by saying that they are immoral women or characterless women and in some cases the mother-in-law tries to arrange another marriage for their already married and HIV+ son. So the HIV/AIDS affected women not only suffer the social stigma, the sorrow of being infected by their husband but above all, they are abused as being characterless or immoral and are chased from the home. If the parents are alive and sympathetic, they take care of their daughters but in some cases parents just avoid them and in some cases parents are dead. Thus the plight of HIV/AIDS affected women is very deplorable, pathetic and beyond any form of justification. They are ill-treated, harassed and tortured.

24. These HIV/AIDS affected women suffer from several diseases at once and their longevity is less than that of men when the disease is discovered at a chronic stage.

25. In spite of so much suffering these women do not feel like committing suicide. We see out of 101 women only 6 felt so. This shows that they do not feel guilty, but among the men 90 of them attempted suicide and 2% of them have committed suicide. Psychology behind this needs research.

26. It is important to note that many women from one particular land-owning community are affected. One has to study this problem from the sociological circumstances particular to that community.

27. Another observation is that from our study is that women from Andhra Pradesh who are infected have husbands who are free from HIV/AIDS. This also needs a special study.

28. Most of these HIV/AIDS infected women in comparison with the HIV/AIDS infected men, are open, cheerful, informative, exchange their views and come freely for discussion and not only do give their ideas but many give the over-all picture since they get feedback from other women patients. This was totally absent in case of men when we interviewed them.

29. We also see that several patients had boils and scabies on their body, and since the wounds were undressed

and open these were feasted by the flies. As few of them are very immobile they are not even able to sit or in a position to cover themselves with a bed sheet. Their hair was matted and unoiled may be it was because of poverty, or perhaps there was no one to apply oil to the hair of the ailing patients who are so ill.

30. 61 of the women have lost faith in god, do not believe in god and have stopped praying to god since they think god has failed to answer their prayers and has made them practically immobile.

31. Most women complained that they were not able to walk or sit for more than half an hour.

32. 18 of the HIV/AIDS patients had grandchildren and they longed to be with them.

33. Women never had any addictive habits like cinema or alcoholism but some used tobacco mainly to curb hunger and to enable them to work for more hours. Most women worked as coolies. Thus they never wasted money or time, they were busy working at home or as coolie. The women looked weak and uncared for and overused.

34. It is unfortunate to note that of the 101 infected 91 were infected by their husbands.

35. Some of the HIV/AIDS patients suffered double discrimination by being infected with by HIV/AIDS and being widows. Thus some were chased away by their husband's family and they end up on the streets with their children.

36. 98% of the woman patients were not careless or spent money lavishly on their paramours or irresponsible in comparison with HIV/AIDS male patients who were careless, spent money on CSWs, used to smoke and drink and had overall irresponsible behaviour.

37. None of the 101 woman consumed alcohol. Of the 101 women interviewed 92 of the women said their husbands were drunkards, smokers and visited CSWs.

38. Of the 101 women, only 3 confessed that they were addicted to cinema.

39. Most women valued education so much that they took much pride in saying the classes up to which they have studied. Even if they have studied only first standard they said so with pride. They were very happy to say the educational qualifications of their husbands. Of the 101 women only two women have entered college, all others were just school dropouts or have never entered school premises. One women did not have TV in her home for she wanted her son to get good marks in plus two who was just in his 11th std.

40. The profession of majority of the HIV/AIDS woman patient's husbands were truck drivers, cab / auto drivers agriculture coolies, construction labourers or petty shop owners or bore pipe workers. None were government servants or teachers or holding a clerical post like a peon or clerk in any private / government concern.

41. Thus from our interview we see the catching of HIV/AIDS by majority of women were only relative due to the careless behaviour of their husbands.

42. Self-help groups must come forward to take care of these women who have no one to support them physically or mentally. When we say physically we mean that some women patients could not sit without someone's help for they suffered from big boils in the legs and hands, which physically prevented them from any movement. Some of them had problem with fingers. These people cannot tie up their hair or comb it. On the top of it, their matted hair and scalp was feasted by lice. The lice had caused some problems in their scalp and head. So even while we were interviewing and observing all these, it made us very unhappy. These women have brought up their family, served their husbands and children and their own parents, but now they die a uncared for death aggravated by flies and lice. Hence some trained self-help groups women with the mission of social service and above all selfless service should take it up to help these patients. Or the hospital could appoint some semi-trained nurses to handle these patients for the ayahs do not assist them unless they pay money or tip them.

When these women patients are themselves deserted where will they go for the money to tip ayahs and ward boys. So at least in the women ward such appointments must be made. After all, most of the women who suffer the disease are selfless women who had become victims of HIV/AIDS by their husbands. At this juncture we feel it is the ardent duty of the nation in general and our state in particular to help these women at least to be in peace and die a peaceful death since their physical sufferings cannot be described in words. This may give them a peace of mind and composed attire, which may lessen their physical sufferings.

43. Another important observation about these women patients is that out of 101 women patients only 4 of them felt like committing suicide. Also we see these women patients in general are not even desperate they are composed only very few (three) of them suffer from depression. On the contrary we observed when we interviewed the HIV/AIDS affected male migrant labourers that for over 89% of them, the first reaction when they were affected by HIV/AIDS had been was to attempt suicide. This was very clear from the way they talked. Also in our interviews from women we found that 3 of the men have committed suicide when they heard that they have been affected by HIV/AIDS. Many of the men (around 40%) have turned to religion for consolation contrary to women (over 60) have lost faith in any religion or god because they feel the power of god is limited and hence they cannot be cured by faith. Only one women was of the opinion that till now people have not invented a special god for HIV/AIDS who can cure it. She compared the goddesses. Mariamma who was a sole goddess to cure smallpox. She only blamed the people for their inability to discover a god for HIV/AIDS who can cure people of this fierceful disease.

44. One of our major observation is that if women are empowered economically, socially then certainly it will not only decrease women being infected by HIV/AIDS and also it will put a very strong hold on the dependent men. So the awareness program would become very powerful.

45. Adolescent boys especially those who are uneducated when they go away for work or stay away from home are at the greatest risk of not only contracting STDs but also the HIV virus. After a year or two of contracting the disease they get married thereby their new-born and their wives become innocent passive victims. Most men before being affected by HIV/AIDS are chronic patients of STD/VD, it is in fact this problems that easily makes them to catch HIV/AIDS.

46. In India, especially in Tamil Nadu, married women are culturally not expected to say 'no' to sex with their husbands. Thus the majority of the infected men do not confide their disease to their wife but without any conscience have unprotected sex with their wives thereby making the women vulnerable to HIV infection. They have practically no voice or control over either absistenence or condom use at home or on their husband's extramarital sexual activities. Here it is pertinent to mention that in India this is the case with all women, whether they are well educated, rich or otherwise. Also the "husband-wife" relationship in India that is inherently patriarchal is best suited for the spread of HIV/AIDS.

47. Most of the women give only the real names of their mother-in-law and father-in-law. However we are not able to understand why they do not give the names of their parents' or native place or village of their parents. They just don't answer this question.

## 6.2 Suggestions given by the HIV/AIDS patients

1.   More than the discovery of medicine a 35-year old woman feels that so far people have not discovered a particular god who can cure HIV/AIDS. She said that as the present gods could not cure or do anything with this disease for they were unaware of how this disease can be cured. So once the god who can cure HIV/AIDS is discovered she hopes everything will work well. However she did not specify who can discover this god.

2.   One may assume that uneducated women may not have inventive and creative ideas but it was an uneducated woman among the interviewed who put forth to us the idea of finding a preventive injection for HIV/AIDS. She says just like polio drops to prevent polio and small-pox vaccination to prevent small pox, why not find some medicine to prevent HIV/AIDS.

3.   Several people from Andhra Pradesh come to take treatment in the Tambaram Sanatorium. They say they get the best treatment in South India for HIV/AIDS only in Tambaram Sanatorium.

4.   One of the major observation is that out of 101 HIV/AIDS affected women patients only one woman was very sympathetic for the CSWs. In our opinion she was very rationalistic and reasonable about CSWs because she said that no women would willingly do the job of CSW, it is acute poverty alone which has forced them to take up such a profession. Driven by poverty and insecurity they have become victims of this profession. She added on. "It is like a ragpicker taking up the work of trash collecting". Thus a woman taking up commercial sex work is like men/women opting for ragpicking. Her talks were still worse but certainly it has a mission, a higher philosophy and a better reason so we felt it very essential to record it.

5. Many women said that the hospital was well maintained and the food was rich and good in the time of the doctor Deivanayagam. Also they said he used to personally visit all the patients. But now the doctors don't even come near and in most cases only nurses give the medicine and injections. They said the rice is foul smelling because they feed the HIV/AIDS patients only with insect ridden rice. Some of them begged us to recommend to government to give good rice. Earlier milk, egg and fruits were given now these are given in lesser and poor quality and irregularly.

6. Some of the women said the men who are inpatients in this hospital smoke secretly without the knowledge of the doctors and nurses some drink and some go out in the nights to visit CSWs, HIV/AIDS spreads by these men's 'loose character'. Even after they are infected they do not correct themselves.

7. Women supported the proposed idea of free HIV/AIDS test by the government since it would be a boon to poor women. But at the same time some of the patients put forth the question "how many people will accept such tests, and come forward to take them?". Now health centers on the advice of TANSACS give free test to woman who come for treatment.

8. Most women said that self-control and morality alone would help prevent the spread of HIV/AIDS.

9. All women uniformly feel that the awareness program about HIV/AIDS given in villages is not sufficient because

a. It does not target the uneducated adults like them.

b. It is very scarce, they forget or no traces about awareness of HIV/AIDS ever makes any form of impact on their mind.

A woman patient felt that awareness and advertisement must be given on Sundays for half an hour or so in the form of songs, dramas or skits or talk by experts on TV, which will make an impact.

10. Media, Cinema and TV made men loose their control and seek sex which was one of the cause of HIV / AIDS spread.

## 6.3 Conclusions from our Mathematical Analysis

1. As over 91 of the women patients from the villages are poor, uneducated and infected only by their husband it is high time we take some strong steps to prevent women becoming passive innocent victims. For this we give the following suggestions. Every remote village must have a health care center within an easily accessible distance specially for women and doctors there should be given training in HIV/AIDS and help these women to know more about HIV/AIDS.

2. The marriage age of a girl in the villages beings at 12 onwards. Thus it is very essential that the government ensures that such child marriages are punished or those who indulge in such acts must be jailed for not less than five years with a lump amount of fine. In particular a special rule must be made to utilize this fine for the rehabilitation of HIV/AIDS affected women patients who are mainly victims of child marriage also.

3. From our analysis we see that all these uneducated women value education greatly which is very evident from their facial expression and the pride they take in saying that they have studied up to $1^{st}$ standard or $2^{nd}$ standard and so on. Not only this we see that they also equally take pride in saying that their husbands have studied up to 1st standard, they disclose it so happily and proudly. Thus from this point we want to make it clear that to be educated is so important to them and the majority of them were only forcefully stopped from their school to look after their younger brother or sister or in some cases their mother was sick so they shouldered the responsibility by cooking, washing vessels and running the family, some female children, were stopped their schooling for they were taken as child labourers.

Thus we see to stop all these malpractice against women, especially female children the only means is they must educate their female children up to 12th standard a stipend money as incentive will be paid to their parents for 12 years apart from making the education completely free for them free books and free uniform will be provided. The government ought to prohibit marriage that interferes or leads to discontinuing of school education. Strictly this concession and opportunity will be provided only to those girls whose parents are uneducated and live in remote villages. Certainly this will not only stop people practice female infanticide but will also give them at least education (for the sake of monetary incentive and above all prevent child marriage). If discontinue all the money spent on the child must be returned to government with interest.

Here it is very important to note that child marriage does not mean both the bride and bridegroom are children but here only the girl is just 12/13/14/15 but in all cases the men with whom these girls get married are twice or thrice the age of these girls. In several cases the girls are married as the second wives. What should this marriage be called? Child-old man marriage or child-widower marriage? The major reason they attribute to this atrocity is that they are very poor and the bridegroom is the child's relative. So is it the right of a child's relative to marry her even before she reaches puberty? Is all these legal, or rational or in any way justifiable? All these can be prevented by the suggested scheme.

Further in most cases these girls loose their husband within 5 years time and they become young widows when they are just in their twenties. What steps does the government take to prevent all these? Why no women's organization has publicized it? Why is the government a passive on looker? In some cases we see the son-in-law is several years older than the mother-in-law. What is all these? What is the physical status of such a female child? What is their mental status? What are the health-related problems they face? Who is the cause for all these? Invariably over 90% of the men who marry such young girls are perverts, drunkards and smokers. Above all they are the ones who frequent the CSWs. Because of the age difference the girl is frightened to ask him anything! Several women expressed their sex-related problems, like they suffered like natural abortions, STD/VD, continuous periods and other gynecological disorders. Some girls have also lost their mental balance for 3 or 4 years, some looked as if they were possessed. We saw out of 101 there were 87 who were victims of one or the other type. Just at the age of fifteen several of them were mothers.

How can a child bring up another child? Who will come to the rescue of these children? How to stop this? Who is going to voice for these victims? Only when solutions to these question are found and solved the HIV/AIDS affected among uneducated rural women can be reduced. Without solving these no solution can be derived and the spread of HIV/AIDS among this strata cannot be even decreased.

It is still unfortunate to see that these women are so uncaring about themselves for they are not slaves of any habits like chewing tobacco, they do not even take care of themselves, even in their late twenties they look like forties, in forties they are so worn out they look like sixties. Their nails are untrimmed, their hair left uncared and they cut a very pathetic figure. How are we justified in harassing these young girls? What is the right of parents? What are their limitations? What is the right of their relatives? What is the duty of the nation, which has so miserably failed to safeguard its women?

4. What is the role played by the women organizations that are based in cities. These interested women organizations/government/NGOs/some service-minded public/NCC cadres/NSS cadres should reach the problem from the grass-roots level. As a group of 10-20 students depending on the village population should camp in the villages for at least 30 days and mingle with the women of the villages, study their habits, their problems, their customs, their occupation, etc. and at the same time they ought to provide awareness to the women.

This is not a big problem if each college adopts a village in the summer vacation, spends a month camping in a village and improves it. Such student-people initiatives will bring about notable changes, they should repeat this process for at least 3 years. However 5 years would be ideal. They can build schools, provide them with health centers and above all empower the women groups in that village, by which a strong women association can be built and 5 years of intensive training can certainly will make them achieve self sufficiency and fight for their own right and ultimately empower them. This is the best way one can empower, educate and eradicate the fear in these women who become easy HIV/AIDS victims. Also this would add sincerity and service mindness in college students and make them more responsible!

We can also channelize the energy of these younger generations and thereby they would live more in reality! This suggestion is not only for women, also for the men in remote villages because they become easy victims of HIV/AIDS and infect (innocently or deliberately) their wife and consequently, their children.

5. UGC and AICTE should formulate and include this in their educational program, which should stipulate a compulsory NCC / NSS programs which makes a compulsory camp in a remote village and the status of the village after 5 years must be submitted to the government. This alone can inculcate the awareness and also help women build self help groups. This might also make men not waste money but productively save the money they spend on smoke, alcohol and CSWs, this money will be a source of support at the time of crisis. All villages must be made to have a self help group and an association to empower women to question and also establish their right by which they can root out evils like child marriage, child-widower marriage, child-old man marriage, female infanticide, infection of HIV/AIDS by HIV/AIDS affected husband and above all education to all females and strict enforcement of Sharda Act that recommends marriage at above 21 years. They can also be trained in some handicrafts.

6. Women organizations should work more to empower women in villages. They can use their energy in empowering women in remote villages where no education is the basic qualification. When we say 'no education' they have not stepped into school campus or a classroom.

7. A part of adult education, HIV/AIDS awareness education can be given to adults. A doctor should give a lecture on what are the basic symptoms they may suffer when infected by HIV/AIDS. Basic indicators of the disease and symptoms, like prolonged headache, skin ailments, frequent fever, big boils, scabies, loss of hair, all types of STD/VD, incurable recurring cold, T.B., indigestion, lack of appetite, loss of weight etc. so that the rural people would be aware of and go for counseling to the nearby health care centers.

Even if they are affected by HIV/AIDS, they can take treatment before the disease becomes chronic for we see women lead a normal life for over fourteen to fifteen years even after being infected by HIV/AIDS. Further the program must not alarm them saying that there is no cure for HIV/AIDS but it has been proved to be a livable disease if medicine and proper living is taken at an appropriate time. By making statement that there is no cure people shun, fear and ultimately associate a social stigma, which is nothing but the fear psychosis in people. If a positive note is proclaimed the HIV/AIDS affected men will not hide it from their wives due to the main fear of social segregation. This in turn makes the women contract the disease from their husbands.

It is important to mention that when men come to know that their wife was infected with HIV/AIDS they never had sex with her but on the other hand when the husband knows fully well that he is a confirmed HIV/AIDS patient he did not inform her of his disease, and also ruthlessly without any social conscience or even with no concern had unprotected sex which made majority of women, that too uneducated, poor, rural women a victim of HIV/AIDS. If this situation is to be changed, only the solution is that women must be empowered. They must be socially and economically empowered then only they can at least save themselves. We see that out of 101 patients, 91 were infected by their husbands. What justification can the men give for these heinous acts?

8.    Just like constant advertisements are made for polio drops by the government in rural areas, the government can make advertisement for people to take up HIV/AIDS test if they suffer from certain symptoms. Confidentiality of the test results must be preserved at all stages.

9.    By giving special preference to jobs in villages, in health centres or in co operative / ration shops run by the government to HIV/AIDS affected persons and also by encouraging those patients by providing them alternative livelihood, the amount of social stigma associated with the disease will be lessened for no one in the village would dare to say that they are going to forgo their ration products since it is sold by the HIV/AIDS patients. Such bold acts alone can wipe out the fear and above all the social stigma associated with HIV/AIDS. When such open encouragement is given to HIV/AIDS patients by the government it would certainly lessen their social stigma. To achieve this the villages must empower women. Generally, women take up a social reform easily than men. This we wish to say by sheer observation.

10. Another important factor to stop the spread of AIDS and above all to help the infected especially in the villages among the poor and uneducated would be to pay some money as subsistence by the government to the affected family. This would serve double purpose:

   a.   People will not feel the stigma. They would feel the government is helping them, so they will come out in open to get the help.
   b.   The poverty stricken family would certainly be benefited by the monetary help.

11. Government can also announce a package for the children, in case of parents who are both HIV/AIDS infected. They can bear the cost of education of these children if they are young and in school, in case they are grown ups, the government can give some jobs to their children so that the family does not suffer both ways.

One important impact of this program would be the educated and urban elite affected with HIV/AIDS will slowly come out in public and acknowledge the disease they suffer openly. At the same time the awareness among public would become more. This will in turn make the doctors not to shun the patients.

12. Above all, free medicine must be given to poor uneducated women from rural first, for they are the priority group in our observation, out of the 101 women at least 76 of them lived like orphans no one was by their side even for the physical support leave alone the material support! So these women, when given fully free treatment with good food may become better and feel better than dying a slow and steady death with an added frustration of living a life of complete suffering. But the government should also take steps to see that the free medicine reaches only the deserving and is not sold out by touts and others in the hospital. We are forced to add this because for some of the patients in the hospital said that there are some touts from Andhra Pradesh who take bribes from the patients in Tambaram sanatorium to get admitted as inpatients. This is a very recent occurrence. Earlier patients were admitted as an inpatient without any such middleman. The government should take stern steps to stop the functioning of these middlemen in Tambaram Sanatorium. Three or four newly admitted women and men spoke about this. It is high time the government curtails such indecent touts who have the heart to exploit the poor uneducated HIV/AIDS patients. This sort of middlemen work in admitting patients is very recent i.e. for the past two years, and these touts are mainly from Andhra Pradesh. Several old patients complained of this but were very frighten to talk out openly, but they said if this stops it would be better for them.

13. By raising the taxes on risk-related products such as alcohol and through emphasise on the policy of obligatory condom use while visiting CSWs since it is the only way to protect the CSWs as the well as men who visit them. For majority of the rural men are in the drunken state when they visit the CSWs so they often fail to use condom.

14. The vulnerability of rural women to HIV/AIDS can be lessened through increased education (awareness education about STD / HIV/AIDS etc and above all sex related problem) access to credit, and equitable rights in the event of divorce, certainly such action could simultaneously diminish women's vulnerability to other ills with similar socioeconomic roots including violence and unwanted pregnancy.

15. However it has become pertinent to mention here that unlike in western countries where some people view the HIV/AIDS epidemic as a failure to respect religious codes but we purely view it as a failure to respect the moral codes for we being secular, an outcome of differentiation in sexual behaviour and sexual decision-making between men and women, a human rights issue or another tragic correlation of poverty and deprivation, to name but a few of the paradigms that have evolved since the start of the epidemic.

16. Men should be educated to treat women especially their wife's properly and with conscience so that the wife is not infected by the erring husband.

17. Unless men are counseled and take in their minds as a social issue no change can be made in the health of their wife. It is an unestablished truth that poverty and illiteracy, are associated with higher disease rates, lack of prenatal and preventive care, lack of health services, and the seeking of medical care only for acute illness. Thus it is among the impoverished rural population that the highest levels of HIV risk behaviour occur.

In rural areas of T.N about 90% of women with HIV/AIDS have contracted it only from their husbands (We don't include CSWs in this list). This is in keeping with the gender inequality that is prevalent in India; for urban men in Tamil Nadu are relatively free to engage in premarital and extra marital sex but women are restricted from doing so. Thus the major risk for most women is from sex with their husbands. A substantial number of these rural men engage sex with CSWs especially when traveling and their wives are forced to accept their behaviour due to economic dependence and in few cases they are ignorant of their husband's habits. Providing these rural uneducated women with greater access to health services and education is crucial to control the HIV/AIDS epidemic.

The government's prominent role in AIDS crisis should be a broad developmental approach that stresses the importance of poverty reduction and that concentrates resources on such basic human needs as primary health care, education, human right and mainly rights of women. As rates of HIV infection have risen in women of child bearing age, HIV/AIDS has evolved into a disease of families with children. The devastating effects of adults HIV infection are reflected in the increasing number of children, adolescents and young adults who are losing one or both parents due to AIDS. For this the women in rural areas should be empowered both with economy and education without which no program can have any effect on their lives. Women must question the presumption of mutual monogamy and deal with issues of trust in relationships.

Spread of HIV/AIDS due to blood transmission, injections or frequent abortion is not properly studied in women HIV/AIDS patients. Women often lack control over their sexual relationship. In marriage the pervasive threat of physical violence or divorce without legal recourse or legal rights to property may totally disempower women even if she knows about AIDS, condoms are available and she knows her husband is HIV infected.

Clearly therefore the central and major issue is the inferior role and subordinate status of women and the disadvantages created by society cannot be redressed through individually focused information and education of HIV specific and health services. However another aspect of patriarchy is the expectation that the man will protect his partner and provide for the family while women are supposed to give priority to the well being of their children, husband and her husband's family.

Women were either ignored in prevention approaches or those targeted were mostly the CSWs or potential child bearers. These approaches mainly instructed women to use condoms, often ignoring the socio-cultural context of their lives. A sustained awareness campaign for the uneducated rural poor must be chalked out, every specific thrust to women.

18. All T.V. and radio channel must talk or give 2 minutes awareness program before news everytime everyday so that this will certainly work in the minds of one and all.

19. Steps must be taken to monitor that the subsidized medicines given by the government reaches the patients without any middleman.

20. A new law must be made to punish the doctor who promise to cure HIV/AIDS patients and swindle them taking large amount of money legally.

21. Top directors can produce good movie about HIV/AIDS in all Indian languages and shows to prevent and create, awareness about HIV/AIDS. Similarly songs in pop styles can deal with AIDS prevention.

22. Most women talk the truth about their habits more openly than men.

23. Wasteland in India, if properly utilized, can certainly eradicate poverty. Allocation of land to poor women will not only empower them but shall also reduce the rural-to-urban migration of men.

24. The awareness of HIV/AIDS can be induced in the rural women along the propaganda of family planning and family welfare schemes.

25. The absence of immunity among rural women and men could be due to the use of insecticides, pesticides which has killed all most all the natural resources like vegetation, edible roots, certain natural edible snails, fishes etc which had rich nutrient and protein values. Thus serious research on this aspect is suggested to environmentalists and agriculture scientists. What is the reason for the death and annihilation of such natural resources? When they annihilate such natural living beings, will these pesticides not have drastic effect on the poor rural people, which may reduce the immune system?

26. The government must carry a compulsory test on every adult patient in rural areas who visit the health center or the hospital for any form of treatment. If they are found to be infected with HIV/AIDS they must be given proper counseling and then treatment by the doctors.

27. It is observed that women were very frank and honest and gave several suggestions. When the husbands become infected, in almost all cases, they infected their wives without any conscience. This needs a psychological analysis of HIV/AIDS affected husbands and their mental make-up in relation to their wives and CSWs.

28. The HIV/AIDS affected women should be financially supported by the state.

29. It has come as a rude shock that several quacks from Kerala are exploiting the poor uneducated rural persons by promising them of cure and they are there by draining all their paltry resources. This should be stopped and the quacks must be imprisoned for exploiting them and cheating the poor HIV/AIDS patients and making still poor and their disease chronic.

30. It is very unfortunate to learn that now brokers from Andhra Pradesh have taken up the admission to beds of HIV/AIDS patients in Tambaram Sanatorium. So it is pertinent that a special cell is set up to monitor these brokers and thus cleanse the hospital premises.

31. One of the major reasons for the rural uneducated poor women becoming the victims of HIV/AIDS in our nations traditional set up which has always made women as an object, an object mainly created to serve men without questioning. So even when the wife is well aware of the fact that her husband has HIV/AIDS or suspects he has HIV/AIDS she has to submit for sex. This is the culture of the traditional Hindu, India. Unless such foolish and slavish barriers are broken, it is not easy to save these women. Indian culture that too specially we mean the South Indian, all the more. Rural Tamil Nadu culture who has so much imbibed the Aryans Manusmirithi, forces women to lead a life of suffering and sacrifice. Now it may be the sufferings of the social stigma of being an HIV/AIDS patient. Unless this social set up is changed rural India would become another Africa, in HIV/ AIDS.

32. The movies / cinema always project men, even men with all bad qualities are portrayed as gods and women with all goodness and original nature is subdued by such bad / ill natured men. If this is on the screen from the very young age male think they are unparalled demigods. As especially in Tamil Nadu the only and the main pleasure or pass time or even full time job is spend their time in cinema theaters or before the Television. Thus the media is one of the major reasons for the increase in HIV/AIDS women patients in rural Tamil Nadu.

33. Unless Women empowerment in India especially in south India is encouraged and established that too among the rural uneducated poor women, women will continue to be passive HIV/AIDS victims.

(i)    The empowerment of women is possible by making property rights, that is the equal right with men to have both mobile and immobile property.

(ii)   Education, by this we do not mean degree or a pass or a school final certificate; education in the true sense should be given to women to first take care of themselves.

(iii)  Women should be policy makers at least for themselves. How can man made policies be favouring women!

(iv)   Special drive must be made to give rural women all support by government at least till the age of 21, marriages done below 21 years should be severely punished. For in India female children at the age of 12 or even less are married to men three times their age. In no nation this can occur. For this the only reason given is that they follow their tradition. That is why we have said in the name of tradition all evil are done, the worst among them is the evil done to the women.

(v)    The labour of women must be duly given credit to, when we say this we bring out this main problem. Suppose a women and a man of equal age are employed even say as a cooli the man is paid three times the wages they pay to a man. The hours of labour are the same. This difference only shows the male domination and their illegal way of legalizing their act. No women's organization has so far taken up this issue in a very serious way to render justice to them. Unless equality is given from the base, how can one expect equality. That is why, in all better job, all promotions are denied even if women is found 10 times better than men. That is why we see recognition of women's labour is to be made, for this women must be empowered and women's organization should take up this issue and fight for them.

34. It is pertinent to mention that in India even after 50 years of independence, India has so far not appointed a single women as a Chief Justice or has not become the President of India. This is clear to show what amount of justice, women in India can get in this male dominated society.

Unless men realize the importance of women and shed their ego it may not be an easy task to save rural poor uneducated women from HIV/AIDS.

We finally thank these 101 women who patiently gave us the interview because of which this book is possible.

# QUESTIONNAIRE USED TO INTERVIEW PEOPLE LIVING WITH HIV/AIDS

## MANUAL FOR FIELD RESEARCHERS

The points to be had in mind while interviewing the AIDS patients.

❑ Names of the AIDS patients, will not be disclosed and their identities shall be protected.

❑ The interviews of the AIDS patients shall be used by us purely for research purposes.

❑ The patients shall be let to speak out freely and fully. The interviewer shall not intercept when the AIDS patients are speaking (to make them tell what the interviewers have in their mind.)

❑ They shall always move with the AIDS patients with kindness and respect, care and compassion.

❑ When interviewing the patients for the first time, no questions disliked by the patients, shall be asked to them, unless they voluntarily come forward to give the answers.

- ❑ The interviewers shall select their questions according to the circumstances/ persons.

- ❑ While gathering information from the AIDS patients, who reside in the hospitals permanently, and from the out-patients of the hospitals an uniform technique shall not be used.

- ❑ Answers to questions with '*' will be graded.

- ❑ Certain specific AIDS patients shall be interviewed at least once or twice in a month.

- ❑ The behaviour of some selected patients in their homes, offices, buses and with their friends, shall be carefully observed and noted down.

- ❑ When these questions are found not adequate, more questions can be framed as per the requirements or depending upon the situation.

- ❑ The interview shall be opened and concluded depending upon how openly the AIDS patients speak out.

- ❑ The interview with AIDS patients shall begin with the casual remark that "the AIDS is only a disease which causes death just as all other disease which cause death".

- ❑ When situations arise in which the AIDS patients conceal certain information / misrepresent, the field researchers shall adopt some psychological techniques and note them down separately.

Interviewee : 
Date and place of interview : 
Circumstance of the interview : 
    ❑ Bedridden ❑ Could walk
    ❑ Could not speak ❑ Looked well
    ❑ Looked ill ❑ Was very feeble
    ❑ Could only sit ❑ Was depressed
    ❑ Normal ❑ Other observation:
    ❑ Somebody spoke ❑ Both of them spoke

**GENERAL PARTICULARS**

1. Name :
2. Age :
3. Husband(living/dead) :
4. Educational qualification :
5. Native place and Address :
6. Present Address :
7. Caste :
8. Religion :

**JOB-RELATED INFORMATION**

9. Profession :
10. What skills do you know? :
11. Do you do this job
    hereditarily? :
12. Have you any other
    sources of income? :
13. How much do you earn
    as monthly income? :
14. How did you spend
    your income? :

    ❑ On self ❑ On family

☐ For accommodation    ☐ All
☐ None

15. Whichever places your husband went from your native place to earn money? :

16. Is your husband HIV/AIDS infected :

☐ Yes    ☐ No
☐ Do not know    ☐ Not yet tested

17. Have you informed or not about your AIDS disease to :

☐ Friends    ☐ None
☐ Siblings (brothers/sisters)    ☐ In-laws
☐ Children    ☐ Spouse
☐ Other relatives    ☐ Strangers
☐ Boss    ☐ All

18. Do you think that your husband has the right to conceal his AIDS disease from you? :

☐ Yes    ☐ No
☐ Partially    ☐ Silent
☐ No reaction    ☐ Just looks sad
☐ Confused    ☐ Suppresses feeling

19. Do you want to talk about your AIDS disease to your friends?
(or, it is not necessary?) :

## MARRIAGE

20. Are you married?                    :
21. Your age at marriage                :
22. Love or arranged marriage           :
23. To whom were you given in
    marriage?                           :
    ❏ Hindu / Muslim/ Christian   ❏ Old man
    ❏ Middle aged man             ❏ Widower
    ❏ Unreactive                  ❏ Other

24. Are you satisfied with
    your marriage?                      :

    ❏ Fully                       ❏ Partially
    ❏ Silent                      ❏ Other reaction
    ❏ Not                         ❏ No reaction

25. If not fully satisfied, why
    did you agree to marriage?          :
26. Were you aware of AIDS
    before marriage?                    :
27. Some people have more
    than one wives/husbands.
    Do you have so…?                    :

## HUSBAND

28. Name                            :
29. Educational qualification       :
30. Age                             :
31. Occupation                      :
32. Govt.. /private sector job?     :
33. Is husband/wife alive?          :
    If no, did you remarry?         :

**CHILDREN**

34. How many children?　　　　　　:

35. Up to what standard the
    children have studied ?　　　:

36. Age gap between children?　:

37. What jobs the children
    are doing? (If over 18)　　:

38. Do your children know
    that you have this disease?　:

39. How you and your children
    were living before you
    contracted HIV?　　　　　　:

40. How are your children
    behaving towards you after
    they became aware of your
    AIDS disease?　　　　　　　:

    ❑ Cold and distant　　　❑ Never visit
    ❑ Kind　　　　　　　　❑ No concern
    ❑ As usual　　　　　　❑ Other reaction

41. Did you contract AIDS before
    /after your child's birth?　:

42. After which child, did you
    contract AIDS?　　　　　　:

43. Even after knowing that
    you were infected with AIDS
    did you have the child?　　:

44. Does that child have AIDS?　:

45. How you had been/have
    been maintaining the child?　:

46. What you expect the Govt.

to do for your children?   :

47. Whoever, and in whichever
    manner, can AIDS education
    be imparted to children?   :

    ❑ Advertisements in TV    ❑ Advts. in Paper
    ❑ Songs    ❑ Posters
    ❑ Dramas    ❑ Radio-broadcast
    ❑ Teaching in School    ❑ Others

48. What do you think about
    giving education to the
    AIDS–affected children
    together with others?   :

    ❑ No    ❑ Yes
    ❑ No reaction    ❑ Never
    ❑ Sign on the face    ❑ Other reaction

49. Are your children spastic?   :

**PARTICULARS OF PARENTS /IN-LAWS**

50. Name   :
51. Age   :
52. Educational Qualification   :
53. Occupation   :
      Govt. or Pvt. employment?   :
54. Whether alive?   :
55. How many wives for father?   :
56. Can you tell about the
    properties your parents
    and ancestors left for you?   :
57. Financial position: Before/After:

| | |
|---|---|
| ❑ Rich | ❑ Poor |
| ❑ In debt | ❑ Upper middle class |
| ❑ Lower middle class | ❑ Others |

58. When they gave you in marriage what did your parents provide them? :

59. When your parents/in-laws came to know that you have AIDS what was their psychological reaction? :

| | |
|---|---|
| ❑ Cold | ❑ No reaction |
| ❑ Scornful | ❑ Blaming |
| ❑ As usual | ❑ Depressed |
| ❑ Dejected | ❑ Aversion |
| ❑ Sad / anger | ❑ Other |

60. Describe the cruelty, if any, undergone by you due to your parents-in-law :

61. Describe the treatment that you received after having been infected by HIV/AIDS :

**PARTICULARS OF RESIDENCE**

62. Is it own house or rented? :

63. What is the nature of the house you live? :

64. If own house, are there any tenants living in it? :

65. If rented house what rent are you paying? :

66. If apartment whether
    government housing
    or private buildings? :

67. Whether nuclear family
    or a joint family ? :

68. Are you living in your
    ancestral home? :

69. If joint family, with
    whom are you living? :

70. What type of locality
    is your house situated? :

71. Whether any of your
    friends have this disease? :

72. How many of your husbands
    friends have contracted
    AIDS? :

73. What is their condition now? :

74. Have you spoken with your
    friends about AIDS? :

    ❑ Never                    ❑ Sometimes
    ❑ Fully                    ❑ Silent
    ❑ Aggressive               ❑ Pathetic
    ❑ Despair                  ❑ Anger
    ❑ Dejection                ❑ Other
    ❑ Shows changes in facial expression like fear,
      sadness

75. What is the nature of your
    husband's friends? What are
    their characteristics? :

    ❑ Careless                 ❑ Responsible
    ❑ Irresponsible            ❑ Commercial

| ❑ Escapist | ❑ For money |
| ❑ Not honest | ❑ Sadist |

76. Can you tell about husbands friendship in the present context? (if no reply is given we tick under these heads)

| ❑ Smoke | ❑ Drink |
| ❑ Visits CSWs | ❑ Very sensitive |
| ❑ Very Commercial | ❑ Commercial |
| ❑ Just sensitive | ❑ Silent |
| ❑ Happy | ❑ Unhappy |
| ❑ Voice chokes | ❑ Looked in contemplation |

77. Can you tell about the occupation of your friends? :

**COMMERCIAL SEX WORKERS**

78. Can you tell the reasons for your husband seeking CSWs :

| ❑ Urge for sex | ❑ Not satisfied |
| ❑ Male Ego | ❑ Silent |
| ❑ Doesn't like to answer | ❑ No reaction |
| ❑ Movie | ❑ Other reaction |

79. What do you think about the living standard of the CSW? :

| ❑ Forced into the trade | ❑ Destitute |
| ❑ Poverty | ❑ Family burden |
| ❑ Due to choice | ❑ For money |
| ❑ Others | ❑ No answer |

80. Do you think your husband
    and yourself got HIV/AIDS
    because of the CSWs?          :
81. Have you consulted your
    friends about the changes
    occurred in your body after you
    had sex with your husband?  :

    ❑ No                    ❑ Yes
    ❑ Fear                  ❑ Sad look
    ❑ Silent                ❑ Looks vacant

82. Did you note the change of
    behaviour in your husband
    after he came home after
    long stay out of town?        :

    ❑ Sad                   ❑ No reaction
    ❑ Happy                 ❑ Guilty
    ❑ Not guilty            ❑ Unhappy
    ❑ Others                ❑ Dejected

**STATE OF MIND OF THE PERSON**

83. Can you tell how you
    contracted AIDS disease?     :
84. In which year, month and
    day, the symptoms of AIDS
    began to show up?            :
85. After changes in the body
    became visible have you
    realized that it is AIDS?     :
86. Can you tell the gradual
    changes that happened in

the body as symptom of AIDS :

87. How was your physical
    condition before you
    contracted HIV/AIDS?      :
    Changes in mind/emotion

    ❏ Blank           ❏ Sad
    ❏ Unhappy         ❏ Irritated
    ❏ Desperate       ❏ Depressed
    ❏ Suicidal        ❏ Others

88. Can you give the physical
    changes in your husband
    after he contracted AIDS   :

89. To whom did you first
    reveal the news of your
    infection?                 :

90. How the people around
    you moved with before
    you contracted AIDS?       :

91. How did your relations
    behave with you after
    they became aware of
    your AIDS infection?       :

92. To what type of physicians
    did you go for consultation
    (Allopathy/Siddha/
    Ayurveda)                  :

93. What is the reasons for
    going to him?              :

94. Can you point out the
    difference between
    the behaviour of your
    relatives towards you
    at the early stages of your

AIDS disease and their
behaviour after your
disease became chronic?     :

95. Who took you to the
    hospital for treatment?
    self? husband? children?   :

96. Can you tell us about the
    manner in which your
    doctor is taking care?      :

    ❑ Well                      ❑ Not so nicely
    ❑ Unconcerned               ❑ No reaction
    ❑ Very caring               ❑ Committed
    ❑ Discriminatory            ❑ Uptight

97. Whom you doubt as the
    person who caused this
    damage to you?              :

98. Do you want to disclose
    to others about your AIDS
    disease?(or, do you want
    to conceal it?)             :

99. Please tell the reasons?    :

100. After knowing about your
     infection, did you have
     sex with your husband?     :

101. Does anyone else in
     your family have AIDS?     :

102. When you came to know
     that you had AIDS, what
     was your mental state?     :

     ❑ Sad                      ❑ Ashamed
     ❑ Highly depressed         ❑ Suicidal

- ❑ No feeling
- ❑ High heartbeat
- ❑ Unhappy
- ❑ Lost Happiness
- ❑ Angry with your husband

- ❑ Shocked
- ❑ Wept uncontrollably
- ❑ Disgraced
- ❑ Courageous

103. What is your present mental condition after your working in the early stages? :

104. At the time when AIDS was said to be a horrible disease some said not to worry about it. What do you think about this? :

105. How long the impact of the AIDS, talks about it, the reading about it and the thoughts about it, were working upon you? :

- ❑ Never
- ❑ No reaction
- ❑ Silent

- ❑ Always
- ❑ Passive
- ❑ Indeterminate

106. What was your mental condition when you read or heard news on AIDS in newspapers, radio and television after you contracted AIDS? :

- ❑ Sad
- ❑ Highly depressed
- ❑ No feeling

- ❑ Ashamed
- ❑ Suicidal
- ❑ Shocked

| | |
|---|---|
| ❑ High heartbeat | ❑ Wept |
| ❑ Unhappy | ❑ Disgraced |
| ❑ Lost Happiness | ❑ Courageous |
| ❑ Suspected my husband | ❑ Blamed husband |

107. What do you expect the mental state of other AIDS patients when they see AIDS advertisements : (pick from above)

108. How was your daily life pattern was affected after you contracted AIDS?:

| | |
|---|---|
| ❑ Socially | ❑ Economically |
| ❑ Mentally | ❑ Materially |
| ❑ Spiritually | ❑ Physically |
| ❑ None | ❑ All |

109. If you had been educated already about AIDS, could you have escaped from getting this disease from your husband? :

110. If not, what reason can you tell? :

**GENERAL QUESTIONS**

111. Have you consulted your family members regarding AIDS? :

112. Are there AIDS patients in your family/relations? :

113. Do you know any friend

of yours who has HIV? :

114. What do you think is the duty of your husband's family towards you? :

115. What you think about how your family members/ relations/friends should behave in taking care of the AIDS patients? :

☐ Kind                    ☐ Harsh
☐ Open                    ☐ Considerate

116.* What is your view on the opinion that if medicine is discovered to fully cure AIDS, social deterioration will further be worsened? :

117.* What is your opinion on the concept that the government should call all the women AIDS patients for free medical examination/checkup? :

118.* What do you think about those who say: AIDS will not touch me. Don't talk about it to me? :

119. Do you think a wife has the right to say no to sex if she knows her husband carries the HIV infection? :

120.* What is your view on the statement that "it is the

social duty of everyone
to help AIDS patients?"             :

122. What do you think
     about sex education?           :

123. What do you think about
     yellow journalism and
     pornographic/blue films?       :

124. Different views prevail in
     different countries on the
     use of contraceptives.
     What is your opinion?          :

125. Whether one has or has
     not the right to know
     whether his colleagues
     in his workplace has AIDS?     :

126.* It is natural to have a
      fear while speaking,
      reading and commenting
      about AIDS. Is there a way
      to change this altitude?      :

127. Can you differentiate
     between HIV/AIDS/STD
     from one another?              :

128. Do you know how many
     years will it take for HIV
     to convert into AIDS?          :

129.* Have you had a discussion
      with your colleagues and
      friends on AIDS?              :

130. About what you would talk
     during such discussion?        :

131.* What is your mental state
      on seeing an AIDS patient?    :

132. Will you explain to the
     public about AIDS if you
     are given an opportunity
     to do so?                          :

133. Many people speak
     differently about AIDS.
     You might have heard/
     spoken to some of them
     Can you recall any one
     unforgettable incident?            :

134. Is the AIDS patient
     alone responsible for
     contracting the disease?           :

135. "When anyone sees/hears
     an AIDS patient, s/he is
     emotionally surcharged.
     S/he is subjected to
     mental affliction/grief/
     anger/insult/frustration
     and behaves so. Can you
     tell the reason for such
     behaviour?                         :

136. "The foremost duty of the
     government is to ascertain
     how many citizens are
     suffering from AIDS. But
     no such precise assessment
     has been made in India."
     What is your opinion?              :

137. Why no free medical
     test is not being conducted
     on all citizens to find out
     whether anyone has AIDS
     like for smallpox/

diabetes/leprosy?      :

138.* We talk about diabetes,
     leprosy and cancer in all
     places. But, we do not
     talk about AIDS. Why?      :

139. What do you think about
     those who say "Increasing
     of AIDS disease is not my
     problem. It is the problem
     of concerned persons".      :

140.* "When you know your
     friend has AIDS how do you
     move with him" why?      :

141. What expenses you want
     the government to meet
     for AIDS control?      :

142. What further measures do
     you expect the government
     to take to control AIDS?      :

143. Which of the measures
     the government is under-
     taking to prevent AIDS is
     most attractive to you?      :

144. In which medical system
     you have faith: Allopathy,
     Siddha/Ayurveda/Unani?
     The reason for it?      :

145. Do you think that AIDS will
     spread through mosquito
     biting/kissing on hand/lips
     and eating in same plate?      :

146. Which advertisements for
     AIDS awareness is mostly

liked by you? :

147. What do you think when
you see in private the
advertisements on AIDS? :

❑ Fear ❑ Guilt
❑ Unconcerned ❑ Sadness
❑ Absence of any feelings ❑ Indeterminate

148. What do you feel when you
see advertisements in the
company of your family? :

❑ Fear ❑ Guilt
❑ Unconcerned ❑ Sadness

149. What do you think of AIDS
awareness advertisements? :

150. Can you suggest good
advertisement techniques
for AIDS awareness? :

151. What further steps the
governments must take to
eradicate HIV transmission
to women by their husbands?:

152. Must all rural women agree
to undergo AIDS test if the
government comes forward
to do the test freely? :

153. What kind of expense you
want the government to
bear for the AIDS control? :

154. Do you think that the
uneducated women from

AIDS rural areas need
financial help and care. :

155. Do you think it is good
initiative if the govt. brings
a special quota in jobs for
people who are HIV+? :

156. "Any disease can be
contracted in the world.
But the AIDS should never
come" what do you think
about this? :

157.* "Women AIDS patients are
to be given top priority in
medical treatment above
all other patients". What is
your view on this opinion? :

158. Do you think that the wife
should also be present when
the doctors discuss the
husband's HIV infection? :

159.* What do you think about
interviewing you by people
like us? :

160. What other questions you
want us to ask you, other
than those we asked you? :

161. Please tell what part of our
questions, you like/don't
like the most? :

**ADVERTISEMENTS**

162.*Through which source you

came to know about AIDS
before you contracted it? :

163. If you came to know it
through advertisements,
what type was the media? :

164. Our govt. has provided
telephone facilities in 77
places throughout the nation
to tell confidentially doubts
questions therein. Have you
used this facility? Do you
know about this facility? :

165. By whatever further steps
the government can take
the advts. to the masses
to prevent the spread of
HIV/AIDS? :

166.* Now you regret that the
immoral behaviour of your
husband which started one
day has now brought dreaded
AIDS disease to you. What do
you think about it now? :

167.*What is your view about
fully eradicating prostitution
from our country? :

168. How can the NGO/social
organizations make the
advertisements on AIDS
reach the rural woman? :

169. In what ways can your
religion publicize AIDS
among rural woman? :

**FOOD HABITS**

170. What kind of food habits
     your family was used to
     vegetarian/meat-eating? :

171. Have you changed your
     food habits after this
     disease contracted you? :

172. If you have changed, how? :

173. How is your physical
     conditions after you changed
     your food habits? :

**RELIGION**

174. Do you profess the same
     religion hereditarily? If not,
     in which generation you
     converted to this religion? :

175. Have you converted into
     any other religion after
     you contracted AIDS? :

176. Did you convert? From
     which religion? :

177. What services you rendered
     to people of that religion
     after your conversion? Or
     was it for personal solace? :

178. What reasons you can
     give for your conversion? :

179. Do you think the religion is
     giving you consolation after

your conversion? Or is it a
testing tool? Is it a escape?   :

180. Did both of your parents
belong to a same religion?   :

181. Did your husband's parents
belong to a same religion?   :

182. Are you continuing to visit
the religious institutions
after you contracted AIDS?   :

183. Do you continue to believe
in God after getting AIDS?   :

# TABLE OF STATISTICS

| No | Age | Age at Marriage | Education | Infected by | Years living with HIV/AIDS |
|---|---|---|---|---|---|
| 1. | 45 | | Nil | Husband (dead) | |
| 2. | Did not say | | Nil | Others | 10 years |
| 3. | Did not say | | Nil | Others | 13 years |
| 4. | 46 | 15 | 8$^{TH}$ | Husband | 8 years |
| 5. | 35 | 25 | 5$^{TH}$ | Husband | Less than a year |
| 6. | 21 | | Nil | Husband | 3 years |
| 7. | 27 | 15 | Nil | Husband | |
| 8. | Did not say | | Nil | Others | 8 years |
| 9. | Did not say | 13 | Did not say | Husband (dead) | Just 6 months |
| 10. | 40 | 18 | Nil | Do not know | 4 years |
| 11. | 40 | 13 | Nil | Husband | 1 ½ years |
| 12. | 30 | 18 | Nil | Husband | - |
| 13. | 22 | 12 | 2$^{ND}$ | Husband | 20 days |

| | | | | |
|---|---|---|---|---|
| 14. | Does not know | 13 | Nil | Husband (dead) | 2 years |
| 15. | 31 | 12+ | 8TH | Husband | |
| 16. | 35 | 15 | NIL | Husband (dead) | |
| 17. | 25 | 19 | 3RD | Due to four abortions | 20 days |
| 18. | 40 | 18 | 6TH | Husband (dead) | 2 years |
| 19. | 30 | | + 2 | Husband | 10 days |
| 20. | 23 | 21 | 6TH | Husband | |
| 21 | 35 | 27 | 8TH | Husband (dead) | |
| 22 | 32 | 15 | 5TH | Husband | 2 years |
| 23. | 45 | 17 | 1ST | Husband | |
| 24. | 40 | | 6TH | Husband | one month |
| 25 | 40 | | 8TH | Husband | Four days |
| 26 | 38 | | Nil | Neighbour | 3 months |
| 27 | 31 | | 6TH | Husband | visited three times |
| 28. | 25 | | Nil | Husband | |
| 29. | Did not say | | Did not say | Husband | |
| 30. | 37 | 15 | Nil | Husband (dead) | |
| 31. | 35 | 14 | Nil | Husband (dead) | |
| 32. | Did not say | | Nil | Husband (dead) | |

| | | | | |
|---|---|---|---|---|
| 33. | 40 | | Nil | Husband | 2 years |
| 34. | 27 | 25 | 12$^{TH}$ | Husband (dead) | few days |
| 35. | 16 | 11 | Nil | Husband (dead) | |
| 36. | 25 | 18 | 6$^{TH}$ | Husband | |
| 37. | 25 | 19 | 5$^{TH}$ | Husband | |
| 38. | 30 | 16 | 6$^{TH}$ | Husband | 1 ½ months |
| 39. | 28 | 17 | 7$^{TH}$ | Husband (dead) | |
| 40. | 35 | | Nil | Husband (dead) | 4 years |
| 41. | 24 | 22 | Nil | Husband | |
| 42. | 45 | 25 | Nil | Husband (doubtful) | |
| 43. | 25 | 21 | 8$^{TH}$ | Husband | |
| 44. | 35 | 16 | Nil | Husband (dead) | 3 years |
| 45. | 25 | 15 | 7$^{TH}$ | Husband | |
| 46. | 22 | 17 | 5$^{TH}$ | Husband (dead) | 3 years |
| 47. | 35 | 18 | 7$^{TH}$ | Husband (dead) | 1 ½ year |
| 48. | 26 | | 10$^{TH}$ | Husband | |
| 49. | 30 | 11 | Nil | Did not say | |
| 50. | 30 | | 8$^{TH}$ | Husband | 2 ½ years |
| 51. | 32 | 13 | Nil | Husband (dead) | |
| 52. | 26 | 13 | 4$^{TH}$ | Husband | 2 years |

| 53. | 35 | 15 | Nil | Husband (dead) | |
| 54. | 35 | 13 | Nil | Husband | |
| 55. | 30 | 17 | 6TH | Husband (dead) | |
| 56. | Did not say | | - | Husband | 1½ months |
| 57. | Did not say | | Nil | Husband | |
| 58. | Did not know | | 5TH | Husband (dead) | |
| 59. | 45 | | 7TH | Unknown | |
| 60. | Did not say | | Nil | Husband (dead) | |
| 61. | 50 | | Nil | Husband (dead) | few days |
| 62. | 42 | | Nil | Husband | |
| 63. | 22 | 16 | SSLC | Husband | |
| 64. | Did not say | | Nil | Husband | |
| 65. | 40 | | Nil | Husband (dead) | |
| 66. | 26 | 17 | 10TH | Husband (dead) | |
| 67. | 50 | 19 | Nil | Husband (dead) | 2 years |
| 68. | 24 | 19 | 5TH | Husband (dead) | 1 year |
| 69. | 35 | 17 | 7TH | Husband (dead) | 2 years |
| 70. | 27 | 15 | Nil | Husband (dead) | |
| 71. | 41 | 16 | Nil | Husband (dead) | |
| 72. | 31 | 25 | Nil | Husband (dead) | |
| 73 | 22 | 21 | 10$^{TH}$ | unknown | |
| 74. | 37 | 25 | 4$^{TH}$ | unknown, abortions | |
| 75. | 28 | 18 | Nil | Husband | 2 years |
| 76. | 32 | 22 | Nil | Husband | 2 years |
| 77. | 25 | 21 | 5$^{TH}$ | Husband | 3 years |
| 78. | 28 | 15 | 7$^{TH}$ | Husband (dead) | |
| 79. | 38 | 20 | 8$^{TH}$ | Husband (dead) | 3 ½ years |
| 80. | 38 | 13 | 5$^{TH}$ | Husband (dead) | 7 years |
| 81. | 25 | 19 | +2 | Husband | |
| 82. | 38 | - | 10$^{TH}$ | Husband | 1 year |
| 83. | 38 | 13 | 8$^{TH}$ | Husband | 1½ months |
| 84. | Did not say | 21 | 4$^{TH}$ | Husband | 4 months |
| 85. | Did not say | Did not say | Nil | Husband (dead) | |
| 86. | 21 | 16 | 10$^{TH}$ | Husband | |
| 87. | 40 | 15 | Nil | Husband (not infected) | |
| 88. | 36 | 13 | Nil | Husband (dead) | |
| 89. | 27 | 23 | 5th Std. | Husband | |
| 90 | 27 | 17 | 5th Std. | Husband | |
| 91 | 23 | 16 | Nursing | Husband (dead) | 4 years |

| 92 | 32 | 12 | Nil | Husband | |
| 93 | 27 | 14 | 9th Std. | Husband | 1 ½ years |
| 94. | 43 | 17 | 8th Std. | Husband | 15 years |
| 95. | 30 | 19 | 8th Std. | Husband | |
| 96. | 31 | 20 | Nil | Husband | 1½ years |
| 97. | 45 | - | 3rd Std. | Husband | 2 years |
| 98 | 22 | 15 | 6th Std. | Husband | 6 month |
| 99. | 35 | | Nil | Husband | 1 year |
| 100. | 31 | | 10th Std. | Husband | |
| 101 | 28 | | 12th Std. | Husband (dead) | |

# NEUTROSOPHY

The neutrosophics are based on neutrosophy, which is an extension of dialectics. Neutrosophy is a theory developed by Florentin Smarandache as a generalization of dialectics.

For more on neutrosophy, please refer:

http://arxiv.org/ftp/math/papers/0010/0010099.pdf

This theory considers every notion or idea <A> together with its opposite or negation <Anti-A> and the spectrum of "neutralities" <Neut-A> (i.e. notions or ideas located between the two extremes, supporting neither <A> nor <Anti-A>). The <Neut-A> and <Anti-A> ideas together are referred to as <Non-A>. The theory proves that every idea <A> tends to be neutralized and balanced by <Anti-A> and <Non-A> ideas—as a state of equilibrium. It was K.Atanassov (1986) who first introduced the Intuitionistic Fuzzy Set(IFS), describing a degree of membership of an element to a set , and a degree of non-membership of that element to the set, what's left was considered as indeterminacy. This was generalized to Neutrosophic Set (NS) in 1995 by Florentin Smarandache.

**Distinctions between Intutionistic Fuzzy Logic (IFL) and Neutrosophic Logic (NL)**

The distinctions between IFL and NL {plus the corresponding intuitionistic fuzzy set (IFS) and neutrosophic set (NS)} are the following.

a) Neutrosophic Logic can distinguish between *absolute truth* (truth in all possible worlds, according to Leibniz) and *relative truth* (truth in at least one world), because NL(absolute truth)=$1^+$ while NL(relative truth)=1. This has application in philosophy (see the neutrosophy). That's why the unitary standard interval [0, 1] used in IFL has been extended to the unitary non-standard interval ]$^-$0, $1^+$[ in NL. Similar distinctions for absolute or relative falsehood, and absolute or relative indeterminacy are allowed in NL.

b) In NL there is no restriction on T, I, F other than they are subsets of ]$^-$0, $1^+$[, thus: $^-0 \leq \inf T + \inf I + \inf F \leq \sup T + \sup I + \sup F \leq 3^+$. This non-restriction allows paraconsistent, dialetheist, and incomplete information to be characterized in NL {i.e. the sum of all three components if they are defined as points, or sum of superior limits of all three components if they are defined as subsets can be >1 (for paraconsistent information coming from different sources) or < 1 for incomplete information}, while that information can not be described in IFL because in IFL the components T (truth), I (indeterminacy), F (falsehood) are restricted either to t + i + f=1 or to $t^2 + f^2 \leq 1$, if T, I, F are all reduced to the points t, i, f respectively, or to sup T + sup I + sup F = 1 if T, I, F are subsets of [0, 1].

c) In NL the components T, I, F can also be *non-standard* subsets included in the unitary non-

standard interval $]^-0, 1^+[$, not only *standard* subsets included   in the unitary standard interval $[0, 1]$ as in IFL.

d)  NL, like dialetheism, can describe     paradoxes, NL(paradox) = (1, I, 1), while IFL can not describe a paradox because the sum of components should be 1 in IFL

## Distinctions between Neutrosophic Set (NS) and Intuitionistic Fuzzy Set (IFS).

a)  Neutrosophic Set can distinguish between *absolute membership* (i.e. membership in all possible worlds; we have extended Leibniz's absolute truth to absolute membership) and *relative membership* (membership in at least one world but not in all), because NS(absolute membership element)=$1^+$ while NS(relative membership element)=1.   This has application in philosophy (see the neutrosophy). That's why the unitary standard interval $[0, 1]$ used in IFS has been extended to the unitary non-standard interval $]^-0, 1^+[$ in NS. Similar distinctions for *absolute or relative non-membership*, and *absolute or relative indeterminant appurtenance* are allowed in NS.

b)  In NS there is no restriction on T, I, F other than they are subsets of $]^-0, 1^+[$, thus: $^-0 [ \inf T + \inf I + \inf F [ \sup T + \sup I + \sup F [ 3^+$.
The inequalities (2.1) and (2.4) of IFS are relaxed.

This non-restriction allows paraconsistent, dialetheist, and incomplete information to be characterized in NS {i.e. the sum of all three components if they are defined as points, or sum of superior limits of all three components if they are defined as subsets can be >1 (for paraconsistent information coming from different sources), or < 1 for incomplete information}, while that information can not be described in IFS because in IFS the components T (membership), I (indeterminacy), F (non-membership) are restricted either to $t + i + f = 1$ or to $t^2 + f^2 [ 1$, if T, I, F are all reduced to the points t, i, f respectively, or to sup $T +$ sup $I +$ sup $F = 1$ if T, I, F are subsets of [0, 1].

Of course, there are cases when paraconsistent and incomplete informations can be normalized to 1, but this procedure is not always suitable.

c) Relation (2.3) from interval-valued intuitionistic fuzzy set is relaxed in NS, i.e. the intervals do not necessarily belong to Int[0,1] but to [0,1], even more general to ]-0, 1+[.

d) In NS the components T, I, F can also be *non-standard* subsets included in the unitary non-standard interval $]^-0, 1^+[$, not only *standard* subsets included in the unitary standard interval [0, 1] as in IFS.

e) NS, like dialetheism, can describe paradoxist elements, NS(paradoxist element) = (1, I, 1), while IFL can not describe a paradox because the sum of components should be 1 in IFS.

f) The connectors in IFS are defined with respect to T and F, i.e. membership and non-membership only (hence the Indeterminacy is what's left from 1), while in NS they can be defined with respect to any of them (no restriction).

g) Component "I", indeterminacy, can be split into more subcomponents in order to better catch the vague information we work with, and such, for example, one can get more accurate answers to the *Question-Answering Systems* initiated by Zadeh (2003). {In Belnap's four-valued logic (1977) indeterminacy is split into Uncertainty (U) and Contradiction (C), but they were inter-related.}

# Index

**U**

Unit hypercube, 58
Unsupervised, 16

**W**

WORLD (Women Organized to Respond to Life threatening Diseases), 10

# About the Authors

**Dr.W.B.Vasantha Kandasamy** is an Associate Professor in the Department of Mathematics, Indian Institute of Technology Madras, Chennai, where she lives with her husband Dr.K.Kandasamy and daughters Meena and Kama. Her current interests include Smarandache algebraic structures, fuzzy theory, coding/ communication theory. In the past decade she has guided eight Ph.D. scholars in the different fields of non-associative algebras, algebraic coding theory, transportation theory, fuzzy groups, and applications of fuzzy theory of the problems faced in chemical industries and cement industries. Currently, six Ph.D. scholars are working under her guidance. She has to her credit 285 research papers of which 209 are individually authored. Apart from this, she and her students have presented around 329 papers in national and international conferences. She teaches both undergraduate and post-graduate students and has guided over 45 M.Sc. and M.Tech. projects. She has worked in collaboration projects with the Indian Space Research Organization and with the Tamil Nadu State AIDS Control Society. She has authored a Book Series, consisting of ten research books on the topic of Smarandache Algebraic Structures which were published by the American Research Press.

She can be contacted at vasantha@itm.ac.in
You can visit her work on the web at: http://mat.iitm.ac.in/~wbv

---

**Dr.Florentin Smarandache** is an Associate Professor of Mathematics at the University of New Mexico, Gallup Campus, USA. He published over 60 books and 80 papers and notes in mathematics, philosophy, literature, rebus. In mathematics his research papers are in number theory, non-Euclidean geometry, synthetic geometry, algebraic structures, statistics, and multiple valued logic (fuzzy logic and fuzzy set, neutrosophic logic and neutrosophic set, neutrosophic probability). He contributed with proposed problems and solutions to the Student Mathematical Competitions. His latest interest is in information fusion were he works with Dr.Jean Dezert from ONERA (French National Establishment for Aerospace Research in Paris) in creasing a new theory of plausible and paradoxical reasoning (DSmT).

He can be contacted at smarand@unm.edu